\documentclass[multiarticle,nochapter,download]{msjmemoirs}
\usepackage{amsmath,amsthm,amssymb}
\usepackage{bbm}

\usepackage{makeidx}
\makeindex

\usepackage[nottoc]{tocbibind}

\usepackage{etoolbox}
\patchcmd{\theindex}{\thispagestyle{plain}}{}{}{}
\usepackage{epic,eepic}
\usepackage[all]{xy}
\usepackage{tikz}
\usetikzlibrary{intersections,calc}

\allowdisplaybreaks
 
\usepackage[T1]{fontenc}
\usepackage{listings}

\makeatletter
\newcommand\BeraMonottfamily{%
  \def\fvm@Scale{0.85}
  \fontfamily{fvm}\selectfont
}
\makeatother

\newtheorem{thm}{Theorem}[section]
\newtheorem{prop}[thm]{Proposition}

\newtheorem{cor}[thm]{Corollary}
\newtheorem{lem}[thm]{Lemma}

\theoremstyle{definition}

\numberwithin{equation}{section}
\newtheorem{rem}[thm]{\bf Remark}
\newtheorem{ex}[thm]{\bf Example}

\newtheorem{defn}[thm]{\bf Definition}

\newtheorem{const}[thm]{\bf Construction}

\newtheorem{algo}[thm]{\bf Algorithm}

\def\notes{\medskip \par\noindent{\centerline{\bf Notes}}\par\nopagebreak\medskip\par\nopagebreak}

\def\bbk{\mathbbm{k}}

\def\bbQ{\mathbb{Q}}
\def\bbR{\mathbb{R}}
\def\bbT{\mathbb{T}}

\def\bbZ{\mathbb{Z}}
\def\bbvarphi{\boldsymbol{\varphi}}
\def\bbpsi{\boldsymbol{\psi}}

\def\rmin{\mathrm{in}}

\def\bfG{\mathbf{G}}

\def\bfa{\mathbf{a}}

\def\bfe{\mathbf{e}}
\def\bfg{\mathbf{g}}
\def\bfh{\mathbf{h}}
\def\bfn{\mathbf{n}}

\def\bfv{\mathbf{v}}
\def\bfw{\mathbf{w}}
\def\bfx{\mathbf{x}}
\def\bfy{\mathbf{y}}
\def\bfz{\mathbf{z}}

\def\rmid{\mathrm{id}}

\def\prin{\mathrm{prin}}
\def\Ker{\mathrm{Ker}}

\def\rmpr{\mathrm{pr}}

\def\calB{\mathcal{B}}
\def\calC{\mathcal{C}}

\def\calH{\mathcal{H}}

\def\calR{\mathcal{R}}

\def\frakC{\mathfrak{C}}
\def\frakc{\mathfrak{c}}
\def\frakd{\mathfrak{d}}
\def\frakD{\mathfrak{D}}
\def\frakg{\mathfrak{g}}
\def\frakh{\mathfrak{h}}
\def\frakj{\mathfrak{j}}
\def\frakp{\mathfrak{p}}

\def\fraks{\mathfrak{s}}

\def\d{\delta}

\def\l{\ell}
\def\ss{\scriptstyle}
\def\kP{\bbk[[P]]}

\def\kQ{\bbk[[Q]]}

\newcommand\mapsfrom{\mathrel{\reflectbox{\ensuremath{\mapsto}}}}

\usepackage{etoolbox}
\makeatletter
\pretocmd{\@startsection}{\gdef\thesectiontype{#1}}{}{}
\pretocmd{\@sect}{\@namedef{the\thesectiontype title}{#8}}{}{}
\pretocmd{\@ssect}{\@namedef{the\thesectiontype title}{#5}}{}{}
\makeatother

\usepackage{fancyhdr}
\begin{document}
\bibliographystyle{amsalpha}
\fontsize{11pt}{12.5pt}\selectfont

\mainmatter[page=1]
\title[\ ]{Cluster Algebras and Scattering Diagrams\\
\ \\
Part III\\
Cluster Scattering Diagrams\thanks{
This is the final manuscript of
 Part III of the monograph ``Cluster Algebras and Scattering Diagrams'',
 MSJ Mem. 41 (2023) by the author. 
} }
\author[\ ]{Tomoki Nakanishi\\
\small Graduate School of Mathematics, Nagoya University}
\maketitle

\begin{abstract}
This is a self-contained exposition of several fundamental properties of  cluster scattering diagrams  introduced and studied by Gross, Hacking, Keel, and Kontsevich.
In particular,   detailed proofs are presented for the construction,  the   mutation invariance, and the positivity of theta functions of cluster scattering diagrams.
 Throughout the text we highlight the fundamental roles of the  dilogarithm elements and the  pentagon relation in cluster scattering diagrams.
\end{abstract}

\newpage
\tableofcontents[depth=2]

\newpage


\setcounter{section}{-1}
\pagestyle{fancy}
\fancyhead[CO]{\thesection.\ \thesectiontitle}

\section{Introduction to Part III}
\label{3sec:intro1}

This is a self-contained exposition of several fundamental properties of
\emph{cluster scattering diagrams} (CSDs for short) 
introduced and studied  by Gross, Hacking, Keel, and Kontsevich \cite{Gross14}.

 Scattering diagrams originated
 in the study of homological mirror symmetry,
 and they are developed
 in the works of
 Kontsevich-Soibelman \cite{Kontsevich06, Kontsevich08, Kontsevich13},
  Gross-Siebert \cite{Gross07},
Gross-Pandharipande-Siebert \cite{Gross09},
  Gross \cite{Gross11},
  Carl-Pumperla-Siebert \cite{Carl10}, etc.
The connection to cluster algebras was established
   by  Gross-Hacking-Keel-Kontsevich \cite{Gross14},
  where CSDs were introduced.
  Then, with a scattering diagram technique,
    several important conjectures on   cluster algebras
          were proved.

Here is a brief description of the contents.
  \begin{itemize}
  \item
  In the first half (Sections 1--3)
  we present a construction of  CSDs
  without assuming any knowledge on scattering diagrams.
The main result here is the existence and the uniqueness of 
a CSD for a given seed (Theorem \ref{3thm:exist2}).
\item
  In the second half (Sections 4--7)
  we present  the fundamental properties of CSDs,
  which are necessary to prove several conjectures in cluster algebras,
including
  \begin{itemize}
  \item
  the mutation invariance (Theorem \ref{3thm:mut1}),
  \item
  the positivity of theta functions  (Theorems  \ref{3thm:pos1} and \ref{3thm:positivetheta1}).
\end{itemize}
   \item
 We present detailed  proofs for all propositions
 in a self-contained way.
  \end{itemize}
  
 The connection between CSDs and cluster patterns are given in Part II
based on the above results.
 Thus, this part is regarded as a supplement to Part II.
 Also, it may be read independently as an introductory text
 to CSDs.
 
 While the most results are taken from \cite{Gross14,Kontsevich13},
 there are several  added or modified proofs and new results due to us.
Most notably, we employ the following approach, which adds
 novelty to our presentation:
 \begin{itemize}
 \item
 We clearly separate the structure group of walls
 and its faithful representation (the \emph{principal $x$-representation}) 
 employed for the construction of CSDs in \cite{Gross14}.
 In the first half we only use the structure group itself.
 This enables us to lift the assumption of  the nondegeneracy for the fixed data
 in \cite{Gross14},
 thereby makes the construction of CSDs more transparent.
In the second half we apply the  principal $x$-representation to 
the wall elements when needed.
 \item
 Throughout 
 we highlight the fundamental role of the \emph{dilogarithm elements\/} and the \emph{pentagon relation\/}
 in CSDs, which is implicit in \cite{Gross14}.
 \end{itemize}

More recently, the \emph{quantum CSDs\/} 
 were studied by Davison-Mandel \cite{Davison19},
 where the counterparts of many results presented here were  given.

\newpage
\fancyhead[CO]{\thesubsection.\ \thesubsectiontitle}

\section{Underlying algebraic structure}
\label{3sec:scat1}

In this section we introduce
some underlying algebraic structure
behind scattering diagrams
we are going to study.

\subsection{Fixed data and seed}
\label{3subsec:fixed1}

\begin{defn}[Fixed data]
\label{3defn:fixed1}
A \emph{fixed data\/} \index{fixed data}
 $\Gamma$ consists of the following:
\begin{itemize}
\item
A lattice $N\simeq \bbZ^r$ of rank $r$ with
a skew-symmetric  bilinear form
\begin{align}
\label{3eq:bili1}
\{\cdot, \cdot \}: N \times N \rightarrow \bbQ.
\end{align}
\item
A sublattice $N^{\circ}\subset N$ of finite index (equivalently, of rank $r$) such that
\begin{align}
\label{3eq:bili2}
\{N^{\circ}, N\}\subset \bbZ.
\end{align}
\item
Positive integers $\d_1,\dots, \d_r$
such that
there is a basis $(e_1,\dots,e_r)$ of $N$, where
 $(\d_1e_1, \dots, \d_r e_r)$ is a basis of $N^{\circ}$.
\item
$M=\mathrm{Hom}(N,\bbZ)$, $M^{\circ}=\mathrm{Hom}(N^{\circ},\bbZ)$.
\end{itemize}
We say that $\Gamma$ is \emph{nondegenerate\/} (resp. \emph{degenerate\/})
if $\{\cdot, \cdot\}$ is  nondegenerate (resp. degenerate). \index{fixed data!nondegenerate}
\end{defn}
For a given $N^{\circ}\subset N$ of finite index,
the above $\d_1,\dots, \d_r$ always exist; for example, take the
 elementary divisors of the embedding $N^{\circ}\hookrightarrow N$.
 However, other choices are equally good.
 
Let $M_{\bbR}=M\otimes_{\bbZ} \bbR\simeq \bbR^r$,
and we regard
\begin{align}
M\subset M^{\circ} \subset M_{\bbR}.
\end{align}
 Let $\langle \cdot, \cdot \rangle: N^{\circ}\times M^{\circ} \rightarrow \bbZ$ be the
 canonical paring.
 We also write its linear extension 
 $N\times M_{\bbR} \rightarrow \bbR$
by the same symbol.

\begin{defn}[Seed]
\label{3defn:seed2}\index{seed!for a fixed data}
A \emph{seed\/} $\fraks=(e_1,\dots,e_r)$ for a fixed data $\Gamma$
is a basis of $N$ such that
$(\d_1e_1, \dots, \d_r e_r)$ is a basis of $N^{\circ}$.
\end{defn}

For a seed $\fraks$ for $\Gamma$, 
we have a basis $(e_1^*$, \dots, $e_r^*)$ of $M$ that is  dual to $(e_1,\dots,e_r)$.
Let $f_i=e^*_i/\d_i$.
Then, 
we also have a basis $(f_1,\dots,f_r)$ of $ M^{\circ}$ that is dual to the basis $(\d_1 e_1,\dots, \d_r e_r)$
of $N^{\circ}$.

Let $\Omega=\Omega_{\Gamma, \fraks}=(\omega_{ij})$ be the $r\times r$ skew-symmetric rational matrix
defined by $\omega_{ij}=\{e_i,e_j\}$, and
let $\Delta=\Delta_{\Gamma, \fraks}=\mathrm{diag}(\d_1, \dots, \d_r)$ be the diagonal  matrix
with  diagonal entries $\d_1, \dots,  \d_r$.
Let $B=B_{\Gamma, \fraks}=(b_{ij})$ be the $r\times r$ integer matrix defined by
\begin{align}
\label{3eq:Bmat1}
B=\Delta\Omega,
\quad
b_{ij}=\{\d_i e_i, e_j\}.
\end{align}
Since $\Omega=\Delta^{-1}B$ is skew-symmetric,
 $B$ is \emph{skew-symmetrizable\/} 
with a (left) skew-symmetrizer  $D=\Delta^{-1}$
in the  sense of \cite{Fomin02}.
 \index{matrix!skew-symmetrizable}
 
 Conversely, for any  $r\times r$ skew-symmetrizable integer matrix $B=(b_{ij})$,
 there is a (not unique) decomposition
 \begin{align}
 \label{3eq:BDO1}
 B=\Delta \Omega,
 \end{align}
 where $\Delta=\mathrm{diag}(\d_1, \dots, \d_r)$ 
with  positive integers $\d_1, \dots,  \d_r$
and $\Omega=(\omega_{ij})$ is an $r\times r$ skew-symmetric rational matrix.
Then, we set
\begin{align}
N=\bigoplus_{i=1}^r \bbZ e_i,
\quad
\{e_i,e_j\}=\omega_{ij},
\quad
N^{\circ}=\bigoplus_{i=1}^r \bbZ \d_i e_i.
\end{align}
This yields
 a  (not unique)  pair of a fixed data $\Gamma$ and a seed $\fraks$ such that $B$
is given by \eqref{3eq:Bmat1}.

The seed $\fraks$ here is identified with the initial seed $(\bfx,\bfy, B)$
for a cluster pattern or a cluster algebra
 in the sense of
 \cite{Fomin07}  with
\begin{align}
\label{3eq:cseed1}
x_i = x^{f_i},
\quad
y_i = y^{e_i},
\quad
b_{ij}=\{\d_i e_i, e_j\},
\end{align}
where $\bfx=(x_i)_{i=1}^r$,  $\bfy=(y_i)_{i=1}^r$, and $B=(b_{ij})_{i,j=1}^r$
 are \emph{cluster variables ($x$-variables)}, \index{$x$-variable (cluster variable)}
 \emph{coefficients ($y$-variables)}, \index{$y$-variable (coefficient)}
 and an \emph{exchange matrix\/} in \cite{Fomin07}, \index{matrix!exchange}
 and $x$ and $y$ are symbols for the formal exponential.
 See Part II for more details.

 \begin{rem}
In \eqref{3eq:bili2} we  took the transpose of the bilinear form $\{\cdot, \cdot\}$  in \cite{Gross14}
to match the convention of exchange matrices in
\cite{Fomin07} in the above correspondence.
As a consequence, every expression concerning the bilinear form 
is the transpose of the one in \cite{Gross14} throughout  Part III.
\end{rem}

Let us introduce a group homomorphism
\begin{align}
\label{3eq:p*3}
\begin{matrix}
p^* : &N &\rightarrow &M^{\circ}\\
& n & \mapsto & \{ \cdot, n\}.
\end{matrix}
\end{align}
By definition, we have
\begin{align}
\label{3eq:p*1}
\langle n', p^*(n)\rangle
=\{n',n\} \in \bbZ
\quad
(n'\in N^{\circ}, n\in N).
\end{align}

\begin{lem}
\label{3lem:p*1}
 Let $e_i$, $f_i$, $B=(b_{ij})$ be as above.
\par
(a). We have
\begin{align}
\label{3eq:pej1}
p^*(e_j) = \sum_{i=1}^r b_{ij} f_i.
\end{align}
Thus, the matrix representation of $p^*$ with respect to the above bases  is given by $B$.
\par
(b). $\Gamma$ is nondegenerate if and only if the matrix $B$ is nonsingular.
\par
(c). The map $p^*$ is injective if and only if $\Gamma$ is nondegenerate.
\end{lem}
\begin{proof}
(a). One can verify it as
\begin{align}
\label{3eq:ppr1}
\langle \d_i e_i , p^*(e_j)\rangle
=\{\d_i e_i ,e_j\}=b_{ij} =\biggl\langle \d_i e_i , \sum_{k=1}^r b_{kj} f_k\biggr\rangle.
\end{align}
The properties (b) and (c) follow from (a).
\end{proof}

Throughout the section, 
we fix a seed $\fraks=(e_1,\dots,e_r)$ of a given fixed data $\Gamma$
unless otherwise mentioned.

Let $N^+$ be the set of \emph{positive elements} in $N$ with respect to $\fraks$, namely, \index{positive element}
\begin{align}
\label{3eq:N+1}
N^+ = N_{\fraks}^+:=
\biggl\{ \sum_{i=1}^r a_i e_i \biggm | a_i \in \bbZ_{\geq 0},
\ \sum_{i=1}^r a_i > 0\biggr\}.
\end{align}
We say that $n\in N$ is \emph{primitive\/} if \index{primitive}
there is no pair $j \in \bbZ_{>1}$ and  $n'\in N$
such that $n=j n'$.
Let $N_{\mathrm{pr}}^+$ be the set of all primitive elements
in $N^+$.

We have a natural \emph{degree function\/} $\deg:N^+\rightarrow \bbZ_{>0}$ given by \index{degree function}
\begin{align}
\deg: \sum_{i=1}^r a_i e_i \mapsto
\sum_{i=1}^r a_i.
\end{align}
For each  integer $\l>0$, we set
\begin{align}
(N^+)^{\leq \l}&= \{ n\in N^+ \mid \deg(n) \leq \l\},
\quad
(N_{\mathrm{pr}}^+)^{\leq \l}=(N^+)^{\leq \l}\cap N_{\mathrm{pr}}^+,
\\
(N^+)^{> \l}&= \{ n\in N^+ \mid \deg(n) > \l\},
\quad
(N_{\mathrm{pr}}^+)^{> \l}=(N^+)^{> \l}\cap N_{\mathrm{pr}}^+,
\\
(N^+)^{\l}&= \{ n\in N^+ \mid \deg(n) = \l\},
\quad
(N_{\mathrm{pr}}^+)^{ \l}=(N^+)^{ \l}\cap N_{\mathrm{pr}}^+.
\end{align}
In particular, $(N^+)^{\leq \l}$ is a finite set.

\subsection{Structure group $G$}
\label{3subsec:structure1}

Let us fix a field  $\bbk$ of characteristic zero arbitrarily, 
e.g., $\bbk=\bbQ$.

For a given fixed data $\Gamma$ and a seed $\fraks$,
we introduce an $N^+$-graded Lie algebra $\frakg$ over $\bbk$
with generators $X_n$ ($n\in N^+$) as follows:

\begin{gather}
\label{3eq:lieg1}
\frakg=\frakg_{\fraks}:=\bigoplus_{n\in N^+}
\frakg_n,
\quad
\frakg_n=\bbk X_n,\\
\label{3eq:Xbra1}
[X_n, X_{n'}]=\{n, n'\} X_{n+n'}.
\end{gather}
The Jacobi identity is easily seen by
the following cyclic expression:
\begin{align}
[X_{n_1},[ X_{n_2}, X_{n_3}]]
=(\{ n_1,n_2\}\{ n_2,n_3\}
-\{ n_2,n_3\}\{ n_3,n_1\})
X_{n_1+n_2+n_3}.
\end{align}

For each integer $\l> 0$,
we define  a Lie algebra ideal and  a quotient of $\frakg$,
\begin{align}
\label{3eq:ggk1}
\frakg^{>\l}:=
\bigoplus_{n\in (N^+)^{>\l}}
\frakg_n,
\quad
\frakg^{\leq \l}:=\frakg/\frakg^{>\l}
\simeq
\bigoplus_{n\in (N^+)^{\leq \l}}
\frakg_n
\text{\ (as a vector space)}.
\end{align}
It is important that
$\frakg^{\leq \l}$ is a \emph{nilpotent\/} Lie algebra.
We have the canonical projections
\begin{align}
\label{3eq:kkl1}
\pi_{\l',\l}:\frakg^{\leq \l'}\rightarrow \frakg^{\leq \l} \quad (\l'>\l).
\end{align}
Thus, we have the completion of $\frakg$,
\begin{align}
\label{3eq:g1}
\widehat{\frakg}=\widehat{\frakg}_{\fraks}:=
\lim_{\longleftarrow}
\frakg^{\leq \l}
\simeq
\prod_{n\in N^+}
\frakg_n
\quad
\text{(as a vector space)}
\end{align}
with the canonical projections $\pi_{\l}:\widehat\frakg \rightarrow \frakg^{\leq \l}$ ($\ell>0$).

For any Lie subalgebra $\frakh$ of  $\frakg$,
 its completion $\widehat\frakh$ is defined by
\begin{align}
\label{3eq:hcomp1}
\frakh^{>\l} = \frakh \cap \frakg^{>\l},
\quad
\frakh^{\leq \l} =\frakh /\frakh^{>\l},
\quad
\widehat\frakh=
\lim_{\longleftarrow}
\frakh^{\leq \l}.
\end{align}

Next, we define the exponential group $G$
for the Lie algebra $\widehat\frakg$
 in \eqref{3eq:g1}.
In fact, there are two equivalent ways to define it.

\par
(a). For each nilpotent Lie algebra $\frakg^{\leq \l}$
in \eqref{3eq:ggk1},
we define the associated  exponential group
\begin{align}
G^{\leq \l}=\exp(\frakg^{\leq \l})
\end{align}
as follows:
As a set, $G^{\leq \l}$ is identified with $\frakg^{\leq \l}$
with $\exp(X)\leftrightarrow X$,
and its product is defined by the \emph{Baker-Campbell-Hausdorff (BCH) formula\/}
\index{Baker-Campbell-Hausdorff (BCH) formula}
(e.g., \cite[\S V.5]{Jacobson79}),
\begin{align}
\begin{split}
\label{3eq:BCH1}
&\ \exp(X) \exp(Y)
\\
=&\ \exp \biggl(X+Y+\frac{1}{2}[X,Y]
+\frac{1}{12}[X,[X,Y]]
-\frac{1}{12}[Y,[X,Y]]
+\cdots
\biggr)
,
\end{split}
\end{align}
where we do not need the explicit form of the higher order
commutators here.
Originally, the BCH formula is a product formula
for the exponential element $e^X = \sum_{k=0}^{\infty} X^k/k!$ in the situation where $e^X$ makes sense.
Here, the product is well-defined due to the nilpotency  of $\frakg^{\leq \l}$.
Moreover, it is compatible with the canonical projections
\begin{align}
\label{3eq:kk1}
\pi_{\l',\l}:G^{\leq \l'}\rightarrow G^{\leq \l} \quad (\l'>\l).
\end{align}
Thus, by the inverse limit we obtain a group
\begin{align}
\label{3eq:G1}
G=G_{\fraks}=G_{\Gamma,\fraks}:=\lim_{\longleftarrow} G^{\leq \l}
\end{align}
with the canonical projections $\pi_{\l}:G \rightarrow G^{\leq \l}$ ($\ell>0$).
Let $ G^{>\l}:=\Ker\, \pi_\l$.
When we discuss $G^{\leq \l}$,  we  often work with elements in $G$
modulo $G^{>\l}$.

(b).
Alternatively, we may directly define
\begin{align}
G:=\exp( \widehat\frakg),
\end{align}
where the infinite sum in the BCH formula \eqref{3eq:BCH1} makes sense
for $\widehat\frakg$ again due to the nilpotency  of $\frakg^{\leq \l}$.
An element $\exp(X)$ ($X\in \widehat\frakg$)  is identified
with an element $(\exp(\pi_\l(X)))_{\l =1}^{\infty}$ of $G$ in \eqref{3eq:G1}.

Here we call $G$ the \emph{structure group\/} for the forthcoming scattering diagrams.
\index{structure group}
We remark that the group $G$  depends only on $N$,  $\{\cdot, \cdot\}$,
and $e_1,\dots, e_r$ among the data from $\Gamma$ and $\fraks$.
In view of the decomposition \eqref{3eq:BDO1}, they corresponds to the matrix $\Omega$.
Thus, one can take a common $G$ for  those $(\Gamma, \fraks)$
where $N$,  $\{\cdot, \cdot\}$,
and $e_1,\dots, e_r$ (and $\Omega$) are fixed but
 $N^{\circ}$ and $\d_1,\dots,\d_r$ (and $\Delta$) are taken differently.

For any $g=\exp(X)\in G$ and any 
 $c\in \bbk$,
a  power $g^c$ of $g$ is defined by
\begin{align}
g^{c}:=\exp(cX).
\end{align}
In particular, for any $c\in \bbQ$, the fractional power $g^c$ of $g$ is unambiguously defined.

We may consider an \emph{infinite} product in $G$ as follows.
Suppose that there is  a countable and linearly ordered sequence $(g_i)_{i\in I}$ ($g_i\in G$).
Suppose that for any $\ell\in \bbZ_{>0}$, there are only finitely many $i\in I$
such that $\pi_{\ell}(g_i)\neq \rmid$.
Let $p_{\ell}$ be the product of those $g_i$ such that 
$\pi_{\ell}(g_i)\neq \rmid$ along its linear order, say, from left to right.
Then, $\lim_{\ell\rightarrow \infty}p_{\ell} $ yields an element in $G$,
which is the infinite product ${\prod}_{i\in I} g_i$.

Let  $\frakh$ be any Lie subalgebra of $\frakg$,
and let $\frakh^{\leq \l}$ and $\widehat\frakh$ be the ones
  in \eqref{3eq:hcomp1}.
Then, the corresponding subgroup $H$ of $G$ is defined by
\begin{align}
\label{3eq:H1}
H^{\leq \l}:=\exp(\frakh^{\leq \l}),
\quad
H:=
\lim_{\longleftarrow} H^{\leq \l}
=
\exp(\widehat\frakh).
\end{align}

The following subgroups of $G$ are especially important.
\begin{defn}[Parallel subgroup] \index{parallel subgroup}
For each $n\in N^+_{\rmpr}$, we define  a  Lie subalgebra of $\frakg$,
\begin{align}
\label{3eq:gpar1}
\frakg_{n}^{\parallel}:=
\bigoplus_{j>0}
\frakg_{jn}.
\end{align}
This is abelian by \eqref{3eq:Xbra1} and the skew-symmetry of 
\eqref{3eq:bili1}.
Then,
we have an abelian subgroup
$G_{n}^{\parallel}$ of $G$
corresponding to $\widehat\frakg_{n}^{\parallel}$
as defined in \eqref{3eq:H1}.
We call $G_n^{\parallel}$ the \emph{parallel subgroup\/} of $n$.
\end{defn}

From the BCH formula  \eqref{3eq:BCH1}, 
we obtain the \emph{Zassenhaus formula\/} \cite{Magnus54}, which is regarded as the dual of
the BCH formula: \index{Zassenhaus formula}
\begin{align}
\label{3eq:Zass1}
\begin{split}
& \exp(X+Y)
\\
=&
\exp(X)\exp(Y)
\exp\biggl(-\frac{1}{2}[X,Y]\biggr)
\exp\biggl( \frac{1}{3}[Y,[X,Y]]+\frac{1}{6}[X,[X,Y]]\biggr)
\cdots,
\end{split}
\end{align}
where again the explicit form of higher commutators are not necessary,
and the infinite product makes sense in $G$.
Together with the BCH formula, 
it also implies the useful formulas
\begin{align}
\label{3eq:Zass2}
\exp(X)\exp(Y)
&=\exp(Y)\exp(X)
\exp([X,Y]) C
\\
\label{3eq:Zass3}
&=\exp(Y)\exp([X,Y]) C' \exp(X)
.
\end{align}
where $C$ and $C'$ are  (possibly infinite) products of the exponentials of higher commutators
of $X$ and $Y$.

\subsection{Representations of $G$}
\label{3subsec:faithful1}

Let us  introduce two kinds of  representations of $G$,
which build a bridge between scattering diagrams and cluster algebras.

Let us recall some basic facts on derivations \cite[\S I.2]{Jacobson79}.
Let $A$ be a (not necessarily associative) algebra over $\bbk$.
A linear map $D: A \rightarrow A$ is a \emph{derivation} of $A$
if $D(ab)=D(a)b+aD(b)$ holds for any $a,b\in A$.
For any derivations $D_1$ and  $D_2$ of $A$,
$[D_1,D_2]:=D_1D_2 - D_2 D_1$ is also a derivation.
Thus, we have the Lie algebra $\mathrm{Der}(A)$ consisting of all derivations of  $A$
called the \emph{derivation (Lie) algebra}.
For a derivation $D$,
suppose that the infinite sum $\exp(D)=\sum_{k=0}^{\infty} D^k/k!$ is
well-defined as a  linear map on $A$. Then, $\exp(D)$ is an algebra automorphism 
of $A$.

(a). \emph{$y$-representation}. \index{$y$-representation}
We consider  a  monoid
 \begin{align}
\label{3eq:Q1}
Q=Q_{\fraks}:=&
 \biggl\{
 \sum_{i=1}^r a_i e_i
 \biggm |
 a_i\in \bbZ_{\geq 0}
  \biggr\}
  =
N^+ \sqcup \{0\}
\subset N.
 \end{align}
 Let $\bbk[Q]$ be the monoid algebra of $Q$  over $\bbk$.
Let $J$ be the maximal ideal of $\bbk[Q]$
generated by $N^+=Q \setminus \{0\}$,
and let 
\begin{align}
\label{3eq:comp1}
\bbk[[Q]]=\lim_{\longleftarrow} \bbk[Q]/J^{\l}
\end{align}
be the completion with respect to $J$.
We express any element $f\in \kQ$
as a formal power series in a symbol $y$ as 
\begin{align}
\label{3eq:fx1}
f=
\sum_{n\in Q} c_n y^n.
\end{align}
Let $\mathrm{Der}(\bbk[[Q]])$ be the derivation   algebra of $\kQ$.

For each $n \in N^+$, 
we define a derivation $\tilde{X}_n$ on $\kQ$
by
\begin{align}
\label{3eq:tildeX1}
 \tilde{X}_n (y^{n'}):=\{ n,n'\} y^{n'+n}
 \quad
 (n'\in Q).
\end{align}
\begin{prop}[{cf.\ \cite[\S 1.1]{Gross14}}]
\label{3prop:gaction1}
The map 
\begin{align}
\label{3eq:Xn1}
\begin{matrix}
\rho_{y}:&  \widehat{\frakg} &\rightarrow &\mathrm{Der}(\bbk[[Q]]) \\
&
\displaystyle
X=  \sum_{n\in N^+} c_n X_n & 
  \mapsto & 
  \displaystyle
 \tilde X= \sum_{n\in N^+} c_n  \tilde{X}_n,
\end{matrix}
\end{align}
yields a Lie algebra homomorphism,
where the sum could be infinite.
Moreover, if $\Gamma$ is nondegenerate,
it is injective.
\end{prop}
\begin{proof}
First, we regard $\rho_{y}$ in \eqref{3eq:Xn1} with finite sum as a map $\frakg \rightarrow\mathrm{Der}(\bbk[Q])$.
Then, we have
\begin{align}
\begin{split}
\tilde{X}_{n_1} \tilde{X}_{n_2} (y^{n})
&=\tilde{X}_{n_1}(\{n_2,n\} y^{n+n_2})
\\
&=(\{n_1, n\} + \{n_1, n_2\}) \{n_2,n\} y^{n+n_2+n_1}.
\end{split}
\end{align}
Thus, we obtain
\begin{align}
\begin{split}
[\tilde{X}_{n_1} , \tilde{X}_{n_2}](y^{n})
&= \{n_1, n_2\}  \{n_1+n_2, n\} y^{n+n_1+n_2}\\
&= \{n_1, n_2\}  \tilde{X}_{n_1+n_2}(y^{n}).
\end{split}
\end{align}
Therefore, $\rho_{y}:\frakg \rightarrow\mathrm{Der}(\bbk[Q])$ is a Lie algebra homomorphism.
We show that it is extended to the one in \eqref{3eq:Xn1}.
First, we note that
$\frakg(J^{\l+1})\subset J^{\l+2}$.
Thus, 
 $\frakg$ acts on $\bbk[Q]^{\leq \l}:=\bbk[Q]/J^{\l+1}$.
 Next, 
we  have
$\frakg^{>\l}(\bbk[Q])\subset J^{\l+1}$.
Thus, 
it is extended to the action of
 $\widehat{\frakg}$.
 Moreover, this action is compatible with
projections
$\pi_{\l',\l}: \bbk[Q]^{\leq \l'}\rightarrow \bbk[Q]^{\leq \l}$
for $\l' >\l$.
Therefore, 
 $\widehat\frakg$ acts on $\kQ$.
 If $\{\cdot, \cdot\}$ is non-degenerate, 
 $ \tilde{X}_n$ ($n\in N^+$) in \eqref{3eq:tildeX1}   are linearly independent.
 Thus,  $\rho_{y}$ is injective.
\end{proof}

\begin{rem}
The above representation  $\rho_{y}$
is equivalent to the direct sum of the adjoint  and the trivial representations of $\widehat{\frakg}$.
\end{rem}

For any $X \in\widehat\frakg$,
let $\tilde{X}:={\rho}_{y}(X)$.
We define 
\begin{align}
\label{3eq:exp1}
\mathrm{Exp}(\tilde{X}):=\sum_{k=0}^{\infty} \frac{1}{k!}\tilde{X}^k
\in \mathrm{GL}(\bbk[[Q]]),
\end{align}
which is  well-defined by the same reason in the proof of Proposition  \ref{3prop:gaction1}.
Since $\tilde{X}$ is a derivation, $\mathrm{Exp}(\tilde{X})$
is an algebra automorphism of $\bbk[[Q]]$
as we mentioned.
Then, we have the  following group homomorphism,
where we abuse the symbol $\rho_y$, for simplicity.
\begin{prop}
\label{3prop:gaction2}
We have a  group homomorphism
\begin{align}
\label{3eq:Xn2}
\begin{matrix}
\rho_{y}:&  G &\rightarrow & \mathrm{Aut}(\bbk[[Q]]) \\
& \exp(X) &  \mapsto &\mathrm{Exp} (\tilde{X}).
\end{matrix}
\end{align}
Moreover, if $\Gamma$ is nondegenerate,
it is injective.
\end{prop}
\begin{proof}
The map $\rho_{y}$ is a group homomorphism,
because the exponential \eqref{3eq:exp1} also satisfies the 
same
BCH formula as mentioned.
Suppose that $\Gamma$ is nondegenerate.
Let us show that $\rho_{y}$ is injective.
Assume that $\mathrm{Exp}(\tilde{X})=\rmid$.
Then, 
by the action \eqref{3eq:tildeX1},
we have $\tilde{X}=0$.
Then, $X=0$ in $\widehat\frakg$ by the injectivity of
$\rho_{y}$ in Proposition  \ref{3prop:gaction1}.
Thus, $\exp(X)=\rmid$.
\end{proof}

The representation $\rho_{y}$ is
 closely related with the mutations of $y$-variables 
in 
\eqref{3eq:cseed1}
as we will see soon.
Thus, we call it the \emph{$y$-represent\-ation\/} of $G$.

\smallskip
(b). \emph{$x$-representation}. \index{$x$-representation}
Here we assume that $\Gamma$ is nondegenerate
from the beginning.
Then, thanks to Lemma \ref{3lem:p*1}, the map $p^*$  in \eqref{3eq:p*3}
is injective.
In particular, $p^*(e_1)$, \dots, $p^*(e_r)$ are $\bbZ$-linearly independent.
We define
a monoid
\begin{align}
 \label{3eq:P1}   
 P=P_{\fraks}:=&
 \biggl\{
 \sum_{i=1}^r a_i p^*(e_i)
 \biggm |
 a_i\in \bbZ_{\geq 0}
  \biggr\}
  =
p^*(Q)
\subset M^{\circ}.
\end{align}
Let $\bbk[P]$ the monoid algebra of $P$ over $\bbk$,
and $\bbk[[P]]$ be the completion by the maximal ideal of 
$\bbk[P]$ generated by $P \setminus \{0\}$.
(Note that, if $\Gamma$ is degenerate, we cannot define $\bbk[[P]]$.)
We express any element $f\in \kP$
as a formal power series in a symbol $x$ as 
\begin{align}
\label{3eq:fx3}
f=
\sum_{m\in P} c_m x^m.
\end{align}

Then,
in parallel to \eqref{3eq:tildeX1},
for each $n\in N^+$,
we define a derivation  $ \tilde{X}_n$ on $\bbk[[P]]$ by
\begin{align}
\label{3eq:tildeX2}
 \tilde{X}_n (x^{m}):=\langle n,m\rangle x^{m+p^*(n)}
 \quad (m\in P).
\end{align}
Then, repeating the same argument for $\rho_{y}$, we obtain
the following group homomorphism.
\begin{prop}
\label{3prop:gaction3}
Suppose that $\Gamma$ is nondegenerate.
Then,
we have an injective group homomorphism
\begin{align}
\label{3eq:Xn4}
\begin{matrix}
\rho_{x}:&  G &\rightarrow & \mathrm{Aut}(\bbk[[P]]) \\
& \exp(X) &  \mapsto &\mathrm{Exp} (\tilde{X}).
\end{matrix}
\end{align}
\end{prop}

The representation $\rho_{x}$ is
 closely related with the mutations of $x$-variables 
in 
\eqref{3eq:cseed1}
as we will see soon.
Thus, we call it the \emph{$x$-representation\/} of $G$.

Even though  a faithful representation of $G$
is desirable in some situation,
it is too restrictive to assume that $\Gamma$ is nondegenerate
in view of cluster algebra theory.
We will treat this problem  later  in Section
\ref{3sec:principal1}.

\subsection{Dilogarithm elements and pentagon relation}
\label{3subsec:pentagon1}

As an application of the actions of $G$ in the previous subsection,
we show  a remarkable relation  among certain distinguished elements in $G$.

Let $\mathrm{Li}_2(x)$ be the \emph{Euler dilogarithm\/} \cite{Lewin81} defined by
\index{Euler dilogarithm} \index{dilogarithm|see{Euler dilogarithm}}
\begin{align}
\mathrm{Li}_2(x):=
\sum_{j=1}^{\infty} \frac{x^j}{j^2}
=-\int_0^x \frac{\log (1-y) }{y} \, dy.
\end{align}
It has the following key property:
\begin{align}
\label{3eq:dilogc2}
x \frac{d}{dx}
\left(
-\mathrm{Li}_2(-x)
\right)
&=
\sum_{j=1}^{\infty} \frac{(-1)^{j+1}}{j}x^j
=
\log (1+x).
\end{align}

Having this in mind, we introduce the following elements of $G$.

\begin{defn}[Dilogarithm element] \index{dilogarithm element}
\label{3defn:diloge1}
For any $n\in N^+$,
let
\begin{align}
\label{3eq:gei1}
\Psi[n]:=\exp
\Biggl(\,
\sum_{j=1}^{\infty} \frac{(-1)^{j+1}}{j^2} X_{j n}
\Biggr)
\in 
G_{n_0}^{\parallel},
\end{align}
where $n_0\in N^+_{\rmpr}$ is the unique one such that $n=t n_0$ for some integer $t>0$.
We call $\Psi[n]$ the \emph{dilogarithm element\/} for $n$.
\end{defn}

\begin{rem}
It is ``only formally''  expressed as
\begin{align}
\Psi[n]=\exp(-\mathrm{Li}_2(-X_n)),
\end{align}
if we interpret $X_n^j$ as $X_{j n}$.
(Caution: $X_n^j$ and $X_{jn}$ act differently on $\kQ$ and $\kP$.)
\end{rem}

The elements $\Psi[n]$ are introduced because of the following property.

\begin{prop}[{Cf.\ \cite[Lemma 1.3]{Gross14}}]
\label{3prop:gn1}
(a).
Under the $y$-represent\-ation $\rho_{y}$, $\Psi[n]$ acts on $\bbk[[Q]]$ as
\begin{align}
\label{3eq:gei2y}
\Psi[n](y^{n'})=y^{n'} (1+y^{n})^{\{ n, n'\}}
\quad (n'\in Q).
\end{align}

\par
(b). Suppose that $\Gamma$ is nondegenerate.
Then,  under the $x$-representation $\rho_{x}$, $\Psi[n]$ acts on $\bbk[[P]]$ as
\begin{align}
\label{3eq:gei2}
\Psi[n](x^m)=x^m (1+x^{p^*(n)})^{\langle n, m\rangle}
\quad (m\in P).
\end{align}
\end{prop}
\begin{proof}
Since the calculations are parallel, we concentrate on the case (b).
\begin{align}
\label{3eq:dilogc1}
\begin{split}
\Psi[n](x^m)
&=
x^m
\exp
\Biggl(
\sum_{j>0}  \frac{(-1)^{j+1}}{j^2} x^{jp^*(n)}\langle jn, m \rangle
\Biggr)
\\
&=
x^m
\exp
\Biggl(\,
 \sum_{j>0}   \frac{(-1)^{j+1}}{j}  x^{jp^*(n)}
\Biggr)
^{\langle n, m \rangle}
\\
&=x^m (1+x^{p^*(n)})^{\langle n, m\rangle}.
\end{split}
\end{align}
\end{proof}

Observe that the calculation in  \eqref{3eq:dilogc1} 
is parallel to  the derivation of
the formula 
\eqref{3eq:dilogc2} by the power series calculation.
The above automorphisms \eqref{3eq:gei2y}
and \eqref{3eq:gei2} are identified with  the automorphism  part of the Fock-Goncharov decomposition of mutations 
of \emph{$y$-variables  (coefficients)\/} and \emph{$x$-variables (cluster variables)\/} in cluster algebras, respectively, in \cite[\S2.1]{Fock03}.

The dilogarithm elements are generators of $G$
in the following sense.
\begin{prop} 
\label{3prop:ggen1}
The group $G$ is generated by $\Psi[tn]^c$ $(n\in N_{\rmpr}^+, t\in \bbZ_{>0}, c\in \bbk)$
admitting the infinite product.
\end{prop}
\begin{proof}
By inverting the expression \eqref{3eq:gei1},
one can express $\exp(X_n)$ ($n\in N^+$) as an infinite product of 
 $\Psi[tn]^{c}$ ($t\in \bbZ_{>0}$, $c\in \bbQ$).
\end{proof}

The dilogarithm elements satisfy the following  remarkable relations.

\begin{prop}
\label{3prop:pent1}
Let $n,n'\in N^+$.
The following relations hold in $G$.
\par
(a). If $\{n',n\}=0$,
for any $c,c'\in \bbk$, 
\begin{align}
\label{3eq:com1}
\Psi[n' ]^{c'} \Psi[ n]^c =\Psi[  n ]^c  \Psi[ n']^{c'}.
\end{align}

(b). (Pentagon relation \cite[Example 1.14]{Gross14} for $c=1$.) \index{pentagon relation}
If $\{n',n\}=c$ $(c\in \bbQ, c\neq 0)$,
\begin{align}
\label{3eq:pent1}
\Psi[n' ]^{1/c} \Psi[ n]^{1/c}=
\Psi[  n ]^{1/c} \Psi[n+n']^{1/c} \Psi[ n']^{1/c}.
\end{align}
\end{prop}
\begin{proof}
(a).
This is  clear from \eqref{3eq:Xbra1}.
(b).
Let us consider the rank 2 sublattice $N'$ of $N$ generated by $n$ and $n'$.
By the assumption $\{n', n \}\neq 0$, the form $\{\cdot, \cdot \}$  restricted on $N'$ is nondegenerate.
Let $G'$ be the subgroup of $G$ corresponding to $N'$.
Similarly, we consider a monoid
\begin{align}
 \label{3eq:P3}   
 P'=&
\{
 a p^*(n) + a'p^*(n') 
 \mid
 a, a'\in \bbZ_{\geq 0}
  \}
\subset M^{\circ}.
\end{align}
Since the relation \eqref{3eq:pent1} involves only elements in $G'$,
one can prove it
by the $x$-representation of $G'$ on $\bbk[[P']]$, which is faithful.
The left hand side is given by
\begin{align}
\begin{split}
&\quad\ \Psi[n']^{1/c} \Psi[ n]^{1/c} (x^m)\\
&= \Psi[n']^{1/c}(x^m (1+x^{p^*(n)})^{\langle n, m\rangle/c})
\\
&= x^m  (1+x^{p^*(n')})^{ \langle  n', m\rangle/c}
 (1+x^{p^*(n)}  (1+x^{p^*(n')})
)^{\langle n, m\rangle/c}.
\end{split}
\end{align}
The right hand side is given by
\begin{align}
\begin{split}
&\quad\ \Psi[ n ]^{1/c} \Psi[ n+ n']^{1/c} \Psi[n']^{1/c}(x^m)\\
&=  \Psi[n ]^{1/c} \Psi[ n+n']^{1/c} (x^m (1+x^{p^*(n')})^{\langle n', m\rangle/c})
\\
&=  \Psi[ n ]^{1/c}( x^m  (1+x^{p^*(n+n')})^{\langle  n+n', m\rangle/c}\\
&\qquad \times
 (1+x^{p^*(n')}  (1+x^{p^*(n+n')})^{-1}
)^{\langle   n', m\rangle/c})
\\
&=  \Psi[ n ]^{1/c} (x^m  (1+x^{p^*(n+n')})^{\langle  n, m\rangle/c}
 (1+x^{p^*(n+n')} +x^{p^*(n')})^{\langle   n', m\rangle/c})
\\
&=  x^m (1+x^{p^*(n)})^{\langle n, m\rangle/c}
  (1+x^{p^*(n+n')} (1+x^{p^*(n)})^{-1})^{\langle  n, m\rangle/c}\\
&\qquad \times
 (1+x^{p^*(n+n')}(1+x^{p^*(n)})^{-1}+x^{p^*(n')}
 (1+x^{p^*(n)})^{-1}
)^{\langle   n', m\rangle/c}
\\
&=  x^m 
  ( 1+x^{p^*(n)}+x^{p^*(n+n')})^{\langle  n, m\rangle/c}
 (1+x^{p^*(n')})^{\langle   n', m\rangle/c}.
\end{split}
\end{align}
Thus, two expressions agree.
\end{proof}

\begin{rem}
(a). 
The above calculation for $c=1$  is essentially  the same for the \emph{pentagon periodicity\/}
of $x$-variables for a  cluster algebra of type $A_2$ in \cite{Fomin02}.
Alternatively,
one can also work with the $y$-representation.
Then, the calculation for $c=1$ is  essentially
  the same for the pentagon periodicity
of $y$-variables for a cluster algebra  of type $A_2$ in \cite{Fomin02}.
\par
(b). The case $c\neq 1$ follows from the case $c=1$ by Proposition
\ref{3prop:Giso1}.
\end{rem}

\subsection{Equivalence of fixed data and normalization factor}
\label{3subsec:equivalence1}

There is some redundancy for choosing a fixed data.
Let us start with a simple situation.

\begin{defn}[Rescaling of fixed data] \index{rescaling (of a fixed data)}
\label{3defn:rescale1}
Let $\Gamma$ and $\Gamma'$
be a pair of fixed data.
We say that $\Gamma'$ is a \emph{rescaling of  $\Gamma$} 
if  there is some  rational number $\lambda>0$ such that
the following relations hold:
\begin{align}
N'&= N,\\
\label{3eq:scale1}
\{n_1, n_2 \}'&=\lambda \{n_1, n_2 \},\\
\label{3eq:scale2}
(N^{\circ})'&= \lambda^{-1} N^{\circ},\\
\label{3eq:scale3}
\d'_i &= \lambda^{-1} \d_i,\\
\label{3eq:scale4}
(M^{\circ})'&= \lambda M^{\circ}.
\end{align}
\end{defn}

The exchange matrix $B$ associated with $(\Gamma, \fraks)$
in  \eqref{3eq:Bmat1}
is  invariant under the above rescaling, namely,
\begin{align}
b'_{ij}
=\{\d'_i e_i, e_j\}'
=\lambda \{ \lambda^{-1}\d_i e_i, e_j \}=
b_{ij}.
\end{align}

Let $\fraks=(e_1,\dots, e_r)$  be a seed for $\Gamma$.
Then, $\fraks$ is also regarded as a seed for $\Gamma'$.
Let $ \widehat{\frakg}'$ and $G'$ be the corresponding Lie algebra and group
for $\Gamma'$ and $\fraks$.

\begin{prop}
\label{3prop:Giso1}
We have a Lie algebra isomorphism
\begin{align}
\label{3eq:XX1}
\begin{matrix}
\phi_{\lambda}: & \widehat{\frakg}'  &\rightarrow& \widehat{\frakg}\\
& X'_n & \mapsto & \lambda X_n
\end{matrix}
\end{align}
and a group isomorphism
\begin{align}
\label{3eq:Giso1}
\begin{matrix}
\phi_{\lambda}: &G'  &\rightarrow&G\\
& \exp(X'_n) & \mapsto & \exp( \lambda X_n).
\end{matrix}
\end{align}
\end{prop}
\begin{proof}
This is 
due to the relation \eqref{3eq:scale1} and the commutation relation \eqref{3eq:Xbra1}.
\end{proof}

\begin{defn}[Normalization factor] \index{normalization factor}
\label{3defn:normal1}
For any $n\in N^+$,
let   $\d(n)$ be the smallest positive rational number
such that $\d(n)  n\in N^{\circ}$. 
We call $\d(n)$ the \emph{normalization factor\/} of $n$ with respect to $N^{\circ}$.
\end{defn}

\begin{ex}
 (a). For the seed $\fraks=(e_1,\dots, e_r)$  under consideration,
we have $\d(e_i)=\d_i$.
\par
(b). For any $n\in N^+_{\rmpr}$ and any  integer $t>0$,
we have 
\begin{align}
\label{3eq:d01}
\d(t n)=\d(n)/t,
\end{align}
because $\d(tn)tn=\d(n)n$.
Also, $\d(n)$ is an integer.
Otherwise, let $\d(n)=r/s$ be the irreducible rational expression with $s\neq 1$.
Then, we have $n, (r/s)n \in N$.
Thus, $(1/s)n\in N$.
This contradicts
that
$n$ is  primitive.
\end{ex}

\begin{prop}
\label{3prop:inv1}
Let $G$ and $G'$ be as above.
For any $n\in N^+_{\rmpr}$,
let $g\in G_{n}^{\parallel}$ and $g'\in (G')_{n}^{\parallel}$ 
with
\begin{align}
 \label{3eq:gc2}
g=\exp\Biggl( \d(n) \sum_{j>0} c_j X_{jn}
\Biggr),
\quad
g'=\exp\Biggl( \d'(n) \sum_{j>0} c'_j X'_{jn}
\Biggr),
\end{align}
where $\d(n)$ and $\d'(n)$ are the normalization factors
with respect to $N^{\circ}$ and $(N^{\circ})'$.
Then,
\begin{align}
\label{3eq:inv2}
\phi_{\lambda}
(g')
=g.
\end{align}
if and only $c_j=c'_j$ for any $j>0$.
\end{prop}
\begin{proof}
This follows from \eqref{3eq:XX1}
and the transformation of the normalization factors
$\d' (n) = \lambda^{-1} \d(n)$ by \eqref{3eq:scale2}.
\end{proof}

In view of Proposition \ref{3prop:inv1},
it is   natural to express an element $g \in G_n^{\parallel}$ as
\begin{align}
\label{3eq:inv3}
g=g_{0}^{\d(n)},
\end{align}
where $g_0=g^{1/\d(n)}$.
We call \eqref{3eq:inv3} the \emph{normalized form\/} of $g$. \index{normalized form}
Then, the coefficients $c_j$ of the expansion
$ g_{0}=\exp(\, \sum_{j>0} c_j X_{jn} )
$
are \emph{invariant\/} under the rescaling of $\Gamma$.

\begin{ex}
\label{3ex:res1}
For $\Psi[n]$ in \eqref{3eq:gei1} with $n=tn_0$ therein, we have
\begin{align}
\phi_{\lambda}
(\Psi'[n]^{\d'(n_0)})
=\Psi[n]^{\d(n_0)}.
\end{align}
Thus, we also have
\begin{align}
\label{3eq:inv1}
\phi_{\lambda}
(\Psi'[n]^{\d'(n)})
=\Psi[n]^{\d(n)}.
\end{align}
These elements play the main role in the forthcoming CSDs.
\end{ex}

Now we consider a more  general situation.

\begin{defn}
Let $(\Gamma, \fraks)$ be a pair of a fixed data and a seed,
and let $B_{\Gamma, \fraks}$
be the associated matrix in \eqref{3eq:Bmat1}.
We say that $(\Gamma, \fraks)$ and $(\Gamma', \fraks')$
are \emph{$B$-equivalent\/}
if $B_{\Gamma, \fraks}=B_{\Gamma', \fraks'}$. \index{$B$-equivalent}
\end{defn}

In Section \ref{3subsec:fixed1}
we have seen that, for a given skew-symmetrizable matrix $B$,
there is a (not unique) pair $(\Gamma,\fraks)$
such that $B=B_{\Gamma, \fraks}$.
Let us show that the associated group $G_{\Gamma, \fraks}$
is uniquely determined from $B$ up to isomorphism.
This fact is natural in view of
the correspondence to cluster algebras.

\begin{prop}
\label{3prop:equiv1}
If $(\Gamma, \fraks)$ and $(\Gamma', \fraks')$ are $B$-equivalent,
then the associated groups $G_{\Gamma, \fraks}$ and $G_{\Gamma', \fraks'}$
are isomorphic.
\end{prop}
\begin{proof}
Let $B=B_{\Gamma, \fraks}=B_{\Gamma', \fraks'}$.
We say that an $n\times n$ matrix $A$ is \emph{indecomposable\/}
if there is no pair $I, J\subset \{1,\dots, r\}$,  $I\sqcup J=\{1,\dots, r\}$,
$I, J \neq \emptyset$ such that $a_{ij}=0$ for any $i\in I$ and $j\in J$.
Suppose that $B$ is decomposed into 
$B_1\oplus \cdots  \oplus B_p$
up to simultaneous permutation of rows and columns,
where $B_i$ are indecomposable submatrices of $B$.
Accordingly, the groups $G=G_{\Gamma, \fraks}$
and $G'=G_{\Gamma', \fraks'}$ are factorized into the corresponding
subgroups $G=\prod_{i=1}^p G_i$ and $G'=\prod_{i=1}^p G'_i$.
Thus, we only need to prove the claim when $B$ is indecomposable.
Then, it is well known and easy to show that a (left) skew-symmetrizer of $B$ is 
unique up to a positive rational number.
This implies that $\d_i$ and $\d'_i$ are related by
$\d'_i = \lambda^{-1}\d_i$  with a common multiple $\lambda^{-1}$.
Let us identify $e_i$ with $e'_i$ ($i=1,\dots, r$) and $N$ with $N'$.
Now we are  in the situation of Definition \ref{3defn:rescale1}.
Therefore, $G$ and $G'$ are isomorphic 
by Proposition \ref{3prop:Giso1}.
\end{proof}

The above proof  also tells that
the parallel results to
Proposition \ref{3prop:inv1}  and Example \ref{3ex:res1}  hold in the situation 
of Proposition \ref{3prop:equiv1}.

\subsection{Formulas for dilogarithm elements}
Let us  present some useful formulas for dilogarithm elements.
\begin{prop}
\label{3prop:bound1}
Let $n, n'\in N^+$.
\par(a).
({\cite[Lemma C.8]{Gross14}}.)
If $\{\d(n') n', n\}=\pm \alpha$ $(\alpha\in \bbZ_{>0})$, then
the following equality holds for
the $x$-representation
on $\bbk[[P']]$ for $P'$ in \eqref{3eq:P3}:
\begin{align}
\label{3eq:gxg1}
\Psi[n']^{\pm \d(n')} &X_{n} \Psi[n']^{\mp \d(n')}
=
\sum_{j=0}^{\alpha}
\begin{binom}
\alpha j
\end{binom}
X_{n+ j n'}.
\end{align}
\par
(b). If $\{\d(n') n', n\}= \alpha$ and $\{\d(n) n, n'\}=- \beta$ $(\alpha, \beta\in \bbZ_{>0})$,
then
\begin{align}
\label{3eq:gxg2}
\Psi[n']^{ \d(n')} \Psi[n]^{ \d(n)}
&=
\Psi[n]^{ \d(n)}
\Psi[\beta n+ n']^{ \d(n')}
C
\Psi[n+\alpha n']^{ \d(n)} \Psi[n']^{ \d(n')},
\end{align}
where $C$ is a (possibly infinite) product of $\Psi[n'']^{c''}$ $(c''\in \bbQ)$
such that  $n''=j_1 n + j_2 n'$ $(j_1, j_2 \in \bbZ_{>0})$ with
$\alpha^{-1}< j_1/j_2<\beta$.
\end{prop}
\begin{proof}
(a).
We have
\begin{align}
\begin{split}
\label{3eq:gxg5}
& \quad \ X_{n} \Psi[n']^{\mp \d(n')} (x^m)
\\
&=
X_{n}( x^{m}(1+x^{p^*(n')})^{\mp\langle  \d(n') n',m\rangle})
\\
&=
\langle n, m\rangle x^{m+p^*(n)}(1+x^{p^*(n')})^{\mp\langle  \d(n') n',m\rangle}
\\
&\quad
+\alpha \langle  n',m\rangle  x^{m+p^*(n'+n)}(1+x^{p^*(n')})^{\mp\langle  \d(n') n',m\rangle-1},
\end{split}
\end{align}
\begin{align}
\begin{split}
\label{3eq:gxg4}
&\quad\
\Psi[n']^{\mp \d(n')} 
\sum_{j=0}^{\alpha}
\begin{binom}
\alpha j
\end{binom}
X_{n+ j n'}
(x^m)
\\
&=
\Psi[n']^{\mp \d(n')} 
\sum_{j=0}^{\alpha}
\begin{binom}
\alpha j
\end{binom}
\langle n+jn' ,m\rangle x^{m+p^*(n+jn')}
\\
&=
\sum_{j=0}^{\alpha}
\begin{binom}
\alpha j
\end{binom}
\langle n+jn' ,m\rangle x^{m+p^*(n+jn')}(1+x^{p^*(n')})^{\mp\langle \d(n')n',m\rangle-\alpha}
\\&=
\sum_{j=0}^{\alpha}
\begin{binom}
\alpha j
\end{binom}
\langle n ,m\rangle x^{m+p^*(n+jn')}(1+x^{p^*(n')})^{\mp\langle \d(n')n',m\rangle-\alpha}
\\
&\quad +
\sum_{j=1}^{\alpha}
\alpha
\begin{binom}
{\alpha-1 }{ j-1}
\end{binom}
\langle n' ,m\rangle x^{m+p^*(n+jn')}(1+x^{p^*(n')})^{\mp\langle \d(n')n',m\rangle-\alpha}
\\
&=
\langle n, m\rangle x^{m+p^*(n)}(1+x^{p^*(n')})^{\mp\langle  \d(n') n',m\rangle}
\\
&\quad
+\alpha \langle  n',m\rangle  x^{m+p^*(n'+n)}(1+x^{p^*(n')})^{\mp\langle  \d(n') n',m\rangle-1}.
\end{split}
\end{align}
Thus, they coincide.
\par
(b). 
As in the proof of Proposition \ref{3prop:pent1}, it is enough to prove it for
 the $x$-representation
on $\bbk[[P']]$ for $P'$ in \eqref{3eq:P3}.
For any $t\in \bbZ_{>0}$, we obtain from \eqref{3eq:gxg1} that
\begin{align}
\label{3eq:gxg6}
\Psi[n']^{ \d(n')} &X_{tn} \Psi[n']^{- \d(n')}
=
X_{tn}
+
\sum_{j=1}^{t \alpha-1}
\begin{binom}
{t \alpha} j
\end{binom}
X_{tn+ j n'}
+
X_{t(n+\alpha n')}.
\end{align}
Taking the sum and exponential in accordance with \eqref{3eq:gei1},
then using the Zassenhaus formula \eqref{3eq:Zass1} and Proposition \ref{3prop:ggen1},
we obtain
\begin{align}
\label{3eq:gxg7}
\Psi[n']^{ \d(n')} &\Psi[n]^{\d(n)} \Psi[n']^{- \d(n')}
=
\Psi[n]^{\d(n)}
C_1
\Psi[n+\alpha n']^{\d(n)},
\end{align}
where $C_1$ is a (possibly infinite) product of $\Psi[n'']^{c''}$ $(c''\in \bbQ)$
such that  $n''=j_1 n + j_2 n'$ $(j_1, j_2\in \bbZ_{>0})$ with
$0< j_2/j_1<\alpha$.
Similarly, by exchange the role of $n$ and $n'$,
we obtain
\begin{align}
\label{3eq:gxg8}
\Psi[n]^{- \d(n)} &\Psi[n']^{\d(n')} \Psi[n]^{\d(n)}
=
\Psi[\beta n+n']^{\d(n')}
C_2
\Psi[n']^{\d(n')},
\end{align}
where $C_2$ is a (possibly infinite) product of $\Psi[n'']^{c''}$ $(c''\in \bbQ)$
such that  $n''=j_1 n + j_2 n'$ $(j_1, j_2\in \bbZ_{>0})$ with
$0< j_1/j_2<\beta$.
Combining \eqref{3eq:gxg7} and \eqref{3eq:gxg8},
 we obtain \eqref{3eq:gxg2}.
\end{proof}

\notes

The contents are mostly taken from \cite[\S1.1]{Gross14}
and partly from \cite[\S2.1]{Kontsevich13}.
The $y$-representation of $G$
in Section \ref{3subsec:pentagon1} appeared in \cite[\S2]{Gross13}.
The pentagon relation  \eqref{3eq:pent1} generalizes the
one in \cite[Example 1.14]{Gross14} for $c=1$.
The equivalence of fixed data in Section \ref{3subsec:equivalence1}
is distilled from the \emph{change of lattice trick\/} in \cite[Appendix C.3]{Gross14}.
Proposition \ref{3prop:bound1} (b) is due to us.

\newpage
\section{Consistent scattering diagrams}

In this section we introduce and study consistent scattering diagrams.
In particular, we show that there is  a natural one-to-one correspondence  between
consistent scattering diagrams and elements in the structure group $G$.
\subsection{Scattering diagrams}

We continue to fix a seed $\fraks$ for a given (possibly  degenerate) fixed data $\Gamma$.

For any $n\in N$, $n\neq 0$, the  rational hyperplane
$n^{\perp}$ in $M_{\bbR}$ that is orthogonal to $n$ is defined by
\begin{align}
n^{\perp} :=\{ z \in M_{\bbR} \mid \langle n, z\rangle=0\}.
\end{align}

\begin{defn}[Cone]
For $m_1,\dots,m_s\in M^{\circ}$, let
\begin{align}
\sigma=\sigma(m_1,\dots,m_s):=\bbR_{\geq 0} m_1 + \cdots  +\bbR_{\geq 0} m_s
\subset M_{\bbR}.
\end{align}
We also set $\sigma(\emptyset)=\{0 \}$.
We call $\sigma$ a \emph{convex rational polyhedral cone}, or simply, a \emph{cone\/} \index{cone} \index{cone!convex rational polyhedral}
in $M_{\bbR}$.
We say that a cone $\sigma$ is \emph{strongly convex\/} if $\sigma\cap (-\sigma)=\{0\}$. \index{cone!strongly convex}
\end{defn}

\begin{defn}[Wall] \index{wall}
\label{3defn:wall1}
A \emph{wall\/}  for a seed $\fraks$
 is a triplet $\bfw=(\frakd, g)_{n}$,
where
\begin{itemize}
\item
$n \in N_{\rmpr}^+$,
\item
$\frakd\subset n^{\perp} $ is a
(not necessarily strongly convex) cone
of codimension 1 in $M_{\bbR}$,
\item
$g\in G_{n}^{\parallel}$.  
\end{itemize}
We call  $n$, $\frakd$, $g$
the \emph{normal vector}, \index{normal vector (of a wall)}
the \emph{support}, \index{support!of a wall}
the \emph{wall element\/} \index{wall element}
of a wall $\bfw$, respectively.
\end{defn}

Informally,  the support $\frakd$  is also  referred to as a \emph{wall\/};
for example, we say ``A curve  in $M_{\bbR}$ crosses a wall.''
The normal vector $n$ is uniquely specified from the data $\frakd$.
So, we may omit it as $\bfw=(\frakd,g)$ when it is not necessary to indicate it explicitly.

\begin{defn}[Scattering diagram] \index{scattering diagram}
A \emph{scattering diagram\/} $\frakD$ for  a seed $\fraks$
is a  collection of walls
$\{\bfw_{\lambda}=(\frakd_{\lambda}, g_{{\lambda}})_{n_{\lambda}}\}_{\lambda\in \Lambda}$
for  $\fraks$,
where $\Lambda$ is a
finite or countably infinite index set, satisfying the following 
\emph{finiteness condition\/}: \index{finiteness condition}
\begin{itemize}
\item
For each integer $\l>0$, there are only finitely many walls
such that $\pi_\l(g_{{\lambda}})\neq \rmid$,
where $\pi_{\l}: G \rightarrow G^{\leq \l}$ is the canonical projection.
\end{itemize}
For any scattering diagram $\frakD$ and any integer $\l>0$,
let  $\frakD_\l$ be the collection of walls 
obtained from  $\frakD$
by removing all walls such that $\pi_\l(g_{{\lambda}})=\rmid$.
We call $\frakD_\l$ the  \emph{reduction\/} of $\frakD$ at degree $\l$. \index{reduction (of a scattering diagram)}
By the finiteness condition, $ \frakD_\l$ has  only finitely many walls.
\end{defn}

We give some related notions.

\begin{defn}[Support/Singular locus]
For a scattering diagram $\frakD=\{ \bfw_{\lambda}= (\frakd_{\lambda}, g_{{\lambda}})_{n_{\lambda}}\}_{\lambda\in \Lambda}$,
the \emph{support\/} and the \emph{singular locus\/} of $\frakD$ are defined by
\index{support!of a scattering diagram} \index{singular locus (of a scattering diagram)}
\begin{align}
\mathrm{Supp}(\frakD)&=\bigcup_{\lambda\in \Lambda} \frakd_{\lambda},\\ 
\mathrm{Sing}(\frakD)&=\bigcup_{\lambda\in \Lambda} \partial\frakd_{\lambda} 
\cup
\bigcup_{\ss \lambda, \lambda'\in \Lambda \atop \ss \dim 
\frakd_{\lambda}\cap \frakd_{\lambda'}\,=\,r-2
} \frakd_{\lambda}\cap \frakd_{\lambda'}.
\end{align}
\end{defn}
We note that, if $\dim 
\frakd_{\lambda}\cap \frakd_{\lambda'}<r-2$,
then it is contained in $\partial \frakd_{\lambda}\cup \partial\frakd_{\lambda'}$.

\begin{defn}[Admissible curve] \index{admissible curve}
\label{3defn:adm1}
A curve $\gamma:[0,1]\rightarrow M_{\bbR}$ is
\emph{admissible\/} for a scattering diagram $\frakD$
if it satisfies the following conditions:
\begin{itemize}
\item
The end points of $\gamma$ are in $M_{\bbR}\setminus \mathrm{Supp}(\frakD)$.
\item
It is a smooth curve, and it intersects $\mathrm{Supp}(\frakD)$ transversally.
\item
$\gamma$ does not intersect  $\mathrm{Sing}(\frakD)$.
\end{itemize}
\end{defn}

\begin{defn}[Path-ordered product] \index{path-ordered product}
\label{3defn:pop1}
Let $\frakD$ be any scattering diagram.
For any admissible curve $\gamma$,
we define an element $\frakp_{\gamma,\frakD}\in G$ as follows:
For each integer $\l>0$, the reduction $\frakD_\l$ at $\l$ has only finitely many  walls.
Suppose that
$\gamma=\gamma(t)$ crosses   walls $\bfw_i=(\frakd_i, g_{i})_{n_i}$
($i=1,\dots,s$)
of $\frakD_\l$
in this order at $t=t_i$ with 
\begin{align}
0<t_1\leq t_2\leq \cdots \leq t_s<1.
\end{align}
Since $\gamma(t_i)\notin \mathrm{Sing}(\frakD)$,
 when $\gamma$ crosses multiple walls  at a time,
these walls have a common normal vector.
We define the \emph{intersection sign\/}  $\epsilon_i$ ($i=1,\dots,s$)  by \index{intersection sign}
\begin{align}
\label{3eq:factor1}
\epsilon_i
=
\begin{cases}
1 & \langle n_{i}, \gamma'(t_i)\rangle<0,\\
-1 & \langle n_{i}, \gamma' (t_i)\rangle>0,\\
\end{cases}
\end{align}
where  $\gamma'(t_i)$ is the velocity vector of $\gamma(t)$ at $t_i$.
Now, we define
\begin{align}
\label{3eq:pdef1}
\frakp_{\gamma,\frakD_\l}:=&\
g_{s}^{\epsilon_s}
\cdots
g_{1}^{\epsilon_1}
, \\
\label{3eq:pdef2}
\frakp_{\gamma,\frakD}:=&\
\lim_{\l \rightarrow\infty}\frakp_{\gamma,\frakD_\l}
\in G.
\end{align}
We call $\frakp_{\gamma,\frakD}$
the \emph{path-ordered product\/} (of wall elements in $\frakD$)
along $\gamma$.
\end{defn}

Note that $\frakp_{\gamma,\frakD}$ only depends on the homotopy class of $\gamma$
in $M_{\bbR}\setminus \mathrm{Sing}(\frakD)$.

\begin{defn}[Equivalence] \index{scattering diagram!equivalent}
Two scattering diagrams $\frakD$ and $\frakD'$ for   a common  seed $\fraks$ are
\emph{equivalent\/} if, for any curve $\gamma$ that is admissible for both $\frakD$ and $\frakD'$, the equality 
$\frakp_{\gamma,\frakD}=\frakp_{\gamma,\frakD'}$ holds.
\end{defn}

\begin{defn}[Generality] \index{general point (in $M_{\bbR}$)}
\label{3defn:general1}
We say that  $z\in M_{\bbR}$ is \emph{general\/}
if there is at most one rational hyperplane $n^{\perp}$ ($n \in N$, $n\neq 0$) such that
$z\in n^{\perp}$.
\end{defn}

\begin{lem}
\label{3lem:gen1}
(a). The set of all general points is dense in $M_{\bbR}$.
\par
(b). For any scattering diagram $\frakD$,
any point $z\in \mathrm{Sing}(\frakD)$ is not general.
\end{lem}
\begin{proof}
(a). This is clear.
(b). Suppose that $z\in \partial\frakd_{\lambda}$ for some  $\frakd_{\lambda}$.
Since $\frakd_{\lambda}$ is a cone of codimension 1, $z$ belongs to a face of  
$\frakd_{\lambda}$ of codimension 2,
which lies in the intersection of two rational hyperplanes.
Thus, $z$ is not general.
For $z\in\frakd_{\lambda}\cap \frakd_{\lambda'}$ with
$\dim \frakd_{\lambda}\cap \frakd_{\lambda'}=r-2$,
the claim is clear.
\end{proof}

\begin{defn}[Total wall element at $z$] \index{total wall element}
For any general $z\in M_{\bbR}$, we define
\begin{align}
\label{3eq:gxD1}
g_z(\frakD):=\prod_{\lambda:\, z\in \frakd_{\lambda}}
g_{{\lambda}}
\in 
\begin{cases}
G_{n_0}^{\parallel}
&
 \text{$z\in n^{\perp}_0$ for some unique $n_0\in N_{\rmpr}^+$},\\
\{\rmid\}
& \text{otherwise}.
\end{cases}
\end{align}
The (possibly infinite) product is well-defined due to the finiteness condition
of $\frakD$.
We call it the \emph{total wall element of $\frakD$ at $z$}.
\end{defn}

The condition for the equivalence
reduces to the following  local condition.

\begin{lem}[{\cite[Lemma 1.9]{Gross14}}]
\label{3lem:gx1}
Two scattering diagrams $\frakD$ and $\frakD'$ are equivalent
if and only if $g_z(\frakD)=g_z(\frakD')$ holds for any general $z \in M_{\bbR}$.
\end{lem}
\begin{proof}
Consider the reductions $\frakD_\l$ and $\frakD'_\l$ at $\l$.
Note that both have only finitely many walls.
Then, thanks to Lemma \ref{3lem:gen1},
the equivalence condition  reduces to
 the following  condition:
 \begin{itemize}
 \item
 For any general  $z \in M_{\bbR}$
 and  for  any  admissible curve $\gamma$
that intersects walls of $\frakD_\l$ and $\frakD'_\l$ only at $z$,
the following equality holds:
\begin{align}
\label{3eq:pp1}
 \frakp_{\gamma,\frakD_\l}\equiv \frakp_{\gamma,\frakD_\l'}
 \mod G^{>\l}.
 \end{align}
\end{itemize}
Under the situation, we have $\frakp_{\gamma,\frakD_\l}
 = g_z(\frakD_\l)^{\epsilon}$ with the intersection sign $\epsilon$.
 Therefore, the above equality is written as 
 $g_z(\frakD_\l)\equiv g_z(\frakD'_\l)$ mod $G^{>\l}$.
 Since $\l$ is arbitrary, 
we have $g_z(\frakD)=g_z(\frakD')$.
\end{proof}

\begin{defn}[Trivial scattering diagram] \index{scattering diagram!trivial}
\label{3defn:trivial1}
A scattering diagram $\frakD$ is \emph{trivial\/} if $g_z(\frakD)=\rmid$
for any general $z\in M_{\bbR}$.
\end{defn}

By Lemma \ref{3lem:gx1}, a scattering diagram is trivial if and only if
it is equivalent to the empty scattering diagram $\emptyset$.

Finally, we introduce
the most crucial notion for scattering diagrams.
\begin{defn}[Consistency] \index{scattering diagram!consistent}
A scattering diagram $\frakD$ is \emph{consistent\/} if for any admissible curve $\gamma$
for $\frakD$,
the associated path-ordered product $\frakp_{\gamma,\frakD}$  depends
only on
the end points of $\gamma$.
\end{defn}

The following restatement is useful.
\begin{lem}
\label{3lem:consist2}
A scattering diagram $\frakD$ is consistent if and only
if
\begin{align}
\label{3eq:pgfi1}
\frakp_{\gamma, \frakD}=\rmid
\end{align}
holds for any admissible loop (i.e., closed curve) $\gamma$ for $\frakD$.
\end{lem}
\begin{proof}
This is proved by the standard argument for fundamental groups.
Assume that $\frakD$ is consistent.
Let $\gamma$ be any admissible loop for $\frakD$.
Split $\gamma$ as $\gamma=(\gamma_2)^{-1} \circ \gamma_1$,
where $\gamma_1$ and $\gamma_2$ are admissible curves having the common ending points.
Then, by the consistency, we have $\frakp_{\gamma_1, \frakD}=
\frakp_{\gamma_2, \frakD}$.
Therefore, $\frakp_{\gamma, \frakD}=
(\frakp_{\gamma_2, \frakD})^{-1}\circ \frakp_{\gamma_1, \frakD}=
\rmid$ holds. The converse is similar.
\end{proof}

\subsection{Rank 2 examples: finite type}
\label{3subsec:rank2}

Based on the pentagon relation in Proposition \ref{3prop:pent1},
we construct some prototypical examples of consistent scattering diagrams of rank 2.
Let $e_1^*, e_2^*\in M$ be the dual basis of $e_1, e_2\in N$.
Accordingly,  let $f_1= e_1^*/\d_1, f_2= e_2^*/\d_2 \in M^{\circ}$ be the dual basis of $
\d_1 e_1, \d_2 e_2\in N^{\circ}$.
Let $\bfe_1, \bfe_2 \in \bbZ^2$ be the unit vectors.
We identify $N \simeq \bbZ^2$, $e_i \mapsto \bfe_i$,
and $M_{\bbR} \simeq \bbR^2$, $f_i \mapsto \bfe_i$.
Then, 
the canonical pairing
$\langle n, z \rangle:
 N\times M_{\bbR}\rightarrow \bbR$
   is given 
   by the corresponding vectors $\bfn$ and $\bfz$ as
\begin{align}
\label{3eq:can1}
\langle \bfn, \bfz \rangle
=
\bfn^T
\begin{pmatrix}
\d_1^{-1} & 0 \\
0 & \d_2^{-1} 
\end{pmatrix}
\bfz.
\end{align}
Let 
$B=(b_{ij})$
be the $2 \times 2$ integer matrix  introduced in \eqref{3eq:cseed1},  namely,
\begin{align}
B=
\begin{pmatrix}
0 & -\d_1 a \\
\d_2 a & 0
\end{pmatrix},
\quad
a=\{e_2,e_1\}.
\end{align}

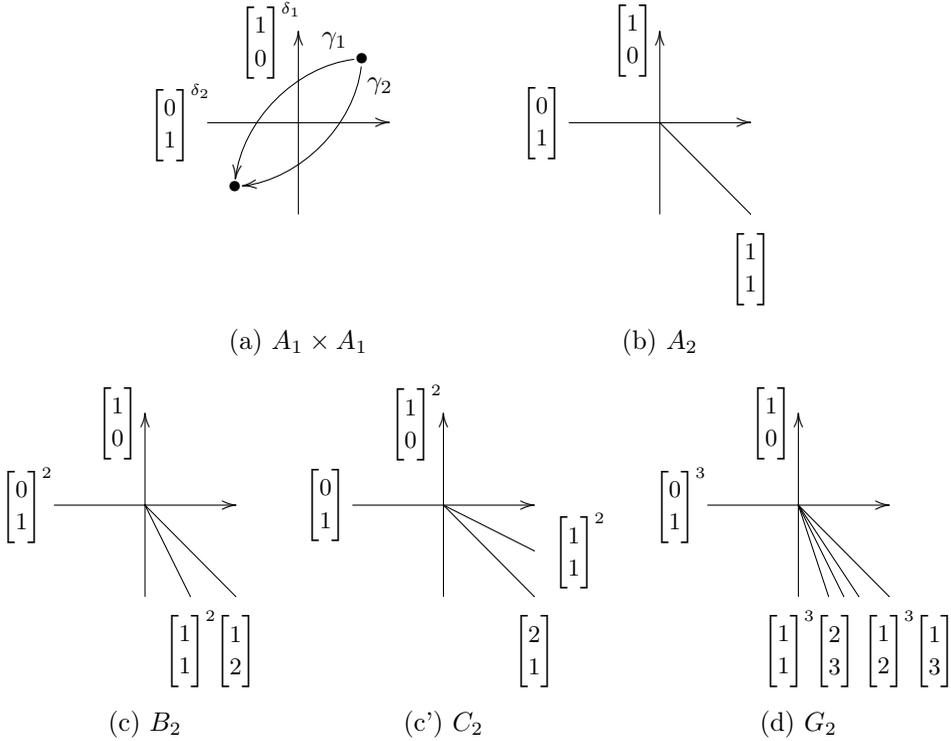
\begin{figure}
\centering
\leavevmode
\begin{xy}
0;/r1.2mm/:,
(0,-24)*{\text{(a) $A_1 \times A_1$}},
(4,9)*{\gamma_1};
(9,4)*{\gamma_2};
(7,7)*+{\bullet};(-7,-7)*+{\bullet};
(7,7)*+{};(-7,-7)*+{}
   **\crv{(0,6.7)&(-6.7,0)}
          ?>*\dir{>};
       (7,7)*+{};(-7,-7)*+{}
       **\crv{(6.7,0)&(0,-6.7)}
       ?>*\dir{>};
(-3, 9)*{\text{\small $\begin{bmatrix}1\\0\end{bmatrix}^{\d_1}$}},
(-13, 0)*{\text{\small $\begin{bmatrix}0\\1\end{bmatrix}^{\d_2}$}},
(0,0)="A",
\ar "A"+(0,0); "A"+(10,0)
\ar "A"+(0,0); "A"+(0,10)
\ar@{-} "A"+(0,0); "A"+(-10,0)
\ar@{-} "A"+(0,0); "A"+(0,-10)
\end{xy}
\hskip45pt
\begin{xy}
0;/r1.2mm/:,
(0,-24)*{\text{(b) $A_2$}},
(-3, 9)*{\text{\small $\begin{bmatrix}1\\0\end{bmatrix}$}},
(-13, 0)*{\text{\small $\begin{bmatrix}0\\1\end{bmatrix}$}},
(10, -16)*{\text{\small $\begin{bmatrix}1\\1\end{bmatrix}$}},
(0,0)="A",
\ar "A"+(0,0); "A"+(10,0)
\ar "A"+(0,0); "A"+(0,10)
\ar@{-} "A"+(0,0); "A"+(-10,0)
\ar@{-} "A"+(0,0); "A"+(0,-10)
\ar@{-} "A"+(0,0); "A"+(10,-10)
\end{xy}
\hskip10pt
\
\vskip10pt
\leavevmode
\begin{xy}
0;/r1.2mm/:,
(0,-24)*{\text{(c) $B_2$}},
(-3, 9)*{\text{\small $\begin{bmatrix}1\\0\end{bmatrix}$}},
(-13, 0)*{\text{\small $\begin{bmatrix}0\\1\end{bmatrix}^2$}},
(5, -15.75)*{\text{\small $\begin{bmatrix}1\\1\end{bmatrix}^2$}},
(10, -16)*{\text{\small $\begin{bmatrix}1\\2\end{bmatrix}$}},
(0,0)="A"
\ar "A"+(0,0); "A"+(10,0)
\ar "A"+(0,0); "A"+(0,10)
\ar@{-} "A"+(0,0); "A"+(-10,0)
\ar@{-} "A"+(0,0); "A"+(0,-10)
\ar@{-} "A"+(0,0); "A"+(5,-10)
\ar@{-} "A"+(0,0); "A"+(10,-10)
\end{xy}
\hskip15pt
\begin{xy}
0;/r1.2mm/:,
(0,-24)*{\text{(c') $C_2$}},
(-3, 9)*{\text{\small $\begin{bmatrix}1\\0\end{bmatrix}^2$}},
(-13, 0)*{\text{\small $\begin{bmatrix}0\\1\end{bmatrix}$}},
(15, -5)*{\text{\small $\begin{bmatrix}1\\1\end{bmatrix}^2$}},
(10, -16)*{\text{\small $\begin{bmatrix}2\\1\end{bmatrix}$}},
(0,0)="A"
\ar "A"+(0,0); "A"+(10,0)
\ar "A"+(0,0); "A"+(0,10)
\ar@{-} "A"+(0,0); "A"+(-10,0)
\ar@{-} "A"+(0,0); "A"+(0,-10)
\ar@{-} "A"+(0,0); "A"+(10,-5)
\ar@{-} "A"+(0,0); "A"+(10,-10)
\end{xy}
\hskip15pt
\begin{xy}
0;/r1.2mm/:,
(0,-24)*{\text{(d) $G_2$}},
(-3, 9)*{\text{\small $\begin{bmatrix}1\\0\end{bmatrix}$}},
(-13, 0)*{\text{\small $\begin{bmatrix}0\\1\end{bmatrix}^3$}},
(-1, -15.75)*{\text{\small $\begin{bmatrix}1\\1\end{bmatrix}^3$}},
(4, -16)*{\text{\small $\begin{bmatrix}2\\3\end{bmatrix}$}},
(10, -15.75)*{\text{\small $\begin{bmatrix}1\\2\end{bmatrix}^3$}},
(15, -16)*{\text{\small $\begin{bmatrix}1\\3\end{bmatrix}$}},
(0,0)="A"
\ar "A"+(0,0); "A"+(10,0)
\ar "A"+(0,0); "A"+(0,10)
\ar@{-} "A"+(0,0); "A"+(-10,0)
\ar@{-} "A"+(0,0); "A"+(0,-10)
\ar@{-} "A"+(0,0); "A"+(3.33,-10)
\ar@{-} "A"+(0,0); "A"+(5,-10)
\ar@{-} "A"+(0,0); "A"+(6.66,-10)
\ar@{-} "A"+(0,0); "A"+(10,-10)
\end{xy}
\caption{Examples of consistent scattering diagrams of rank 2.
Column vectors with exponents represent the wall elements  $\Psi[(n_1,n_2)]^c$.}
\label{3fig:scat1}
\end{figure}

(a). Type $A_1 \times A_1$.
We start with the case $a=0$, so that $\Gamma$ is degenerate.
Let $\d_1$ and $\d_2$ be arbitrary.
By \eqref{3eq:com1}, we have
\begin{align}
\label{3eq:pent0}
\Psi[\bfe_2 ]^{\d_2} \Psi[ \bfe_1]^{\d_1}=
\Psi[  \bfe_1 ]^{\d_1}  \Psi[ \bfe_2]^{\d_2}.
\end{align}
This equality is naturally interpreted as
a consistent scattering diagram of rank 2 in Figure \ref{3fig:scat1} (a).
Namely, it consists of two walls
\begin{align}
(\bfe_1^{\perp}, \Psi[\bfe_1]^{\d_1})_{\bfe_1},
\quad
(\bfe_2^{\perp}, \Psi[\bfe_2]^{\d_2})_{\bfe_2}.
\end{align}
The left hand side of the equality \eqref{3eq:pent0}
is the path-ordered product $\frakp_{\gamma_1, \frakD}$
along $\gamma_1$,
while the right hand side is the one $\frakp_{\gamma_2, \frakD}$ along $\gamma_2$.

Now assume $a\neq 0$, so that $\Gamma$ is nondegenerate.
Using the rescaling of  fixed data
in Section \ref{3subsec:equivalence1}
and interchanging $e_1$ and $e_2$, if necessary,
we may assume
$\{e_2,e_1\}=1$
without loosing generality.
Thus, we have
\begin{align}
\label{3eq:B1}
B=
\begin{pmatrix}
0 & -\d_1 \\
\d_2 & 0
\end{pmatrix}
=
\begin{pmatrix}
\d_1 & 0\\
0 & \d_2
\end{pmatrix}
\begin{pmatrix}
0 & -1 \\
1 & 0
\end{pmatrix}
=\Delta\Omega.
\end{align}
Let $\bfn=(n_1, n_2)\in \bbZ^2$ be the  corresponding vector to  $n\in N^+$.
Then, we have $\{\bfn',\bfn \}=n'_2 n_1 - n'_1 n_2$.

Here, we concentrate on the case $\d_1\d_2\leq 3$.

\par
(b). Type $A_2$. Let $\d_1=\d_2=1$.
Since $\{\bfe_2, \bfe_1\}=1$,
we apply
the pentagon relation \eqref{3eq:pent1} with $c=1$,
and we have
\begin{align}
\label{3eq:pent2}
\Psi[\bfe_2 ] \Psi[ \bfe_1]=
\Psi[  \bfe_1 ] \Psi[(1,1)] \Psi[ \bfe_2].
\end{align}
This equality is naturally interpreted as
a consistent scattering diagram of rank 2 in Figure \ref{3fig:scat1} (b).
Namely, it consists of three walls
\begin{align}
(\bfe_1^{\perp}, \Psi[\bfe_1])_{\bfe_1},
\quad
(\bfe_2^{\perp}, \Psi[\bfe_2])_{\bfe_2},
\quad
(\bbR_{\geq 0} (1,-1), \Psi[(1,1)])_{(1,1)}.
\end{align}
All normal vectors $\bfe_1$, $\bfe_2$, $(1,1)$ exhaust the positive roots of the root system of type
$A_2$.

\par
(c). Type $B_2$ and $C_2$. 
Here we consider two cases,
$\d_1=1$, $\d_2=2$ (type $B_2$) and
$\d_1=2$, $\d_2=1$ (type $C_2$).
They  are essentially the same.
Nevertheless,
we present both because they are useful    to construct the rank 3 scattering diagrams
later
in Section \ref{3subsec:rank31}.
First, we consider the   type $B_2$ case.
We consider an element $\Psi[\bfe_2 ]^2 \Psi[ \bfe_1]$.
For simplicity,
let us write $\Psi[(n_1,n_2)]$ as $[n_1,n_2]$.
Then, 
by applying
the pentagon relation \eqref{3eq:pent1}  with $c=1$ repeatedly for adjacent pairs
$\bfn'$, $\bfn$
with $\{\bfn', \bfn\}=1$, we have
\begin{align}
\label{3eq:pent3}
\begin{split}
\begin{bmatrix}
0\\
1
\end{bmatrix}
^2
\begin{bmatrix}
1\\
0
\end{bmatrix}
&=
\begin{bmatrix}
0\\
1
\end{bmatrix}
\begin{bmatrix}
1\\
0
\end{bmatrix}
\begin{bmatrix}
1\\
1
\end{bmatrix}
\begin{bmatrix}
0\\
1
\end{bmatrix}
=
\begin{bmatrix}
1\\
0
\end{bmatrix}
\begin{bmatrix}
1\\
1
\end{bmatrix}
\begin{bmatrix}
0\\
1
\end{bmatrix}
\begin{bmatrix}
1\\
1
\end{bmatrix}
\begin{bmatrix}
0\\
1
\end{bmatrix}
\\
&=
\begin{bmatrix}
1\\
0
\end{bmatrix}
\begin{bmatrix}
1\\
1
\end{bmatrix}
^
2
\begin{bmatrix}
1\\
2
\end{bmatrix}
\begin{bmatrix}
0\\
1
\end{bmatrix}
^2
.
\end{split}
\end{align}
This equality is naturally interpreted as
a consistent scattering diagram in Figure \ref{3fig:scat1} (c),
which consists of four walls
\begin{gather}
(\bfe_1^{\perp}, \Psi[\bfe_1])_{\bfe_1},
\quad
(\bfe_2^{\perp}, \Psi[\bfe_2]^2)_{\bfe_2},
\\
(\bbR_{\geq 0} (1,-2), \Psi[(1,1)]^2)_{(1,1)},
\quad
(\bbR_{\geq 0} (1,-1), \Psi[(1,2)])_{(1,2)}.
\end{gather}
Some remarks are in order.
Firstly, by \eqref{3eq:can1},
we have
\begin{align}
\bbR_{\geq 0} (1,-2)\subset (1,1)^{\perp},
\quad
\bbR_{\geq 0} (1,-1)\subset (1,2)^{\perp}.
\end{align}
Secondly,
the exponent $c=1$ or 2 of  each wall element $\Psi[\bfn]^c$ is the normalization factor $\d(\bfn)$
in Definition \ref{3defn:normal1}, where $N^{\circ}$ is identified with the sublattice of $\bbZ^2$ generated by
$(1,0)$ and $(0,2)$.
Thirdly, the result is compatible with the formula \eqref{3eq:gxg2}.
All normal vectors $\bfe_1$, $\bfe_2$, $(1,1)$, $(1,2)$ exhaust the positive roots of the root system of type
$B_2$. 
The  type $C_2$ case is similar,
and the corresponding scattering diagram is presented in Figure \ref{3fig:scat1} (c').

\par
(d). Type $G_2$. Let $\d_1=1$, $\d_2=3$.
We consider an element  $\Psi[\bfe_2 ]^3 \Psi[ \bfe_1]$.
Again,
by applying
the pentagon relation \eqref{3eq:pent1}  with $c=1$ repeatedly, we have
\begin{align}
\label{3eq:pent4}
\begin{split}
\begin{bmatrix}
0\\
1
\end{bmatrix}
^3
\begin{bmatrix}
1\\
0
\end{bmatrix}
&=
\begin{bmatrix}
0\\
1
\end{bmatrix}
\biggl(
\begin{bmatrix}
1\\
0
\end{bmatrix}
\begin{bmatrix}
1\\
1
\end{bmatrix}
^
2
\begin{bmatrix}
1\\
2
\end{bmatrix}
\begin{bmatrix}
0\\
1
\end{bmatrix}
^2
\biggr)
\\
&=
\begin{bmatrix}
1\\
0
\end{bmatrix}
\begin{bmatrix}
1\\
1
\end{bmatrix}
^2
\begin{bmatrix}
1\\
2
\end{bmatrix}
\begin{bmatrix}
1\\
1
\end{bmatrix}
\begin{bmatrix}
1\\
2
\end{bmatrix}
^2
\begin{bmatrix}
1\\
3
\end{bmatrix}
\begin{bmatrix}
0\\
1
\end{bmatrix}
^3
\\
&=
\begin{bmatrix}
1\\
0
\end{bmatrix}
\begin{bmatrix}
1\\
1
\end{bmatrix}
^3
\begin{bmatrix}
2\\
3
\end{bmatrix}
\begin{bmatrix}
1\\
2
\end{bmatrix}
^3
\begin{bmatrix}
1\\
3
\end{bmatrix}
\begin{bmatrix}
0\\
1
\end{bmatrix}
^3
.
\end{split}
\end{align}
This equality is naturally interpreted as
a consistent scattering diagram in Figure \ref{3fig:scat1} (d),
which consists of six walls
\begin{gather}
(\bfe_1^{\perp}, \Psi[\bfe_1])_{\bfe_1},
\quad
(\bfe_2^{\perp}, \Psi[\bfe_2]^3)_{\bfe_2},
\\
(\bbR_{\geq 0} (1,-3), \Psi[(1,1)]^3)_{(1,1)},
\quad
(\bbR_{\geq 0} (1,-2), \Psi[(2,3)])_{(2,3)},
\\
(\bbR_{\geq 0} (2,-3), \Psi[(1,2)]^3)_{(1,2)},
\quad
(\bbR_{\geq 0} (1,-1), \Psi[(1,3)])_{(1,3)}.
\end{gather}
The same remarks in the case (c) are applied as before.
The case $\d_1=3$, $\d_2=1$ is similar.

These examples are particular cases of
\emph{cluster scattering diagrams}, which we are going to study.

\subsection{Decompositions of $G$ at general point}

Let us begin with a general fact on nilpotent Lie algebras.

\begin{lem}
\label{3lem:decom1}
Let $\frakg$ be a nilpotent Lie algebra that has
a decomposition (as a vector space)
\begin{align}
\label{3eq:gdecom1}
\frakg=\frakg_1 \oplus \frakg_2
\end{align}
 by its Lie subalgebras $\frakg_1$ and $\frakg_2$.
 Let $G=\exp(\frakg)$ and $G_i=\exp(\frakg_i)$ $(i=1,2)$
 be the corresponding exponential groups 
 whose products are defined
 by the BCH formula \eqref{3eq:BCH1}.
 Then, we have the following decomposition of $G$:
 \begin{align}
 G=G_1G_2,
 \quad
 G_1\cap G_2=\{ \mathrm{id}\}.
 \end{align}
In other words, any element $g\in G$  is uniquely factorized as
 $g=g_1g_2$ $(g_i\in G_i)$.
 \end{lem}
\begin{proof}
By the assumption \eqref{3eq:gdecom1},
any element $g\in G$ is expressed as $g=\exp(X+Y)$, $X\in \frakg_1$, $Y\in \frakg_2$.
Then, apply the Zassenhaus formula  \eqref{3eq:Zass1}. If the right hand side of \eqref{3eq:Zass1}
is already of the form $g_1g_2$ $(g_i\in G_i)$, we are done.
Otherwise, 
expand the terms $\exp(Z)$ $(Z\not\in \frakg_1 \cup \frakg_2)$ therein
by \eqref{3eq:Zass1} again,
and transpose the left-most pair  $\exp(X_2) \exp(X_1)$ 
($X_i\in \frakg_i$) in the opposite order by \eqref{3eq:Zass2}.
Continue the process until we end up with the expression
$g_1g_2$, $g_i\in G_i$.
Due to the nilpotency of $\frakg$, the process completes 
in finitely many steps.
Thus, $G=G_1G_2$.
If $\exp(X_1)=\exp(X_2)$ ($X_i\in \frakg_i$),
we have $X_1=X_2$.
Then, we have $X_1=X_2=0$ by \eqref{3eq:gdecom1}.
Therefore, $ G_1\cap G_2=\{ \mathrm{id}\}$.
\end{proof}

We apply the above lemma to the following situation.
For a given  seed $\fraks$ of a fixed data $\Gamma$,
let $\frakg$, $\widehat\frakg$, and $G$ be the ones in Section \ref{3subsec:structure1}.
Let $z\in M_{\bbR}$ be any general point.
We introduce the \emph{decomposition of $\frakg$ at $z$} as
\begin{align}
\label{3eq:decom1}
\frakg=
\frakg_+^{z}
\oplus
\frakg_0^{z}
\oplus
\frakg_-^{z},
\end{align}
where 
\begin{align}
\frakg_+^{z}
=\bigoplus_{\ss n\in N^+ \atop \ss \langle n,z \rangle>0}\frakg_n,
\quad
\frakg_0^{z}
=\bigoplus_{\ss n\in N^+ \atop\ss \langle n,z \rangle=0}\frakg_n,
\quad
\frakg_-^{z}
=\bigoplus_{\ss n\in N^+ \atop \ss \langle n,z \rangle<0}\frakg_n.
\end{align}
By \eqref{3eq:Xbra1}, 
$\frakg_+^{z}$,
$\frakg_0^{z}$,
$\frakg_-^{z}$
are Lie subagebras of $\frakg$.
By the generality assumption of $z$,
we have
\begin{align}
\label{3eq:g0x1}
\frakg_0^{z}
=
\begin{cases}
\frakg_{n_0}^{\parallel}
& \text{$z\in n^{\perp}_0$ for some unique $n_0\in N_{\rmpr}^+$},
\\
\{0 \}
& \text{otherwise},
\end{cases}
\end{align}
where 
$\frakg_{n}^{\parallel}$ is the one in
\eqref{3eq:gpar1}.

Let
$\widehat\frakg_+^{z}$,
$\widehat\frakg_0^{z}$,
$\widehat\frakg_-^{z}$
be their completions,
and
let $G_+^{z}$, $G_0^{z}$, $G_-^{z}$
be the corresponding subgroups of $G$.
In particular,
\begin{align}
\label{3eq:g0x2}
G_0^{z}
=
\begin{cases}
G_{n_0}^{\parallel}
& \text{$z\in n^{\perp}_0$ for some unique $n_0\in N_{\rmpr}^+$},
\\
\{0 \}
& \text{otherwise}.
\end{cases}
\end{align}
Correspondingly, we have the \emph{decomposition of $G$ at $z$}
as follows. \index{decomposition of $G$!at $z$}

\begin{prop}[{\cite[Proof~of~Theorem~1.17]{Gross14}}]
\label{3prop:Gdecom1}
For any general $z\in M_{\bbR}$,
we have
\begin{align}
\label{3eq:Gdecom1}
G=G_+^{z}G_0^{z}G_-^{z}
\end{align}
so that
any element $g$ of $G$ is
uniquely factorized as
\begin{align}
\label{3eq:fac1}
g=g_+^{z}g_0^{z}g_-^{z}
\quad
(g_+^{z}\in G_+^{z},
\
g_0^{z}\in G_0^{z},
\
g_-^{z}\in G_-^{z}).
\end{align}
\end{prop}
\begin{proof}
Apply Lemma \ref{3lem:decom1}
to the decomposition
\begin{align}
\label{3eq:decom5}
\frakg^{\leq \l}=
(\frakg_+^{z})^{\leq \l}
\oplus
(\frakg_0^{z})^{\leq \l}
\oplus
(\frakg_-^{z})^{\leq \l},
\end{align}
and we obtain the decomposition of $G^{\leq \l}$
\begin{align}
\label{3eq:Gdecom2}
G^{\leq \l}
=(G_+^{z})^{\leq \l}
(G_0^{z})^{\leq \l}
(G_-^{z})^{\leq \l}
\end{align}
with the unique factorization property.
Moreover, the decomposition is compatible with the projection
$\pi_{\l',\l}$
in \eqref{3eq:kk1}.
Then, taking the inverse limit, we obtain
\eqref{3eq:Gdecom1}.
The unique factorization property is preserved in 
the inverse limit.
\end{proof}

By  Proposition \ref{3prop:Gdecom1}
and \eqref{3eq:g0x2}, for any $g\in G$ and any general $z\in M_{\bbR}$,
we have a unique $g_0^z$ such that
\begin{align}
\label{3eq:gx2}
g_0^z
\in 
\begin{cases}
G_{n_0}^{\parallel}
&
 \text{$z\in n^{\perp}_0$ for some unique $n_0\in N_{\rmpr}^+$},\\
\{\rmid\}
& \text{otherwise}.
\end{cases}
\end{align}
Compare it with \eqref{3eq:gxD1}.

\begin{figure}
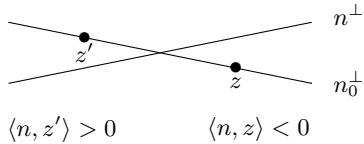

\centering
\leavevmode
\xy
(13,-10)*{\text{\small $ \langle n, z \rangle <0$}};
(-13,-10)*{\text{\small $ \langle n, z' \rangle >0$}};
(25,-4)*{\text{\small $ n_0^{\perp}$}};
(25,5)*{\text{\small $ n^{\perp}$}};
(-10,0)*{\text{\small $z'$}};
(10,-4)*{\text{\small $z$}};
(-10,2)*+{\bullet};
(10,-2)*+{\bullet};
\ar@{-} (-20,4); (20,-4)
\ar@{-} (-20,-4); (20,4)
\endxy
\caption{Decomposition at $z$ is not constant on $n_0^{\perp}$.}
\label{3fig:decomx1}
\end{figure}

What is important here is that, for a general  $z \in n_0^{\perp}$,
$g_0^{z}$
varies discontinuously  on the whole hyperplane  $n_0^{\perp}$,
in general.
When $z$ crosses the intersection $n_0^{\perp}\cap n^{\perp}$  with
 $n\in N_{\mathrm{pr}}^+$ other than $n_0$,
the sign  of $\langle n, z\rangle $ changes.
See Figure \ref{3fig:decomx1}.
Thus, the decomposition 
\eqref{3eq:Gdecom1}  changes.

\begin{figure}
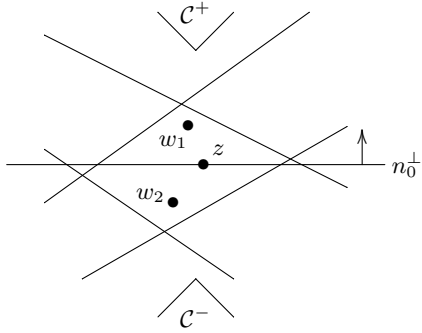

\centering
\leavevmode
\xy
(28,0)*{\text{\small $ n_0^{\perp}$}};
(3,2)*{\text{\small $z$}};
(-3,3)*{\text{\small $w_1$}};
(-6,-4)*{\text{\small $w_2$}};
(0,20)*{\text{\small $\calC^+$}};
(0,-20)*{\text{\small $\calC^-$}};
(,0)*+{\bullet};
(-1,5)*+{\bullet};
(-3,-5)*+{\bullet};
\ar@{-} (-25,0); (25,0)
\ar@{-} (-20,-5); (15,20)
\ar@{-} (-20,17); (20,-3)
\ar@{-} (-20,2); (5,-15)
\ar@{-} (-15,-15); (20,5)
\ar@{-} (0,15); (-5,20)
\ar@{-} (0,15); (5,20)
\ar@{-} (0,-15); (-5,-20)
\ar@{-} (0,-15); (5,-20)
\ar@{->} (22,0); (22,4.5)
\endxy
\caption{Configuration of $z$, $w_1$, $w_2$ 
for Lemma \ref{3lem:local1}.}
\label{3fig:local1}
\end{figure}

We are going to relate
 the  decomposition of $g\in G$ at $z$
with a consistent scattering diagram.
To begin with,
let us fix an integer $\l>0$.
There are only finitely many hyperplanes 
$n^{\perp}$ with
$n\in (N_{\mathrm{pr}}^+)^{\leq \l}$.
They intersect each other, so that each hyperplane
is subdivided into cones.
Let $z$, $w_1$, and $w_2$ be general points in $M_{\bbR}$
as follows: (See Figure \ref{3fig:local1}.)

\begin{itemize}
\item
$z$ belongs to $n_0^{\perp}$ for some  unique $n_0\in (N_{\mathrm{pr}}^+)^{\leq \l}$.
\item
$w_1$ and $w_2$ do not belong to  $n^{\perp}$ for any 
$n \in (N_{\mathrm{pr}}^+)^{\leq \l}$.
\item
$\langle n_0, w_1 \rangle >0$ and $\langle n_0, w_2 \rangle <0$.
\item
For  any $n \in (N_{\mathrm{pr}}^+)^{\leq \l}$ other than  $n_0$,
 $\langle n, w_i \rangle >0$ if and only if  $\langle n, z \rangle >0$.
\end{itemize}

We have a key lemma.

\begin{lem}[{\cite[Lemma~3.2]{Bridgeland16}}]
\label{3lem:local1}
For any $g\in G$ and
 the above $z,w_1,w_2\in M_{\bbR}$,
 we have
\begin{align}
\label{3eq:ggg1}
g_0^{z}\equiv
(g_+^{w_2})^{-1} g_+^{w_1}
\mod G^{>\l}.
\end{align}
Also, when the configuration of $w_1$ and $w_2$ are interchanged,
we have
\begin{align}
\label{3eq:ggg2}
(g_0^{z})^{-1}\equiv
(g_+^{w_2})^{-1} g_+^{w_1}
\mod G^{>\l}.
\end{align}
\end{lem}
\begin{proof}
Consider the decomposition $g= g_+^{z} g_0^{z}g_-^{z}$ at $z$.
Under the condition of $z$, $w_1$, and $w_2$,
the equality $g\equiv g_+^{z} g_0^{z}g_-^{z}$ mod $G^{>\l}$ is interpreted also as
the decompositions $g\equiv (g_+^{z} g_0^{z} )g_-^{z}$ at $w_1$
and  $g\equiv g_+^{z}(g_0^{z}g_-^{z})$ at $w_2$  mod $G^{>\l}$ as well,
where the middle terms  $g_0^{w_1}$ and $g_0^{w_2}$
are the identity  mod $G^{>\l}$.
Therefore,  $g_+^{w_1}\equiv g_+^{z}g_0^{z}$ and $g_+^{w_2}\equiv g_+^{z}$.
Thus, we have \eqref{3eq:ggg1}.
The formula  \eqref{3eq:ggg2} is obtained from  \eqref{3eq:ggg1} by taking the inverse.
\end{proof}

\subsection{Construction of consistent scattering diagrams}

For a given seed $\fraks$, we define
the \emph{positive and negative orthants\/}  $\calC^{\pm}$ in $M_{\bbR}$
by \index{cone!positive/negative}
\begin{align}
\label{3eq:chamber1}
\calC^+&=\calC_{\fraks}^+ :=\{ z \in M_{\bbR} \mid 
\langle e_i, z\rangle \geq 0\ (i=1,\dots,r)
\},
\\
\calC^-&=\calC_{\fraks}^- :=\{ z \in M_{\bbR} \mid 
\langle e_i, z\rangle \leq 0\ (i=1,\dots,r)
\}.
\end{align}
\begin{lem}
\label{3lem:cham1}
For any wall for $\fraks$,
its support $\frakd$ intersects
$\calC^{\pm}$ only in the boundary of $\calC^{\pm}$.
\end{lem}
\begin{proof}
We have $\frakd \subset n^{\perp}$ for its normal vector $n\in N_{\rmpr}^+$.
On the other hand,
for any $n\in N_{\rmpr}^+$,
$\langle n, z\rangle> 0$ if $z\in \mathrm{Int}(\calC^{+})$,
and $\langle n, z\rangle<0$ if $z\in \mathrm{Int}(\calC^{-})$.
Thus, $\frakd\cap \mathrm{Int}(\calC^{\pm})=\emptyset$.
\end{proof}

Let $\frakD$ be a consistent scattering diagram for $\fraks$.
Let $\gamma_{+-}$ be any admissible curve for $\frakD$ with the initial point in 
$\mathrm{Int}(\calC^+)$ and the final point in $\mathrm{Int}(\calC^-)$.
By the consistency and Lemma \ref{3lem:cham1}, 
the element
\begin{align}
g(\frakD):=\frakp_{\gamma_{+-}, \frakD}\in G
\end{align}
 is independent of the choice of such $\gamma_{+-}$.

Surprisingly, the single element $g(\frakD)$ contains all essential information
of the whole consistent scattering diagram $\frakD$
up to equivalence.

\begin{prop}[{\cite[Proof of Theorem 1.17]{Gross14}}]
\label{3prop:gxg1}
Let $\frakD$ be a consistent scattering diagram.
Let $g=g(\frakD)\in G$.
Then, for any general $z\in M_{\bbR}$, we have
\begin{align}
\label{3eq:gD1}
g_z(\frakD)=g_0^z.
\end{align}
\end{prop}

\begin{figure}
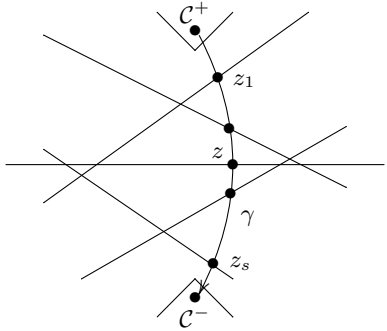

\centering
\leavevmode
\xy
(3,1.5)*{\text{\small $z$}};
(6.5,11)*{\text{\small $z_1$}};
(6,-13)*{\text{\small $z_s$}};
(7,-7)*{\text{\small $\gamma$}};
(0,20)*{\text{\small $\calC^+$}};
(0,-20)*{\text{\small $\calC^-$}};
(5,0)*+{\bullet};
(0,17.5)*+{\bullet};
(0,-17.5)*+{\bullet};
(3,11.4)*+{\bullet};
(4.5,4.7)*+{\bullet};
(4.7,-3.8)*+{\bullet};
(2.4,-13.1)*+{\bullet};
(0,18)*+{};(0,-18)*+{}
   **\crv{(10, 0)}
          ?>*\dir{>};
\ar@{-} (-25,0); (25,0)
\ar@{-} (-20,-5); (15,20)
\ar@{-} (-20,17); (20,-3)
\ar@{-} (-20,2); (5,-15)
\ar@{-} (-15,-15); (20,5)
\ar@{-} (0,15); (-5,20)
\ar@{-} (0,15); (5,20)
\ar@{-} (0,-15); (-5,-20)
\ar@{-} (0,-15); (5,-20)
\endxy
\caption{A path $\gamma$ for the proof of
Proposition \ref{3prop:gxg1}.
}
\label{3fig:gamma1}
\end{figure}

\begin{proof}
Let $z\in M_{\bbR}$ be a given general point.
One can find an admissible curve $\gamma$ 
for $\frakD$ satisfying the following conditions:
\begin{itemize}
\item
$\gamma$ intersects $z$.
\item
$\gamma$ has the  initial point  in 
$\mathrm{Int}(\calC^+)$ and the  final point  in $\mathrm{Int}(\calC^-)$.
\item
The velocity vector of $\gamma$ at any point
satisfies $\langle e_i, \gamma'\rangle <0$ for any $i$.
\end{itemize}
For example,  take the straight line passing through $z$
  from $z+a \sum_{i=1}^r e_i^*\in \mathrm{Int}(\calC^+)$ to $z-a \sum_{i=1}^r e_i^*
\in \mathrm{Int}(\calC^-)$
for sufficiently large $a>0$.
Then, slightly deform it  to avoid  intersecting with $\mathrm{Sing}(\frakD)$.

Consider the reduction $\frakD_\l$ at $\l$.
Let $z_1$, \dots, $z_s$ be the intersections of  $\gamma$ with walls of $\frakD_\l$
in this order such that $z_{j}=z$.
(If $z$ is not on any wall of $\frakD_\l$,
we regard  it  on a wall with the trivial wall element $\rmid$.)
See Figure \ref{3fig:gamma1}.
 Thanks to Lemma \ref{3lem:gen1},
by slightly deforming $\gamma$ if necessary,
we may assume that all $z_1$, \dots, $z_s$ are general.
The  intersection signs $\epsilon_i$ are all 1
due to the velocity condition.
Then, we have
\begin{align}
g
\equiv
\frakp_{\gamma,\frakD_\l}
=
g_{z_s}(\frakD)
\cdots
g_{z_{j+1}}(\frakD)
g_{z}(\frakD)
g_{z_{j-1}}(\frakD)
\cdots
g_{z_1}(\frakD)
\mod
G^{>\l}.
\end{align}
Let $n_i\in (N_{\rmpr}^+)^{\leq \l}$ be the
normal vector of the wall crossed by $\gamma$ at $z_i$.
Then,
we have $\langle n_i, z\rangle <0$
 for $i< j$,
 and $\langle n_i, z\rangle >0$
 for $i> j$.
This implies that
\begin{align}
g_z(\frakD) \equiv  g_0^z
\mod 
G^{>\l}.
\end{align}
Thus, we have $g_z(\frakD)=g_0^z$.
\end{proof}

The above result, together with Lemma
\ref{3lem:gx1}, shows that a consistent scattering diagram $\frakD$ 
is recovered from
a single element $g(\frakD)$ up to equivalence.
In fact, we have the following stronger result.

\begin{thm}[{\cite[Theorem 2.1.6]{Kontsevich13}, \cite[Theorem 1.17]{Gross14}}]
\label{3thm:bij1}
The assignment $\bbvarphi: \frakD \mapsto g(\frakD)$ gives a one-to-one correspondence between equivalence classes of 
 consistent scattering diagrams and elements in $G$.
\end{thm}

\begin{proof}
{\bf Injectivity of $\bbvarphi$.}
Let $\frakD$ and $\frakD'$ be consistent scattering diagrams,
and suppose that $g(\frakD)=g(\frakD')=g$.
Then, by \eqref{3eq:gD1},
\begin{align}
g_z(\frakD)=g_z(\frakD')=g_0^z
\end{align}
holds for any general $z\in M_{\bbR}$.
Then, by  Lemma
\ref{3lem:gx1}, $\frakD$ and $\frakD'$ are equivalent.

{\bf Surjectivity of $\bbvarphi$.}
Step 1
(Construction of $\bbvarphi^{-1}: g \mapsto \frakD(g)$).
Take
an arbitrary element $g\in G$.
Let $g_0^z$ be the one in
 Proposition \ref{3prop:Gdecom1}.
First, we construct  a scattering diagram $\frakD=\frakD(g)$ such that 
\begin{align}
g_z(\frakD)=g_0^z
\end{align}
holds for any general  $z\in M_{\bbR}$.
Let us fix any $n\in N_{\mathrm{pr}}^+$.
For each $\l>0$,
the hyperplane $n^{\perp}$ is divided  into
finitely many cones of codimension 1
by $n'{}^{\perp}$
with 
$n'\in (N_{\mathrm{pr}}^+)^{\leq \l}$,
$n'\neq n$.
Then, for general $z\in n^{\perp}$,
$\pi_\l(g_0^z)$ is constant with respect to $z$ on each cone.
Let $\frakd$ be any of such cones
and $z\in \frakd$ be any general point.
By \eqref{3eq:gx2}, $g_0^z$ is uniquely written as
\begin{align}
\label{3eq:expand1}
g_0^z=
\exp\biggl(\, \sum_{j>0}  c_{z,jn} X_{jn}\biggr)
\quad
(c_{z,jn}\in \bbk).
\end{align}
Then, with $\frakd$,
 we associate a wall $(\frakd, g_{\frakd})_{n}$,
where
\begin{align}
\label{3eq:gd1}
g_{\frakd}=
\begin{cases}
\exp( c_{z, jn} X_{j n})
& \text{$\l= j \deg(n)$ for some unique integer $j>0$},\\
\rmid
& \text{otherwise}.
\end{cases}
\end{align}
Let  $\frakD$  be the collection of all such walls for all
 $n\in N_{\mathrm{pr}}^+$ and $\ell$.
Then, $\frakD$ satisfies
the finiteness condition by construction.
Therefore,  $\frakD$ is a scattering diagram.
Moreover, 
 for any general  $z\in M_{\bbR}$,
 we have
\begin{align}
\label{3eq:gxDg1}
\begin{split}
g_z(\frakD)&=
\begin{cases}
\displaystyle
\prod_{j>0} \exp(c_{z,jn}X_{jn})
& \text{$z\in n^{\perp}$ for some unique $n\in N_{\rmpr}^+$},
\\
\rmid
& \text{otherwise}
\end{cases}\\
&=g_0^z
\end{split}
\end{align}
as desired.

Step 2. Let us show that $\frakD$ is consistent.
Consider the reduction $\frakD_\l$ at $\l$.
Let $\gamma$ be any admissible curve with the initial point $w$
 and and the final point $w'$.
Let $z_1$, \dots, $z_s$ be the intersections of  $\gamma$ with walls of $\frakD_\l$
in this order.
By slightly deforming $\gamma$, if necessary,
we may assume that all $z_1$, \dots, $z_s$ are general.
Then, by \eqref{3eq:gxDg1},
we have
\begin{align}
\label{3eq:pgamma1}
\frakp_{\gamma, \frakD_\l}=
(g_0^{z_s})^{\epsilon_s}
\cdots
(g_0^{z_1})^{\epsilon_1},
\end{align}
where $\epsilon_i$ is the intersection sign at $z_i$.
Let $w_1=w$, $w_2$, \dots, $w_{s+1}=w'$,
where $w_{i}$ ($i=2,\dots,s$)
is any general points on $\gamma$ between $z_{i-1}$ and $z_{i}$.
Note that, when $\epsilon_i=1$ (resp.~$\epsilon_i=-1$), 
$z_i$, $w_i$, $w_{i+1}$ (resp.~$z_i$, $w_{i+1}$, $w_{i}$)
are in the same configuration as
$z$, $w_1$, $w_2$ in Lemma \ref{3lem:local1}.
Therefore, by \eqref{3eq:ggg1} and \eqref{3eq:ggg2},
we have
\begin{align}
(g_0^{z_i})^{\epsilon_i}\equiv (g_+^{w_{i+1}})^{-1} g_+^{w_i}
\mod G^{>\l}.
\end{align}
Putting it into \eqref{3eq:pgamma1}, 
\begin{align}
\label{3eq:pgg1}
\frakp_{\gamma, \frakD_\l} \equiv
(g_+^{w'})^{-1} g_+^{w}
\mod G^{>\l}.
\end{align}
Thus, we conclude that
\begin{align}
\label{3eq:pgg2}
\frakp_{\gamma, \frakD} =
(g_+^{w'})^{-1} g_+^{w},
\end{align}
which only depends on the end points $w$ and $w'$ of $\gamma$.

Step 3. Finally, we show that $g(\frakD)=g$.
Let $\gamma_{+-}$ be any admissible curve with the initial point  $w$ in 
$\mathrm{Int}(\calC^+)$ and the final point $w'$ in $\mathrm{Int}(\calC^-)$.
Then,  
 we have the equality \eqref{3eq:pgg2}.
Meanwhile,
for any $n\in N^+$, we have
$\langle n, w\rangle >0$ and $\langle n, w'\rangle <0$.
Thus,  the decompositions of $g$ with respect to $w$ and $w'$ 
yield $g_+^{w}= g$ and $g_+^{w'}= \rmid$.
Therefore, we have $\frakp_{\gamma_{+-}, \frakD}=g$.
Thus, we conclude that $\bbvarphi$ is surjective.
\end{proof}

\begin{rem}
Theorem \ref{3thm:bij1} and its proof only depend on the facts that $\frakg$
is $N^+$-graded and $\frakg_n^{\parallel}$ is abelian.
Thus, it is applicable to a much wider class of $\frakg$.
\end{rem}

To summarize the results in this section,
all essential information of any consistent scattering diagram $\frakD$ is
encoded in the corresponding element $g(\frakD)=\frakp_{\gamma_{+-}, \frakD}$.
Conversely,
for any element $g$ of $G$,
one can construct a consistent scattering diagram $\frakD(g)$,
and admissible curves $\gamma_{+-}$ therein
correspond to  various factorizations of $g$
with $g_0^z$'s.

\begin{rem}
\label{3rem:face1}
The construction of a consistent scattering diagram in
Theorem \ref{3thm:bij1} also tells that
one may assume that, for the support $\frakd\subset n^{\perp}$ of any wall,
each face of $\frakd$ of codimension 2 
lies in the intersection with a  hyperplane $n'^{\perp}$
with a \emph{positive\/} normal vector $n'\in N_{\rmpr}^{+}$.
\end{rem}

\subsection{Minimal support}
\label{3subsec:appl1}

Often it is convenient to consider a consistent
scattering diagram whose support is \emph{minimal\/}
among all equivalent diagrams.

\begin{defn}[Minimal support] \index{scattering diagram!with minimal support}
Let $\frakD$ be a consistent scattering diagram and $g=g(\frakD)$ be
the corresponding element in $G$.
We say that $\frakD$ is \emph{with minimal support\/}
if it satisfies the following condition:
\begin{itemize}
\item
For any general $z\in M_{\bbR}$ such that $g_0^z=\rmid$,
there is no wall $(\frakd,g_{\frakd})_n$ of $\frakD$ such that $z\in \frakd$.
 \end{itemize}
 It is clear that the support of such $\frakD$ is the  minimal set among the supports of all consistent scattering diagrams
that are equivalent to $\frakD$.
\end{defn}

\begin{ex}
All  examples of rank 2 consistent scattering diagrams in Section \ref{3subsec:rank2} are 
with minimal support.
\end{ex}

\begin{prop}
\label{3prop:minimal1}
For any consistent scattering diagram $\frakD$, there is 
a (not unique)  consistent scattering diagram with minimal support
that is equivalent to $\frakD$.
\end{prop}
\begin{proof}
By Theorem \ref{3thm:bij1},
 $\frakD$  is equivalent to the consistent scattering diagram $\frakD(g)$ for some $g\in G$ constructed in the proof
of the surjectivity of $\bbvarphi$ in Theorem  \ref{3thm:bij1}.
Let  $z\in M_{\bbR}$ be a general point such that $g_0^z=\rmid$.
Then, for any wall $(\frakd, g_{\frakd})_n$ in $\frakD(g)$ such that $z\in \frakd$,
we have $g_{\sigma}=\rmid$ by
\eqref{3eq:gd1}.
Therefore, any such  wall is removed form $\frakD(g)$ up to equivalence.
The resulting diagram has the desired property.
\end{proof}

\notes

The contents are  taken mostly from \cite[Sect.~1.2]{Gross14},
\cite[Sect.~2]{Kontsevich13}
and partly from \cite{Bridgeland16}
with added/modified proofs.
The consistent scattering diagrams here are a special case of
consistent scattering diagrams
 in \cite{Gross07} and
  \emph{wall-crossing structures\/} (WCS)
  in
 \cite{Kontsevich13}.
The contents in Section \ref{3subsec:appl1} is added by us.

\newpage

\section{Cluster scattering diagrams}
\label{3sec:const1}
Here we introduce a special class of consistent scattering diagrams
 called \emph{cluster scattering diagrams\/} (CSDs).
Then, we show the construction/existence of CSDs.
This is the first fundamental result in the  application of scattering diagrams
 to cluster algebras.

\subsection{Parametrization of $G$ by parallel subgroups}

For a given   seed $\fraks$ of a (possibly degenerate) fixed data $\Gamma$,
let $\frakg$, $\widehat\frakg$, and $G$ be the ones in Section \ref{3subsec:structure1}.
We introduce  another decomposition of $G$.

For any  $n \in N_{\mathrm{pr}}^+$,
we introduce the  \emph{decomposition of $\frakg$ by $n$} as
\begin{align}
\label{3eq:decom6}
\frakg=
\frakg_+^{n}
\oplus
\frakg_0^{n}
\oplus
\frakg_-^{n}
\quad
(n \in N_{\mathrm{pr}}^+),
\end{align}
where 
\begin{align}
\frakg_+^{n}
=\bigoplus_{\ss n'\in N^+ \atop\ss  \{n',n\}>0}\frakg_{n'},
\quad
\frakg_0^{n}
=\bigoplus_{\ss n'\in N^+ \atop\ss  \{n',n\}=0}\frakg_{n'},
\quad
\frakg_-^{n}
=\bigoplus_{\ss n'\in N^+ \atop\ss  \{n',n\}<0}\frakg_{n'}.
\end{align}
By \eqref{3eq:Xbra1}, 
$\frakg_+^{n}$,
$\frakg_0^{n}$,
$\frakg_-^{n}$
are Lie subalgebras of $\frakg$.
Let 
$\widehat\frakg_+^{n}$,
$\widehat\frakg_0^{n}$,
$\widehat\frakg_-^{n}$
be their completions,
and let
$G_+^{n}$, $G_0^{n}$, $G_-^{n}$
be the corresponding subgroups of $G$.

\begin{ex}
In the extreme case of degeneracy $\{\cdot, \cdot\}=0$,
$G=G_0^n$
for any $n$.
\end{ex}

By  Lemma \ref{3lem:decom1},
we have the {\em decomposition of $G$ by $n$}. \index{decomposition of $G$!by $n$}
\begin{lem}[{\cite[\S1.2]{Gross14}}]
\label{3lem:Gdecom2}
For any
$n \in N_{\mathrm{pr}}^+$,
we have
\begin{align}
G=G_+^{n}G_0^{n}G_-^{n}
\end{align}
so that
any element $g$ of $G$ is
uniquely factorized as
\begin{align}
\label{3eq:fac2}
g=g_+^{n}g_0^{n}g_-^{n}
\quad
(g_+^{n}\in G_+^{n},
\
g_0^{n}\in G_0^{n},
\
g_-^{n}\in G_-^{n}).
\end{align}
\end{lem}

In contrast to the decomposition in 
Proposition \ref{3prop:Gdecom1},
the above decomposition directly depends on $n \in N_{\mathrm{pr}}^+$,
not each point $z$ on $n^{\perp}$.

One can further decompose 
$\frakg_0^{n}$ as
\begin{align}
\label{3eq:decom7}
\frakg_0^{n}
=
\frakg_{n}^{\parallel}
\oplus
\frakg_{n}^{\perp}
\quad
(n \in N_{\mathrm{pr}}^+),
\end{align}
where $\frakg_{n}^{\parallel}$ is
the Lie subalgebra of $\frakg_0^{n}$ in \eqref{3eq:gpar1}, and
\begin{align}
\label{3eq:gperp1}
\frakg_{n}^{\perp}
:=\bigoplus_{\ss n'\in N^+\atop \ss
\{n',n\}=0,\,
 n' \not\in \bbZ n }\frakg_{n'}.
\end{align}

\begin{lem}
\label{3lem:gideal1}
 $\frakg_{n}^{\perp}$ is
 a Lie algebra ideal of $\frakg_0^{n}$.
 \end{lem}
\begin{proof}
It is enough to show that
\begin{align}
[\frakg_0^{n}, \frakg_0^{n}]
\subset \frakg_{n}^{\perp},
\end{align}
or equivalently,  $[\frakg_0^{n}, \frakg_0^{n}]\cap 
 \frakg_{n}^{\parallel}=\{0 \}$.
 Suppose that $X_{n_1}, X_{n_2}\in \frakg_0^n$
 and $[X_{n_1}, X_{n_2}]\in \frakg_n^{\parallel}$.
 Then, we have
 $\{n_i,n\}=0$, and $n_1+n_2=j n$ for some $j>0$.
It follows that $ \{n_1,n_2\}=\{ n_1, jn\} =0$.
 Thus, by \eqref{3eq:Xbra1},
 $[X_{n_1}, X_{n_2}]=\{n_1,n_2\}X_{n_1+n_2}=0$.
\end{proof}

Let 
$\widehat\frakg_{n}^{\parallel}$ and
$\widehat\frakg_{n}^{\perp}$
be the completions of $\frakg_{n}^{\parallel}$ and
$\frakg_{n}^{\perp}$,
and let $G_{n}^{\parallel}$ and $G_{n}^{\perp}$
be the corresponding subgroups of $G$, respectively.

\begin{lem}[{\cite[\S1.2]{Gross14}}]
\label{3lem:Gdecom3}
For any $n \in N_{\mathrm{pr}}^+$, the following facts hold.
\par
(a). $G_{n}^{\perp}$ is a normal subgroup of $G_0^{n}$.
\par
(b).
We have
\begin{align}
G_0^n=G_{n}^{\parallel}G_{n}^{\perp}
\end{align}
so that any element $g_0^{n}$ of $G_0^{n}$ is
uniquely factorized as
\begin{align}
g_0^{n}=g_{n}^{\parallel} g_{n}^{\perp}
\quad
(g_{n}^{\parallel}\in G_{n}^{\parallel},
\
g_{n}^{\perp}\in G_{n}^{\perp}).
\end{align}
\par
(c).
The map 
\begin{align}
\label{3eq:Gdecom3}
G_0^{n}\rightarrow G_{n}^{\parallel},
\quad
g_0^{n}=g_{n}^{\parallel} g_{n}^{\perp}
\mapsto
g_{n}^{\parallel}
\end{align}
is a group homomorphism with  the kernel $G_{n}^{\perp}$.
\end{lem}
\begin{proof}
(a).
Let  $X\in (\frakg_{n}^{\perp})^{\leq \l}$ and $Y \in 
(\frakg_0^{n})^{\leq \l}$.
Then, 
by \eqref{3eq:Zass2},
the product
$\exp(Y)\exp(X)\exp(Y)^{-1}
$
is rewritten as a product of $\exp(X)$ and exponentials of
the commutators and
the higher commutators of $X$ and $Y$.
By Lemma \ref{3lem:gideal1},
it belongs to $(G_{n}^{\perp})^{\leq \l}$.
(b).
This follows from Lemma \ref{3lem:decom1}.
(c).
This follows from (a).
\end{proof}

Combining Lemmas \ref{3lem:Gdecom2} and \ref{3lem:Gdecom3},
we obtain  a  map
\begin{align}
\label{3eq:pro1}
\psi_{n}:G \rightarrow G_0^{n} \rightarrow G_{n}^{\parallel},
\quad g \mapsto g_0^{n} \mapsto g_{n}^{\parallel}
\quad
(n \in N_{\mathrm{pr}}^+).
\end{align}
(The map $\psi_n$ is not  a group homomorphism except for $n=e_i$.
See Lemma \ref{3lem:Gdecom4}.)
Consider a set (not a group)
\begin{align}
\label{3eq:S1}
S^{\parallel}=S^{\parallel}_G
:=\prod_{n\in N_{\mathrm{pr}}^+}
 G_{n}^{\parallel}.
\end{align}
We then have a map  (not a group homomorphism),
\begin{align}
\bbpsi: G \rightarrow
S^{\parallel},
 \quad
 g \mapsto \bfg^{\parallel} = (g_{n}^{\parallel})_{n\in N_{\mathrm{pr}}^+}.
\end{align}
We note that
\begin{align}
\label{3eq:Gid1}
 (G_{n}^{\parallel})^{\leq \l}
 =\{ \rmid \}
\quad
(n \in  (N_{\mathrm{pr}}^+)^{> \l}).
\end{align}
In particular, we have
\begin{align}
\label{3eq:pgn1}
\pi_{\l}(g)=\rmid
\quad
(g\in G_n^{\parallel},\
n\in (N_{\rmpr}^+)^{>\l}).
\end{align}

\begin{thm}[{\cite[Prop.~1.20]{Gross14}}]
\label{3thm:bij2}
The map $\bbpsi:G \rightarrow
S^{\parallel}$ is a bijection.
\end{thm}

Our proof is based on the following lemma.
\begin{lem}[{cf. \cite[Appendix C.1]{Gross14}}]
\label{3lem:center1}
\par
(a).
For the projection $\pi_{\l+1,\l}: G^{\leq \l+1} \rightarrow G^{\leq \l}$ in \eqref{3eq:kk1},
$\Ker\, \pi_{\l +1,\l}$ is in the center of $G^{\leq \l+1}$,
and  it is isomorphic   (as a group) to the direct product of groups,
\begin{align}
\label{3eq:abel1}
\Ker\, \pi_{\l+1,\l}
\simeq
\prod_{n\in (N^+)^{ \l+1}}
G_n,
\quad
G_n=\exp (\frakg_n),
\end{align}
where the isomorphism is given by $\pi_{\l+1}(\prod_{n} g_n)\mapsfrom (g_n)$.
\par
(b). $G^{\leq 1}$
is an abelian group, and it   is isomorphic   (as a group) to the direct product
of groups,
\begin{align}
\label{3eq:abel2}
G^{\leq 1}
\simeq
\prod_{n\in (N^+)^{\leq 1}}
G_n,
\end{align}
where the isomorphism is given by $\pi_{1}(\prod_{n} g_n)\mapsfrom (g_n)$.
\end{lem}
\begin{proof}
We work with the corresponding Lie algebras.
Let $\pi_{\l+1,\l}: \frakg^{\leq \l+1} \rightarrow \frakg^{\leq \l}$ be the 
canonical projection.
Then, as  a vector space, we have
\begin{align}
\label{3eq:ker1}
\Ker\, \pi_{\l+1,\l}
\simeq
\bigoplus_{n\in (N^+)^{ \l+1}}
\frakg_n.
\end{align}
Moreover,
$\Ker\, \pi_{\l+1,\l }$ is in the center of $\frakg^{\leq \l+1} $ by 
\eqref{3eq:Xbra1}.
Therefore, the isomorphism \eqref{3eq:ker1} holds also as a  Lie algebra,
where the right hand side is viewed as the direct sum of Lie algebras.
The proof of (b) is similar.
\end{proof}

\begin{proof}[Proof of Theorem \ref{3thm:bij2}]
{\bf Surjectivity of $\bbpsi$.}
Let $ \bfh=(h_{n})_{n\in N_{\mathrm{pr}}^+}$ be a given  element in $S^{\parallel}$.
Let us construct an element $g\in G$ 
  inductively with respect to degree $\l$
 such that $\bbpsi(g)=\bfh$,
  that is,
 $g_{n}^{\parallel}=h_{n}$ ($n\in N_{\rmpr}^+$) holds.
 To show it, we recursively construct elements $g_{\l}\in G$ ($\l=1,2,\dots$)
 such that
 \begin{align}
 \label{3eq:psiconst1}
 \pi_\l((g_\l)_{n}^{\parallel})&=\pi_\l(h_{n})
\quad
(n\in (N_{\mathrm{pr}}^+)^{\leq \l}),
\\
 \label{3eq:psiconst2}
 \pi_\l(g_{\l+1})
 &=
  \pi_\l(g_{\l}).
\end{align}
By \eqref{3eq:psiconst2},
the limit $g=\lim_{\l \rightarrow \infty} g_{\l}\in G$ is well-defined.
Moreover, by \eqref{3eq:psiconst1}, $g$ has the desired property.

 Step 1.  We first consider the case $\l=1$.
We define
\begin{align}
g_1= \prod_{n \in  (N_{\mathrm{pr}}^+)^{\leq 1}}
h_{n} \in G,
\end{align}
where the order of the (finite) product is arbitrarily chosen,
and
the element $g_1$ depends on the order of the product.
Then, we have
\begin{align}
\label{3eq:pig1}
\pi_1(g_1)= \prod_{n \in  (N_{\mathrm{pr}}^+)^{\leq 1}}
\pi_1(h_{n}),
\end{align}
where factors $\pi_1(h_{n})$ commute with each other,
because $ G^{\leq 1}$ is abelian by Lemma \ref{3lem:center1} (b).
It follows that
\begin{align}
\pi_1((g_1)_{n}^{\parallel})=\pi_1(h_{n})
\quad
(n \in  (N_{\mathrm{pr}}^+)^{\leq 1}).
\end{align}

Step 2.
Next, for some $\l>0$, suppose that we have $g_\l\in G$
such that
\begin{align}
\label{3eq:pig2}
\pi_\l((g_\l)_{n}^{\parallel})=\pi_\l(h_{n})
\quad
(n\in (N_{\mathrm{pr}}^+)^{\leq \l})
\end{align}
holds.
For any $n \in  (N_{\mathrm{pr}}^+)^{\leq \l+1}$,
let $h'_{n}\in  G_{n}^{\parallel}$
be the one such that
\begin{align}
(g_\l)_{n}^{\parallel}
h'_{n}
=h_{n}.
\end{align}
We note that
\begin{align}
\label{3eq:pik1}
\pi_{\l+1}(h'_{n})\in \Ker\, \pi_{\l+1,\l}.
\end{align}
This is true for $n\in (N_{\mathrm{pr}}^+)^{\leq \l}$
by \eqref{3eq:pig2},
and  for $n\in (N_{\mathrm{pr}}^+)^{ \l+1}$
by \eqref{3eq:pgn1}.
Now, we define
\begin{align}
\label{3eq:gk1}
g_{\l+1}= g_\l
\prod_{n \in  (N_{\mathrm{pr}}^+)^{\leq \l+1}}
h'_{n},
\end{align}
where the order of the (finite) product is arbitrarily chosen,
and
the element $g_{\l+1}$ depends on the order of the product.
Then,  we have
\begin{align}
\label{3eq:gk2}
\pi_{\l+1}(g_{\l+1})=\pi_{\l+1}( g_\l)
\prod_{n \in  (N_{\mathrm{pr}}^+)^{\leq \l+1}}
\pi_{\l+1}(h'_{n}).
\end{align}
By \eqref{3eq:pik1} and Lemma \ref{3lem:center1} (a), $\pi_{\l+1}(h'_{n})$ commutes with
any elements in  $G^{\leq \l+1}$.
It follows that
\begin{align}
\pi_{\l+1}((g_{\l+1})_{n}^{\parallel})=
\pi_{\l+1}((g_\l)_{n}^{\parallel})
\pi_{\l+1}(h'_{n})
=
\pi_{\l+1}(h_{n})
\quad
(n \in  (N_{\mathrm{pr}}^+)^{\leq  \l+1}).
\end{align}
Moreover, by \eqref{3eq:pik1} and \eqref{3eq:gk2},
we have
\begin{align}
\pi_\l(g_{\l+1})
 &=
  \pi_\l(g_{\l}).
\end{align}

{\bf Injectivity of $\bbpsi$.}
For   given $g, g' \in G$, assume that
 $\bbpsi(g)=\bbpsi(g')=\bfg^{\parallel}$.
Let us show that $g=g'$; namely, we have,
 for any $\l>0$,
\begin{align}
\label{3eq:pgpg1}
\pi_\l(g)=\pi_\l (g').
\end{align}
By Lemma \ref{3lem:center1} (b) and the assumption, we have
\begin{align}
\pi_1(g)
= \prod_{n \in  (N_{\mathrm{pr}}^+)^{\leq 1}}
\pi_1(g_{n}^{\parallel})
= \prod_{n \in  (N_{\mathrm{pr}}^+)^{\leq 1}}
\pi_1(g'{}_{n}^{\parallel})
=\pi_1(g').
\end{align}
Thus, the equality \eqref{3eq:pgpg1} holds for $\l=1$.
Suppose that the equality holds  for some $\l$.
Then, $\pi_{\l+1}(g) \pi_{\l+1}(g')^{-1}\in \Ker\, \pi_{\l+1,\l}$.
So, by Lemma \ref{3lem:center1} (a),
there are some $h_{n}\in G_n$ such that
\begin{align}
\label{3eq:pik2}
\pi_{\l+1}(g)=
\pi_{\l+1}(g')
\prod_{n\in (N^+)^{ \l+1}}
\pi_{\l+1}(h_{n})
.
\end{align}
Since
 $\pi_{\l+1}(h_{n})$ are central elements in $G^{\leq \l +1}$,
 for $n\in (N^+)_{\rmpr}^{\leq \l}$ 
and a positive integer $t$ such that
$tn \in (N^+)^{ \l+1}$,
  we have
\begin{align}
\pi_{\l+1}(g_{n}^{\parallel})
=
\pi_{\l+1}(g'{}_{n}^{\parallel})
\pi_{\l+1}(h_{tn}),
\end{align}
Meanwhile,
by  the assumption,
we have
$g_{n}^{\parallel}
=
(g')_{n}^{\parallel}
$.
Therefore,
$\pi_{\l+1}(g_{n}^{\parallel})
=
\pi_{\l+1}((g')_{n}^{\parallel})$.
It follows that $\pi_{\l+1}(h_{n})=\rmid$
for any $n \in (N^+)^{ \l+1}$.
Putting it back in \eqref{3eq:pik2}, we conclude that $\pi_{\l+1}(g)=
\pi_{\l+1}(g')$.
\end{proof}

\subsection{Incoming and outgoing walls}

Let 
$ \bfh = (h_{n})_{n\in N_{\mathrm{pr}}^+}$
be an arbitrary element of 
the set
$S^{\parallel}$ in \eqref{3eq:S1}.
The bijection $\bbpsi$ in Theorem \ref{3thm:bij2}
uniquely determines the element $g(\bfh)\in G$
such that $\bbpsi(g(\bfh))=\bfh$.
In turn,
the bijection $\bbvarphi$ in Theorem \ref{3thm:bij1}
uniquely determines a consistent scattering diagram
$\frakD(g(\bfh))$ up to equivalence
such that $\bbvarphi( \frakD(g(\bfh)))=g(\bfh)$.
Let us write $\frakD(g(\bfh))$ as $\frakD(\bfh)$ for simplicity.
Below we will give an alternative description
of $\frakD(\bfh)$ directly by $\bfh$.

Let $p^*$ be the map in \eqref{3eq:p*3}.
For  $n\in N_{\rmpr}^+$, we have
\begin{align}
p^*(n)\in n^{\perp},
\end{align}
because $\langle n, p^*(n)\rangle = \{n,n\}=0$.

\begin{rem}
When $\Gamma$ is degenerate, $p^*(n)$ is not injective.
However, this fact is irrelevant to all results in Section \ref{3sec:const1}.
\end{rem}

For any  $g\in G$,
we have the following relation between
$g_0^z$ ($z\in n^{\perp}$) and $g_n^{\parallel}$
around $p^*(n)$.

\begin{lem}[{\cite[Theorem~1.21]{Gross14}}]
\label{3lem:gg2}
Let $g\in G$ and $n\in N_{\mathrm{pr}}^+$.
Then, for any  $\l>0$,
there is an open  neighborhood $U_\l\subset n^{\perp}$ of $p^* (n)$
such that 
\begin{align}
\label{3eq:g01}
g_{0}^{z}\equiv g_{n}^{\parallel}
\mod G^{>\l}
\quad
(z\in U_\l,\ \text{$z$: general}\,).
\end{align}
In other words,
\begin{align}
\lim_{\ss z\rightarrow p^*(n)
\atop \ss
z\in n^{\perp},\, z:\, \rm{general}
}
g_0^z = g_{n}^{\parallel}.
\end{align}
\end{lem}
\begin{proof}
\emph{Case 1}. Suppose that  $p^*(n)$ is general.
This occurs if and only if 
the fixed data $\Gamma$ is of rank 2 and nondegenerate.
Thus, any hyperplane $n^{\perp}$ intersects
other hyperplane $n'{}^{\perp}$
with 
$n'\in N_{\mathrm{pr}}^+$,
$n'\neq n$
only at the origin.
Then, $n^{\perp}$ is divided into two rays,
on each of which $g_0^{z}$ is constant.
It follows that,  
if  a general $z\in n^{\perp}$ and $p^*(n)\neq 0$ belong to
the same ray,
\begin{align}
g_0^{z}= g_0^{p^*(n)}
\end{align}
holds.
Therefore, it is enough to prove
\begin{align}
\label{3eq:gg1}
g_0^{p^*(n)}=  g_{n}^{\parallel}.
\end{align}
By comparing \eqref{3eq:decom1} and \eqref{3eq:decom6},
we see  that
two decompositions
\begin{align}
g&=g_+^{p^*(n)} g_0^{p^*(n)} 
g_-^{p^*(n)},\\
g&=g_+^{n} g_0^{n} 
g_-^{n}
\end{align}
coincide.
Thus, we have
\begin{align}
g_0^{n} =g_0^{p^*(n)} \in G_{n}^{\parallel}.
\end{align}
This further implies that
\begin{align}
g_0^{n} =g_{n}^{\parallel}.
\end{align}
Thus, we have the equality \eqref{3eq:gg1}.
\par
\begin{figure}
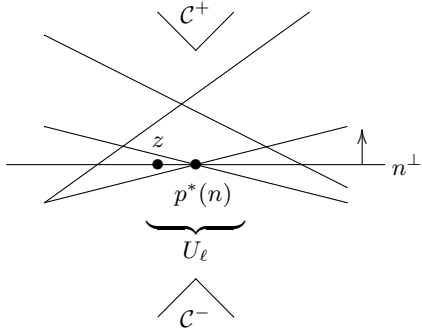

\centering
\leavevmode
\xy
(28,0)*{\text{\small $ n^{\perp}$}};
(1,-4)*{\text{\small $p^*(n)$}};
(-5,3)*{\text{\small $z$}};
(0,20)*{\text{\small $\calC^+$}};
(0,-20)*{\text{\small $\calC^-$}};
(0,-10)*{\text{\small $\underbrace{\mbox{\qquad\qquad}}_{\displaystyle U_\l}$}};
(0,0)*+{\bullet};
(-5,0)*+{\bullet};
\ar@{-} (-25,0); (25,0)
\ar@{-} (-20,-5); (15,20)
\ar@{-} (-20,17); (20,-3)
\ar@{-} (-20,5); (20,-5)
\ar@{-} (-20,-5); (20,5)
\ar@{-} (0,15); (-5,20)
\ar@{-} (0,15); (5,20)
\ar@{-} (0,-15); (-5,-20)
\ar@{-} (0,-15); (5,-20)
\ar@{->} (22,0); (22,4.5)
\endxy
\caption{Configuration
for Case 2 of the proof of Lemma \ref{3lem:gg2}.}
\label{3fig:gdecoms1}
\end{figure}

\emph{Case 2}.
Suppose that  $p^*(n)$  is not general.
Fix an integer $\ell>0$.
Then, there  might be  some (possibly no) $n'\in (N_{\mathrm{pr}}^+)^{\leq \l}$,
$n'\neq n$
such that $\langle n', p^*(n)\rangle=\{n', n\}=0$.
We choose an open neighborhood of 
$U_\l\subset n^{\perp}$
of 
$p^*(n)$ that does not intersect any $n''{}^{\perp}$ ($n''
\in (N_{\mathrm{pr}}^+)^{\leq \l}$) such that
$\langle n'', p^*(n)\rangle = \{n'',n\}\neq 0$.
See Figure \ref{3fig:gdecoms1}.
Consider the decomposition of $g$ by $n$
\begin{align}
g&=g_+^{n} g_0^{n} 
g_-^{n}.
\end{align}
We further decompose $g_0^{n} $ at general $z\in U_\l$,
\begin{align}
\label{3eq:ghhh1}
g_0^{n}=h_+^z h_0^z h_-^z
\quad
(h_+^z\in G_+^z,\
h_0^z\in G_0^z= G_{n}^{\parallel} ,\
h_-^z\in G_-^z).
\end{align}
Then, $h_{\pm}^z \in G_{n}^{\perp}$ and $h_{0}^z \in G_{n}^{\parallel}$.
 Therefore, by applying the group homomorphism 
 in \eqref{3eq:Gdecom3} to \eqref{3eq:ghhh1},
we obtain 
\begin{align}
\label{3eq:hg1}
h_0^z=g_{n}^{\parallel}.
\end{align}
On the other hand, by  the assumption on $U_\l$,
we have $\pi_{\ell}(g_{\pm}^{n})\in (G_{\pm}^z)^{\leq \ell}$.
Thus,
\begin{align}
g&=(g_+^{n} h_+^z) h_0^z (h_-^z
g_-^{n})
\end{align}
is regarded as the decomposition of $g$ at $z$ modulo $G^{>\l}$.
Therefore, we have
\begin{align}
\label{3eq:hg2}
h_0^z \equiv  g_0^z
\mod G^{>\l}.
\end{align}
Then, by \eqref{3eq:hg1} and \eqref{3eq:hg2},   we obtain the equality \eqref{3eq:g01}.
\end{proof}

\begin{defn} \index{wall!incoming} \index{wall!outgoing} 
A wall $\bfw=(\frakd, g)_n$ 
is \emph{incoming\/} (\emph{resp.~outgoing\/}) if
\begin{align}
p^*(n)\in \frakd
\end{align}
 holds (resp.~otherwise)
for its  normal vector $n$.
\end{defn}

Now we come back to our problem.
For any $\bfh\in S^{\parallel}$,
we define a set of incoming walls,
\begin{align}
\label{3eq:in1}
\mathrm{In}(\bfh)
:=\{  (n^{\perp}, h_n)_n
\mid n\in N_{\rmpr}^+
\}.
\end{align}
For a scattering diagram $\frakD$, let $\frakD_{\rmin}$
denote the set of all incoming walls of $\frakD$.

The following theorem characterizes
a consistent scattering diagram $\frakD(\bfh)$
corresponding to $\bfh\in S^{\parallel}$  up to equivalence.
\begin{thm}[{\cite[Theorem~1.21]{Gross14}}]
\label{3thm:exist1}
For any $\bfh\in S^{\parallel}$, 
let $\frakD(\bfh)$ be a corresponding consistent scattering diagram
up to equivalence.
Then, the following facts hold.
\par
(a). There is a consistent scattering diagram $\frakD$
with $\frakD_{\rmin}=\mathrm{In}(\bfh)$ such that
$\frakD$ is equivalent to $\frakD(\bfh)$.
\par
(b).
Any  consistent scattering diagram $\frakD$
with $\frakD_{\rmin}=\mathrm{In}(\bfh)$
is equivalent to $\frakD(\bfh)$.
\end{thm}

\begin{proof}
(a).
As a representative of $\frakD(\bfh)$,
we choose  $\frakD(g)$
 in
Step 1 in the proof of the surjectivity of $\bbvarphi$
in Theorem \ref{3thm:bij1}
with $g=g(\bfh)$.
For each $n\in N_{\rmpr}^+$, up to the equivalence we rearrange
the walls of $\frakD(g)$ in $n^{\perp}$ 
 to have the desired property as follows.
For each integer $\l>0$,
the hyperplane $n^{\perp}$ is divided into
cones, say, $\frakd_1$, \dots, $\frakd_p$
such that  $\pi_{\l}(g_0^z)$ is constant
with respect to general $z$ on each cone.
Let us write the associated walls therein
as $(\frakd_1,g_1)_n$, \dots, $(\frakd_p,g_p)_n$.
Let $U_\l$ be the open neighborhood of $p^*(n)$ in $n^{\perp}$ in
Lemma \ref{3lem:gg2}.
When $p^*(n)$ is in the interior of some cone, say, $\frakd_1$,
we may assume that
$p^*(n) \in \frakd_1$
and  $p^*(n) \not\in\frakd_i$ ($i\neq 1$).
When $p^*(n)$ is on the boundary of some cone,
there are several cones that
intersect  $U_\l$.
By the property \eqref{3eq:g01},
$\pi_{\l}(g_0^z)$ is constant on these cones so that
\begin{align}
\pi_{\l}(g_0^z)\equiv g_n^{\parallel}= h_n
\quad
\mod G^{>\ell}.
\end{align}
Therefore, the attached wall element in \eqref{3eq:gd1}
is also constant.
Thus, one can join these walls into a single wall.
After renumbering cones,
  we may assume that $p^*(n)\in \frakd_1$
and $p^*(n)\not\in \frakd_i $ ($i\neq 1$).
Then, we do the following:
\begin{itemize}
\item
Remove the incoming wall $(\frakd_1,g_1)_n$.
\item
Replace the outgoing walls $(\frakd_i,g_i)_n$ with $(\frakd_i,g_ig_1^{-1})_n$
for $i\neq 1$.
\end{itemize}
We do this for all $n\in N_{\rmpr}^+$ and $\l>0$,
so that all walls 
that are left
are outgoing.
Then, finally, we add walls of $\mathrm{In}(\bfh)$. 
The resulting scattering diagram $\frakD$ is equivalent to $\frakD(\bfh)$,
and it has the desired property.
\par
(b).
Suppose that $\frakD$ is a  consistent scattering diagram
with $\frakD_{\rmin}=\mathrm{In}(\bfh)$.
Let $g=g(\frakD)$. 
By Theorems 
\ref{3thm:bij1} and \ref{3thm:bij2}, it is enough to prove
\begin{align}
\label{3eq:in4}
g_{n}^{\parallel}=h_n
\quad
(n  \in N_{\rmpr}^+).
\end{align}
Let us fix $n  \in N_{\rmpr}^+$.
There is an incoming wall $\bfw_n=(n^{\perp}, h_n)_n$,
and all other walls with the normal vector $n$ do not contain $p^*(n)$.
Then,  for each integer $\l>0$,
there is a neighborhood $V_{\l}$ of $p^*(n)$ in $n^{\perp}$
that does not intersect any outgoing walls of $\frakD_{\l}$
with the normal vector $n$.
Thus, for any general $z\in V_{\l}$,
\begin{align}
g_0^z=g_z(\frakD)
\equiv 
h_n
\mod
G^{>\l}.
\end{align}
It follows that
\begin{align}
\lim_{\ss z\rightarrow p^*(n)
\atop
\ss
z\in n^{\perp},\ z :\, \rm{general}
}
g_0^z 
= h_n.
\end{align}
Meanwhile,
by Lemma \ref{3lem:gg2},
 the left hand side coincides with
$g_n^{\parallel}$.
Thus, the equality \eqref{3eq:in4}  holds.
\end{proof}

\subsection{Cluster scattering diagrams}

Now we are ready to introduce a special class of consistent scattering diagrams,
which is the main subject of Part III.

\begin{defn}[Cluster scattering diagram (CSD)]
\index{cluster scattering diagram (CSD)}\index{scattering diagram!cluster}
\label{3defn:CSD1}
For a given seed $\fraks$, 
we define a set of incoming walls
\begin{align}
\label{3eq:in2}
\mathrm{In}_{\fraks}:=&\ \{  \bfw_{e_i}:=(e_i^{\perp}, \Psi[e_i]^{\d_i})_{e_i}
\mid i=1,\dots,r\}.
\end{align}
Then, a \emph{cluster scattering diagram\/} (CSD for short)
for $\fraks$
is  a consistent scattering diagram $\frakD_{\fraks}$ 
with  $(\frakD_{\fraks})_{\rmin}=\mathrm{In}_{\fraks}$.
 \end{defn}
 
 Note that the wall elements in \eqref{3eq:in2}
 is in the normalized form in \eqref{3eq:inv1}
 so that this definition is invariant under the rescaling
 of the fixed data $\Gamma$.

The following theorem is the first fundamental result on CSDs.

\begin{thm}[Existence/uniqueness of CSD {\cite[Theorem~1.21]{Gross14}}]
\label{3thm:exist2}
For any seed $\fraks$, a CSD $\frakD_{\fraks}$ for $\fraks$ uniquely exists
up to equivalence.
\end{thm}
\begin{proof}
For a given seed $\fraks$, let $\bfh_{\fraks}=(h_{\fraks,n})\in S^{\parallel}$
be the one with
\begin{align}
\label{3eq:hs1}
h_{\fraks,n}=
\begin{cases}
\Psi[e_i]^{\d_i}& n =e_i\ (i=1,\dots,r),\\
\rmid & \text{other $n \in N_{\rmpr}^+.$}
\end{cases}
\end{align}
Then, by Theorem \ref{3thm:exist1},
there is a unique  consistent scattering diagram $\frakD$
with $\frakD_{\rmin}=\mathrm{In}(\bfh_{\fraks})$
up to equivalence.
Then, by removing all trivial incoming walls $(n^{\perp}, \rmid)_n$
($n\neq e_1,\dots, e_r$),
we obtain $\frakD_{\fraks}$.
\end{proof}

\begin{defn}[Cluster element] \index{cluster element}
For the above $\bfh_{\fraks}$,
the corresponding element $g_{\fraks}:=g(\bfh_{\fraks})\in G$ 
is called the \emph{cluster element\/} for a seed $\fraks$.
\end{defn}

\begin{ex}
All examples of consistent scattering diagrams
in Section \ref{3subsec:rank2} are CSDs $\frakD_{\fraks}$
with the cluster element  $g_{\fraks}=\Psi[\bfe_2]^{\d_2}\Psi[\bfe_1]^{\d_1}$.
By \eqref{3eq:pej1}, we have
\begin{align}
p^*(e_j) =\sum_{i=1}^2 b_{ij} f_i.
\end{align}
Thus, in the identification $N\simeq \bbZ^2$ and $M_{\bbR}\simeq \bbR^2$ therein, $p^*(\bfn)=B\bfn = (-\d_1 n_1, \d_2 n_2)$,
which is in the second quadrant of $M_{\bbR}\simeq \bbR^2$.
It follows that the rays in the  fourth quadrant are all outgoing walls.
They are of \emph{finite type}, that is, one can choose a representative of $\frakD_{\fraks}$
 with only  finitely many walls.
In fact, they are the only  rank 2 CSDs
of finite type
up to equivalence.
\end{ex}

\subsection{Constancy of  total wall-element on $e_i^{\perp}$}
\label{3subsec:constancy1}

Here we present some  property of walls with the normal vector $e_i$
for any consistent scattering diagram.
For this purpose, we introduce  yet  another decomposition of $G$.

For each  $e_i $ ($i=1,\dots,r$),
we introduce the  \emph{decomposition of $\frakg$ by $e_i$} as
\begin{align}
\label{3eq:decom11}
\frakg=
\frakg_{e_i}^{\parallel}
\oplus
\frakg_{e_i}^{\nparallel}
,
\end{align}
where $\frakg_{e_i}^{\nparallel}$
is a Lie subalgebra
of $\frakg$ defined by
\begin{align}
\frakg_{e_i}^{\nparallel}
:=\bigoplus_{\ss n\in N^+ \atop  \ss n\notin \bbZ e_i  }\frakg_{n}
=
\frakg_+^{e_i}
\oplus
\frakg_-^{e_i}
\oplus
\frakg_{e_i}^{\perp}.
\end{align}

\begin{lem}
\label{3lem:gideal2}
 $\frakg_{e_i}^{\nparallel}$ is
 a Lie algebra ideal of $\frakg$.
 \end{lem}
 \begin{proof}
 This follows from the following simple fact:
For $n,n'\in N^+$, we have $n+n'\in \bbZ e_i$ if and only if $n, n'\in \bbZ e_i$.
 \end{proof}
 
 \begin{rem}
 On the contrary, for  $n\in N_{\mathrm{pr}}^+$
 other than $e_1$, \dots, $e_r$,
 the Lie subalgebra
 $\frakg_{n}^{\nparallel}$ defined similarly
 is not a  Lie algebra ideal of $\frakg$.
 \end{rem}

Let 
$\widehat\frakg_{e_i}^{\nparallel}$
be the completion of $\frakg_{e_i}^{\nparallel}$,
and let $G_{e_i}^{\nparallel}$ be the corresponding subgroup of $G$.
We have
the {\em decomposition of $G$ by $e_i$} and 
 a parallel result to
Lemma \ref{3lem:Gdecom3}. \index{decomposition of $G$!by $e_i$}

\begin{lem}
\label{3lem:Gdecom4}
The following facts hold.
\par
(a). $G_{e_i}^{\nparallel}$ is a normal subgroup of $G$.
\par
(b).
We have
\begin{align}
G=G_{e_i}^{\parallel}G_{e_i}^{\nparallel}
\end{align}
so that
any element $g$ of $G$ is
uniquely factorized as
\begin{align}
g=g_{e_i}^{\parallel} g_{e_i}^{\nparallel}
\quad
(g_{e_i}^{\parallel}\in G_{e_i}^{\parallel},
\
g_{e_i}^{\nparallel}\in G_{e_i}^{\nparallel}).
\end{align}
\par
(c).
The map 
\begin{align}
\label{3eq:Gdecom4}
\chi_{e_i}:G\rightarrow G_{e_i}^{\parallel},
\quad
g=g_{e_i}^{\parallel} g_{e_i}^{\nparallel}
\mapsto
g_{e_i}^{\parallel}
\end{align}
is a group homomorphism with the  kernel $G_{e_i}^{\nparallel}$.
\par
(d). The map $\chi_{e_i} $ in \eqref{3eq:Gdecom4}
coincides with the map $\psi_{e_i}$  in \eqref{3eq:pro1}.

\end{lem}
\begin{proof}
The proofs of (a)--(c) are the same as Lemma \ref{3lem:Gdecom3}.
Let us prove (d).
Let $g=g_+^{e_i}g_{e_i}^{\parallel}g_{e_i}^{\perp}g_-^{e_i}$
be the decomposition in Lemmas \ref{3lem:Gdecom2}
and \ref{3lem:Gdecom3}.
By applying the map $\chi_{e_i}$ to this equality, we
obtain
$\chi_{e_i}(g)=g_{e_i}^{\parallel}=\psi_{e_i}(g)$.
\end{proof}

We have the following  consequence of Lemma \ref{3lem:Gdecom4}.
\begin{prop}
\label{3prop:ei1}
(a).
Let  $g$ be any element of  $G$.
Then, 
for any general $z\in e_i^{\perp}$,
the equality
\begin{align}
\label{3eq:g0xg1}
g_0^z=g_{e_i}^{\parallel}
\end{align}
holds. In particular, $g_0^z$ is constant with respect to general  $z\in e_i^{\perp}$.
\par
(b).
Let $\frakD$ be  any consistent scattering diagram,
and let $g=g(\frakD)$ be the corresponding element in $G$.
Then,
for any general $z\in e_i^{\perp}$,
the equality
\begin{align}
g_z(\frakD)=g_{e_i}^{\parallel}
\end{align}
holds.
In particular, $g_z({\frakD})$ is constant with respect to general  $z\in e_i^{\perp}$.
\end{prop}
\begin{proof}
(a).
Let $g$ be any element of $G$.
For any general $z\in e_i^{\perp}$,
we have the decomposition of $g$ at $z$ in \eqref{3eq:fac1},
\begin{align}
g=g_+^z g_0^z 
g_-^z
\quad
(g_+^z\in G_+^z,\ g_0^z \in G_{e_i}^{\parallel},\ g_-^z\in G_-^z).
\end{align}
Then,  applying the map $\chi_{e_i}$ to this equality, we
obtain
$\chi_{e_i}(g)=g_{0}^{z}$.
\
Therefore, $g_{0}^{z}=g_{e_i}^{\parallel}$.
\par
(b). This follows from (a) and Proposition \ref{3prop:gxg1}.
\end{proof}

In particular, we obtain the following property of a CSD,
which will be used later.

\begin{prop}[{\cite[Theorems~1.28, Remark~1.29]{Gross14}}]
\label{3prop:exist3}
There is a CSD $\frakD_{\fraks}$ 
such that the normal vector of any outgoing wall is other than $e_1$, \dots, $e_r$.
\end{prop}
\begin{proof}
Take any CSD $\frakD_{\fraks}$.
The total wall element $g_z(\frakD_{\fraks})$ on $e_i^{\perp}$ 
of $\frakD_{\fraks}$ is constant by Proposition  \ref{3prop:ei1}.
Since the wall $\bfw_{e_i}$ in  \eqref{3eq:in2} is the only incoming wall
with the normal vector $e_i$,
$g_z(\frakD_{\fraks})$ on $e_i^{\perp}$  coincides with $\Psi[e_i]^{\d_i}$
by Lemma \ref{3lem:gg2}.
This means that
the total contribution to $g_z(\frakD_{\fraks})$ on $e_i^{\perp}$
 from all outgoing walls is trivial.
 Therefore, one can remove all outgoing walls whose normal vectors
are $e_1, \dots, e_r$ up to equivalence.
\end{proof}

\subsection{Rank 2 examples: infinite type}
\label{3subsec:rank2inf}

We have already seen
all rank 2 CSDs of finite type 
in Section \ref{3subsec:rank2}.
Here we  consider rank 2  CSDs of infinite type,
where $\d_1 \d_2 \geq 4$ for $B$ in \eqref{3eq:B1}.

\smallskip
\par\noindent 
 {\bf Affine type.}
The case $\d_1\d_2=4$ is exceptional, and it is said to be  \emph{of affine type}.
There are two cases.

\begin{figure}
\centering
\leavevmode
\begin{xy}
0;/r1.2mm/:,
(0,-14)*{\text{(a) $A_1^{(1)}$}},
(8.7,-10)*{\cdot},
(9.2,-10)*{\cdot},
(10,-8.7)*{\cdot},
(10,-9.2)*{\cdot},
(0,0)="A"
\ar "A"+(0,0); "A"+(10,0)
\ar "A"+(0,0); "A"+(0,10)
\ar@{-} "A"+(0,0); "A"+(-10,0)
\ar@{-} "A"+(0,0); "A"+(0,-10)
\ar@{-} "A"+(0,0); "A"+(5,-10)
\ar@{-} "A"+(0,0); "A"+(6.66,-10)
\ar@{-} "A"+(0,0); "A"+(7.5,-10)
\ar@{-} "A"+(0,0); "A"+(8,-10)
\ar@{-} "A"+(0,0); "A"+(10,-5)
\ar@{-} "A"+(0,0); "A"+(10,-6.66)
\ar@{-} "A"+(0,0); "A"+(10,-7.5)
\ar@{-} "A"+(0,0); "A"+(10,-8)
\ar@{-} "A"+(0,0); "A"+(10,-10)
\end{xy}
\hskip40pt
\begin{xy}
0;/r1.2mm/:,
(0,-14)*{\text{(b) $A_2^{(2)}$}},
(4.1,-10)*{\cdot},
(4.6,-10)*{\cdot},
(5.5,-10)*{\cdot},
(6,-10)*{\cdot},
(0,0)="A"
\ar "A"+(0,0); "A"+(10,0)
\ar "A"+(0,0); "A"+(0,10)
\ar@{-} "A"+(0,0); "A"+(-10,0)
\ar@{-} "A"+(0,0); "A"+(0,-10)
\ar@{-} "A"+(0,0); "A"+(2.5,-10)
\ar@{-} "A"+(0,0); "A"+(3.3,-10)
\ar@{-} "A"+(0,0); "A"+(3.75,-10)
\ar@{-} "A"+(0,0); "A"+(10,-10)
\ar@{-} "A"+(0,0); "A"+(7.5,-10)
\ar@{-} "A"+(0,0); "A"+(6.6,-10)
\ar@{-} "A"+(0,0); "A"+(5,-10)
\end{xy}
\hskip35pt
\begin{xy}
0;/r1.2mm/:,
(0,-14)*{\text{(c) non-affine type}},
(0,0)="A"
\ar@{-} (5,-10); (6.3,-8.7)
\ar@{-} (3,-10); (5.5,-7.5)
\ar@{-} (2.4,-8.6); (4.5,-6.5)
\ar@{-} (1.8,-7.2); (3.8,-5.2)
\ar@{-} (1.6,-5.4); (3,-4)
\ar "A"+(0,0); "A"+(10,0)
\ar "A"+(0,0); "A"+(0,10)
\ar@{-} "A"+(0,0); "A"+(-10,0)
\ar@{-} "A"+(0,0); "A"+(0,-10)
\ar@{-} "A"+(0,0); "A"+(2,-10)
\ar@{-} "A"+(0,0); "A"+(2.5,-10)
\ar@{-} "A"+(0,0); "A"+(2.66,-10)
\ar@{-} "A"+(0,0); "A"+(2.72,-10)
\ar@{-} "A"+(0,0); "A"+(2.76,-10)
\ar@{-} "A"+(0,0); "A"+(10,-10)
\ar@{-} "A"+(0,0); "A"+(8,-10)
\ar@{-} "A"+(0,0); "A"+(7.5,-10)
\ar@{-} "A"+(0,0); "A"+(7.33,-10)
\ar@{-} "A"+(0,0); "A"+(7.27,-10)
\ar@{-} "A"+(0,0); "A"+(7.23,-10)
\end{xy}
\caption{Rank 2 CSD: infinite type}
\label{3fig:scat3}
\end{figure}
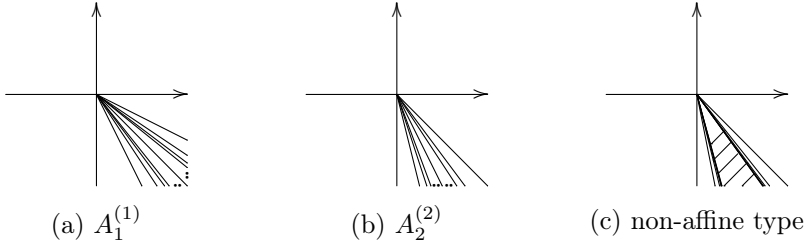

(a). Type $A_1^{(1)}$. Let $\d_1=\d_2=2$.
Note that, for any  primitive $n \in N_{\rmpr}^+$,
$\d(n)=2$.
The following description of $\frakD_{\fraks}$ is well known \cite{Reineke08,Gross09,Reading18}:
The walls of $\frakD_{\fraks}$ are given by
\begin{gather}
\label{3eq:a111}
(\bfe_1^{\perp}, \Psi[\bfe_1]^2)_{\bfe_1},
\quad
(\bfe_2^{\perp}, \Psi[\bfe_2]^2)_{\bfe_2},
\\
\label{3eq:a112}
(\bbR_{\geq 0} (p,-p-1), \Psi[(p+1,p)]^2)_{(p+1,p)}
\quad
(p\in \bbZ_{>0}),
\\
\label{3eq:a113}
(\bbR_{\geq 0} (p+1,-p), \Psi[(p,p+1)]^2)_{(p,p+1)}
\quad
(p\in \bbZ_{>0}),
\\
\label{3eq:a114}
(\bbR_{\geq 0} (1,-1), \prod_{j=0}^{\infty} \Psi[2^j\bfn_0]^{2^{2-j}})_{\bfn_0},
\end{gather}
where  $\bfn_0=(1,1)$.
Note that $\d(2^j \bfn_0)=2^{1-j}$
by \eqref{3eq:d01}.
See Figure \ref{3fig:scat3} (a).

Here we demonstrate  to derive the above walls, up to some degree, by repeatedly applying the pentagon relation
\eqref{3eq:pent1}
in a similar manner
for finite type.
We use the same notation $[n_1,n_2]^c$
for $\Psi[(n_1,n_2)]^c$
as in  Section \ref{3subsec:rank2}.
We apply the pentagon relation 
to the cluster element $g_{\fraks}=\Psi[\bfe_2]^2 \Psi[\bfe_1]^2$
repeatedly
and
obtain
the equality
\begin{align}
\label{3eq:com2}
\begin{bmatrix}
0\\
1
\end{bmatrix}
^2
\begin{bmatrix}
1\\
0
\end{bmatrix}
^2
&=
\begin{bmatrix}
1\\
0
\end{bmatrix}
\begin{bmatrix}
1\\
1
\end{bmatrix}
^2
\begin{bmatrix}
1\\
2
\end{bmatrix}
\begin{bmatrix}
1\\
0
\end{bmatrix}
\begin{bmatrix}
1\\
1
\end{bmatrix}
^2
\begin{bmatrix}
1\\
2
\end{bmatrix}
\begin{bmatrix}
0\\
1
\end{bmatrix}
^2.
\end{align}
At this moment, the expression in the right hand side
is not yet ``ordered'' to be presented by
a scattering diagram.
Here, we say that 
the above product is \emph{ordered\/} (resp. \emph{anti-ordered\/}) \index{ordered/anti-ordered}
if, for  any adjacent pair $[\bfn']^{c'}[\bfn]^{c}$,
$\{\bfn',\bfn\}=n'_2 n_1- n'_1 n_2\leq 0$ (resp. $\{\bfn',\bfn\}\geq 0$) holds.
Equivalently, 
if we view $[n_1,n_2]$ as a fraction
$n_1/n_2$,
then, the numbers should be aligned in the decreasing order
(resp. increasing order)
form left to right.
The left hand side of \eqref{3eq:com2} is anti-ordered.
To make the right hand side of \eqref{3eq:com2} ordered, we need to interchange
$[1,2]$ and $[1,0]$ in the middle,
where
$\{ \bfn',\bfn\}=2$ for $\bfn'=(1,2)$ and $\bfn=(1,0)$.
As the lowest approximation, we consider modulo
$G^{>3}$.
Then, $[1,2]$ and $[1,0]$ commute, and we have,
modulo $G^{>3}$,
\begin{align}
\begin{bmatrix}
0\\
1
\end{bmatrix}
^2
\begin{bmatrix}
1\\
0
\end{bmatrix}
^2
&
\equiv
\begin{bmatrix}
1\\
0
\end{bmatrix}
^2
\begin{bmatrix}
2\\
1
\end{bmatrix}
^2
\begin{bmatrix}
1\\
1
\end{bmatrix}
^4
\begin{bmatrix}
1\\
2
\end{bmatrix}
^2
\begin{bmatrix}
0\\
1
\end{bmatrix}
^2.
\end{align}
This certainly agrees with 
the data in \eqref{3eq:a111}--\eqref{3eq:a114}.
This is also compatible with the formula \eqref{3eq:gxg2}.

To proceed to higher degree, we
now apply the pentagon relation  \eqref{3eq:pent1}
with $c=2$ to the pair $[1,2]^{1/2}$ and $[1,0]^{1/2}$
in \eqref{3eq:com2},
which are obtained by factorizing the elements $[1,2]$ and $[1,0]$.
To clarify the structure among fractional powers,
we use the notation $[n_1,n_2]_2:=[n_1,n_2]^{1/2}$.
Then, we have
\begin{align}
\label{3eq:com3}
\begin{bmatrix}
1\\
2
\end{bmatrix}
\begin{bmatrix}
1\\
0
\end{bmatrix}
=
\begin{bmatrix}
1\\
2
\end{bmatrix}
_2^2
\begin{bmatrix}
1\\
0
\end{bmatrix}
_2^2
&=
\begin{bmatrix}
1\\
0
\end{bmatrix}
_2
\begin{bmatrix}
2\\
2
\end{bmatrix}
_2
^2
\begin{bmatrix}
3\\
4
\end{bmatrix}
_2
\begin{bmatrix}
1\\
0
\end{bmatrix}
_2
\begin{bmatrix}
2\\
2
\end{bmatrix}
_2
^2
\begin{bmatrix}
3\\
4
\end{bmatrix}
_2
\begin{bmatrix}
1\\
2
\end{bmatrix}
_2
^2.
\end{align}
This is parallel to \eqref{3eq:com2}.
As the next approximation, we consider modulo
$G^{>7}$.
Then, $[3,4]$ and $[1,0]$ commute, and we have,
modulo $G^{>7}$,
\begin{align}
\label{3eq:fact1}
\begin{bmatrix}
1\\
2
\end{bmatrix}
\begin{bmatrix}
1\\
0
\end{bmatrix}
&
\equiv
\begin{bmatrix}
1\\
0
\end{bmatrix}
_2
^2
\begin{bmatrix}
3\\
2
\end{bmatrix}
_2
^2
\begin{bmatrix}
2\\
2
\end{bmatrix}
_2
^4
\begin{bmatrix}
3\\
4
\end{bmatrix}
_2
^2
\begin{bmatrix}
1\\
2
\end{bmatrix}
_2
^2
=
\begin{bmatrix}
1\\
0
\end{bmatrix}
\begin{bmatrix}
3\\
2
\end{bmatrix}
\begin{bmatrix}
2\\
2
\end{bmatrix}
^2
\begin{bmatrix}
3\\
4
\end{bmatrix}
\begin{bmatrix}
1\\
2
\end{bmatrix}
.
\end{align}
Then, we plug it into
\eqref{3eq:com2}
and apply the pentagon relation,
and
we have,
modulo $G^{>7}$,
\begin{align}
\label{3eq:a11mod7}
\begin{bmatrix}
0\\
1
\end{bmatrix}
^2
\begin{bmatrix}
1\\
0
\end{bmatrix}
^2
&
\equiv
\begin{bmatrix}
1\\
0
\end{bmatrix}
^2
\begin{bmatrix}
2\\
1
\end{bmatrix}
^2
\begin{bmatrix}
3\\
2
\end{bmatrix}
^2
\begin{bmatrix}
4\\
3
\end{bmatrix}
^2
\begin{bmatrix}
1\\
1
\end{bmatrix}
^4
\begin{bmatrix}
2\\
2
\end{bmatrix}
^2
\begin{bmatrix}
3\\
4
\end{bmatrix}
^2
\begin{bmatrix}
2\\
3
\end{bmatrix}
^2
\begin{bmatrix}
1\\
2
\end{bmatrix}
^2
\begin{bmatrix}
0\\
1
\end{bmatrix}
^2.
\end{align}
Again,
this  agrees with 
the data in \eqref{3eq:a111}--\eqref{3eq:a114}.
By repeating this procedure modulo  $G^{>2^{\l}-1}$,
we obtain the equality
\begin{align}
\label{3eq:a115}
\begin{bmatrix}
0\\
1
\end{bmatrix}
^2
\begin{bmatrix}
1\\
0
\end{bmatrix}
^2
&
=
\begin{bmatrix}
1\\
0
\end{bmatrix}
^2
\begin{bmatrix}
2\\
1
\end{bmatrix}
^2
\begin{bmatrix}
3\\
2
\end{bmatrix}
^2
\cdots
\prod_{j=0}^{\infty}
\begin{bmatrix}
2^j\\
2^j
\end{bmatrix}
^{2^{2-j}}
\cdots
\begin{bmatrix}
2\\
3
\end{bmatrix}
^2
\begin{bmatrix}
1\\
2
\end{bmatrix}
^2
\begin{bmatrix}
0\\
1
\end{bmatrix}
^2,
\end{align}
which
is equivalent to the result \eqref{3eq:a111}--\eqref{3eq:a114}.
 A complete proof of the equality \eqref{3eq:a115} 
 based on the pentagon relation is found in 
  \cite{Matsushita21}.

\smallskip

(b). Type $A_2^{(2)}$. Let $\d_1=1$, $\d_2=4$.
The situation is parallel to $A_1^{(1)}$ though a little more complicated.
So, we  present results with less explanation.
The following description of $\frakD_{\fraks}$ is  known \cite{Reading18}:
The walls of $\frakD_{\fraks}$ are given by
\begin{gather}
\label{3eq:a221}
(\bfe_1^{\perp}, \Psi[\bfe_1])_{\bfe_1},
\quad
(\bfe_2^{\perp}, \Psi[\bfe_2]^4)_{\bfe_2},
\\
\label{3eq:a222}
(\bbR_{\geq 0} (p,-2p-1), \Psi[(2p+1,4p)])_{(2p+1,4p)}
\quad
(p\in \bbZ_{>0}),
\\
\label{3eq:a223}
(\bbR_{\geq 0} (2p-1,-4p), \Psi[(p,2p-1)]^4)_{(p,2p-1)}
\quad
(p\in \bbZ_{>0}),
\\
\label{3eq:a224}
(\bbR_{\geq 0} (p,-2p+1), \Psi[(2p-1,4p)])_{(2p-1,4p)}
\quad
(p\in \bbZ_{>0}),
\\
\label{3eq:a225}
(\bbR_{\geq 0} (2p+1,-4p), \Psi[(p,2p+1)]^4)_{(p,2p+1)}
\quad
(p\in \bbZ_{>0}),
\\
\label{3eq:a226}
\begin{split}
&\quad\ (\bbR_{\geq 0} (1,-2), \Psi[\bfn_0]^6\prod_{j=1}^{\infty} \Psi[2^j \bfn_0]^{2^{2-j}})_{\bfn_0},
\end{split}
\end{gather}
where  $\bfn_0=(1,2)$.
Note that $\d(2^j \bfn_0)=2^{1-j}$.
See Figure \ref{3fig:scat3} (b).

We apply the pentagon relation 
to the cluster element $g_{\fraks}=\Psi[\bfe_2]^4 \Psi[\bfe_1]$
repeatedly
and
obtain
the equality
\begin{align}
\label{3eq:com4}
\begin{bmatrix}
0\\
1
\end{bmatrix}
^4
\begin{bmatrix}
1\\
0
\end{bmatrix}
&=
\begin{bmatrix}
1\\
0
\end{bmatrix}
\begin{bmatrix}
1\\
1
\end{bmatrix}
^3
\begin{bmatrix}
2\\
3
\end{bmatrix}
\begin{bmatrix}
1\\
2
\end{bmatrix}
^3
\begin{bmatrix}
1\\
3
\end{bmatrix}
\begin{bmatrix}
1\\
1
\end{bmatrix}
\begin{bmatrix}
1\\
2
\end{bmatrix}
^3
\begin{bmatrix}
2\\
5
\end{bmatrix}
\begin{bmatrix}
1\\
3
\end{bmatrix}
^3
\begin{bmatrix}
1\\
4
\end{bmatrix}
\begin{bmatrix}
0\\
1
\end{bmatrix}
^4
.
\end{align}
Again, this
is not yet ordered,
and
to order it, we need to interchange
$[1,3]$ and $[1,1]$ in the middle.

As the lowest approximation, we consider modulo
$G^{>5}$.
Then, $[1,3]$ and $[1,1]$ commute, and we have,
modulo $G^{>5}$,
\begin{align}
\begin{bmatrix}
0\\
1
\end{bmatrix}
^4
\begin{bmatrix}
1\\
0
\end{bmatrix}
&
\equiv
\begin{bmatrix}
1\\
0
\end{bmatrix}
\begin{bmatrix}
1\\
1
\end{bmatrix}
^4
\begin{bmatrix}
2\\
3
\end{bmatrix}
^4
\begin{bmatrix}
1\\
2
\end{bmatrix}
^6
\begin{bmatrix}
1\\
3
\end{bmatrix}
^4
\begin{bmatrix}
1\\
4
\end{bmatrix}
\begin{bmatrix}
0\\
1
\end{bmatrix}
^4
.
\end{align}
This certainly agrees with 
the data in \eqref{3eq:a221}--\eqref{3eq:a226}.

To proceed to higher degree,
we do
in the same way as before.
By apply the pentagon relation  \eqref{3eq:pent1} with $c=2$, 
we have
\begin{align}
\label{3eq:com5}
\begin{bmatrix}
1\\
3
\end{bmatrix}
\begin{bmatrix}
1\\
1
\end{bmatrix}
=
\begin{bmatrix}
1\\
3
\end{bmatrix}
_2^2
\begin{bmatrix}
1\\
1
\end{bmatrix}
_2^2
&=
\begin{bmatrix}
1\\
1
\end{bmatrix}
_2
\begin{bmatrix}
2\\
4
\end{bmatrix}
_2
^2
\begin{bmatrix}
3\\
7
\end{bmatrix}
_2
\begin{bmatrix}
1\\
1
\end{bmatrix}
_2
\begin{bmatrix}
2\\
4
\end{bmatrix}
_2
^2
\begin{bmatrix}
3\\
7
\end{bmatrix}
_2
\begin{bmatrix}
1\\
3
\end{bmatrix}
_2
^2,
\end{align}
which has a common  structure with
 \eqref{3eq:com2}.
As the next approximation, we consider modulo
$G^{>11}$.
Then, $[3,7]$ and $[1,1]$ commute, and we have,
modulo $G^{>11}$,
\begin{align}
\begin{bmatrix}
1\\
3
\end{bmatrix}
\begin{bmatrix}
1\\
1
\end{bmatrix}
&
\equiv
\begin{bmatrix}
1\\
1
\end{bmatrix}
_2
^2
\begin{bmatrix}
3\\
5
\end{bmatrix}
_2
^2
\begin{bmatrix}
2\\
4
\end{bmatrix}
_2
^4
\begin{bmatrix}
3\\
7
\end{bmatrix}
_2
^2
\begin{bmatrix}
1\\
3
\end{bmatrix}
_2
^2
=
\begin{bmatrix}
1\\
1
\end{bmatrix}
\begin{bmatrix}
3\\
5
\end{bmatrix}
\begin{bmatrix}
2\\
4
\end{bmatrix}
^2
\begin{bmatrix}
3\\
7
\end{bmatrix}
\begin{bmatrix}
1\\
3
\end{bmatrix}
.
\end{align}
Then, we plug it into
\eqref{3eq:com4}
and apply the pentagon relation,
and
we have,
modulo $G^{>11}$,
\begin{align}
\begin{split}
\begin{bmatrix}
0\\
1
\end{bmatrix}
^4
\begin{bmatrix}
1\\
0
\end{bmatrix}
&
\equiv
\begin{bmatrix}
1\\
0
\end{bmatrix}
\begin{bmatrix}
1\\
1
\end{bmatrix}
^4
\begin{bmatrix}
3\\
4
\end{bmatrix}
\begin{bmatrix}
2\\
3
\end{bmatrix}
^4
\begin{bmatrix}
3\\
5
\end{bmatrix}
^4
\begin{bmatrix}
4\\
7
\end{bmatrix}
^4
\begin{bmatrix}
1\\
2
\end{bmatrix}
^6
\begin{bmatrix}
2\\
4
\end{bmatrix}
^2
\\
&
\qquad
\times
\begin{bmatrix}
3\\
7
\end{bmatrix}
^4
\begin{bmatrix}
2\\
5
\end{bmatrix}
^4
\begin{bmatrix}
3\\
8
\end{bmatrix}
\begin{bmatrix}
1\\
3
\end{bmatrix}
^4
\begin{bmatrix}
1\\
4
\end{bmatrix}
\begin{bmatrix}
0\\
1
\end{bmatrix}
^4.
\end{split}
\end{align}
Again,
this  agrees with 
the data in \eqref{3eq:a221}--\eqref{3eq:a226}.
By repeating this procedure modulo  $G^{>3\cdot 2^{\l}-1}$,
we  obtain the equality
\begin{align}
\label{3eq:a227}
\begin{split}
&
\begin{bmatrix}
0\\
1
\end{bmatrix}
^4
\begin{bmatrix}
1\\
0
\end{bmatrix}
=
\begin{bmatrix}
1\\
0
\end{bmatrix}
\begin{bmatrix}
1\\
1
\end{bmatrix}
^4
\begin{bmatrix}
3\\
4
\end{bmatrix}
\begin{bmatrix}
2\\
3
\end{bmatrix}
^4
\begin{bmatrix}
5\\
8
\end{bmatrix}
\begin{bmatrix}
3\\
5
\end{bmatrix}
^4
\cdots
\\
&
\hskip60pt
\times
\begin{bmatrix}
1\\
2
\end{bmatrix}
^6
\prod_{j=1}^{\infty}
\begin{bmatrix}
 2^j\\
 2^{j+1}
\end{bmatrix}
^
{2^{2-j}}
\cdots
\begin{bmatrix}
5\\
12
\end{bmatrix}
\begin{bmatrix}
2\\
5
\end{bmatrix}
^4
\begin{bmatrix}
3\\
8
\end{bmatrix}
\begin{bmatrix}
1\\
3
\end{bmatrix}
^4
\begin{bmatrix}
1\\
4
\end{bmatrix}
\begin{bmatrix}
0\\
1
\end{bmatrix}
^4
,
\end{split}
\end{align}
which is equivalent to the result \eqref{3eq:a111}--\eqref{3eq:a114}.
 A complete proof of the equality \eqref{3eq:a227} 
 based on the pentagon relation is found in 
  \cite{Matsushita21}.
The case $\d_1=4$, $\d_2=1$ is similar.

\smallskip
\par\noindent 
 {\bf Non-affine type.}
The case $\d_1\d_2>4$.
The entire structure of a CSD is not yet known.
Let 
$\sigma(\bfv_+, \bfv_-)$ be the irrational
cone spanned
by two
vectors 
 \begin{align}
 \label{3eq:ds1}
 \bfv_{\pm}
 =
 \begin{pmatrix}
\d_1\d_2\pm\sqrt{\d_1\d_2(\d_1\d_2-4)} \\
- 2\d_2\\
\end{pmatrix}.
\end{align}
Under the canonical pairing \eqref{3eq:can1},
an outgoing wall $\bfw$ is in the region $\sigma(\bfv_+,\bfv_-)$
if and only if its normal vector $\bfn=(n_1,n_2)$ satisfies
the condition
\begin{align}
\label{3eq:ds2}
\frac{2}{\d_2+ \sqrt{\d_2(\d_2-4/\d_1)}}
<
\frac{n_1}{n_2}
<
\frac{2}{\d_2- \sqrt{\d_2(\d_2-4/\d_1)}}.
\end{align}
See Figure \ref{3fig:scat3} (c),
where the  cone $\sigma(\bfv, \bfv')$ is depicted as a hatched region.
The region outside $\sigma(\bfv_+, \bfv_-)$
is decomposed into infinitely many chambers
by the  walls converging
to  $\bfv_+$ and $\bfv_-$
as in the affine case.
They are the \emph{$G$-chambers\/}
in Section \ref{3subsec:chamber1},
and its behavior is well-known in  cluster algebra theory
(e.g., \cite[\S3]{Reading18}).
On the other hand,
in the region $\sigma(\bfv_+, \bfv_-)$,
 which is informally called the \emph{Badlands}, \index{Badlands}
it is expected that every ray with the normal vector
$n\in N_{\rmpr}^+$
 in the  cone $\sigma(\bfv, \bfv')$ appears
 as a nontrivial wall \cite[Example 1.15]{Gross14}.
For the skew-symmetric case $\d_1=\d_2$,
this is proved in \cite[Example 7.10]{Davison19}.
However, the wall elements are not known yet.
In Section \ref{3ex:nonaffine1},
we present    results calculated by computer,
which support the above expectation.

\notes
Most contents are  taken  from \cite[\S1.2]{Gross14}
with added/modified proofs.
In particular,
 an alternative and simple proof of
Proposition
\ref{3prop:exist3}
is given
based on  
Proposition \ref{3prop:ei1}.
The presentation of rank 2 CSDs of affine type
based on the pentagon relation  \eqref{3eq:pent1} is also due to us.
The procedure demonstrated above is
the implementation of the \emph{perturbation trick\/} and
the \emph{change of lattice  trick\/}
 described in \cite[Appendix C.3]{Gross14},
and it appears later as Algorithm \ref{3algo:ord1}
 in more generality.

\newpage
\section{Principal extension method}
\label{3sec:principal1}

In this short section we give a faithful representation of
the structure group $G$ for any (possibly degenerate) fixed data
based on the principal extension method.

\subsection{Principal extension of fixed data}

In the previous sections we have successfully constructed 
a CSD  for any (possibly degenerate) fixed data.
  On the other hand, to study further properties of CSDs,
 there is a situation  where
a faithful representation of $G$ is desired even
 for \emph{degenerate\/} fixed data.
  In  cluster algebra theory, 
this corresponds to the situation where
the algebraically independence of  \emph{$\hat{y}$-variables\/} in \cite{Fomin07}
is desirable.
To resolve it, the \emph{principal coefficients},
or equivalently,  the \emph{principal extensions\/} of exchange matrices were introduced in \cite{Fomin07}.
  Here we present a parallel method 
 in scattering diagrams by \cite[Construction 2.1]{Gross13},
 which in particular provides a faithful representation of $G$
 for any (possibly degenerate) fixed data.

\begin{defn}[Principal extension of fixed data] \index{principal extension!of a fixed data}
Let $\Gamma=\{ N, \{\cdot, \cdot\}, N^{\circ},\allowbreak (\d_i), \allowbreak M, M^{\circ}\}$
be a (possibly degenerate) fixed data of rank $r$.
The \emph{principal extension  $\Gamma_{\prin}$ of   $\Gamma$}
is a fixed data 
 of rank $2r$ consisting of the following data:
\begin{itemize}
\item
A lattice $\tilde{N}=N\oplus M^{\circ}$ of rank $2r$ with
a skew-symmetric bilinear form
\begin{gather}
\{\cdot, \cdot \}: \tilde N \times \tilde  N \rightarrow \bbQ,
\\
\label{3eq:bili3}
\{(n,m), (n',m')\}:=\{n,n'\}+\langle n',m\rangle - \langle n,m'\rangle.
\end{gather}
\item
A sublattice $\tilde{N}^{\circ}=N^{\circ}\oplus M \subset \tilde N$.
It is of finite index and
\begin{align}
\label{3eq:bili4}
\{\tilde N^{\circ}, \tilde N\}\subset \bbZ.
\end{align}
\item
A $2r$-tuple of
positive integers $(\tilde{\d}_1,\dots, \tilde{\d}_{2r})=(\d_1,\dots, \d_r, \d_1,\dots, \d_r)$.
\item
$\tilde M=\mathrm{Hom}(\tilde N,\bbZ)\simeq M\oplus N^{\circ}$, 
$\tilde M^{\circ}=\mathrm{Hom}(\tilde N^{\circ},\bbZ)\simeq M^{\circ}\oplus N$.
\end{itemize}
The bilinear form \eqref{3eq:bili3} is always nondegenerate,
so that
$\Gamma_{\prin}$ is an nondegenerate fixed data.
\end{defn}

The canonical pairing 
 $\langle \cdot, \cdot \rangle: \tilde N^{\circ}\times \tilde M^{\circ} \rightarrow \bbZ$ is given by
 \begin{align}
 \label{3eq:can2}
 \langle (n,m), (m',n')\rangle
 = \langle n,m' \rangle+ \langle n',m\rangle.
 \end{align}

\begin{ex}
\label{3ex:pseed1}
For any seed $\fraks=(e_1,\dots,e_r)$ for $\Gamma$, let $(f_1,\dots,f_r)$ be 
the basis of $M^{\circ}$ in Section \ref{3subsec:fixed1}.
Then, $\tilde\fraks=((e_1,0), \dots, (e_r, 0), (0, f_1), \dots, (0, f_r))$ is a seed
for $\Gamma_{\prin}$.
Also, $((0, f_1), \dots, (0, f_r), (e_1,0), \dots, (e_r, 0))$ is a basis of
$\tilde M^{\circ}$.

\end{ex}

\begin{rem}
In the above example,
by applying the correspondence in \eqref{3eq:cseed1},
and
using $x_i$ and $y_i$ in \eqref{3eq:cseed1},
we have
\begin{align}
\label{3eq:cseed2}
\tilde \bfx=(x_1,\dots,x_r, y_1,\dots, y_r).
\end{align}
They are regarded as $x$-variables $x_i$ with \emph{principal coefficients\/} $y_i$
in \cite{Fomin07}.
\end{rem}

\subsection{Principal $x$-representation of $G$}
\label{3subsec:principal1}

To clarify the background,
let us briefly describe 
the scheme of the construction of a CSD in \cite{Gross14}
for degenerate $\Gamma$.
For a seed $\tilde\fraks$ for $\Gamma_{\prin}$ in Example \ref{3ex:pseed1}, every result in the previous 
sections are applicable,
where the structure group for scattering diagrams
is the group $\tilde{G}$ associated with $\tilde\fraks$,
and the ambient space 
for their supports is $\tilde M_{\bbR}\simeq \bbR^{2r}$.
Moreover, since $\Gamma_{\prin}$ is nondegenerate,
$\tilde{G}$ has the faithful $x$-representation.
Then, one can construct a CSD for $\fraks$ by reducing the structure group $\tilde{G}$ to
the subgroup $G$  associated with $\fraks$.
Moreover, $G$ has  a faithful representation
induced from the faithful $x$-representation of $\tilde{G}$.
The resulting CSD is still in the ambient space $\tilde M_{\bbR}\simeq M_{\bbR}\oplus \bbR^r$,
and it is invariant along the fiber $\bbR^r$.
Then,
 we apply the projection $\pi: \tilde M_{\bbR}\rightarrow M_{\bbR}$ to the
support, and we obtain the CSD constructed here.

Since  we have already  constructed a CSD with the ambient space $M_{\bbR}\simeq \bbR^{r}$
in the previous section
 directly
 for any (possibly degenerate) fixed data,
 we do not need to  repeat the above construction of a CSD.
 Thus, our strategy here is just borrowing 
 the faithful representation of $G$
 induced from the one of $\tilde{G}$,
 and applying it to the wall elements when needed.

Let $\fraks$ be a seed for a given (possibly degenerate) fixed data $\Gamma$.
From now on, we also use the internal-sum notation for 
 $\tilde{N}=N\oplus M^{\circ}$ and $\tilde M^{\circ}=\mathrm{Hom}(\tilde N^{\circ},\bbZ)= M^{\circ}\oplus N$
 such as $\tilde n= n+m$, $\tilde m = m+n$, ($n\in N$, $m\in M^{\circ}$).
 This notation is compatible with the notation in 
 \eqref{3eq:bili3} because
$\{(n,0), (n',0)\}=\{n,n'\}$.
Also, for $n\in N$ and $\tilde m =m+n'$ ($m\in M^{\circ}$, $n'\in N$),
we have
\begin{align}
\label{3eq:bibi1}
\langle  n,\tilde m\rangle = 
\langle  n, m\rangle.
\end{align}

We introduce a group homomorphism (cf.\ \eqref{3eq:p*3}).
\begin{align}
\label{3eq:p*4}
\begin{matrix}
\tilde p^* : & N &\rightarrow &\tilde{M}^{\circ}\\
&  n& \mapsto & \{ \cdot, n \}.
\end{matrix}
\end{align}
By the definition of $\tilde p^*$ and \eqref{3eq:bili3}, we have
\begin{align}
\label{3eq:p*5}
\langle n' + m', \tilde p^*(n)\rangle
=\{n',n\} +\langle n,m'\rangle \in \bbZ
\quad
(n'\in N^{\circ}, n\in N, m'\in M).
\end{align}
\begin{lem}
\label{3lem:injp2}
The map $\tilde p^*$ is injective.
In particular, $\tilde p^*(e_1), \cdots, \tilde p^*(e_r)\in \tilde{M}^{\circ}$ are $\bbZ$-linearly independent.
\end{lem}
\begin{proof}
We take a basis 
$\fraks=(e_1,\dots, e_r)$
 of $N$ 
and  a basis
$(f_1, \dots, \allowbreak f_r, \allowbreak e_1, \dots, e_r)$
 of $\tilde M^{\circ}$.
Then,  as \eqref{3eq:ppr1},
we can easily verify that
\begin{align}
\label{3eq:pej2}
\tilde p^*(e_j) = \sum_{i=1}^r b_{ij} f_i + e_j.
\end{align}
Thus, the matrix representation of $\tilde p^*$ is given by a $2r\times r$ matrix
\begin{align}
\label{3eq:Bprin1}
\tilde B=
\begin{pmatrix}
B\\
I
\end{pmatrix},
\end{align}
where $B=(b_{ij})_{i,j=1}^r$, $b_{ij}=\{\d_ie_i, e_j\}$ as before.
The matrix $\tilde B$ has rank $r$.
\end{proof}

\begin{rem}
The matrix $\tilde B$ is the principal extension 
of an exchange matrix $B$ in \cite{Fomin07}.
\end{rem}

In place of $P$  in \eqref{3eq:P1},
we define a monoid $\tilde P\subset \tilde M^{\circ}$ as follows:
\begin{itemize}
\item[(i).]
$\tilde P=\sigma \cap \tilde M^{\circ}$,
where $\sigma $ is a $2r$-dimensional strongly convex cone in $\tilde M_{\bbR}= \tilde M^{\circ} \otimes_{\bbR}\bbR$.
\item[(ii).]
$\tilde p^*(e_1), \cdots, \tilde p^*(e_r)\in \tilde P$.
\end{itemize}
Such $\tilde P$ is not unique at all, and we choose one arbitrarily.
The results we will present here do not depend on the choice of $\tilde P$.
For example, we may take the monoid  generated by
$f_1$, \dots, $ f_r$, $\tilde p^*(e_1)$, \dots, $\tilde p^*(e_r)$,
which are a basis of $\tilde M^{\circ}$ by \eqref{3eq:pej2}.

\begin{lem}
\label{3lem:Mnd1}
For any $n\in N$ with $n\neq 0$,  there is some $\tilde m\in \tilde P$ such that
$\langle n, \tilde m \rangle \neq 0$.
\end{lem}
\begin{proof}
By the  condition (i), 
if $\langle n, \tilde m \rangle = 0$ for any $\tilde m \in \tilde P$,
$\langle n,  \tilde m \rangle =0$ for any $\tilde m \in \tilde M^{\circ}$
by linearity.
This is a contradiction.
\end{proof}
 Let $\bbk[\tilde P]$ be the monoid algebra of $\tilde P$  over $\bbk$.
Let $J$ be the maximal ideals of $\bbk[\tilde P]$
generated by $\tilde P \setminus \{0\}$,
and let 
$\bbk[[\tilde P]]$
be the completion by $J$.
We express an arbitrary element  $f\in \bbk[[\tilde P]]$
as a formal power series in $x$ as 
\begin{align}
\label{3eq:fx2}
f=
\sum_{\tilde m\in \tilde P} c_{\tilde m} x^{\tilde m}.
\end{align}
In parallel to \eqref{3eq:tildeX2},
we define the action of $X_n \in\frakg $ ($n\in N^+$) on $\bbk[[\tilde P]]$ by
 \begin{align}
\label{3eq:tildeX4}
 \tilde{X}_n (x^{\tilde m}):=
 &\ \langle n, m  \rangle x^{\tilde m+\tilde p^*(n)}
\quad
(\tilde m  \in \tilde P).
\end{align}
Here and in the related formulas below,
  $\tilde m = m+n'$ ($m\in M^{\circ}$, $n'\in N$).

\begin{prop}
\label{3prop:gaction4}
We have an injective group homomorphism
\begin{align}
\label{3eq:Xn5}
\begin{matrix}
\tilde \rho_x:&  G &\rightarrow & \mathrm{Aut}(\bbk[[\tilde P]]) \\
& \exp(X) &  \mapsto &\mathrm{Exp} (\tilde{X}).
\end{matrix}
\end{align}
\end{prop}
\begin{proof}
We only need to prove the injectivity, 
since the rest is the same with Propositions \ref{3prop:gaction2}
and \ref{3prop:gaction3}.
To show the injectivity, it is enough to prove the linear independence
of $\tilde X_n$ ($n\in N^+$) in \eqref{3eq:tildeX4}.
This follows from the injectivity of $\tilde p^*$ and Lemma 
\ref{3lem:Mnd1}.
\end{proof}

The representation $\tilde \rho_x$ is
 closely related with the mutations of $x$-variables 
 with principal coefficients
in 
\eqref{3eq:cseed2}.
Thus, we call it the \emph{principal $x$-represent\-at\-ion\/} of $G$.
 \index{principal $x$-representation}\index{$x$-representation!principal}

\subsection{Faithful representation of wall elements}
\label{3susec:faithful1}

Let $\fraks$ be any seed for
a given (possibly degenerate) fixed data.

\begin{lem}[{\cite[Lemma 1.3]{Gross14}}]
\label{3lem:gn2}
For any $n\in N^+_{\rmpr}$,
 let
\begin{align}
g=\exp\Biggl(\, \sum_{j>0} c_j X_{jn}
\Biggr)
\in G_{n}^{\parallel}
\end{align}
be an arbitrary element in $G_{n}^{\parallel}$.
Correspondingly, let
\begin{align}
f_{g}:=\exp
\Biggl(\, \sum_{j>0} j c_j x^{j\tilde p^*(n)}
\Biggr)
\in 
\bbk[[x^{\tilde p^*(n)}]]
\subset
 \bbk[[\tilde P]].
\end{align}
Then, under the principal $x$-representation $\tilde \rho_x$,   $g$ acts on   $ \bbk[[\tilde P]]$ by 
\begin{align}
\label{3eq:gp2}
g(x^{\tilde m}) = x^{\tilde m} f_g^{\langle  n, m\rangle }
\quad (\tilde m\in \tilde P).
\end{align}
\end{lem}
\begin{proof}
We have
\begin{align}
\begin{split}
g(x^{\tilde m})
&=
x^{\tilde m}
\exp
\Biggl(
\sum_{j>0} c_j x^{j\tilde p^*(n)}\langle jn,  m \rangle 
\Biggr)
\\
&=
x^{\tilde m}
\exp
\Biggl(\,
 \sum_{j>0}  j c_j x^{j\tilde p^*(n)}
\Biggr)
^{\langle n,  m \rangle }
=x^{\tilde m} f_g^{\langle  n, m\rangle }.
\end{split}
\end{align}
\end{proof}

Thanks to the injectivity of $\tilde p^*$ and $\tilde \rho_x$,
$g\in G$ is uniquely specified by $f_g$ and its action
\eqref{3eq:gp2}.
Also, if we consider the normalized form of $g$ in \eqref{3eq:inv3},
then $f_{g_0}$ corresponding to $g_0=g^{1/\d(n)}$  is
also invariant under the rescaling of $\Gamma$.
In view of this,
we introduce the following notion.
 
\begin{defn}[Normalized automorphism $\frakp_{n,f}$] \index{normalized automorphism}
\label{3defn:auto2}
For any $n\in N^+_{\rmpr}$ and 
\begin{align}
f=1+\sum_{j>0} c_j x^{j\tilde p^*(n)}
\in 
\bbk[[x^{\tilde p^*(n)}]]
\subset
\bbk[[\tilde P]],
 \end{align}
we define an algebra automorphism $\frakp_{n,f}$ of  $\bbk[[\tilde P]]$
 by
\begin{align}
\label{3eq:wca1}
\frakp_{n,f}(x^{\tilde m}):=x^{\tilde m} f^{\langle \d(n)  n, m\rangle }
\quad (\tilde m \in \tilde P)
\end{align}
We call $\frakp_{n,f}$ the \emph{normalized automorphism\/} by $(n,f)$.
\end{defn}

\begin{defn}[Wall function] \index{wall function}
\label{3defn:wallf1}
Let  $\bfw=(\frakd, g)_n$ ($n\in N_{\rmpr}^+$) be any wall for $\fraks$.
Let $g_0:=g^{1/\d(n)}$,
and  let $f=f_{g_0}$
be the  formal power series
in $x^{\tilde p^*(n)}$ in Lemma \ref{3lem:gn2}
corresponding to $g_0$.
Then, $g_{\frakd}$ acts on $\bbk[[\tilde P]]$ by $\frakp_{n,f}$.
Thus, we can specify the wall $\bfw$ also by $\bfw=[\frakd, f]_n$
without ambiguity,
where we use the square bracket to avoid confusion
with the usual notation.
We call $f$  the \emph{wall function\/} of  a wall $\bfw$.
\end{defn}

The above specification of a wall is used   in \cite{Gross14}.

The following example is particularly important.

\begin{ex}[Cf. Proposition \ref{3prop:gn1}]
\label{3ex:gn2}
For any $n\in N^+_{\rmpr}$ and $t\in \bbZ_{>0}$,
let $g=\Psi[tn]^{\d(tn)}=\Psi[tn]^{\d(n)/t}$,
where
$\Psi[tn]$ is the dilogarithm element  
in \eqref{3eq:gei1}.
Let
\begin{align}
\label{3eq:gnorm1}
g_0 :=g^{1/\d(n)} =\Psi[tn]^{1/t}=
 \exp\Biggl( \frac{1}{t} \sum_{j>0} \frac{(-1)^{j+1}}{j^2}X_{jtn}
\Biggr).
\end{align}
Then, the corresponding wall function  $f=f_{g_0}$  is given by
$f[tn]$, where
\begin{align}
\label{3eq:fei2}
f[tn]:=
\exp
\Biggl(\,
\sum_{j>0} \frac{(-1)^{j+1}}{j} x^{j\tilde p^*(tn)}
\Biggr)
=1+x^{\tilde p^*(tn)}.
\end{align}
Note that the factor $t^{-1}$ in \eqref{3eq:gnorm1} cancels the factor $t$ from $X_{jtn}$
in \eqref{3eq:fei2}.
Therefore, we have
\begin{align}
\label{3eq:gei5}
\Psi[tn]^{\d(tn)}(x^{\tilde m} )
=\frakp_{n,f[tn]}
=x^{\tilde m} (1+x^{\tilde p^*(tn)})^{\langle {\d(n)} n,  m\rangle }
.
\end{align}
This automorphism is identified with  the automorphism part of the Fock-Goncharov decomposition of mutations 
of $x$-variables with principal coefficients in cluster algebras
\cite[\S2.1]{Fock03}.
In terms of $f[tn]$,
the set of incoming walls in \eqref{3eq:in2} is  represented as
\begin{align}
\label{3eq:in3}
\mathrm{In}_{\fraks}=&\ \{  \bfw_{e_i}=[e_i^{\perp}, 1+x^{\tilde p^*(e_i)}]_{e_i}
\mid i=1,\dots,r\}.
\end{align}
This is  the description of a CSD given in \cite{Gross14}.
\end{ex}
 
 \notes
All results in this section are taken from
\cite[\S1, Appendix B]{Gross14} and \cite[\S2]{Gross13},
and adapted for the present context.

\newpage
\section{Positive realization of  CSD
and pentagon relation}

In this section we
present an alternative construction of CSDs
following \cite{Gross14}.
The exponents of wall elements for the
resulting CSDs are all positive;
thus we call them \emph{positive realizations} of CSDs.
This realization is a key to prove the Laurent positivity
of the corresponding cluster algebra.
The construction also reveals the fundamental role 
of the pentagon relation for CSDs.

\subsection{Positive realization of  CSD}

We start with a general fact on a consistent scattering diagram,
which easily follows from
Proposition   \ref{3prop:ggen1} and  the proof of Theorem \ref{3thm:bij1}.
Recall the correspondence between $\Psi[tn]^{\d(tn)}$ and $f[tn]$ in Example \ref{3ex:gn2}.

\begin{prop}
Let $\frakD$ be any consistent scattering diagram for a seed $\fraks$.
Then, there is a consistent scattering diagram $\frakD'$
that is equivalent to  $\frakD$
such that
the wall element/function of any wall $\bfw=(\frakd, g)_{n}=[\frakd, f]_{n}$
of $\frakD'$
has the following form
\begin{align}
\label{3eq:ggen1}
g=\Psi[tn]^{s \d(tn) },
\quad
f=(1+x^{\tilde p^*(tn)})^s
\quad
(n\in N_{\rmpr}^+, \  t \in \bbZ_{>0},\ s\in \bbk).
\end{align}
\end{prop}

\begin{proof}
Let $\frakD=\frakD(g)$ be a consistent scattering diagram corresponding to $g\in G$.
We construct $\frakD'$ by slightly modifying
the construction of  $\frakD(g)$ in the proof of the surjectivity of $\bbvarphi$ in Theorem \ref{3thm:bij1}.
As in the proof of Proposition \ref{3prop:ggen1}, we uniquely expand $g_0^z$ in \eqref{3eq:expand1} as 
\begin{align}
g_0^z 
=
\prod_{j>0}   \Psi[jn]^{c'_{z,jn}}
\quad
(c'_{z,jn}\in \bbk).
\end{align}
Accordingly, we place $g_{\frakd}$ in \eqref{3eq:gd1} with
\begin{align}
\label{3eq:gd3}
g_{\frakd}=
\begin{cases}
\Psi[jn]^{c'_{z,jn}}
&  \text{$\l= j \deg(n)$ for some unique integer $j>0$},\\
\rmid
& \text{otherwise}.
\end{cases}
\end{align}
Then, the resulting consistent scattering diagram has
 the property \eqref{3eq:ggen1}.
\end{proof}

The goal of this section  is to prove the following  
positivity result
on a CSD,
where the fact $s\in \bbZ_{> 0}$ is crucial.

\begin{thm}
[{\cite[Theorem~1.13]{Gross14}}]
\label{3thm:pos1}
Let $\frakD_{\fraks}$ be a CSD for a seed $\fraks$.
Then, there is a consistent scattering diagram $\frakD$
that is equivalent to  $\frakD_{\fraks}$
such that
the wall element/function of any wall $\bfw=(\frakd, g)_{n}=[\frakd, f]_{n}$
of $\frakD$
has the following form
\begin{align}
\label{3eq:gpos1}
g=\Psi[tn]^{s \d(tn) },
\quad
f=(1+x^{\tilde p^*(tn)})^s
\quad
(n\in N_{\rmpr}^+, \ s,t \in \bbZ_{>0}).
\end{align}
In particular, the wall function $f$ is a polynomial in $x$ with  positive coefficients.
\end{thm}
We call the above CSD a \emph{positive realization} of $\frakD_{\fraks}$. \index{positive realization}
By decomposing the wall in \eqref{3eq:gpos1} into $s$ walls with $g=\Psi[tn]^{ \d(tn) }$, we may assume that $s=1$ for any wall whenever necessary.

\begin{ex}
\label{3ex:wallaffine1}
One can directly confirm  the property \eqref{3eq:gpos1} for finite and affine types of rank 2
by the results presented in Sections \ref{3subsec:rank2} and  \ref{3subsec:rank2inf},
respectively.
For walls other than  \eqref{3eq:a114} and \eqref{3eq:a226},
we have $s=1$.
For the wall $ \eqref{3eq:a114}$, one can split it into infinitely many walls
so that the wall elements are given by
\begin{align}
g=\Psi[2^j \bfn_0]^{2 \d(2^j \bfn_0) }
\quad
(j=0,1,2,\dots).
\end{align}
Thus, only the case $s=2$ in \eqref{3eq:gpos1} appears.
Similarly, for the wall  \eqref{3eq:a226},
we have
\begin{align}
g=\Psi[ \bfn_0]^{3 \d(  \bfn_0)},\  \Psi[2^j  \bfn_0]^{2 \d(2^j  \bfn_0) }
\quad
(j=1,2,\dots).
\end{align}
Thus, the cases $s=2$ and $3$  in \eqref{3eq:gpos1} appear.
\end{ex}
\subsection{Ordering Algorithm}

Here we concentrate on  a fixed data  $\Gamma$ of rank 2.
The following result on dilogarithm elements 
is a key to prove   Theorem \ref{3thm:pos1}.
Let us
recall the notion of ordered and anti-ordered products
in Section \ref{3subsec:rank2inf}.
Namely, a (possibly infinite) product of
$\Psi[n]^c$'s $(n\in N^+,\ c\in \bbQ_{>0})$ is 
\emph{ordered\/} (resp. \emph{anti-ordered\/}),  \index{ordered/anti-ordered}
if, for  any adjacent pair $\Psi[n']^{c'}\Psi[n]^{c}$,
$\{n',n\}\leq 0$ (resp. $\{n',n\}\geq 0$) holds.

\begin{prop}[Ordering Lemma] \index{Ordering Lemma}
\label{3prop:rank21}
Let $\fraks$ be a seed for a fixed data $\Gamma$ of rank 2.
Let 
\begin{align}
\label{3eq:init1}
C^{\mathrm{in}}=
\Psi[t'_j n'_j]^{s'_j  \d(t'_j n'_j)}
\cdots
\Psi[t'_1n'_1]^{s'_1\d(t'_1 n'_1)}
\quad
(n'_i\in N_{\rmpr}^+, \ s'_i,t'_i \in \bbZ_{>0})
\end{align}
be any finite anti-ordered product.
 Then, $C^{\mathrm{in}}$
equals to a (possibly infinite) ordered product $C^{\mathrm{out}}$ of 
factors  of the same form
\begin{align}
\label{3eq:gpos2}
\Psi[t_i n_i]^{s_i \d(t_i n_i) }
\quad
(n_i\in N_{\rmpr}^+, \ s_i,t_i \in \bbZ_{>0}).
\end{align}
Moreover, 
after gathering all powers of a common $\Psi[t_i n_i]$,
the above product is unique up to 
the reordering of the commuting adjacent pairs with $n_i=n_{i+1}$.
\end{prop}

When $\Gamma$ is degenerate, the claim holds trivially by Proposition \ref{3prop:pent1} (a).
Thus, we  may concentrate on the case where $\Gamma$ is nondegenerate.

The uniqueness is easily shown by a routine argument with Lemma \ref{3lem:center1}.
Namely, let $\pi_{\ell}(C^{\mathrm{out}})\in G^{\leq \ell}$, which is also regarded as an ordered product.
The uniqueness of the ordered product $\pi_{1}(C^{\mathrm{out}})$ holds by Lemma \ref{3lem:center1} (b).
Suppose that the uniqueness holds for the ordered product $\pi_{\ell}(C^{\mathrm{out}})$.
Then, the uniqueness of  the ordered product $\pi_{\ell+1}(C^{\mathrm{out}})$
follows from the assumption and  Lemma \ref{3lem:center1} (a).

To prove the existance, we introduce an algorithm based on the pentagon relation \eqref{3eq:pent1},
which is a systematic generalization of the one used in Section \ref{3subsec:rank2inf}.
In the algorithm,
without loss of generality, we may assume that $\{e_2,e_1\}=1$
as before.
Let us write $\Psi[n]$ as $[n]$, for simplicity, as in Section  \ref{3subsec:rank2inf}.

\begin{defn}[$p$-exchangeable pair] \index{$p$-exchangeable pair}
We say that, for any integer $p>0$, an anti-ordered adjacent pair  $[n']^{c'}[n]^{c}$ is \emph{$p$-exchangeable\/}
if $\{n',n\}=p$ and $c'=j'/p$ and $c=j/p$ for some $j',j\in \bbZ_{>0}$.
Note that we can apply the pentagon relation \eqref{3eq:pent1} for such a pair
by decomposing them into factors $[n']^{1/p}$ and $[n]^{1/p}$.
\end{defn}

\begin{ex}
\label{3ex:padm1}
Let  $[t_1n_1]^{s_1/t_1}[t_2n_2]^{s_2/t_2}$
$(n_i\in N_{\rmpr}^+,\ s_i,t_i\in \bbZ_{>0})$
be an anti-ordered adjacent pair with $\{t_1n_1,t_2n_2\}=p$.
Then, $p=t_1t_2\{n_1,n_2\}$ is a multiple of both $t_1$ and $t_2$.
Thus, it is $p$-exchangable.
\end{ex}

Here is the algorithm with  annotation.
\begin{algo}[Ordering Algorithm] \index{Ordering Algorithm}
\label{3algo:ord1}
Fix any degree $\ell>0$.
For any finite anti-ordered  product $C^{\mathrm{in}}$
with factors of the form $[tn]^{s/t}$ $(n\in N_{\rmpr}^+, \ s,t \in \bbZ_{>0})$,
obtain a finite ordered  product $C^{\mathrm{out}}({\ell})$ of
factors of the form 
$[tn]^{s/t}$ $(n\in N_{\rmpr}^+, \ s, t \in \bbZ_{>0})$
 by the following procedure.
 (The fact $s\in \bbQ_{>0}$ is obvious from the algorithm,
 but the fact $s\in \bbZ_{>0}$ is not so,
 and it will be shown in the proof of Lemma
\ref{3lem:nofail1}.)
\begin{itemize}
\item[($p_0$).]
This is the main routine. Below $p_0:=1$.
\begin{itemize}  
\item[1.]
Set  $C_0=C^{\mathrm{in}}$. (initial data)

\item[2.]
Decompose every factor $[n]^c$ $(c\in \bbZ_{>0})$
in $C_0$ into a product of a single factor $[n]$.
(Here we distinguish a single factor $[n]^c$ and a product $[n]\cdots [n]$ of a single factor $[n]$.)
\item[3.]
Repeat the following operation to $C_0$ in this order until we do nothing.
\begin{itemize}
\item[i.]
Pick up any adjacent pair $[n']^{c'}[n]^c$ in $C_0$.
($c'$ and $c$ are 1 or non-integers.)
\item[ii.]
If $\{n', n\}<0$ (ordered), we do nothing.
\item[iii.]
If $\{n', n\}=0$ (parallel),
\begin{itemize}
\item[a.]
If $\deg(n')>\deg(n)$, we exchange it as $[n]^c [n']^{c'}$.
(Align parallel vectors from left to right in the increasing order of degree.)
\item[b.]
Otherwise, we do nothing.
\end{itemize}
\item[iv.]
If $\{n', n\}>0$  (anti-ordered),
\begin{itemize}
\item[a.]
If $\deg(n'+n)> \l$, we exchange it as $[n]^c [n']^{c'}$.
(They are commutative modulo $G^{>\ell}$.)
\item[b.] If $\{n', n\}=p_0$ and $c'=c=1$ ($p_0$-exchangeable), 
replace  it with
$[n][n'+n][n']$.
(Apply the pentagon relation  \eqref{3eq:pent1} with $c=p_0$.)
\item[c.] Otherwise, we do nothing.
\end{itemize}
\end{itemize}
\item[4.]
If the resulting $C_0$ is ordered, 
we join all factors $[n]^{c_1}, \dots, [n]^{c_j}$ with common $n$ 
 (which are now adjacent to each other) into a single factor
$[n]^{c_1+\cdots +c_j}$.
Then, we set $C^{\mathrm{out}}({\ell}):=C_0$,
and the process completes.
\item[5.]
Otherwise, 
pick up  any $p_1$-exchangeable  adjacent pair in $C_0$,
\begin{align}
\label{3eq:Cin1}
C_1^{\mathrm{in}}=[n_1']^{c'}[n_1]^{c},
\quad
\{n'_1,n_1\}=p_1
\quad
(c', c \in (1/p_1)\bbZ_{>0}) 
\end{align}
for some integer $p_1>p_0$, and go to the subroutine ($p_1$).
\end{itemize}
\item[($p_1$).] 
\begin{itemize}
\item[1.]
Set $C_1=C_1^{\mathrm{in}}$.
\item[2.]
Decompose every factor $[n]^{c}$ $(c\in (1/p_1)\bbZ_{>0})$
in $C_1$ into a product of a single factor $[n]_{p_1}:=[n]^{1/p_1}$.
\item[3.]
Repeat the same process (3) of ($p_0$)  for $C_1$,
where 
\begin{itemize}
\item $[n]$ is replaced with $[n]_{p_1}$.
\item $p_0$ is replaced with $p_1$.
\end{itemize}
(In particular, the pentagon relation \eqref{3eq:pent1}  with $c=p_1$
is applied to $p_1$-exchangeable anti-ordered adjacent pairs
in $C_1$.)
\item[4.]
If the resulting  $C_1$ is ordered,
we join all factors $[n]^{c_1}_{p_1}, \dots, [n]^{c_j}_{p_1}$ with common $n$ 
 into a single factor $[n]^{(c_1+\cdots +c_j)/{p_1}}$.
 (The case $c_i\neq 1$ may happen by the result of the forthcoming subroutine ($p_2$).)
Then, we set $C_1^{\mathrm{out}}:=C_1$,
 replace $C_1^{\mathrm{in}}$ in $C_0$ with $C_1^{\mathrm{out}}$, and we come back to
 the process (2) of ($p_0$).
 (The ``compatibility'' of the exponent $(c_1+\cdots +c_j)/p_1$ in the parental routine ($p_0$) is the issue we are going to examine.)
 \item[5.]
Otherwise, 
pick up  any $p_2$-exchangeable  adjacent pair 
\begin{align}
C_2^{\mathrm{in}}=[n_2']^{c'}_{p_1} 
[n_2]^{c}_{p_1}
\quad
\{n'_2,n_2\}=p_2
\quad
(c', c \in (p_1/p_2) \bbZ_{>0}) 
\end{align}
in $C_1$
for some integer $p_2 \neq p_1$, and go to the subroutine ($p_2$).
(Since $n_2'$ and $n_2$ are $\bbZ$-linear combinations of $n_1'$ and $n_1$,
$p_2>p_1$ and $p_2$ is a multiple of $p_1$.)
\end{itemize}

\item[($p_2$).]
\begin{itemize}
\item[1--3.]
Repeat the same process (1)--(3) of $(p_1)$
by replacing $p_1$ with $p_2$.
\item[4.]
If the resulting  $C_2$ is ordered, 
join all factors $[n]^{c_1}_{p_2}, \dots, [n]^{c_j}_{p_2}$ with common $n$ 
 into a  factor
$[n]_{p_1}^{(c_1+\cdots +c_j)p_1 /{p_2}}$.
Then, we set $C_2^{\mathrm{out}}:=C_2$,
 replace $C_2^{\mathrm{in}}$ in $C_1$ with $C_2^{\mathrm{out}}$, and 
we come back to the process (2) of ($p_1$).
(Again, the ``compatibility'' of the exponent $(c_1+\cdots +c_j)p_1 /{p_2}$ in the parental routine ($p_1$) is the issue.)
 \item[5.]
Otherwise, go to the next subroutine ($p_3$) for some $p_3>p_2$.
\end{itemize}
\end{itemize}
\end{algo}

\begin{figure}
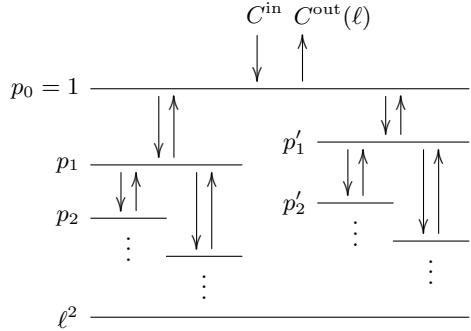

\centering
\leavevmode
\xy
(-2,40)*{\text{\small $C^{\mathrm{in}}$}};
(7,39.4)*{\text{\small $C^{\mathrm{out}}(\ell)$}};
(-31,30)*{\text{\small $p_0=1$}};
(-28,20)*{\text{\small $p_1$}};
(-28,13)*{\text{\small $p_2$}};
(2,23)*{\text{\small $p'_1$}};
(2,15)*{\text{\small $p'_2$}};
(-28,0)*{\text{\small $\ell^2$}};
(-20,10)*{\text{\small $\vdots$}};
(-10,5)*{\text{\small $\vdots$}};
(10,12)*{\text{\small $\vdots$}};
(20,7)*{\text{\small $\vdots$}};
\ar@{->} (-3,37); (-3,31)
\ar@{->} (3,31); (3,37)
\ar@{-} (-25,30); (25,30)
\ar@{-} (-25,0); (25,0)
\ar@{-} (-25,20); (-5,20)
\ar@{-} (5,23); (25,23)
\ar@{-} (-25,13); (-15,13)
\ar@{-} (-15,8); (-5,8)
\ar@{-} (5,15); (15,15)
\ar@{-} (15,10); (25,10)
\ar@{->} (-16,29); (-16,21)
\ar@{->} (-14,21); (-14,29)
\ar@{->} (-21,19); (-21,14)
\ar@{->} (-19,14); (-19,19)
\ar@{->} (-11,19); (-11,9)
\ar@{->} (-9,9); (-9,19)
\ar@{->} (14,29); (14,24)
\ar@{->} (16,24); (16,29)
\ar@{->} (9,22); (9,16)
\ar@{->} (11,16); (11,22)
\ar@{->} (19,22); (19,11)
\ar@{->} (21,11); (21,22)
\endxy
\caption{Flow of Ordering Algorithm.}
\label{3fig:algorithm1}
\end{figure}

See Figure \ref{3fig:algorithm1} for a schematic diagram of
the flow of the algorithm.
Note that each sequence of  subroutines  $p_0=1 < p_1< p_2< \cdots$ has a finite depth,
because  $p_i \leq \ell^2$.

We say that the above algorithm \emph{fails\/} if, at some stage,
$C_i$ is not ordered and all remaining anti-ordered adjacent pair
therein is not $p$-exchangeable for any integer $p>0$.
As mentioned in Algorithm \ref{3algo:ord1},
 some fractional power created by a subroutine
may  cause the failure of the process.
For example, in \eqref{3eq:fact1}, 
if the power of $[3,2]$ is $1/2$, instead of $1$ therein,
one cannot carry out the remaining ordering in \eqref{3eq:a11mod7},
because the pentagon relation with $c=1$ is not applicable.
We are going to show in Lemma 
\ref{3lem:nofail1}
 that such a situation never occurs.

If the algorithm does not fail, then
the algorithm completes (successfully) in finite steps
by the following reason:
\begin{itemize}
\item
Due to the structure of the pentagon relation \eqref{3eq:pent1},
all factors created in the subroutine ($p_i$) are
of higher degree than the factors of the initial $p_i$-exchangeable adjacent pair.
\item
In particular, for the lowest degree $\ell_0$ of 
 the factors  in $C^{\rm{in}}$,
 there is no new factor of degree $\ell_0$ that is
 created during the ordering.
\item
Then, by the induction on $\ell'$ with $\ell_0\leq \ell' \leq \ell$,
we can show that
 the number of the factors of degree $\ell'$
that are created during the ordering has some upper bound.
\end{itemize}

Let us go back to the situation in Proposition \ref{3prop:rank21}.
For a given $C^{\mathrm{in}}$ therein, suppose that the algorithm does not fail for any $\ell$.
Then, we have 
\begin{align}
C^{\mathrm{in}}\equiv C^{\mathrm{out}}({\ell}) \mod G^{>\ell},
\quad \pi_{\ell',\ell}(C^{\mathrm{out}}({\ell'}))=C^{\mathrm{out}}({\ell}).
\end{align}
Let $C^{\mathrm{out}}:=\lim_{\ell\rightarrow \infty} C^{\mathrm{out}}({\ell})$.
Then, we have
\begin{align}
\label{3eq:ee1}
C^{\mathrm{in}}
=
C^{\mathrm{out}},
\end{align}
and  $C^{\mathrm{out}}$ is a  (possibly infinite) ordered product of 
factors 
\begin{align}
\label{3eq:gwall1}
\Psi[t_i n_i]^{s_i \d(t_i n_i) }
\quad
(n_i\in N_{\rmpr}^+,\
t_i \in \bbZ_{>0},\
s_i\in \bbQ_{>0}
),
\end{align}
where we may assume that all $\Psi[t_i n_i]$ are mutually distinct
by joining the factors with  common $\Psi[t_i n_i]$.
The positivity of $s_i$ follows from the fact that the pentagon relation
does not involve negative factors.

To prove that $s_i$ is an integer,
 we use the following property of the  principal $x$-representation $\tilde \rho_x$.
\begin{lem}
\label{3lem:int0}
Under the principal $x$-representation $\tilde \rho_x$,
$\Psi[tn]^{s\d(tn)}$ acts on $\bbZ[[\tilde P]]$ if and only if
$s$ is an integer.
\end{lem}
\begin{proof}
Recall (see Example \ref{3ex:gn2}) that the action of
$\Psi[tn]^{s\d(tn)}$ ($n\in N_{\rmpr}^+,\
s\in \bbQ,\
t \in \bbZ_{>0})$
on $\bbk[[\tilde P]]$
is given by
\begin{align}
\label{3eq:gei3}
\Psi[tn]^{s\d(tn)}(x^{\tilde m})=x^{\tilde m} (1+x^{\tilde p^*(tn)})^{s \langle \d(n) n,   m\rangle }
\quad
(\tilde m \in \tilde P).
\end{align}
Thus,
the if-part is clear.
Suppose that $s$ is not an integer.
Note that
the coefficient of the  term  $x^{\tilde m}x^{\tilde p^*(tn)}$
 in \eqref{3eq:gei3}
is $c=s \langle \d(n) n,  m\rangle$.
Then,  it is enough to show that $c$ is not an integer for some $\tilde m\in  \tilde P$.
Let $a_0>0$ be the smallest positive integer such that
 $\langle \d(n) n,  m\rangle \in a_0 \bbZ$ for any $\tilde m\in \tilde P$.
By the condition (i) of the definition of $\tilde P$ in 
Section \ref{3subsec:principal1},
$\tilde P$ contains a $\bbZ$-basis of $\tilde{M}^{\circ}$,
for example $\tilde m + f_i, \tilde m + e_i$ ($i=1,\dots,r$)
for some $\tilde m\in \tilde P$.
Then, by  $\bbZ$-linearity,
we have $\langle \d(n) n,  m\rangle \in a_0 \bbZ$ for any $\tilde m\in \tilde M^{\circ}$.
If $a_0\neq 1$,
we have $\d(n) n/a_0\in N^{\circ}$,
which contradicts the definition of $\d(n)$.
Thus, $a_0=1$. Then, $c=s$ is not an integer for $\tilde m\in \tilde P$
attaining the minimum $a_0=1$.
\end{proof}

\begin{lem}
\label{3lem:int1}
 The factor $s_i$ in \eqref{3eq:gwall1} is an integer.
\end{lem}
\begin{proof}
By Lemma \ref{3lem:int0}, the product $C^{\mathrm{in}}$   acts on $\bbZ[[\tilde P]]$.
Then,  again by Lemma \ref{3lem:int0},
we can prove that
all $s_i$ in \eqref{3eq:gwall1} are   integers
by the induction on $\deg(t_in_i)$.
\end{proof}

It remains to show that Algorithm \ref{3algo:ord1} never fails.
This is proved by the same method for Lemma \ref{3lem:int1}.

\begin{lem}
\label{3lem:nofail1}
For any $C^{\mathrm{in}}$,
Algorithm \ref{3algo:ord1} never fails,
so that it completes in finitely many steps.
\end{lem}
\begin{proof}
We first note a key observation to prove the lemma.

For the first subroutine ($p_{1}$) entered from the main
 routine $(p_{0})$,
the following facts hold:
\begin{itemize}
\item[({\bf i}).]
The initial product $C^{\mathrm{in}}_1$ in \eqref{3eq:Cin1} has the form
\begin{align}
\label{3eq:form1}
C^{\mathrm{in}}_1=[t_1 n_1]^{s_1/t_1}[t_2 n_2]^{s_2/t_2}
\quad (n_i\in N_{\rmpr}^+,\ s_i, t_i\in \bbZ_{>0}).
\end{align}
\item[({\bf o}).]
If the algorithm does not fail,
any factor in  the final product $C^{\mathrm{out}}_1$ also has the form
$[t n]^{s/t}$
with the same condition in \eqref{3eq:form1}.
\end{itemize}
The property ($\text{\bf o}$) follows from the property ($\text{\bf i}$) by applying the same proof of
Lemma \ref{3lem:int1}, where we  set $N^{\circ}=N$, $M^{\circ}=M$
and consider the action on $\bbZ[[\tilde P]]$.

Similarly, but slightly differently, suppose that a subroutine $(p_i)$ $(i\geq 1)$ starts with
the initial product
\begin{align}
C^{\mathrm{in}}_i
=
[\tilde n_1]^{s_1}_{p_i}  [\tilde n_2]^{s_2}_{p_i},
\quad
\{\tilde n_1,\tilde n_2\}=p_i
\quad 
(\tilde n_i\in N^+,\ s_i \in \bbZ_{>0}).
\end{align}
Let
\begin{align}
 N_{p_i}^+ :=\{  a_1 \tilde n_1 + a_2 \tilde n_2 \in N^+\mid a_i \in \bbZ_{\geq 0},\ a_1+a_2 >0\}.
\end{align}
Let $ ( N_{p_i}^+)_{\rmpr}$ be the set of all primitive vectors in $ N_{p_i}^+$
with respect to  $ N_{p_i}^+$ (not with respect to $N^+$).
Note that, for any $n_1, n_2\in   N_{p_i}^+$,
$\{n_1,n_2\}$ is a multiple of $p_i$.
Then, for a subroutine ($p_{i+1}$) $(i\geq 1)$  entered from a
 subroutine $(p_{i})$,
the following facts hold:
\begin{itemize}
\item[({\bf i}).]
The initial product $C^{\mathrm{in}}_{i+1}$ has the form
\begin{align}
\label{3eq:form2}
C^{\mathrm{in}}_{i+1}=
[t_1 n_1]^{s_1/t_1}_{p_i}
[t_2 n_2]^{s_2/t_2}_{p_i}
\quad (n_i\in ( N_{p_i}^+)_{\rmpr},\ s_i, t_i\in \bbZ_{>0}).
\end{align}
\item[({\bf o}).]
If the algorithm does not fail,
any factor in  the final product $C^{\mathrm{out}}_{i+1}$ also has the form
$[t n]^{s/t}_{p_i}$
with the same condition in \eqref{3eq:form2}.
\end{itemize}
Again, the property ($\text{\bf o}$) follows from the property ($\text{\bf i}$) by
the same argument as above,
where we consider the action on $(1/p_i)\bbZ[[\tilde P]]$,
and $N$ and $M^{\circ}$ for $\tilde P$ are replaced with
the sublattice $ N_{p_i}$ generated by $\tilde n_1$ and $\tilde n_2$
and its dual lattice $ M_{p_i}$, respectively.

Let us prove  the lemma together with the property ({\bf i}) 
in both cases
along  the flow of the algorithm.
\begin{itemize}
\item[1.]
The algorithm starts from the main routine ($p_0$),
where  the initial product $C^{\mathrm{in}}$ certainly has factors of the
form
in \eqref{3eq:form1} by assumption.
Also,
any factor created by the pentagon relation from a $1$-exchangeable pair  has the form $[tn]^{s}$
with $s\in \bbZ_{>0}$; in particular,
it has the form in \eqref{3eq:form1}.
Suppose that we are in the process (5) of the routine ($p_0$).
Take any remaining anti-ordered pair $[t_1 n_1]^{s_1/t_1}[t_2 n_2]^{s_2/t_2}$,
and set $p_1:=\{t_1n_1, t_2n_2\}$.
Then,  by Example \ref{3ex:padm1},
it is $p_1$-exchangeable.
Thus, we go to the subroutine ($p_1$) with $C_{1}^{\mathrm{in}}$ satisfying the property ({\bf i}).
\item[2.]
In the subroutine $(p_1)$
we decompose $C_{1}^{\mathrm{in}}$  into integer powers of $[n]_{p_1}$,
and after ordering all $p_1$-exchangeable pairs by the pentagon relation,
we only obtain integer powers of $[n]_{p_1}$.
Thus, if we have an anti-ordered pair in the process (5) of the routine ($p_1$),
it has the form in \eqref{3eq:form2}  with $\{t_1n_1,t_2n_2\}=p_2>p_1$
and $s_i/t_i$ are integers.
In particular, 
the pair is $p_2$-exchangeable, and 
we go to the next subroutine ($p_2$) with $C_{2}^{\mathrm{in}}$ satisfying the property ({\bf i}).
(So far, we do not encounter a serious problem.)
\item[3.]
The process continues
until we get to the deepest subroutine $(p_k)$.
Then, by definition, $(p_k)$ returns the ordered product $C^{\mathrm{out}}_k$
satisfying the property ({\bf o}) to
the subroutine $(p_{k-1})$.
We claim that the subroutine $(p_{k-1})$ continues to run.
 (The following is the key argument.)
In fact, suppose that we still have any anti-ordered pair
in $C_{k-1}$. Then, 
it has the form
\begin{align}
 [t_1n_1]_{p_{k-1}}^{s_1/t_1}
  [t_2n_2]_{p_{k-1}}^{s_2/t_2}
  \quad (n_i\in ( N_{p_{k-1}}^+)_{\rmpr},\ s_i, t_i\in \bbZ_{>0})
 \end{align}
 thanks to the property ({\bf o}) just we had.
Let  $p'_k=\{t_1n_1,t_2n_2\}>0$.
 Since $\{n_1,n_2\}$ is  a multiple of $p_{k-1}$ as already mentioned,
 $p'_k$ is a multiple of both $t_1 p_{k-1}$ and $t_2 p_{k-1}$.
 Thus, the pair is $p'_k$-exchangeable.
 If $p'_k=p_{k-1}$, we stay in the subroutine $(p_{k-1})$.
If $p'_k>p_{k-1}$, we go to another subroutine
($p'_k$) again with $C_{k}^{\mathrm{in}}$ satisfying the property ({\bf i}).
The same story continues.
\end{itemize}
Thus, the algorithm never fails.
\end{proof}

By Lemmas \ref{3lem:int1} and \ref{3lem:nofail1},
we obtain Proposition \ref{3prop:rank21}.
Moreover, as a special case of Proposition \ref{3prop:rank21},
we obtain Theorem \ref{3thm:pos1} in the rank 2 case.
\begin{cor}
\label{3cor:pos1}
Theorem \ref{3thm:pos1}
holds for $\mathrm{rank}\, N=2$.
\end{cor}
\begin{proof}
Apply Proposition \ref{3prop:rank21} to
the product $C^{\mathrm{in}}=\Psi[e_2]^{\d_2}\Psi[e_1]^{\d_1}$.
\end{proof}

\begin{ex}
\label{3ex:nonaffine1}
It is not difficult to translate Algorithm \ref{3algo:ord1}
into a computer program.
The source code is given in Appendix
\ref{3sec:source1}
for the reader's convenience.
Let us exhibit some examples
of  non-affine infinite type
calculated by computer.
Under the same assumption in 
Section \ref{3subsec:rank2inf},
we have the following results modulo $G^{>7}$:
\begin{align}
\label{3eq:algres1}
\begin{split}
\begin{bmatrix}
0\\
1
\end{bmatrix}
^5
\begin{bmatrix}
1\\
0
\end{bmatrix}
^1
&
\equiv
\begin{bmatrix}
1\\
0
\end{bmatrix}
^1
\begin{bmatrix}
1\\
1
\end{bmatrix}
^5
\begin{bmatrix}
3\\
4
\end{bmatrix}
^5
\biggl(
\begin{bmatrix}
2\\
3
\end{bmatrix}
^{10}
\begin{bmatrix}
1\\
2
\end{bmatrix}
^{10}
\begin{bmatrix}
2\\
4
\end{bmatrix}
^{10}
\begin{bmatrix}
2\\
5
\end{bmatrix}
^{27}
\begin{bmatrix}
1\\
3
\end{bmatrix}
^{10}
\biggr)
\\
&\qquad\ \times
\begin{bmatrix}
1\\
4
\end{bmatrix}
^{5}
\begin{bmatrix}
1\\
5
\end{bmatrix}
^1
\begin{bmatrix}
0\\
1
\end{bmatrix}
^5
,
\end{split}
\\
\label{3eq:algres2}
\begin{split}
\begin{bmatrix}
0\\
1
\end{bmatrix}
^6
\begin{bmatrix}
1\\
0
\end{bmatrix}
^1
&
\equiv
\begin{bmatrix}
1\\
0
\end{bmatrix}
^1
\begin{bmatrix}
1\\
1
\end{bmatrix}
^6
\biggl(
\begin{bmatrix}
3\\
4
\end{bmatrix}
^{15}
\begin{bmatrix}
2\\
3
\end{bmatrix}
^{20}
\begin{bmatrix}
1\\
2
\end{bmatrix}
^{15}
\begin{bmatrix}
2\\
4
\end{bmatrix}
^{30}
\begin{bmatrix}
2\\
5
\end{bmatrix}
^{102}
\begin{bmatrix}
1\\
3
\end{bmatrix}
^{20}
\begin{bmatrix}
1\\
4
\end{bmatrix}
^{15}
\biggr)
\\
&\qquad\ \times
\begin{bmatrix}
1\\
5
\end{bmatrix}
^6
\begin{bmatrix}
1\\
6
\end{bmatrix}
^1
\begin{bmatrix}
0\\
1
\end{bmatrix}
^6
,
\end{split}
\\
\label{3eq:algres3}
\begin{split}
\begin{bmatrix}
0\\
1
\end{bmatrix}
^3
\begin{bmatrix}
1\\
0
\end{bmatrix}
^2
&
\equiv
\begin{bmatrix}
1\\
0
\end{bmatrix}
^2
\begin{bmatrix}
2\\
1
\end{bmatrix}
^3
\biggl(
\begin{bmatrix}
3\\
2
\end{bmatrix}
^6
\begin{bmatrix}
4\\
3
\end{bmatrix}
^{14}
\begin{bmatrix}
1\\
1
\end{bmatrix}
^{6}
\begin{bmatrix}
2\\
2
\end{bmatrix}
^6
\begin{bmatrix}
3\\
3
\end{bmatrix}
^6
\\
&\qquad\ \times
\begin{bmatrix}
3\\
4
\end{bmatrix}
^{36}
\begin{bmatrix}
2\\
3
\end{bmatrix}
^{14}
\begin{bmatrix}
1\\
2
\end{bmatrix}
^6
\begin{bmatrix}
2\\
4
\end{bmatrix}
^6
\biggr)
\begin{bmatrix}
2\\
5
\end{bmatrix}
^3
\begin{bmatrix}
1\\
3
\end{bmatrix}
^2
\begin{bmatrix}
0\\
1
\end{bmatrix}
^3
,
\end{split}
\\
\label{3eq:algres4}
\begin{split}
\begin{bmatrix}
0\\
1
\end{bmatrix}
^4
\begin{bmatrix}
1\\
0
\end{bmatrix}
^2
&
\equiv
\begin{bmatrix}
1\\
0
\end{bmatrix}
^2
\begin{bmatrix}
2\\
1
\end{bmatrix}
^4
\biggl(
\begin{bmatrix}
3\\
2
\end{bmatrix}
^{12}
\begin{bmatrix}
4\\
3
\end{bmatrix}
^{44}
\begin{bmatrix}
1\\
1
\end{bmatrix}
^{8}
\begin{bmatrix}
2\\
2
\end{bmatrix}
^{12}
\begin{bmatrix}
3\\
3
\end{bmatrix}
^{24}
\\
&\qquad\ \times
\begin{bmatrix}
3\\
4
\end{bmatrix}
^{182}
\begin{bmatrix}
2\\
3
\end{bmatrix}
^{44}
\begin{bmatrix}
1\\
2
\end{bmatrix}
^{12}
\begin{bmatrix}
2\\
4
\end{bmatrix}
^{40}
\begin{bmatrix}
2\\
5
\end{bmatrix}
^{44}
\begin{bmatrix}
1\\
3
\end{bmatrix}
^8
\biggr)
\begin{bmatrix}
1\\
4
\end{bmatrix}
^{2}
\begin{bmatrix}
0\\
1
\end{bmatrix}
^4
,
\end{split}
\\
\label{3eq:algres5}
\begin{split}
\begin{bmatrix}
0\\
1
\end{bmatrix}
^3
\begin{bmatrix}
1\\
0
\end{bmatrix}
^3
&
\equiv
\begin{bmatrix}
1\\
0
\end{bmatrix}
^3
\begin{bmatrix}
3\\
1
\end{bmatrix}
^3
\biggl(
\begin{bmatrix}
5\\
2
\end{bmatrix}
^9
\begin{bmatrix}
2\\
1
\end{bmatrix}
^{9}
\begin{bmatrix}
4\\
2
\end{bmatrix}
^{18}
\begin{bmatrix}
3\\
2
\end{bmatrix}
^{39}
\begin{bmatrix}
4\\
3
\end{bmatrix}
^{204}
\begin{bmatrix}
1\\
1
\end{bmatrix}
^{9}
\begin{bmatrix}
2\\
2
\end{bmatrix}
^{18}
\\
&\qquad\ \times
\begin{bmatrix}
3\\
3
\end{bmatrix}
^{54}
\begin{bmatrix}
3\\
4
\end{bmatrix}
^{204}
\begin{bmatrix}
2\\
3
\end{bmatrix}
^{39}
\begin{bmatrix}
1\\
2
\end{bmatrix}
^9
\begin{bmatrix}
2\\
4
\end{bmatrix}
^{18}
\begin{bmatrix}
2\\
5
\end{bmatrix}
^9
\biggr)
\begin{bmatrix}
1\\
3
\end{bmatrix}
^3
\begin{bmatrix}
0\\
1
\end{bmatrix}
^3
.
\end{split}
\end{align}
Let us examine the behavior
in the Badlands $\sigma(\bfv_+, \bfv_-)$
defined in \eqref{3eq:ds1}.
The condition \eqref{3eq:ds2} for a normal vector  $\bfn=(n_1,n_2)$ is given as follows:
\begin{alignat}{3}
(\d_1,\d_2)&=(1,5), &\quad &
\frac{2}{5+\sqrt{5}}=0.276\dots
&\leq \frac{n_1}{n_2}
\leq
\frac{2}{5-\sqrt{5}}=0.723\dots,
\\
(\d_1,\d_2)&=(1,6), &\quad &
\frac{1}{3+\sqrt{3}}=0.211\dots
&\leq \frac{n_1}{n_2}
\leq
\frac{1}{3-\sqrt{3}}=0.788\dots,
\\
(\d_1,\d_2)&=(2,3), &\quad &
\frac{2}{3+\sqrt{3}}=0.422\dots
&\leq \frac{n_1}{n_2}
\leq
\frac{2}{3-\sqrt{3}}=1.577\dots,
\\
(\d_1,\d_2)&=(2,4), &\quad &
\frac{1}{2+\sqrt{2}}=0.292\dots
&\leq \frac{n_1}{n_2}
\leq
\frac{1}{2-\sqrt{2}}=1.707\dots,
\\
(\d_1,\d_2)&=(3,3), &\quad &
\frac{2}{3+\sqrt{5}}=0.381\dots
&\leq \frac{n_1}{n_2}
\leq
\frac{2}{3-\sqrt{5}}=2.618\dots.
\end{alignat}
The corresponding elements $\Psi[n]$ are inside the
parentheses in \eqref{3eq:algres1}--\eqref{3eq:algres5}.
We observe that, in each case, \emph{every} $n\in N^+$ with $\deg(n)\leq 7$ satisfying the above condition
appears in the product.
This agrees with (and even stronger than) the expectation of \cite{Gross14}
explained in Section \ref{3subsec:rank2inf}.
We expect this is true up to any degree, and we do not have any counterexample so far.

Let us also present an example of the ordering of a more general
anti-ordered product in the left hand side.
We have, modulo $G^{>9}$,
\begin{align}
\begin{split}
&\quad\,
\begin{bmatrix}
0\\
1
\end{bmatrix}
^2
\begin{bmatrix}
2\\
4
\end{bmatrix}
^{1/2}
\begin{bmatrix}
3\\
3
\end{bmatrix}
^{2/3}
\begin{bmatrix}
1\\
0
\end{bmatrix}
^2
\\
&
\equiv
\begin{bmatrix}
1\\
0
\end{bmatrix}
^2
\begin{bmatrix}
2\\
1
\end{bmatrix}
^{2}
\begin{bmatrix}
6\\
3
\end{bmatrix}
^{40/3}
\begin{bmatrix}
5\\
3
\end{bmatrix}
^{10}
\begin{bmatrix}
3\\
2
\end{bmatrix}
^{2}
\begin{bmatrix}
4\\
3
\end{bmatrix}
^{6}
\begin{bmatrix}
5\\
4
\end{bmatrix}
^{146}
\begin{bmatrix}
1\\
1
\end{bmatrix}
^{4}
\begin{bmatrix}
2\\
2
\end{bmatrix}
^{2}
\begin{bmatrix}
3\\
3
\end{bmatrix}
^{2/3}
\begin{bmatrix}
4\\
4
\end{bmatrix}
^{47}
\begin{bmatrix}
4\\
5
\end{bmatrix}
^{230}
\\
&\qquad\ \times
\begin{bmatrix}
3\\
4
\end{bmatrix}
^{10}
\begin{bmatrix}
2\\
3
\end{bmatrix}
^2
\begin{bmatrix}
3\\
5
\end{bmatrix}
^{34}\begin{bmatrix}
1\\
2
\end{bmatrix}
^2\begin{bmatrix}
2\\
4
\end{bmatrix}
^{1/2}\begin{bmatrix}
3\\
6
\end{bmatrix}
^{220/3}
\begin{bmatrix}
2\\
5
\end{bmatrix}
^2\begin{bmatrix}
2\\
6
\end{bmatrix}
^3
\begin{bmatrix}
2\\
7
\end{bmatrix}
^2
\begin{bmatrix}
0\\
1
\end{bmatrix}
^2
.
\end{split}
\end{align}

As we see in  these examples, some exponents of $\Psi[n]$ quickly increases when $\deg(n)$ increases.
In fact, they easily exceed 1,000 when the degree is around  10 even for these tiny examples.
This puts a limitation of the explicit computation of higher degrees even with computer,
unfortunately.
\end{ex}

\subsection{Alternative construction of CSDs}

Here we present an alternative construction of CSDs
based on Proposition \ref{3prop:rank21} (Ordering lemma).
The resulting CSD
is a positive realization in
Theorem \ref{3thm:pos1}.

Let us introduce the following notion following \cite{Gross14}.

\begin{defn}[Joint]
Let $\frakD$ be a scattering diagram.
\par
(a). 
For any pair of walls $\bfw_i=(\frakd_i,g_{i})_{n_i}$ ($i=1,2$), the intersection
of their supports $\frakj=\frakd_1\cap \frakd_2$
is called a \emph{joint\/} of $\frakD$ if $\frakj$
is a cone of codimension 2.  \index{joint}
\par
(b) For any joint $\frakj$, let
\begin{align}
\label{3eq:Njn1}
N_{\frakj}:=
\{ n\in N \mid \langle n, z\rangle=0
\
( z\in \frakj) \}.
\end{align}
Then,
a joint $\frakj$ is \emph{parallel\/} 
(resp.\ \emph{perpendicular\/}) if 
the skew-symmetric form $\{\cdot, \cdot\}$ restricted on $N_{\frakj}$
vanishes (resp. otherwise).
\index{joint!parallel}\index{joint!perpendicular}
Thus,  the joint $\frakj$ in (a)  is parallel if and only if $\{n_1,n_2\}=0$.
\end{defn}

By the topological reason,
any consistency relation of $\frakD$  reduces to the consistency relations for admissible loops around joints.

For a given seed $\fraks$,
we will inductively construct  scattering diagrams $\frakD_{\l}$
($\l=1, 2, \dots$) satisfying the following properties:
\begin{itemize}
 \item[(a).]
 $\frakD_1=\mathrm{In}_{\fraks}$,
 where $\mathrm{In}_{\fraks}$ is  the set of  incoming walls 
in \eqref{3eq:in2}.
\item[(b).] $\frakD_1\subset \frakD_2 \subset \cdots$,
and
all walls of $\frakD_{\l+1}\setminus \frakD_{\l}$ 
are outgoing walls and  have the form
in \eqref{3eq:gpos1}
with $\deg(tn)=\l+1$.
\item[(c).]
Let $\frakj$ be any perpendicular joint  of $\frakD_{\l}$,
and let $\gamma_{\frakj}$ be any sufficiently small loop  around
$\frakj$ such that it only crosses the walls of $\frakD_{\l}$ containing $\frakj$.
Then, we have
\begin{align}
\label{3eq:jid1}
\frakp_{\gamma_{\frakj}, \frakD_{\l}}\equiv\rmid
 \mod
G^{>\l}.
\end{align}
 \end{itemize}

 Note that  we do not directly impose the condition \eqref{3eq:jid1} for
 parallel joints.

\begin{const}[{Cf.\ \cite[Construction C.1, Appendix C.3]{Gross14}}]
\label{3const:CSD1}
\
\par
(1). We start with $\frakD_1=\mathrm{In}_{\fraks}$.

  \par
  (2).
Next, we 
describe the construction of $\frakD_{2}$ 
from $\frakD_1$
to illustrate the idea for the  general case (3).
For any \emph{perpendicular\/} joint
$\frakj=e_i^{\perp}\cap e_j^{\perp}
$
of $\frakD_{1}$,
where
$\{e_i,e_j\}\neq 0$,
let $N_{\frakj}$ be the one in \eqref{3eq:Njn1},
which is the rank 2 sublattice of $N$
generated by $e_i$ and $e_j$.
By exchanging $i$ and $j$ if necessary, we may assume that
 that $c=\{ e_j, e_i\}>0$.
Then, for any $n=n_i e_i + n_je_j\in N_{\frakj}^+:=N^+\cap N_{\frakj}$,
we have
\begin{align}
\langle e_i, p^*(n)\rangle
=\{e_i,n\} \leq 0,
\quad 
\langle e_j, p^*(n)\rangle
=\{e_j, n\} \geq 0.
\end{align}
Thus, $p^*(n)$ is in the ``second quadrant''
in $M_{\bbR}$ with respect to $e_i^{\perp}$ and $e_j^{\perp}$.
Now, we apply Proposition \ref{3prop:rank21}
to the anti-ordered product
$C^{\mathrm{in}}=\Psi[e_j]^{\d_j} \Psi[e_i]^{\d_i}$
modulo $G^{>2}$,
and obtain the ordered product
\begin{align}
 \Psi[e_i]^{\d_i} \Psi[e_i+e_j]^{c \d_i \d_j}  \Psi[e_j]^{\d_j}.
\end{align}
Then,  we add a  new wall
\begin{align}
\label{3eq:add1}
(\sigma(\frakj, -p^*(e_i+e_j)),
\Psi[e_i+e_j]^{c \d_i\d_j})_{e_i+e_j}
\end{align}
to $\frakD_1$,
where $\sigma(\frakj, -p^*(e_i+e_j))$ is the cone in $M_{\bbR}$ 
spanned by $\frakj$ and $-p^*(e_i+e_j)$.
The perpendicular condition implies that $p^*(e_i+e_j)\not\in \frakj$.
Thus,
the cone $\sigma(\frakj, -p^*(e_i+e_j))$ is of codimension 1.
Also,
it is outgoing because it is in the ``fourth quadrant'' in $M_{\bbR}$
with respect to $e_i^{\perp}$ and $e_j^{\perp}$.
We do the procedure for all  perpendicular joints of $\frakD_1$
to obtain $\frakD_2$.

 \begin{figure}
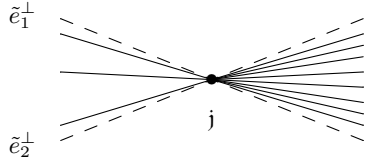

 \centering
\leavevmode
\xy
(0,-5)*{\text{\small $\frakj$}};
(-25,8.5)*{\text{\small $ \tilde{e}_1^{\perp}$}};
(-25,-8.5)*{\text{\small $ \tilde{e}_2^{\perp}$}};
(0,0)*+{\bullet};
\ar@{--} (-20,8); (20,-8)
\ar@{--} (-20,-8); (20,8)
\ar@{-} (-20,6); (20,-6)
\ar@{-} (-20,1); (20,-1)
\ar@{-} (-20,-6); (20,6)
\ar@{-} (0,0); (20,3)
\ar@{-} (0,0); (20,4.5)
\ar@{-} (0,0); (20,1)
\ar@{-} (0,0); (20,-3)
\ar@{-} (0,0); (20,-4.5)
\endxy
\caption{Walls containing a perpendicular joint $\frakj$.}
\label{3fig:2d1}
\end{figure}

\par
(3). Generalizing the procedure in (2), we construct $\frakD_{\l+1}$ from $\frakD_{\l}$ as follows.
By decomposing the supports of walls of $\frakD_{\l}$ if necessary, we may assume
that any pair of  joints $\frakj_1\neq \frakj_2$ may intersect each other only in their boundaries.
For any \emph{perpendicular\/} joint
$\frakj
$
of $\frakD_{\l}$,
let $N_{\frakj}$ be the one in \eqref{3eq:Njn1}.
Let $N_{\frakj}^+ :=N^+\cap N_{\frakj}$
and
 $N_{\frakj,\rmpr}^+ :=N^+_{\rmpr}\cap N_{\frakj}$.
  Then, there exist a unique pair $\tilde{e}_1, \tilde{e}_2\in N_{\frakj,\rmpr}^+$ such that
  \begin{align}
  N_{\frakj}^+ \subset \bbQ_{\geq 0}   \tilde{e}_1+  \bbQ_{\geq 0}\tilde{e}_2,
  \quad
\{ \tilde{e}_2, \tilde{e}_1 \}>0.
  \end{align}
  Then,
 for any $n\in N_{\frakj}^+$, 
\begin{align}
\label{3eq:quad1}
\langle \tilde{e}_1, p^*(n)\rangle
=\{\tilde{e}_1, n\} \leq 0,
\quad 
\langle \tilde{e}_2, p^*(n)\rangle
=\{\tilde{e}_2, n\} \geq 0.
\end{align}
See Figure \ref{3fig:2d1}.
Thus, $p^*(n)$ is in the ``second quadrant''
in $M_{\bbR}$ with respect to $\tilde{e}_1^{\perp}$ and $\tilde{e}_2^{\perp}$.
Then, by the construction of walls as described below,
all wall elements in $\frakD_{\l}$ has the form in  \eqref{3eq:gpos1}.
We take the wall elements of all walls
that  contain $\frakj$
and intersect the second quadrant
 with respect to $\tilde{e}_1^{\perp}$ and $\tilde{e}_2^{\perp}$.
 Then, we consider the anti-ordered product of them
\begin{align}
\label{3eq:C0}
C^{\mathrm{in}}=
\Psi[t'_{j'} n'_{j'}]^{s'_{j'}  \d(t'_{j'} n'_{j'})}
\cdots
\Psi[t'_1n'_1]^{s'_1\d(t'_1 n'_1)}
\quad
(\deg(t'_i n'_i)\leq \ell),
\end{align}
where  $n'_i \in N_{\frakj,\rmpr}^+$ and   $s'_i, t'_i\in \bbZ_{>0}$.
Then, we apply Proposition \ref{3prop:rank21}
to $C^{\mathrm{in}}$
modulo $G^{>{\l+1}}$,
and obtain a finite ordered product
\begin{align}
\label{3eq:C1}
C^{\mathrm{out}}(\ell+1)=
\Psi[t_j n_j]^{s_j  \d(t_j n_j)}
\cdots \Psi[t_1n_1]^{s_1\d(t_1 n_1)}
\quad
(\deg(t_i n_i)\leq \ell+1),
\end{align}
such that $C^{\mathrm{in}}\equiv C^{\mathrm{out}}(\ell+1)$ mod $G^{>\ell+1}$.
By the construction of walls which we describe below,
all factors in \eqref{3eq:C1} with $\deg(t_i n_i)\leq \l$
already appear as wall elements in $\frakD_{\l}$.
Now, we add new walls
\begin{align}
\label{3eq:add2}
(\sigma(\frakj, -p^*(n_i)),
\Psi[t_i n_i ]^{s_i \d(t_i n_i)})_{n_i}
\end{align}
 to $\frakD_{\l}$ for all factors  in \eqref{3eq:C1}
with $\deg(t_i n_i)= \l+1$.
By the perpendicular condition,
the cone $\sigma(\frakj, -p^*(n_i))$ is of codimension 1.
Also,
it is outgoing because it is in the ``fourth quadrant'' 
with respect to $\tilde{e}_1^{\perp}$ and $\tilde{e}_2^{\perp}$.
We call such a wall in \eqref{3eq:add2} \emph{a wall attached to $\frakj$}. \index{wall!attached to a joint}
We do the procedure for all  perpendicular joints of $\frakD_{\l}$
to obtain $\frakD_{\l+1}$.
\end{const}

  Next we  prove the condition (c) for $\frakD_{\l+1}$
  (not by the induction on $\l$).
  It is useful to note the following property.
  
\smallskip
\par\noindent
\begin{lem}
\label{3lem:pcont1}
Let $\bfw$ be any wall of $\frakD_{\l+1}\setminus \frakD_{\l}$
 attached to a joint $\frakj$ in $\frakD_{\l}$
in \eqref{3eq:add2}.
Let $\bfw'$ be any wall of $\frakD_{\l+1}$ other than $\bfw$.
Let $\frakd$ and $\frakd'$ be the supports
of $\bfw$ and $\bfw'$, respectively.
Suppose that $\frakj'=\frakd \cap \frakd'$ is a joint of $\frakD_{\l+1}$
such that $\frakj'\not\subset \frakj$.
\par
(a).
If $\frakj'\subset \partial \frakd$,
then $\frakj'$ is parallel.
\par
(b).
If $\frakj'$ is perpendicular,
  the wall $\bfw$ does not
 contribute to
    $\frakp_{\gamma_{\frakj'}, \frakD_{\l+1}}$ modulo $G^{>\l+1}$
\end{lem}
 \begin{proof}
 (a).
 Let $\frakd=\sigma(\frakj, -p^*(n))$ for some $n\in N_{\frakj,\rmpr}^+$.
 Then,
 by the assumptions that $\frakj'\not\subset \frakj$ and $\frakj'\subset \partial \frakd$,
the joint $\frakj'=\frakd \cap \frakd'$ is contained in a vector space of codimension 2 that  contains $p^*(n)$.
Thus, $p^*(n)\in (n')^{\perp}$, where $n'\in N^+_{\rmpr}$ is the normal vector for $\bfw'$.
 It follows that $0=\langle n', p^*(n)\rangle = \{n',n\}$.
Therefore,  $\frakj'$ is parallel.
 \par
 (b).
 By (a), we have $\frakj' \not\subset \partial\frakd$.
 Then, $\gamma_{\frakj'}$ crosses $\frakd$ twice with opposite intersection signs.
 Since the wall element of $\bfw$ is in the center of  $G^{\leq \l+1}$,
two contributions  of the wall $\bfw$
to $\frakp_{\gamma_{\frakj'}, \frakD_{\l+1}}$
cancel     modulo $G^{>\l+1}$.
 \end{proof}

  Let $\frakj'$ be any perpendicular joint of $\frakD_{\l+1}$.
There are two cases:
\begin{itemize}
\item
Suppose that $\frakj'\not\subset \frakj$ for any perpendicular joint $\frakj$ of $\frakD_{\l }$.
Then, $\frakj'$ is a new joint created by some added wall to $\frakD_{\l}$.
Then,
we have 
 $\frakp_{\gamma_{\frakj'}, \frakD_{\l+1}}\equiv \rmid$ mod $G^{>\l+1}$
 by Lemma \ref{3lem:pcont1} (b).
\item
Suppose that $\frakj'\subset \frakj$ for some perpendicular joint $\frakj$ of $\frakD_{\l }$.
Again by  Lemma \ref{3lem:pcont1} (b), one can ignore the contribution to
$\frakp_{\gamma_{\frakj'}, \frakD_{\l+1}}$ modulo $G^{>\l+1}$
from the added walls other than the ones attached to $\frakj$.
Then, by the construction of the added walls attached to $\frakj$
we have 
$\frakp_{\gamma_{\frakj'}, \frakD_{\l+1}}\equiv 
\frakp_{\gamma_{\frakj}, \frakD_{\l+1}}\equiv \rmid$ mod $G^{>\l+1}$.

\end{itemize}
Thus, the condition (c) holds for $\frakD_{\l+1}$.
 
 Using the diagrams $\frakD_1\subset \frakD_2 \subset \cdots$
 constructed above, 
we define
 \begin{align}
 \frakD=
 \bigcup_{\l=1}^\infty \frakD_{\l},
 \end{align}
 whose reduction
 at $\l$ is $\frakD_{\l}$.
 The scattering diagram $\frakD$ has the following
 properties:
 \begin{itemize}
 \item[(a).]
 $\frakD_{\mathrm{in}}=\mathrm{In}_{\fraks}$.
\item[(b).]
All outgoing walls of $\frakD$ 
 have the form
in \eqref{3eq:add1}
with $\d(tn)>1$.
\item[(c).]
For 
any perpendicular joint   $\frakj$  of $\frakD_{\l}$,
 $\frakp_{\gamma_{\frakj}, \frakD_{\l}}\equiv \rmid$ mod
$G^{>\l}$ holds.
 \end{itemize}
 However, we have not yet shown that $\frakD$ is consistent,
 because the consistency around parallel joints is not examined.

The proof of Theorem  \ref{3thm:pos1} completes
by showing the following  lemma.

\begin{lem}
[{\cite[Lemma~C.7]{Gross14}}]
\label{3lem:unique1}
Let $\frakD$ be the one constructed above.
Let $\frakD_{\fraks}$ be a CSD in
Theorem \ref{3thm:exist2}.
Then, 
$\frakD$ and $\frakD_{\fraks}$ are equivalent.
In particular, $\frakD$ is consistent and it is a CSD for $\fraks$.
\end{lem}
\begin{proof}
We note that any CSD $\frakD_{\fraks}$
also
satisfies  the property (a) by definition
and the property (c) by the consistency.
Let $\frakC=\frakD_{\fraks}- \frakD$
be a scattering diagram obtained as follows.
Take all walls of ${\frakD}$, invert
all wall elements,  add all walls of  $ {\frakD}_{\fraks}$,
and eliminate all incoming walls because they cancel each other.
In particular, $\frakD_{\fraks}$ is equivalent to $\frakD\sqcup \frakC$
and $\frakC$ has no incoming wall.

We show that $\frakC$ is trivial
in the sense of Definition \ref{3defn:trivial1}.
Let $\frakC_{\l}$ be the reduction of $\frakC$ at $\l$.
Then, we prove the following equivalent claim:
\par
\noindent
{\bf Claim.}
For any general $z\in M_{\bbR}$, we have
\begin{align}
g_z(\frakC_{\l})\equiv \rmid
\mod G^{>\l}.
\end{align}

\smallskip
We show the claim by the induction on $\l$.
The claim is true for $\l=1$, because $(\frakD_{\fraks})_1=\frakD_1=\mathrm{In}(\fraks)$,
so that $\frakC_1=\emptyset$.
Now assume that Claim holds for  $\l$.
Take any perpendicular joint $\frakj$ of $\frakD_{\l+1}$.
By the induction assumption,
for any general $z\in M_{\bbR}$,
$g_z(\frakC_{\l+1})$ is central in $G^{\leq \l+1}$.
Thus, we have
\begin{align}
\frakp_{\gamma_{\frakj}, \frakC_{\l+1}}
\equiv 
\frakp_{\gamma_{\frakj}, (\frakD_{\fraks})_{\l+1}}
\frakp_{\gamma_{\frakj}, \frakD_{\l+1}}^{-1}
\mod
G^{>\l+1}.
\end{align}
On the other hand,
by the consistency of $\frakD_{\fraks}$,
we have $\frakp_{\gamma_{\frakj}, (\frakD_{\fraks})_{\l+1}}
\equiv \rmid$ mod $G^{>\l+1}$.
Also, by the property (c) of $\frakD$, 
we have $\frakp_{\gamma_{\frakj}, \frakD_{\l+1}}
\equiv \rmid$ mod $G^{>\l+1}$.
It follows that 
\begin{align}
\label{3eq:pG1}
\frakp_{\gamma_{\frakj}, \frakC_{\l+1}}
\equiv \rmid
\mod
G^{>\l+1}.
\end{align}
Take any $n\in N_{\rmpr}^+$ 
such that $\frakj\subset n^{\perp}$.
Suppose that $\gamma_{\frakj}$ intersects
$n^{\perp}$ at only two general $z_1$ and $z_2$,
and that they locate in the opposite side of $\frakj$ on $n^{\perp}$ each other.
Then, the  congruence \eqref{3eq:pG1} implies that $g_{z_1}(\frakC)\equiv g_{z_2}(\frakC)$
mod $
G^{>\l+1}$.
Thus, $n^{\perp}$ is decomposed into walls
only by parallel joints so that $g_{z}(\frakC)$ is constant mod $
G^{>\l+1}$ on each wall.
Then, each wall should contain the  line $\bbR p^*(n)$.
Thus, 
if $g_z(\frakC_{\l})\not\equiv \rmid$ mod 
$G^{>\l+1}$ on some wall, the wall is incoming.
This is a contradiction.
Therefore, Claim holds for $\l+1$.
\end{proof}

The above construction of a CSD 
also reveals the following fundamental result
on the structure of a CSD,
which is implicit in \cite{Gross14}.
\begin{thm}
[{cf.\ \cite[Appendix C.3]{Gross14}}]
\label{3thm:struct1}
For any admissible loop $\gamma$ for $\frakD_{\fraks}$,
the consistency relation  $\frakp_{\gamma, \frakD_{\fraks}}=\rmid$ is reduced to a trivial one
(i.e., $a=a$)
by applying the  relations in Proposition \ref{3prop:pent1}
possibly infinitely many times.
\end{thm}
\begin{proof}
Consider the reduction $\frakD_{\fraks,\ell+1}$ at $\ell+1$.
We may concentrate on
 the consistency relation for an admissible loop
around a joint in  $\frakD_{\fraks,\ell+1}$.
For a parallel joint $\gamma$,
the consistency relation around $\gamma$ has the form
\begin{align}
\label{3eq:triv1}
\begin{split}
&\quad \
\Psi[t'_{j'} n'_{j'}]^{s'_{j'}  \d(t'_{j'} n'_{j'})}
\cdots
\Psi[t'_1n'_1]^{s'_1\d(t'_1 n'_1)}
\\
&
=
\Psi[t'_1n'_1]^{s'_1\d(t'_1 n'_1)}
\cdots
\Psi[t'_{j'} n'_{j'}]^{s'_{j'}  \d(t'_{j'} n'_{j'})},
\end{split}
\end{align}
where the left hand side is an anti-ordered product.
By applying \eqref{3eq:com1}, it reduces to a trivial relation.
For a perpendicular joint $\gamma$, 
the consistency relation around $\gamma$  has the form
\begin{align}
\label{3eq:cong1}
\begin{split}
&\quad \ \Psi[t'_{j'} n'_{j'}]^{s'_{j'}  \d(t'_{j'} n'_{j'})}
\cdots
\Psi[t'_1n'_1]^{s'_1\d(t'_1 n'_1)}
\\
&\equiv
\Psi[t_j n_j]^{s_j  \d(t_j n_j)}
\cdots \Psi[t_1n_1]^{s_1\d(t_1 n_1)}
\mod
G^{>{\l+1}},
\end{split}
\end{align}
where the left and right hand sides are the ones in \eqref{3eq:C0} and \eqref{3eq:C1}, respectively.
The right hand side was constructed from the left hand side
by Algorithm \ref{3algo:ord1}, which only depends on  the relations \eqref{3eq:com1}
and  \eqref{3eq:pent1}.
Thus, by applying the relations in the reverse way,
the relation \eqref{3eq:cong1} is reduced to a trivial one.
\end{proof}

In summary,
a CSD consists only of the dilogarithm elements of the form $\Psi[tn]^{s \d(tn)}$ ($s\in \bbZ_{>0}$) ,
or $\Psi[tn]^{ \d(tn)}$ by decomposing them,
and
its consistency relations are generated by
the pentagon relation.

We also record the following properties of the CSD constructed here.
\begin {prop}
\label{3prop:CSD1}
The CSD $\frakD_{\fraks}$ constructed in Construction \ref{3const:CSD1}
is with minimal support.
Also, the normal vector of any outgoing wall is not $e_i$ $(i=1,\dots, r)$.
\end{prop}
\begin{proof}
The wall elements in the form \eqref{3eq:gpos1} never cancel each other due to the positivity of
$s$. Thus, for any general $z$  on a wall of $\frakD_{\fraks}$, $g_z(\frakD_{\fraks})\neq \rmid$ holds.
Also, any outgoing wall  of $\frakD_{\fraks}$ created in
\eqref{3eq:C1} has the normal vector other than $e_i$ due  to Algorithm \ref{3algo:ord1}.
(This is true also by Proposition \ref{3prop:exist3}.)
\end{proof}

\subsection{Rank 3 examples: finite type}
\label{3subsec:rank31}
Let us demonstrate how  Construction \ref{3const:CSD1} really works for  rank 3 CSDs.
To represent rank 3 CSDs, we use the projection of $\mathrm{Supp}(\frakD_{\fraks})$ to the unit 2-sphere $S^2$ in $M_{\bbR}\simeq
\bbR^2$,
then apply the standard stereographic projection of $S^2$ to the plane $\bbR^2$ centered at
the point $(1/\sqrt{3})(1,1,1) \in S^2$.
Then, the supports of  incoming walls (resp.\ outgoing walls) in $\frakD_{\fraks}$ are represented by circles (resp.\ arcs), respectively.
Below,  figures are only  topological, that is,  without caring the precise metrics.

Here we concentrate on rank 3 CSDs of \emph{finite type}.

\begin{figure}
\centering
\includegraphics{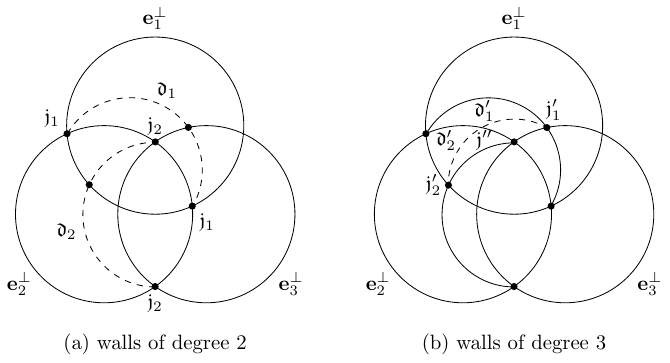}
\vskip-10pt
\caption{Construction of CSD of type $A_3$.
A filled circle $\bullet$ represents a perpendicular joint.
A dashed arc represents an  attached wall of a given degree.
}
\label{3fig:rank31}
\end{figure}

(a). Type $A_3$.
Let $B$ be the matrix in \eqref{3eq:Bmat1}.
Consider the case
\begin{align}
B=
\begin{pmatrix}
0 & -1 & 0\\
1 & 0 & -1\\
0 & 1 &0
\end{pmatrix},
\quad \d_1=\d_2=\d_3=1.
\end{align}
This is the simplest nontrivial example of a rank 3 CSD, where all perpendicular joints are
of type $A_2$ in Section \ref{3subsec:rank2}.
We follow  the procedure of Construction \ref{3const:CSD1} step by step.
See Figure \ref{3fig:rank31},
where only perpendicular joints are marked by $\bullet$ in figures.

At degree $1$, there are three incoming walls $\bfw_{e_i}$ ($i=1,2,3$)  as the initial input.
Correspondingly, we have three circles representing their supports  $\bfe_i^{\perp}$.
Note that the inside and the outside
 of each circle $\bfe_i^{\perp}$ corresponds to the half spaces
 \begin{align}
  \calH_i^{> 0}:=
\{ z \in M_{\bbR}\mid \langle e_i, z\rangle > 0\},
\quad
  \calH_i^{< 0}:=
\{ z \in M_{\bbR}\mid \langle e_i, z\rangle < 0\},
\end{align}
respectively.

Now, we construct  walls of degree $2$.
We have two perpendicular joints $\frakj_1=\bfe_1^{\perp}\cap \bfe_2^{\perp}=\bbR \bfe_3$,
 $\frakj_2=\bfe_2^{\perp}\cap \bfe_3^{\perp}=\bbR \bfe_1$ to take care of,
 and the other joint $\frakj_3=\bfe_1^{\perp}\cap \bfe_3^{\perp}=\bbR \bfe_2$
 is parallel.
Note that (only) these joints intersect the unit circle $S^2$ twice in $M_{\bbR}$.
Therefore, there are two points representing each joint  in the figure.
First, consider the joint $\frakj_1$. Since $\{\bfe_2, \bfe_1\}=1$,
the attached wall is
\begin{gather}
\bfw_1 = (\frakd_1, \Psi[(1,1,0)])_{(1,1,0)},
\\
\frakd_1:=\sigma(\frakj_1,-p^*((1,1,0)))=\sigma((0,0,1),(0,0,-1), (1,-1,-1)).
\end{gather}
where we note that $-p^*(\bfn)=-B\bfn$.
Since $\frakd_1$ is in the ``fourth quadrant'',
it is inside the circle $\bfe_1^{\perp}$
and outside the circle  $\bfe_2^{\perp}$
except for the boundary.
It is depicted in the dashed arc in Figure \ref{3fig:rank31} (a).
Similarly, for  the joint $\frakj_2$,  
we have $\{\bfe_3, \bfe_2\}=1$,
and
the attached wall is
\begin{gather}
\bfw_1 = (\frakd_2, \Psi[(0,1,1)])_{(0,1,1)},
\\
\frakd_2:=\sigma(\frakj_2,-p^*((0,1,1)))=\sigma((1,0,0),(-1,0,0), (1,1,-1)).
\end{gather}
This finishes the degree 2 walls.

Next, we construct  walls of degree $3$.
All joints that  appeared  above are already consistent.
Thus, we concentrate on the newly created perpendicular joints $\frakj'_1=\frakd_1\cap \bfe_3^{\perp}=\bbR_{\geq 0} (1,-1,0)$ 
and $\frakj'_2=\frakd_2 \cap\bfe_1^{\perp}=\bbR_{\geq 0} (0,1,-1)$ in 
Figure \ref{3fig:rank31} (b).
First, consider the joint $\frakj'_1$.
Since $\{\bfe_3, (1,1,0)\}=1$, the attached wall is
\begin{gather}
\bfw'_1 = (\frakd'_1, \Psi[(1,1,1)])_{(1,1,1)},
\\
\frakd'_1:=\sigma(\frakj'_1,-p^*((1,1,1)))=\sigma((1,-1,0), (1,0,-1)).
\end{gather}
Similarly, for the joint $\frakj'_2$,
we have $\{(0,1,1), \bfe_1\}=1$,
and
the attached wall is
\begin{gather}
\bfw'_2 = (\frakd'_2, \Psi[(1,1,1)])_{(1,1,1)},
\\
\frakd'_2:=\sigma(\frakj'_2,-p^*((1,1,1)))=\sigma((0,1,-1), (1,0,-1)).
\end{gather}
The supports $\frakd'_1$ and $\frakd'_2$ are in the same hyperplane $(1,1,1)^{\perp}$.
Moreover, they satisfy
\begin{align}
\frakj'':=
\frakd'_1\cap \frakd'_2 =
 \bbR_{\geq 0} (1,0,-1),
\quad
\frakd'_1\cup \frakd'_2=\sigma(\frakj'_1,\frakj'_2).
\end{align}
Thus,   $\frakd'_1$ and $\frakd'_2$ ``shake hands'' 
at $\frakj''$ to make a cone $\sigma(\frakj'_1,\frakj'_2)$.
(This shake-hands phenomenon is guaranteed to occur  due to the consistency proved in Lemma \ref{3lem:unique1}.)
Note that  the algorithm does not request us to joint $\bfw'_1$ and $\bfw'_2$ into one wall.
The cone $\sigma(\frakj'_1,\frakj'_2)$ is depicted in the dashed arc in
 Figure \ref{3fig:rank31} (b).
 Moreover, we have
 $\sigma(\frakj'_1,\frakj'_2) \cap e_2^{\perp}=\frakj''$,
 and  $\{(1,1,1), \bfe_2\}=0$.
 Thus, the joint $\frakj''$ is parallel (as it should be because $-p^*(1,1,1)\in \frakj''$ by
 construction).
 Therefore, the construction completes.
 Observe that
all normal vectors
\begin{align}
 \bfe_1,\ \bfe_2,\ \bfe_3,\ (1,1,0),\ (0,1,1),\ (1,1,1)
\end{align}
 exhaust the \emph{positive roots\/} of the root system of type $A_3$.

\begin{rem}
 As this is of finite type, $\mathrm{Supp}(\frakD_{\fraks})$ in  Figure  \ref{3fig:rank31} (b) coincides with
 the \emph{cluster complex\/} of a cluster pattern of type $A_3$ \cite{Fomin03a}.
 We give more accounts in Section \ref{3subsec:chamber1}.
 It also coincides with the \emph{picture\/} for the \emph{picture group\/} of type $A_3$ in \cite{Igusa15,Igusa16, Igusa20},
where the picture groups are defined based on the representation theory of (valued) quivers.
Moreover, it also coincides with the \emph{mutation fan\/} in \cite{Reading12, Reading17}.
The same remark applies to the following examples (b) and (c).
The appearance of the positive roots is explained by these theories.
Nevertheless, the phenomenon that the pentagon relation
reproduces root systems is remarkable and important, and it deserves to be  investigated further
from this new perspective.
\end{rem}

\begin{figure}
\centering
\includegraphics{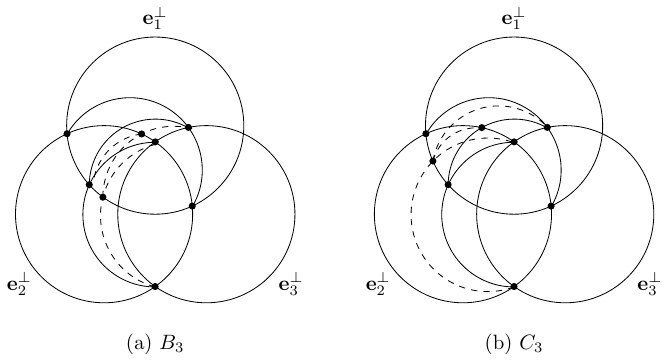}
\vskip-10pt
\caption{Construction of CSDs of types $B_3$ and $C_3$.
The solid arcs are the supports of type $A_3$,
while
the dashed arcs are the additional supports.
}
\label{3fig:rank32}
\end{figure}

(b). Type $B_3$. Consider the case
\begin{align}
B=
\begin{pmatrix}
0 & -1 & 0\\
1 & 0 & -1\\
0 & 2 &0
\end{pmatrix},
\quad
\d_1=\d_2=1,\ \d_3=2.
\end{align}
This is again of finite type.
The construction is similar to  Type $A_3$,
but it involves  joints of type $B_2$ in Section \ref{3subsec:rank2} as well.
For example, to the perpendicular joint $\frakj=\bfe_2^{\perp}\cap \bfe_3^{\perp}$,
 a wall $\bfw=(\frakd, \Psi[(0,1,1)]^2)_{(0,1,1)}$ of degree 2 and
a  wall $\bfw'=(\frakd', \Psi[(0,1,2)])_{(0,1,2)}$  of degree 3 are attached.
The details are left to the reader as an exercise.
The resulting CSD is depicted in Figure \ref{3fig:rank32} (a).
As observed therein, all supports of type $A_3$ in the case (a) appear again,
and they are indicated by solid arcs.
Newly added supports are indicated by dashed arcs,
where they have normal vectors $(0,1,2)$, $(1,1,2)$, $(0,2,2)$
from longer to shorter arcs.
All normal vectors 
\begin{align}
\bfe_1,\ \bfe_2,\ \bfe_3,\ (1,1,0),\ (0,1,1),\ (1,1,1),\
(0,1,2),\ (1,1,2),\ (1,2,2) 
\end{align}
exhaust the  positive roots of the root system of type $B_3$ as expected.

(c). Type $C_3$. Consider the case
\begin{align}
B=
\begin{pmatrix}
0 & -1 & 0\\
1 & 0 & -2\\
0 & 1 &0
\end{pmatrix},
\quad
\d_1=\d_2=2,\ \d_3=1.
\end{align}
This is again of finite type.
The construction is similar to  type $B_3$.
However, it involves joints of type $C_2$,
 instead of  type $B_2$,
in Section  \ref{3subsec:rank2}.
Note that this difference affects the configuration of attached walls.
Also, for joints of type $A_2$ we apply the pentagon relation \eqref{3eq:pent1}
with $c=1/2$.
The resulting CSD is depicted in Figure \ref{3fig:rank32} (b).
Once again, all supports  of type $A_3$ appear,
and they are indicated by solid arcs.
Newly added supports are indicated by dashed arcs,
where they have normal vectors $(0,2,1)$, $(2,2,1)$, $(1,2,1)$
from longer to shorter arcs.
All normal vectors
\begin{align}
 \bfe_1,\ \bfe_2,\ \bfe_3,\ (1,1,0),\ (0,1,1),\ (1,1,1),\
(0,2,1),\ (1,2,1),\ (2,2,1)
\end{align}
exhaust the  positive roots of the root system of type $C_3$. 
Moreover,
we see in Figure \ref{3fig:rank32} (a) and (b)  that 
the CSDs of type $B_3$ and $C_3$ are topologically identical
even though the process of the construction and the assignment
of normal vectors are  different.
This agrees with the well-known \emph{duality\/} between the cluster patterns of types $B_r$ and $C_r$ \cite{Fomin03a,Fomin03b},
and a more general result called the \emph{synchronicity\/} 
on cluster patterns in \cite{Nakanishi19}; namely, any pair of cluster patterns with
the initial exchange matrices $B$ and $-B^{T}$ share the same structure.

The  examples (a)--(c) exhaust all rank 3 CSD of  finite type with indecomposable $B$
up to mutation-equivalence in \cite{Fomin03a}.
We will continue to study examples of rank 3 CSDs of infinite type
in Section \ref{3subsec:rank32}
after clarifying the \emph{chamber structure\/} of CSDs.

\notes

The construction of a CSD
in Construction \ref{3const:CSD1}
is  a hybrid of the construction by
 \cite[\S10]{Kontsevich06} and \cite[\S3]{Gross07}
 presented in
\cite[Appendix C.1]{Gross14}
and the proof of the positivity in Theorem \ref{3thm:pos1}
in \cite[Appendix C.3]{Gross14}.
 Algorithm \ref{3algo:ord1},
  which is based on the pentagon relation 
 in Proposition \ref{3prop:pent1},
 substitutes for the \emph{perturbation trick\/}
 and the \emph{change of lattice trick\/}
 in \cite[Appendix C.3]{Gross14}.
 Theorem \ref{3thm:struct1} is
 implicit 
 in \cite[Appendix C.3]{Gross14}.

\newpage
\section{Mutation invariance of CSDs}
\label{3sec:mutation1}

In this section we study the mutation invariance of CSDs,
which is important to make the connection with cluster patterns.

\subsection{Mutations of seed and CSD}

Throughout the section, 
we fix $k\in \{1,\dots,r\}$ arbitrarily.

For any integer $a$, we define
\begin{align}
[a]_+=
\begin{cases}
a & a \geq 0,\\
0 & a <0.
\end{cases}
\end{align}
We have a useful formula
\begin{align}
\label{3eq:a+1}
[a]_+ - [-a]_+ = a.
\end{align}
The following notion originates in \cite[\S1.2]{Fock03}.
\begin{defn}[Seed mutation] \index{mutation!seed (for a fixed data)}
Let $\fraks=(e_1,\dots,e_r)$ be a seed 
for a given fixed data $\Gamma$,
and let $k\in \{1,\dots, r\}$.
We define a new seed $\mu_k(\fraks)=
\fraks'=(e'_1,\dots,e'_r)$ for $\Gamma$
by
\begin{align}
\label{3eq:mut1}
e'_i=
\begin{cases}
-e_k & i=k,\\
e_i+ [\{ \d_k e_k,e_i\}]_+ e_k
& i \neq k,
\end{cases}
\end{align}
which is called the \emph{mutation of $\fraks$ in direction $k$}.
\end{defn}

\begin{lem}
The above $\fraks'$ is indeed a seed for $\Gamma$.
\end{lem}
\begin{proof}
The transformation \eqref{3eq:mut1}  is written in the following matrix form,
\begin{align}
\label{3eq:etrans1}
(e'_1, \dots, e'_r)=
(e_1, \dots, e_r)
\begin{pmatrix}
1 & \\
& \ddots \\
[b_{k1}]_+ &\cdots & -1 &
\cdots &
[b_{kr}]_+
\\
&&&\ddots\\
&&&&1
\end{pmatrix},
\end{align}
where we recall that $b_{ki}=\{\d_ke_k, e_i\}$.
The matrix is unimodular.
Thus, $\fraks'$ is a basis of $N$.
Also, the basis $(\d_1e_1,\dots, \d_r e_r)$ of $N^{\circ}$ transforms as
\begin{align}
\label{3eq:mut2}
\d_i e'_i=
\begin{cases}
-\d_k e_k & i=k,\\
\d_i e_i+ [\{  e_k, \d_i e_i\}]_+ \d_k e_k
& i \neq k.
\end{cases}
\end{align}
Then, $(\d_1e'_1,\dots, \d_r e'_r)$ is a basis of $N^{\circ}$ by the same reason
as above.
\end{proof}

\begin{rem}
\label{3rem:inv1}
The mutation $\mu_k$ above is not involutive.
Namely, $\mu_k(\fraks')\neq \fraks$ in general.
Indeed, the inverse of \eqref{3eq:mut1} is given by
\begin{align}
\label{3eq:mutinv1}
e_i=
\begin{cases}
-e'_k & i=k,\\
e'_i + [\{ - \d_k e'_k,e'_i\}]_+ e'_k
& i \neq k.
\end{cases}
\end{align}
\end{rem}

Even though $\mu_k$ is not involutive,
it is ``involutive'' in the following sense.
\begin{lem}

Let $\fraks''=\mu_k(\fraks')$. Let  
$B_{\Gamma,\fraks}$ be the matrix in \eqref{3eq:Bmat1}.
Then, we have
$B_{\Gamma,\fraks''}=B_{\Gamma,\fraks}$.
Thus, the associated groups
$G_{\Gamma,\fraks''}$ and $G_{\Gamma,\fraks}$
are isomorphic.
\end{lem}
\begin{proof}
By \eqref{3eq:mut1}, we have
\begin{align}
e''_i=
\begin{cases}
e_k & i=k,\\
e_i+ \{ \d_k e_k,e_i\} e_k
& i \neq k.
\end{cases}
\end{align}
Then, it is easy to check the equality 
$B_{\Gamma,\fraks''}=B_{\Gamma,\fraks}$.
The last claim follows from
Proposition \ref{3prop:equiv1}.
\end{proof}

Under the change of a basis  of $N$ from $\fraks$ to   $\fraks'$,
we have the following coordinate transformation.
\begin{lem}
\label{3lem:ct1}
For $n=\sum_{i=1}^r  a_i e_i \in N$,
let $n=\sum_{i=1}^r  a'_i e'_i$.
Then, we have
\begin{align}
\label{3eq:aa1}
a'_i=
\begin{cases}
\displaystyle
-a_k + \sum_{j=1}^r  [\{\d_k e_k, e_j\}]_+ a_j& i=k,
\\
a_i & i \neq k.
\end{cases}
\end{align}
\end{lem}
\begin{proof}
Let $A$ be the matrix in \eqref{3eq:etrans1}.
Then, $A^2=I$ holds.
Thus, the column vector $\bfa=(a_1,\dots, a_r)$ transforms by $\bfa'=A\bfa$.
\end{proof}

Below we use the principal $x$-representation in
Section \ref{3sec:principal1} and its variation.
Thus,
a wall is expressed by  $\bfw=(\frakd,g)_n=[\frakd, f]_n$ 
interchangeably as specified  in Definition \ref{3defn:wallf1}.
As a general rule to remember, we apply $p^*$  for $\frakd$ and $n$,
while we apply
 $\tilde p^*$  for wall functions $f$.

For any  seed $\fraks=(e_1,\dots,e_r)$ and $k\in \{1,\dots,r\}$,
we define linear maps
\begin{align}
\label{3eq:Sk1}
&
\begin{matrix}
S_k=S_{k,\fraks}:& M_{\bbR}& \rightarrow &M_{\bbR}\\
& z & \mapsto&
z+ \langle \d_k e_k, z\rangle p^*(e_k),
\end{matrix}
\\
\label{3eq:Sk11}
&
\begin{matrix}
\tilde S_k=\tilde S_{k,\fraks}:& \tilde M^{\circ}& \rightarrow &\tilde M^{\circ}\\
& \tilde m & \mapsto&
\tilde m + \langle \d_k e_k,  m\rangle \tilde p^*(e_k),
\end{matrix}
\\
\label{3eq:Sk2}
&
\begin{matrix}
S_k^* =S_{k,\fraks}^*:& N& \rightarrow &N\\
& n & \mapsto&
n+  \{ \d_k e_k, n\} e_k,
\end{matrix}
\end{align}
where 
we continue to use the notation
  $\tilde m = m+n'$ ($m\in M^{\circ}$, $n'\in N$).
They are dual to each other in the following sense.
\begin{lem} 
\label{3lem:S*1}
(a). The following equalities hold.
\begin{align}
\langle S_{k}^*(n), S_{k}(z) \rangle
&=
\langle n, z \rangle,\\
\langle S_{k}^*(n), \tilde S_{k}(\tilde m) \rangle
&=
\langle n, \tilde m \rangle
=
\langle n,  m \rangle,\\
S_{k}(p^*(n))
&=p^*(S_{k}^*(n)),
\\
\tilde S_{k}(\tilde p^*(n))
&=\tilde p^*(S_{k}^*(n)).
\end{align}
\par
(b). We have
\begin{align}
\label{3eq:Sk3}
S_k^*(\d_i e_i)
=
\begin{cases}
\d_k e_k & i =k,\\
\d_i e_i + \{ e_k, \d_i e_i \} \d_k e_k
&i \neq k.
\end{cases}
\end{align}
In particular, $S_{k}^*$ preserves $N^{\circ}$.
\par
(c).
For any $n\in N$,
we have
$n \in N^{\circ}$ if and only if 
$S^*_k(n) \in N^{\circ}$.
\par
(d). Let $\d(n)$ $(n\in N)$ be the normalization factor in Definition \ref{3defn:normal1}. Then, we have
\begin{align}
\d(S^*_k(n) )=\d(n).
\end{align}
\end{lem}
\begin{proof}
(a).
We have
\begin{align}
\begin{split}
\langle S_{k}^*(n), S_{k}(z) \rangle
&=
 \langle n+  \{ \d_k e_k, n\} e_k,
z+ \langle \d_k e_k, z\rangle p^*(e_k)
\rangle
\\
&=\langle n, z \rangle 
+ \langle \d_k e_k, z\rangle\{ n,e_k\} 
+ \{ \d_k e_k, n\}   \langle e_k, z\rangle\\
&=\langle n, z \rangle.
\end{split}
\\
\label{3eq:ss1}
\begin{split}
\langle S_{k}^*(n), \tilde S_{k}(\tilde m ) \rangle
&=
 \langle n+  \{ \d_k e_k, n\} e_k,
\tilde m + \langle \d_k e_k,  m \rangle \tilde p^*(e_k)
\rangle
\\
&=\langle n, \tilde m \rangle
+ \langle \d_k e_k,  m \rangle \{ n,e_k\}
+ \{ \d_k e_k, n\}   \langle e_k, \tilde m \rangle\\
&=\langle n,  \tilde m \rangle =\langle n,   m \rangle.
\end{split}
\\
\begin{split}
S_{k}(p^*(n))
&=
p^*(n) + 
 \langle \d_k e_k, p^*(n) \rangle p^*(e_k)
 \\
&=
p^*(n + 
 \{ \d_k e_k, n\} e_k) 
=
p^*(S_{k}^*(n)).
\end{split}
\\
\label{3eq:ss2}
\begin{split}
\tilde S_{k}(\tilde p^*(n))
&=
\tilde p^*(n) + 
 \langle \d_k e_k, \tilde p^*(n) \rangle \tilde p^*(e_k)
 \\
&=
\tilde p^*(n + 
 \{ \d_k e_k, n\} e_k) 
=
\tilde p^*(S_{k}^*(n)).
\end{split}
\end{align}
\par
(b). This is clear
by \eqref{3eq:Sk2}.
\par
(c). The only-if-part follows from (b). 
The inverse of $S_k^*$ is given by
\begin{align}
\label{3eq:invS1}
S_k^*{}^{-1}:N \rightarrow N,
\quad
n \mapsto n- \{ \d_k e_k, n\} e_k,
\end{align}
which also preserves $N^{\circ}$.
 Thus,  the if-part also holds.
 \par
 (d). This follows from (c).
\end{proof}

We define the half spaces of $M_{\bbR}$
\begin{align}
\label{3eq:Hk1}
 \calH_{k}^+&= \calH_{k,\fraks}^+:=\{ z\in M_{\bbR} \mid \langle  e_k, z\rangle\geq 0\},
 \\
\calH_{k}^-&= \calH_{k,\fraks}^-:=\{ z\in M_{\bbR} \mid \langle  e_k, z\rangle\leq 0\}.
\end{align}
Then, we define a piecewise-linear transformation
\begin{align}
\label{3eq:Tk1}
\begin{matrix}
T_{k}=T_{k,\fraks}:& M_{\bbR}& \rightarrow &M_{\bbR}\\
& z & \mapsto& 
z+ [\langle \d_k e_k, z\rangle]_+ p^*(e_k)
=
\begin{cases}
S_{k}(z) & z \in  \calH_{k}^+,\\
z & z \in  \calH_{k}^-.
\end{cases}
\end{matrix}
\end{align}

\begin{defn}[Mutation of CSD] \index{mutation!of a cluster scattering diagram}
\label{3defn:mutD1}
Let   $\frakD_{\fraks}$ be
the CSD for a seed $\fraks$  in  Construction \ref{3const:CSD1},
thus, satisfying the properties in Proposition \ref{3prop:CSD1}.
Any wall 
of $\frakD_{\fraks}$ has the form
\begin{align}
(\frakd, g)_n=\ (\frakd, \Psi[tn]^{s \d(tn)})_n =[\frakd, (1+x^{\tilde p^*(tn)})^s]_n.
\end{align}
Then, we define
\begin{align}
\label{3eq:Skw1}
S_k(\frakd, g)_n:=&\ (S_k(\frakd), \Psi[S_k^*(tn)]^{s \d(S_k^*(tn))})_{S_k^*(n)}
\\
=&\ [S_k(\frakd), (1+x^{\tilde S_k(\tilde p^*(tn))})^s]_{S_k^*(n)}.
\end{align}
Note that for $n\in N^+$, $S_k^*(n)$ does not necessarily belong to $N^+$.
However, we temporarily ignore this fact.
For each $k \in \{1,\dots, r\}$,
let $T_k(\frakD_{\fraks})$ be  a collection of ``walls''
 obtained from $\frakD_{\fraks}$ as follows:
 \par
(a).
For each wall $\bfw=(\frakd, g)_n$ of $\frakD_{\fraks}$
other than $\bfw_{e_k}=(e_k^{\perp}, \Psi[e_k]^{\d_k})_{e_k}$,
let $\frakd_+=\frakd\cap \calH_{k}^+$,
and $\frakd_-=\frakd\cap \calH_{k}^-$.
We split the wall $\bfw$ into two walls
\begin{align}
\label{3eq:wallpm1}
\bfw_+=(\frakd_+, g)_n,
\quad
\bfw_-=
(\frakd_-, g)_n,
\end{align}
where, if one of  $\frakd_+$ and  $\frakd_-$ is codimension more than $1$,
we throw that part away.
Then,
we replace $\bfw_+$ with a new wall
$S_k(\bfw_+)$ in \eqref{3eq:Skw1},
while we leave $\bfw_-$ as it is.
\par
(b). Replace $\bfw_{e_k}$ with  $\bfw_{e'_k}=(e'_k{}^{\perp}, \Psi[e'_k]^{\d_k})_{e'_k}$. (Clearly,  $e'_k=-e_k \not\in N^+$.)

We call $T_k(\frakD_{\fraks})$
 the \emph{mutation of $\frakD_{\fraks}$ in direction $k$}.
\end{defn}
By construction, we have
\begin{align}
\label{3eq:supp2}
\mathrm{Supp}(T_k(\frakD_{\fraks}))
=T_k(\mathrm{Supp}(\frakD_{\fraks})).
\end{align}

Let 
 $n=\sum_{i=1}^r  a_i e_i \in N_{\fraks}^+$
 be the normal vector of any wall in $\frakD_{\fraks}$.
 Then, by Lemma \ref{3lem:S*1}, the normal vector $n'$ of the corresponding wall
 in $T_k(\frakD_{\fraks})$ 
defined above is given by
\begin{align}
\label{3eq:n'1}
n'=
\sum_{i=1}^r  b_i e_i
=
\begin{cases}
-e_k
& n = e_k
\\
n+  \{ \d_k e_k, n\} e_k.
& \text{$n$ is the normal vector of $\frakd_+$},
\\
n 
& \text{$n$ is the normal vector of $\frakd_-$}.
\end{cases}
\end{align}

\begin{lem}
\label{3lem:n'2}
Let $n'$ be as above.
Then, 
under the change of the coordinate transformation
$
n'=
\sum_{i=1}^r  b_i e_i
=
\sum_{i=1}^r  b'_i e'_i
$,
we have
$b'_i=a_i$ for $i\neq k$,
while
\begin{gather}
\label{3eq:aa2}
b'_k=
\begin{cases}
1
& n = e_k,
\\
\displaystyle
-a_k+\sum_{j=1}^r[\mp \{\d_ke_k, e_j\}]_+a_j
& \text{$n$ is the normal vector of $\frakd_{\pm}$}.
\end{cases}
\end{gather}
\end{lem}
\begin{proof}
This follows from \eqref{3eq:n'1} and Lemma \ref{3lem:ct1},
where we  use
\eqref{3eq:a+1} for the second case.
\end{proof}

The goal of this section to prove the following theorem,
which is the second fundamental result on cluster scattering diagrams.

\begin{thm}[{Mutation invariance of CSD \cite[Theorem 1.24]{Gross14}}] \index{mutation invariance!of a cluster scattering diagram}
\label{3thm:mut1}
The collection of walls $T_k(\frakD_{\fraks})$
is a consistent scattering diagram for the seed
$\fraks'=\mu_k(\fraks)$.
Moreover, $T_k(\frakD_{\fraks})$
is equivalent to a CSD $\frakD_{\fraks'}$
for $\fraks'$.
\end{thm}

We see in \eqref{3eq:aa2} that
$b'_k <0$ may occur.
Therefore,
for $T_k(\frakD_{\fraks})$
to be a scattering diagram for $\fraks'$,
$g_z(\frakD_{\fraks})=\rmid$ should hold
if $z$ is on a wall of $\frakD_{\fraks}$ such that $b'_k<0$.
This property is  nontrivial,
and so far we have no direct proof.
(In the rank 2 case, however, we can directly prove it.
See Example \ref{3ex:prohibit1}.)
Below we follow the strategy 
by \cite{Gross14} to avoid this difficulty.
The main idea of  \cite{Gross14}  is that we temporarily 
extend the notion of scattering diagrams
by accepting
the walls with $b'_k<0$,
then prove the uniqueness property of a CSD $\frakD_{\fraks'}$
\emph{up to equivalence
among such extended scattering diagrams}.

\subsection{Extended structure group}

As mentioned above,
following the strategy by \cite{Gross14}, we temporarily loosen
the condition  $n\in N_{\rmpr}^+$ for the normal vector $n$ of a wall.
We  define
\begin{align}
\label{3eq:N+2}
\overline{N}^{+} = \overline{N}_{k,\fraks}^{+} :=
\biggl\{ \sum_{i=1}^r a_i e_i \mid a_i \in \bbZ_{\geq 0} \text{\ for $i\neq k$},\ 
a_k\in \bbZ,\
\sum_{\ss i=1
\atop \ss i\neq k}^n a_i > 0\biggr\}.
\end{align}
In other words, we allow  negative coefficients of $e_k$ in $n$,
while we exclude $n= j e_k$ ($j\in \bbZ$).
Let $\overline{N}_{\rmpr}^{+}$ be the set of all primitive elements
in $\overline{N}^{+} $.
\begin{lem}
\label{3lem:Nbar1}
The following facts hold:
\par
(a).
$e_k \not\in \overline{N}_{\rmpr}^{+}$.
\par
(b).
For $n\in N_{\rmpr}^+$ with $n\neq e_k$, we have
$n, S^*_k(n) \in \overline{N}_{\rmpr}^{+}$.
\end{lem}
\begin{proof}
The only nontrivial claim is that $S^*_k(n) $ is primitive in (b).
If it is not, $S^*_k(n)=j n' $ for some $j\in \bbZ_{>1}$ and $n'\in N$.
Then, applying $S_k^*{}^{-1}$ in \eqref{3eq:invS1},
we see that $n= S_k^*{}^{-1}(jn')=j S_k^*{}^{-1}(n')$ is not primitive.
\end{proof}

We define a degree function
\begin{align}
\begin{matrix}
\overline{\deg}:&\overline{N}^{+}& \rightarrow &\bbZ_{>0} \\
&
\displaystyle
\sum_{i=1}^r a_i e_i
& \mapsto
&
\displaystyle
\sum_{\ss i=1\atop \ss i\neq k}^r a_i.
\end{matrix}
\end{align}
Define an $\overline{N}^{+}$-graded Lie algebra  $\overline\frakg$,\begin{align}
\overline\frakg=\bigoplus_{n\in \overline{N}^{+}}
\frakg_n,
\end{align}
which  is parallel to $\frakg$ in   \eqref{3eq:lieg1}.
Then, we repeat the construction of 
$\widehat{g}$ and $G$ in Section \ref{3sec:scat1},
where
 $N^+$ and $\deg(n)$ therein are replaced with
 $\overline{N}^{+}$ and $\overline{\deg}(n)$.
Note that  the corresponding Lie algebras $\overline\frakg^{\leq \l}$
 are still nilpotent.
 Therefore, the construction works,
and  we obtain the completion $\widehat{\overline\frakg}$
 and the corresponding exponential group $\overline{G}$.
 The reason to exclude $e_k$ from $\overline{N}^{+}$
 is now clear, because it is incompatible with the completion of $\overline \frakg$.
 We call $\overline{G}$ the  \emph{extended structure group for $\frakD_{\fraks}$ and $k$}.
 \index{structure group!extended}

This, however, is not yet enough for our purpose.
Let   $G(\Psi[e_k])$ be the group freely generated by $\Psi[e_k]$.
Since 
 $\Psi[e_k]\not\in \overline{G}$,
 the product
 of $\Psi[e_k]$ and
elements of $ \overline{G}$ is not defined in  $\overline{G}$.
Instead of further extending $\overline{G}$,
 we  extend the algebra $\bbk[[\tilde P]]$ so that 
 both $G(\Psi[e_k])$ and $\overline{G}$
 act faithfully on it.

First,  we define a monoid $\overline  P\subset \tilde M^{\circ}$ as follows:
\begin{itemize}
\item[(i).]
$\overline  P=\overline \sigma \cap \tilde M^{\circ}$,
where $\overline\sigma $ is a $2r$-dimensional (not strongly convex) cone in $\tilde M_{\bbR}$
such that $\overline\sigma \cap (-\overline\sigma)=\bbR \tilde p^*(e_k)$.
\item[(ii).]
$\tilde p^*(e_i), \cdots, \tilde p^*(e_r)\in \overline P$.
\end{itemize}
We may take $\tilde P$ in Section \ref{3subsec:principal1} as a submonoid of $\overline P$
in the positive side of $\tilde p^*(e_k)$.
 Let $\bbk[\overline{P}]$ be the monoid ring of $\overline{P}$ over $\bbk$.
 Then, the  group of units of $\overline{P}$ is given by $\overline{P}^{\times}=\overline{P}\cap \bbR \tilde p^*(e_k)$.
Let $\overline{J}$
  be the maximal ideal of $\bbk[\overline{P}]$
generated by $\overline{P} \setminus
\overline{P}^{\times}$,
and let ${\bbk[[\overline{P}]]}$
be the completion with respect to $\overline{J}$.
Furthermore, we consider the localization
$R_k:={\bbk[[\overline{P}]]}_{f[e_k]}$
of  ${\bbk[[\overline{P}]]}$ by $f[e_k]=1+x^{\tilde p^*(e_k)}$.
Namely, any element $f\in R_k$
is written as 
\begin{align}
\label{3eq:f2}
f=
\sum_{\tilde m\in \overline{P}} 
b_{\tilde m}
(1+x^{\tilde p^*(e_k)})^{-a_{\tilde m}}
 x^{\tilde m}
\quad
(a_{\tilde m}\in \bbZ_{\geq 0},\
b_{\tilde m} \in \bbk),
\end{align}
where in the infinite  sum
there are only finitely many $j\in \bbZ$ such that $b_{\tilde m + j \tilde p^*(e_k)}\neq 0$
for each $\tilde m\in \overline{P}$.
We define the action of $\Psi[n]\in \overline{G}$ ($n\in \overline{N}^+$) 
and $\Psi[e_k]\in G(\Psi[e_k])$
on $R_k$
in the usual way as
\begin{align}
\begin{split}
&\Psi[n](f)
=\sum_{\tilde m\in \overline{P}} 
b_{\tilde m}
(1+x^{\tilde p^*(e_k)} (1+x^{\tilde p^*(n)})^{\langle n, \tilde p^*(e_k) \rangle} )^{-a_{\tilde m}}
 x^{\tilde m}(1+x^{\tilde p^*(n)})^{\langle n, m \rangle}
\\
&\quad =\sum_{\tilde m\in \overline{P}} 
b_{\tilde m}
(1+x^{\tilde p^*(e_k)})^{-a_{\tilde m}}
  \biggl(
  1+
   x^{\tilde p^*(e_k)}
  \frac
{  (1+x^{\tilde p^*(n)})^{\langle n, \tilde p^*(e_k) \rangle} -1}
{1+x^{\tilde p^*(e_k)}}
  \biggr)^{-a_{\tilde m}}
\\
&\qquad\qquad\quad
\times
 x^{\tilde m}(1+x^{\tilde p^*(n)})^{\langle n, m \rangle},
 \end{split}
\\
&\Psi[e_k](f)
=\sum_{\tilde m\in \overline{P}} 
b_{\tilde m}
(1+x^{\tilde p^*(e_k)})^{-a_{\tilde m}}
 x^{\tilde m}(1+x^{\tilde p^*(e_k)})^{\langle e_k, m \rangle}.
\end{align}
The right hand sides are in $R_k$.
Moreover, the actions are faithful as Proposition
\ref{3prop:gaction2}.
Thus, we have the following result.
\begin{lem}
\label{3lem:prod1}
Both groups $G(\Psi[e_k])$ and $\overline{G}$ 
act faithfully as automorphisms of  $R_k$.
In particular,  the product of $\Psi[e_k]$
and $g\in \overline{G}$ is
defined in 
$\mathrm{Aut}(R_k)$.
\end{lem}

The wall elements $\Psi[e_k]^{\pm1}$ appear
only for the consistency relation around a joint $\frakj\in e_k^{\perp}$.
Moreover, the consistency relation has the form
\begin{align}
g_1 \Psi[e_k]^{\pm \d_k} g_2 \Psi[e_k]^{\mp \d_k}=\rmid
\quad (g_1,g_2\in \overline G_{\frakj}),
\end{align}
where $G_{\frakj}$ is the one in the following lemma.
Then, the lemma tells that the consistency relation essentially
reduces to the one in $\overline{G}$.
\begin{lem}[{\cite[Lemma C.8]{Gross14}}]
\label{3lem:ad1}
Let $n \in \overline N_{\rmpr}^+$
and  $\frakj = e_k^{\perp} \cap n^{\perp}$.
Let $N_{\frakj}$ be the one in \eqref{3eq:Njn1},
and let $\overline N_{\frakj,\rmpr}^+:= N_{\frakj}\cap \overline N_{\rmpr}^+$.
Let $\overline G_{\frakj}$ be the subgroup of $\overline G$ generated
by the parallel subgroups $\overline G_{n}^{\parallel}$ $(n\in \overline N_{\frakj,\rmpr}^+)$.
Then, 
 we have
\begin{align}
\label{3eq:ad1}
&\Psi[e_k]^{\d_k} \overline G_{\frakj} \Psi[e_k]^{-\d_k}\subset \overline G_{\frakj},
\quad \mbox{if $\{e_k, n\}\geq 0$},
\\
&\Psi[e_k]^{-\d_k} \overline G_{\frakj} \Psi[e_k]^{\d_k}\subset \overline G_{\frakj},
\quad \mbox{if $\{e_k, n\}\leq 0$}.
\end{align}
\end{lem}
\begin{proof}
This is clear if $\{e_k, n\}=0$.
Suppose that $\{e_k, n\}\neq 0$.
Then, we obtain the same formula in Proposition \ref{3prop:bound1} (a) by working with the action on $R_k$.
The claim immediately follows from it.
\end{proof}

For the seed $\fraks'=\mu_k(\fraks)=(e'_1,\dots, e'_r)$,
we repeat the same construction as above.
Then, despite the asymmetry of the mutation
  in Remark \ref{3rem:inv1},
  the following result holds,
  where $'$ indicates the corresponding ones for $\fraks'$:
  \begin{lem}
  \label{3lem:nn'1}
  \begin{gather}
  \overline{N}^+=\overline{N}'{}^+,
  \quad
  \overline{\frakg}=  \overline{\frakg}',
  \quad
  \overline{G}=  \overline{G}'.
  \end{gather}
Moreover, $\overline{P}$ serves as $\overline{P}'$ for $\fraks'$
so that
  \begin{gather}
  \label{3eq:pp'1}
   {\bbk[[\overline{P}]]}
   =
   {\bbk[[\overline{P}']]}, 
   \quad
R_k
   =
R'_k.
     \end{gather}
  \end{lem}
  \begin{proof}
For example, the equality $  \overline{N}^+=\overline{N}'{}^+$ follows from
\eqref{3eq:mut1} and \eqref{3eq:N+2}.
The last equality in \eqref{3eq:pp'1} follows from the identity
\begin{align}
1+x^{\tilde p^*(e_k)}
=x^{-\tilde p^*(e'_k)}(1+x^{\tilde p^*(e'_k)}).
\end{align}
  \end{proof}

\subsection{Extended scattering diagrams}

\begin{defn}[Extended scattering diagram] \index{scattering diagram!extended}
For a seed $\fraks$ and $k\in \{1,\dots, r\}$,
an \emph{extended scattering diagram\/}
$\overline{\frakD}$ for $(\fraks,k)$
is a collection of walls satisfying the following properties:
\begin{itemize}
\item A wall is either $\bfw_{e_k}=(e_k^{\perp}, \Psi[e_k]^{\d_k})_{e_k}$, or
$\bfw=(\frakd, g)_n$ with $n\in \overline{N}^{+}_{\rmpr}$,
$\frakd \subset n^{\perp}$, $g\in \overline{G}_{n}^{\parallel}$,
where $\frakd$ is a cone of codimension 1 as usual.
(In particular, $n\neq e_k$.)
\item It satisfies the  finiteness condition
with respect to $\overline{G}$ and $\overline\deg$
for walls other than $\bfw_{e_k}$.
\item There is at most one $\bfw_{e_k}$.
\end{itemize}
Thanks to Lemma \ref{3lem:prod1},
for any admissible curve $\gamma$ in $\overline \frakD$,
the path-ordered product $\frakp_{\gamma, \overline\frakD}$ is defined 
in $\mathrm{Aut}(R_k)$
in the same way as before.
Thus, the equivalence and the consistency of extended scattering diagrams
and other related notions
 are defined in the same way as before.
 (In particular, the consistency condition around a  joint $\frakj\subset e_k^{\perp}$ makes
 sense in $\overline G$  thanks to Lemma \ref{3lem:ad1}.)
\end{defn}

\begin{rem}
In \cite{Gross14} the distinguished wall $\bfw_{e_k}$ is excluded from walls,
and called the \emph{slab}.
\end{rem}

\begin{lem}
\label{3lem:es1}
(a).
Let   $\frakD_{\fraks}$ be
the CSD 
constructed in Construction \ref{3const:CSD1}.
Then,
 $\frakD_{\fraks}$ is a consistent extended scattering diagram for
$(\fraks,k)$.
\par
(b). 
For the same CSD $\frakD_{\fraks}$ as above,
the mutation $T_k(\frakD_{\fraks})$
is an extended scattering diagram for $(\fraks',k)$,
where $\fraks'=\mu_k(\fraks)$.
\end{lem}
\begin{proof}
(a). For any wall of $\frakD_{\fraks}$ other than $\bfw_{e_k}$,
its normal vector $n\in N_{\rmpr}^+$ is not $e_k$
by Proposition \ref{3prop:CSD1}.
Thus, $n$ belongs to $n\in \overline{N}^{+}_{\rmpr}$
by Lemma \ref{3lem:Nbar1} (b).

Next, we show that $\frakD_{\fraks}$
satisfies the finiteness condition with respect to $\overline{G}$.
This is not trivial, and it is equivalent to  the following claim:
\smallskip
\par\noindent
{\bf Claim.}
For each $n\in N^+$ with $n-e_k\not\in N^+$,
 there are only finitely many $j\in \bbZ_{>0}$
such that
$\Psi[n+je_k]^{s\d(n+je_k)}$ is a wall element
 of $\frakD_{\fraks}$.
\par
\smallskip
 Let us prove Claim by the induction on $\ell=\overline\deg(n)$.
 Let $\ell=1$, namely, $n=e_i$ ($i\neq k$).
 Then, by Construction \ref{3const:CSD1}, 
such walls are created only from the perpendicular joint
$\frakj\subset e_k^{\perp}\cap e_i^{\perp}$
with an anti-ordered product of
 $\Psi[e_i]^{\d_i} $ and $\Psi[e_k]^{\d_k}$. 
By Proposition \ref{3prop:bound1} (b),
 $j$ is bounded  as
\begin{align}
\label{3eq:up1}
j\leq |\{\d_k e_k , e_i\}|.
\end{align}
Suppose that Claim holds up to $\ell$. Let 
$\overline\deg(n)={\ell}+1$.
Then,  a wall element of the form $\Psi[n+je_k]^{s\d(n+je_k)}$ is created 
in the following two ways:
\begin{itemize}
\item[(a).] It is created from a  perpendicular joint $\frakj\not \subset e_k^{\perp}$
with an anti-ordered product of $\Psi[n'+j'e_k]$'s with $\overline\deg(n')\leq \ell$.
\item[(b).] It is created from  a perpendicular joint $\frakj\subset e_k^{\perp}\cap n^{\perp}$
with an anti-ordered product of $\Psi[e_k]$ and $\Psi[n+j'e_k]$, where the latter is created in (a).
\end{itemize}
There are only finitely many such walls by the induction assumption and Proposition \ref{3prop:bound1} (b).
Thus, Claim holds for ${\ell}+1$.
Therefore, 
 $\frakD_{\fraks}$ is an extended scattering diagram for
$(\fraks,k)$.

Finally,  we show the consistency of $\frakD_{\fraks}$
as an extended scattering diagram for $(\fraks,k)$.
Recall that any wall element $\Psi[n]^{s \d(n)}$ 
of $\frakD_{\fraks}$ also acts on $R_k$.
Thus, for any admissible loop $\gamma$, the consistency relation $\frakp_{\gamma,\frakD_{\fraks}}=\rmid$
along $\gamma$
holds in $\mathrm{Aut}(R_k)$,
where the path-ordered product 
is defined  in the limit with respect to $\deg$.
 Let us temporally write 
 the path-ordered products along $\gamma$ 
 defined in the limit with respect to $\deg$ and $\overline\deg$
 as $\frakp_{\gamma}$ and $\overline \frakp_{\gamma}$, respectively.
 Then, it is enough to prove the equality
 \begin{align}
 \label{3eq:pcomp1}
 \frakp_{\gamma}=\overline \frakp_{\gamma}.
 \end{align}
For any integer $\ell>0$,
let $C_{\ell}$ be the largest $j$ for $n$ with $\overline\deg(n)=\ell$
in Claim above.
Then,
the following inequalities hold for any $n\in N^+$
(the first inequality is trivial):
\begin{align}
\overline\deg(n)\leq 
\deg(n)\leq
\overline\deg(n)+ C_{\overline\deg(n)}.
\end{align}
Thus,
the limit $\ell=\deg(n) \rightarrow \infty$ implies
the limit $\overline\ell=\overline\deg(n) \rightarrow \infty$,
and \emph{vice versa}.
Thus, we have the equality \eqref{3eq:pcomp1}.
Therefore, $\frakD_{\fraks}$ is consistent as an extended scattering diagram for
$(\fraks,k)$.
\par

(b). 
For any wall of $T_k(\frakD_{\fraks})$ other than $\bfw_{e_k}$,
its normal vector $S^*_k(n)$ belongs to $ \overline{N}^{+}_{\rmpr}$
by Lemma \ref{3lem:Nbar1} (b).
The finiteness condition with respect to 
$\overline{G}$
follows from the  finiteness condition
of $\frakD_{\fraks}$,
because $\overline \deg(S_k^*(n))=\overline \deg(n)$.
Thus,  $T_k(\frakD_{\fraks})$
is an extended scattering diagram for 
$(\fraks',k)$.
\end{proof}

We will soon show that $T_k(\frakD_{\fraks})$
is also consistent as an 
extended scattering diagram for $(\fraks',k)$.
We first prove
 the following uniqueness theorem,
 which is a key of the strategy by \cite{Gross14}.
The claim is parallel to Theorem
\ref{3thm:exist1} (b),
but it will be proved by a different argument.
\begin{prop}[{Cf.\ \cite[Theorem~1.28]{Gross14}}]
\label{3prop:ext1}
Let $\mathrm{In}_{\fraks}$ be the one in  \eqref{3eq:in2}.
Then,
any consistent  extended  scattering diagram $\overline{\frakD}$ for $(\fraks,k)$
with
\begin{align}
\label{3eq:incond1}
\overline\frakD_{\rmin}=\mathrm{In}_{\fraks}
\end{align}
is equivalent to a CSD $\frakD_{\fraks}$ as an extended scattering diagram.
\end{prop}
\begin{proof}
Any CSD $\frakD_{\fraks}$ satisfies the
  condition \eqref{3eq:incond1}
  by Lemma \ref{3lem:es1} and  Definition \ref{3defn:CSD1}.
Thus, we prove the following more general claim.

\par\noindent
{\bf Claim.} Two consistent extended  scattering diagrams $\overline{\frakD}$, $\overline{\frakD}'$
satisfying the condition \eqref{3eq:incond1} are equivalent.

Suppose that
both
 $\overline{\frakD}$ and  $ \overline{\frakD}'$  
 satisfy  the  condition  \eqref{3eq:incond1}.
 Let $\overline\frakC
=\overline{\frakD}' - \overline{\frakD}.
$
be an extended scattering diagram
defined in the same way as in the proof of Lemma \ref{3lem:unique1}.
Namely,
take all walls of $\overline{\frakD}$, invert
all wall elements,  add all walls of  $ \overline{\frakD}'$,
and eliminate all incoming walls because they cancel each other.
Then, we have  the following properties:
\begin{itemize}
\item
 $\overline{\frakD}\sqcup  \overline{\frakC}$ and $\overline{\frakD}'$
 are equivalent.
\item
All walls of $\overline\frakC$ are outgoing.
\item
All   walls of $\overline\frakC$
have normal vectors other than $e_k$.
\end{itemize}
Now we prove that, for any $n\in  \overline{N}_{\rmpr}^{+}$,
\begin{align}
\label{3eq:pigxC1} 
\pi_{\l}(
g_z  (\overline{\frakC}))=\rmid
\end{align}
holds for any general $z\in n^{\perp}$.
We prove it by the induction on  $\l$.
Let $ \overline{\frakC}_\l$ be the reduction of
$ \overline{\frakC}$ at $\l$.
We note the following:
\begin{itemize}
\item
A parallel result to Lemma \ref{3lem:center1} holds  for $\overline{G}$.
\end{itemize}

\begin{figure}
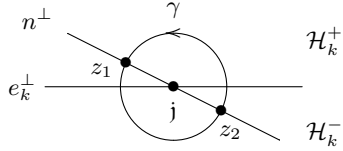

\centering
\leavevmode
\xy
(-20,0)*{\text{\small $e_k^{\perp}$}};
(0,-3)*{\text{\small $ \frakj$}};
(-18,9)*{\text{\small $ n^{\perp}$}};
(-9.5,2.2)*{\text{\small $z_1$}};
(7.5,-6)*{\text{\small $z_2$}};
(20,6)*{\text{\small $\calH_k^+$}};
(20,-6)*{\text{\small $\calH_k^-$}};
(0,10)*{\text{\small $\gamma$}};
(-6.3,3.15)*+{\bullet};
(0,0)*+{\bullet};
(6.3,-3.15)*+{\bullet};
(0,7)*+{};(-2,7)*+{}
          ?>*\dir{>};
(0,0)="A",
(0,0)*\xycircle<20pt>{}
\ar@{-} "A"+(-14,7); "A"+(14,-7)
\ar@{-} "A"+(17,0); "A"+(-17,0)
\endxy
\caption{Admissible loop $\gamma$ in the proof of Proposition 
\ref{3prop:ext1}.}
\label{3fig:adm1}
\end{figure}

First, consider the case $\l=1$.
Let us fix $n\in (\overline{N}_{\rmpr}^{+})^{\leq 1}$,
namely, $n=e_i+j e_k$ ($i\neq k)$ for some $j\in \bbZ$.
We consider the case $\{e_k,n\}\geq 0$. 
(We work with all such $n$ simultaneously.)
We first note that $p^*(n)\in \calH_{k}^+$
and the claim \eqref{3eq:pigxC1} holds for $\ell=1$
in some neighborhood of $p^*(n)$ due to the absence
of any incoming wall.
Let $z_1$ and $z_2$ be any general points in $n^{\perp}$.
We take   an admissible loop $\gamma$ 
 intersecting $n^{\perp}$ transversally only at $z_1$ and $z_2$. 
Let $\frakj=n^{\perp}\cap e_k^{\perp}$.
First, suppose that $\frakj$ is not inside the loop $\gamma$.
Then,  every nontrivial element
appearing in the products $\frakp_{\gamma, \overline\frakD_1}$ 
and $\frakp_{\gamma, \overline\frakD'_1}$ belong to $\overline G$.
Since $\overline{G}^{\leq 1}$ is abelian
and $ \overline{\frakD}$ and  $\overline{\frakD}'$ are consistent,
we have
\begin{align}
\label{3eq:pigxC3}
\pi_1(\frakp_{\gamma, \overline{\frakC}_1})=
\pi_1(\frakp_{\gamma, \overline{\frakD}'_1})
\pi_1(\frakp_{\gamma, \overline{\frakD}_1}^{-1})
=
\rmid.
\end{align}
Then, by the decomposition  of  $\overline{G}^{\leq 1}$,
which is parallel to Lemma   \ref{3lem:center1},
we have
\begin{align}
\label{3eq:pigxC2}
\pi_{1}(
g_{z_1}  (\overline{\frakC}))
= \pi_{1}(
g_{z_2}  (\overline{\frakC})).
\end{align}
It follows that $\pi_{1}(
g_{z}  (\overline{\frakC}))=\rmid$ holds in $\mathrm{Int}(\calH_k^+)$,
while $\pi_{1}(
g_{z}  (\overline{\frakC}))$ is constant in $\mathrm{Int}(\calH_k^-)$.

Next, suppose that $\frakj$ is inside the loop $\gamma$.
One can choose $z_1\in n^{\perp}\cap \mathrm{Int}(\calH_k^+)$ and $z_2\in n^{\perp}\cap \mathrm{Int}(\calH_k^-)$ sufficiently close to each other  so that $\gamma$ intersects
only walls  whose normal vectors are in 
$\overline N_{\frakj,\rmpr}^+$ of Lemma \ref{3lem:ad1}.
See Figure \ref{3fig:adm1}.
It follows from the result just obtained above that $\pi_{1}(
g_{z_1}  (\overline{\frakC}))=\rmid$ on $n^{\perp}\cap \mathrm{Int}(\calH_k^+)$.
Also,  any other $n'{}^{\perp}$ $(n'\in \overline N_{\frakj,\rmpr}^+)$ containing $\frakj$ also satisfies
$\{e_k, n'\}>0$.
It follows that
the factor $\Psi[e_k]^{ \d_k}C \Psi[e_k]^{- \d_k}\in \overline G$ in 
$\frakp_{\gamma, \overline\frakD_1}$ 
and $\frakp_{\gamma, \overline\frakD'_1}$ coming from the region $\calH_k^+$ coincide
modulo $\overline G^{>1}$.
Thus, we have \eqref{3eq:pigxC3} again,
and we obtain \eqref{3eq:pigxC2}.
Therefore, we conclude that
\eqref{3eq:pigxC1}  holds for $\ell=1$.
The other case  $\{e_k,n\}\leq 0$ is similar by exchanging $\calH_k^+$ and $\calH_k^-$.

Next, suppose that \eqref{3eq:pigxC1} holds up to $\l$ for any general $z\in M_{\bbR}$.
Let us fix $n\in (\overline{N}_{\rmpr}^{+})^{ \l+1} $.
We take an admissible loop $\gamma$ in a similar way to the above.
Then, by the induction assumption,
all factors of
$\pi_{\l+1}(\frakp_{\gamma, \overline{\frakC}_{\l+1}})$
belongs to $\Ker\, \pi_{\l+1,\l}$, which is abelian as in  Lemma \ref{3lem:center1}.
Then, by repeating the same argument as $\ell=1$,
 we have
\begin{align}
\pi_{\l+1}(\frakp_{\gamma, \overline{\frakC}_{\l+1}})=
\pi_{\l+1}(\frakp_{\gamma, \overline{\frakD}'_{\l+1}})
\pi_{\l+1}(\frakp_{\gamma, \overline{\frakD}_{\l+1}}^{-1})
=
\rmid.
\end{align}
Then, again by the decomposition  of  $\Ker\, \pi_{\l+1,\l}$,
which is parallel to Lemma   \ref{3lem:center1},
we conclude that 
\eqref{3eq:pigxC1}  holds for $\ell+1$.
\end{proof}

\subsection{Proof of mutation invariance of CSDs}

Let us complete the proof of Theorem \ref{3thm:mut1}.
We  have already seen in Lemma \ref{3lem:es1} that 
$T_k(\frakD_{\fraks})$ is an 
extended scattering diagram for $(\fraks',k)$.

We first prove the consistency of 
$T_k(\frakD_{\fraks})$ as an 
extended scattering diagram for $(\fraks',k)$.
Let us introduce  the following algebra automorphism
of  $R_k=R'_k$,
\begin{align}
\begin{matrix}
\tilde S_k=\tilde S_{k,\fraks}: &R_k
&\rightarrow 
&
R_k\\
& x^{\tilde m} &\mapsto &x^{\tilde S_{k}(\tilde m)}=
x^{\tilde m + \langle \d_k e_k,  m\rangle \tilde p^*(e_k)}
\quad (\tilde m \in \overline P),
\\
& (1+x^{\tilde p^*(e_k)})^{-1} &\mapsto &(1+x^{\tilde p^*(e_k)})^{-1}.
\end{matrix}
\end{align}
\begin{lem}
\label{3lem:ST1}
The following equalities hold as an algebra automorphism of $R_k$:
\begin{align}
\label{3eq:trans1}
\Psi[S^*_{k}(n)] \circ \tilde S_k
&=\tilde S_k \circ \Psi[n]
\quad
(n\in {N}^{+}),
\\
\label{3eq:trans2}
\Psi[e'_k]^{\d_k}
&=
\tilde S_k \circ \Psi[e_k]^{-\d_k}.
\end{align}
\end{lem}
\begin{proof}
Here, we consider the action only on $x^{\tilde m}$.
We have
\begin{align}
\begin{split}
(\Psi[S^*_{k}(n)] \circ \tilde S_k)(x^{\tilde m})
&=
\Psi[S^*_{k}(n)] (x^{\tilde S_{k}(\tilde m)})\\
&=x^{\tilde S_{k}(\tilde m)} \tilde S_{k}(f[n])^{\langle  S^*_{k}(n),{\tilde S_{k}(\tilde m)}\rangle}\\
&=\tilde S_k (x^{\tilde m} f[n]^{\langle  n,   m\rangle})\\
&=(\tilde S_k\circ \Psi[n])(x^{\tilde m}),
\end{split}
\end{align}
where  we used  Lemma \ref{3lem:S*1} (a) and (d).
Also, we have
\begin{align}
\begin{split}
\Psi[e'_k]^{\d_k} (x^{\tilde m})
&=x^{\tilde m} (1+x^{\tilde p^*(e'_k)}  )^{\langle \d_k e'_k,  m\rangle}\\
&=x^{\tilde m}  (1+x^{-\tilde p^*(e_k)}  )^{-\langle  \d_k e_k,   m\rangle}\\
&=x^{\tilde m + \langle \d_k e_k,  m\rangle \tilde p^*(e_k)}  (1+x^{\tilde p^*(e_k)}  )^{-\langle \d_k e_k,   m\rangle}\\
&=\tilde S_k(x^{\tilde m}  (1+x^{\tilde p^*(e_k)}  )^{-\langle  \d_k e_k,   m\rangle})\\
&=(\tilde S_k\circ \Psi[e_k]^{-\d_k})(x^{\tilde m}).
\end{split}
\end{align}
The action on  $(1+x^{\tilde p^*(e_k)})^{-1}$ is similar and easier.
\end{proof}

\begin{lem}
\label{3lem:const1}
$T_k(\frakD_{\fraks})$ is a consistent 
extended scattering diagram for $(\fraks',k)$.
\end{lem}

\begin{figure}
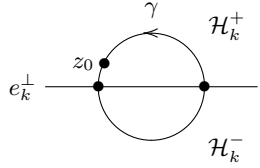

\centering
\leavevmode
\xy
(-17,0)*{\text{\small $ e_k^{\perp}$}};
(10,8)*{\text{\small $ \calH_{k}^+$}};
(10,-8)*{\text{\small $ \calH_{k}^-$}};
(-9,3)*{\text{\small $z_0$}};
(0,10)*{\text{\small $\gamma$}};
(7,0)*+{\bullet};
(-7,0)*+{\bullet};
(-6.2,3)*+{\bullet};
(0,7)*+{};(-2,7)*+{}
          ?>*\dir{>};
(0,0)="A",
(0,0)*\xycircle<20pt>{}
\ar@{-} "A"+(14,0); "A"+(-14,0)
\endxy
\caption{Admissible loop $\gamma$ in the proof of Lemma 
\ref{3lem:const1}.}
\label{3fig:scat2}
\end{figure}

\begin{proof}
It is enough to prove the consistency
as an extended scattering diagram for $(\fraks',k)$.
In view of Lemma \ref{3lem:consist2},
it is enough to  show that $\frakp_{\gamma, T_k(\frakD_{\fraks})}=\rmid $  for any admissible loop $\gamma$.
If $\gamma$ is entirely in $\calH_{k}^-$, the claim follows from the consistency of $\frakD_{\fraks}$
as an extended scattering diagram  for $(\fraks,k)$
in Lemma \ref{3lem:es1} (a).
If $\gamma$ is entirely in $\calH_{k}^+$, the claim follows from the consistency of $\frakD_{\fraks}$
and \eqref{3eq:trans1}.
Therefore, it is enough to consider the case where $\gamma$ crosses
the hyperplane $e_k^{\perp}$ as in Figure \ref{3fig:scat2},
where $z_0$ is the base point of $\gamma$.
Then, from the consistency of $\frakD_{\fraks}$, we have
\begin{align}
\label{3eq:pid2}
\frakp_{\gamma, \frakD_{\fraks}}=
g_+\circ \Psi[e_k]^{-\d_k}\circ g_{-}\circ \Psi[e_k]^{\d_k}
=\rmid,
\end{align}
where $g_{\pm}$ is the contribution from each half space $\calH_{k}^{\pm}$.
Then, by the definition of $T_k(\frakD_{\fraks})$ and  Lemma \ref{3lem:ST1},
\begin{align}
\begin{split}
\frakp_{\gamma, T_k(\frakD_{\fraks})}
&=
(\tilde S_k \circ g_+ \circ \tilde S_k^{-1})\circ \Psi[e'_k]^{\d_k}\circ g_{-}\circ \Psi[e'_k]^{-\d_k}\\
&=
\tilde S_k \circ (g_+\circ \Psi[e_k]^{-\d_k}\circ g_{-}\circ \Psi[e_k]^{\d_k})
\circ \tilde S_k^{-1}
\\
&=\rmid.
\end{split}
\end{align}
\end{proof}

Next we prove the following property.

\begin{lem}
\label{3lem:init1}
$T_k(\frakD_{\fraks})$ is equivalent to
an extended scattering diagram $\overline{\frakD}$ for $(\fraks',k)$
with 
\begin{align}
\label{3eq:incond2}
\overline\frakD_{\rmin}=\mathrm{In}_{\fraks'},
\end{align}
where $\fraks'=\mu_k(\fraks)$.
\end{lem}
\begin{proof}
Let us write $\bfw'_i:=\bfw_{e'_i}$ for simplicity.
For the sake of convenience, we restate here what should be proved.

\par\noindent
{\bf Claim.}
$T_k(\frakD_{\fraks})$ is equivalent to
an extended scattering diagram $\overline{\frakD}$ for $(\fraks',k)$
satisfying the following properties:
\begin{itemize}
\item[(i).]
$\overline{\frakD}$ contains a unique wall $\bfw'_k=(e'_k{}^{\perp}, \Psi[e'_k]^{\d_k})_{e'_k}$.
\item[(ii).]
For any $i\neq k$, 
$\overline{\frakD}$ contains a unique wall $\bfw'_i=(e'_i{}^{\perp}, \Psi[e'_i]^{\d_i})_{e'_i}$.
\item[(iii).]
All other walls in  $\overline{\frakD}$ are outgoing.
\end{itemize}

The property (i) holds for $T_k(\frakD_{\fraks})$ by the construction.
So, let us consider the properties (ii) and (iii).
We note that, by the construction of   $T_k(\frakD_{\fraks})$ in Definition \ref{3defn:mutD1},
a wall $S_k(\bfw_+)$ is incoming if and only
if $\bfw_+$ is incoming.
Thus, any  incoming wall in $T_k(\frakD_{\fraks})$
other than $\bfw'_k$
only comes  from $\bfw_{i,+}=((e_i {}^{\perp})_+, \Psi[e_i]^{\d_i})_{e_i}$ or $\bfw_{i,-}=((e_i {}^{\perp})_-, \Psi[e_i]^{\d_i})_{e_i}$
 $(i\neq k)$.
For each  $i\neq k$,  we consider three cases.
\par
{\bf Case 1.} $ \{\d_k e_k, e_i\} > 0$.
We have $p^*(e_i)\in \calH_{k}^+$.
Thus, $\bfw_{i,+}$ is an incoming wall.
Also, $S_k^*(e_i)=e_i+\{\d_ke_k, e_i\} e_k = e'_i$.
Therefore, we have an incoming wall in $T_k(\frakD_{\fraks})$,
\begin{align}
S_k(\bfw_{i,+})=((e'_i {}^{\perp})_-, \Psi[e'_i]^{\d_i})_{e'_i}
=\bfw'_{i,-},
\end{align}
where the sign $\pm$ of $(e'_i {}^{\perp})_{\pm}$ is with respect to $e'_k = -e_k$. 
Then, up to equivalence, it is extended to 
$\bfw'_i$ by adding outgoing walls.
\par
{\bf Case 2.} $ \{\d_k e_k, e_i\} < 0$.
We have $p^*(e_i)\in \calH_{k}^-$.
Thus, $\bfw_{i,-}$ is an incoming wall.
Also, $e_i= e'_i$.
Therefore, we have an incoming wall in $T_k(\frakD_{\fraks})$,
\begin{align}
\bfw_{i,-} =((e'_i {}^{\perp})_+, \Psi[e'_i]^{\d_i})_{e'_i}=\bfw'_{i,+}.
\end{align}
Again, up to equivalence, it is extended to 
$\bfw'_i$ by adding outgoing walls.
\par
{\bf Case 3.} $ \{\d_k e_k, e_i\} = 0$.
We have $p^*(e_i)\in \calH_{k}^+\cap \calH_{k}^-$.
Thus, $\bfw_{i,+}$
and $\bfw_{i,-}$ are
both  incoming walls, where $p^*(e_i)$ is at the boundary  of both walls
$\bfw_{i,\pm}$.
Also, $S_k^*(e_i) = e'_i$.
Therefore, we have incoming walls $\bfw_{i,+}=\bfw'_{i,-}$ and  $\bfw_{i,-}=\bfw'_{i,+}$
 in $T_k(\frakD_{\fraks})$.
Then, joining them back again, we have
$\bfw'_i$.
\par
It is clear that the resulting extended scattering diagram $\overline{\frakD}$
satisfies all desired properties.
\end{proof}

By 
Lemmas \ref{3lem:const1} and \ref{3lem:init1},
up to equivalence,
$T_k(\frakD_{\fraks})$ satisfies the assumption
of $\overline{\frakD}$ in  Proposition \ref{3prop:ext1} (by replacing $\fraks$ therein with $\fraks'$).
Thus, we conclude that 
$T_k(\frakD_{\fraks})$ is equivalent to $\frakD_{\fraks'}$
as an extended scattering diagram for $(\fraks',k)$.
Moreover,
since our choice of $\frakD_{\fraks}$ is with minimal support by Proposition \ref{3prop:CSD1},
so is $T_k(\frakD_{\fraks})$.
Thus, we have
\begin{align}
\mathrm{Supp}(T_k(\frakD_{\fraks}))=
\mathrm{Supp}(\frakD_{\fraks'}).
\end{align}
This implies that $T_k(\frakD_{\fraks})$ is actually
an  ordinary scattering diagram for $\fraks'$,
and that it is  equivalent to $\frakD_{\fraks'}$
as an ordinary scattering diagram for $\fraks'$.
This completes the proof of Theorem \ref{3thm:mut1}.

\subsection{Admissible region}
In view of the remark after Theorem \ref{3thm:mut1},
we have the following corollary of Theorem \ref{3thm:mut1}.

\begin{cor}
\label{3cor:mi1}
Let $\frakD_{\fraks}$ be a CSD for $\fraks$.
Let $n=\sum_{i=1}^r a_i e_i \in N_{\rmpr}^+$
such that $n\neq e_1,\dots, e_r$,
and let $z$ be any general point in $n^{\perp}$.
For each $k=1,\dots, r$, we define
\begin{gather}
\label{3eq:aa3}
b'_k=
\begin{cases}
-a_k + \sum_{j=1}^r  [- \{\d_k e_k, e_j\}]_+ a_j
& z\in \calH_k^+,
\\
-a_k + \sum_{j=1}^r  [\{\d_k e_k, e_j\}]_+ a_j
&z\in \calH_k^-.
\end{cases}
\end{gather}
Then, we have
\begin{align}
\label{3eq:prohibit1}
g_z(\frakD_{\fraks})\neq \rmid
\quad
\Longrightarrow
\quad
\text{$b'_k \geq 0$ for any $k$.}
\end{align}
\end{cor}
It is an interesting problem to give a direct proof of this property
and also to extend it to   other classes of   consistent scattering diagrams.

Let us call the region of $z\in M_{\bbR}$ that satisfies
the condition $b'_k \geq 0$ for any $k$ the \emph{admissible region\/} for a CSD $\frakD_{\fraks}$, \index{admissible region}
where nontrivial outgoing walls are allowed.

\begin{figure}
\centering
\leavevmode
\begin{xy}
(0,0)="A"
\ar@{-} (5,-10); (10,-5)
\ar@{-} (3,-10); (9.3,-3.7)
\ar@{-} (2,-9); (7.8,-3.2)
\ar@{-} (1.6,-7.4); (6.4,-2.6)
\ar@{-} (1.3,-5.7); (5,-2)
\ar@{-} (1,-4); (3.6,-1.4)
\ar@{-} (0.7,-2.3); (2, -1)
\ar "A"+(0,0); "A"+(10,0)
\ar "A"+(0,0); "A"+(0,10)
\ar@{-} "A"+(0,0); "A"+(-10,0)
\ar@{-} "A"+(0,0); "A"+(0,-10)
\ar@{-} "A"+(0,0); "A"+(2,-10)
\ar@{-} "A"+(0,0); "A"+(10,-4)
\end{xy}
\caption{Admissible region  for a CSD of rank 2.}
\label{3fig:scat4}
\end{figure}
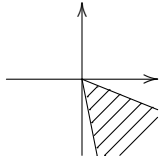

\begin{ex}
\label{3ex:prohibit1}
Let us explicitly write down  the condition \eqref{3eq:prohibit1} 
for the rank 2  CSDs in Sections \ref{3subsec:rank2} and \ref{3subsec:rank2inf}.
For $\{e_2,e_1\}=0$, no $n\neq e_1, e_2$ is allowed.
This agrees with the CSD of type $A_1\times A_1$ therein.
For $\{e_2,e_1\}=1$,
the quantities in  \eqref{3eq:aa3} are given as
\begin{gather}
\label{3eq:aa4}
b'_1=
\begin{cases}
-a_1 + \d_1 a_2 
& z\in \calH_1^+,
\\
-a_1 
&z\in \calH_1^-.
\end{cases}
\quad
b'_2=
\begin{cases}
-a_2 
& z\in \calH_2^+,
\\
-a_2+ \d_2 a_1
&z\in \calH_2^-,
\end{cases}
\end{gather}
where $a_1, a_2\neq 0$.
Then, the condition $b'_k \geq 0$ ($k=1, 2$) is equivalent to the condition
\begin{align}
\label{3eq:allowed1}
z\in \calH_1^+ \cap \calH_2^-,
\quad \frac{1}{\d_1} \leq \frac{a_2}{a_1} \leq \d_2.
\end{align}
The first condition is equivalent to that a wall is outgoing, 
and the second condition follows from the formula
\eqref{3eq:gxg2}.
Therefore, we  directly verified 
the property \eqref{3eq:prohibit1} in the rank 2 case.
Let $\bfz=(z_1,z_2) \in \bbR^2$ be the point corresponding  to $z\in M_{\bbR}$.
Then, by \eqref{3eq:can1},
we have
\begin{align}
\d_1^{-1} a_1z_1 + \d_2^{-1} a_2z_2 =0.
\end{align}
Thus, the condition \eqref{3eq:allowed1} is rewritten as
\begin{align}
z_1 \geq 0,\quad z_2\leq 0,
\quad - \d_2 \leq \frac{z_2}{z_1} \leq - \frac{1}{\d_1}.
\end{align}
The admissible region is depicted in Figure \ref{3fig:scat4}.
\end{ex}

\subsection{Chamber structure, $G$-fan, and cluster complex}
\label{3subsec:chamber1}

There is  an important consequence of the mutation invariance of CSDs.

Suppose that $\frakD_{\fraks}$ is a CSD with minimal support.
We first give few remarks:
\begin{itemize}
\item[
(a).]
By Lemma \ref{3lem:cham1}, 
the interior  $\mathrm{Int}(\calC^+_{\fraks})$ of  the positive orthant
$\calC^+_{\fraks}$ is a \emph{chamber\/} of $\mathrm{Supp}(\frakD_{\fraks})$,
\index{chamber}
namely, it is a connected component
of $M_{\bbR}\setminus \mathrm{Supp}(\frakD_{\fraks})$.
Also,
since $\mathrm{Int}(\calC^+_{\fraks})$ is an open set,
 its boundary belongs to $\mathrm{Supp}(\frakD_{\fraks})$.
\item[
(b).]
For $\fraks'=\mu_k(\fraks)$,
we define a cone in $M_{\bbR}$ of full dimension by
\begin{align}
\calC^{\fraks}_{\fraks'}
:=
T_{k,\fraks}^{-1}(\calC^+_{\fraks'}).
\end{align}
Since $\mathrm{Int}(\calC^+_{\fraks'})$
is a chamber of $\mathrm{Supp}(\frakD_{\fraks'})$,
$\mathrm{Int}(\calC^{\fraks}_{\fraks'})$ is a
chamber of $\mathrm{Supp}(\frakD_{\fraks})$
by \eqref{3eq:supp2}.
More explicitly, $\calC^{\fraks}_{\fraks'}$ is the cone
defined by,
for $i=k$,
\begin{align}
\label{3eq:ek1}
\langle e'_k, T_k(z)\rangle=
\langle -e_k,
 z+\langle \d_ke_k,z\rangle p^*(e_k)\rangle
 =- \langle e_k, z\rangle\geq 0,
\end{align}
and,
for $i\neq k$, by \eqref{3eq:ek1},
\begin{align}
\begin{split}
\langle e'_i, T_k(z)\rangle&=
\langle e_i+ [\{ \d_k e_k,e_i\}]_+e_k,
 z+[\langle \d_ke_k,z\rangle]_+ p^*(e_k)\rangle
 \\
 &=\langle e_i, z\rangle
 +[\{\d_k e_k, e_i\} ]_+ \langle e_k,z\rangle
 \geq 0.
 \end{split}
\end{align}
Thus, the cone $\calC^{\fraks}_{\fraks'}$ is adjacent to $\calC_{\fraks}^{+}$
with a common boundary in $e_k^{\perp}$.

\item[
(c).]
Let $\sigma$ be a face of $\calC^{\fraks}_{\fraks'}$
of codimension 1
with the  normal vector $n$.
Let $z\in \sigma$ be general.
Then, $g_z(\frakD_{\fraks})$ is obtained by
the pull-back of  $g_{T_k(z)}(\frakD_{\fraks'})=\Psi[e'_i]^{\d_i}$
by $T_k$.
Thus, we have
\begin{align}
\label{3eq:ggd1}
g_z(\frakD_{\fraks})=\Psi[n]^{\d(n)}.
\end{align}
Moreover, $\d(n)=\d_i$ for some $i\in \{1,\dots, r\}$.
\end{itemize}

We repeat the above argument for a sequence of mutations of seeds.

\begin{defn}[$G$-cone/$G$-chamber] \index{cone!$G$-}\index{$G$-cone} \index{chamber!$G$-}\index{$G$-chamber}
Suppose that $\frakD_{\fraks}$ is a CSD with minimal support.
We consider any sequence of mutations of seeds
\begin{align}
\fraks=\fraks_0 \buildrel \mu_{k_0} \over \mapsto \fraks_1
\buildrel \mu_{k_2} \over \mapsto
\cdots 
\buildrel \mu_{k_{i-1}} \over \mapsto
\fraks_{i}.
\end{align}
We define
a cone in $M_{\bbR}$ by
\begin{align}
\label{3eq:clusterc1}
\calC^{\fraks}_{\fraks_{i}} :=
(T_{k_0,\fraks_0}^{-1}\circ \cdots \circ T_{k_{i-1},\fraks_{i-1}}^{-1})(\calC^+_{\fraks_{i}}).
\end{align}
As a special case $i=0$, we set $\calC^{\fraks}_{\fraks}:=\calC^+_{\fraks}$.
Then,
$\mathrm{Int}(\calC^{\fraks}_{\fraks_{i}})$ is 
a chamber of $\mathrm{Supp}(\frakD_{\fraks})$.
We call $\calC^{\fraks}_{\fraks_{i}}$ and $\mathrm{Int}(\calC^{\fraks}_{\fraks_{i}})$ a \emph{$G$-cone\/} and a \emph{$G$-chamber\/} of
$\frakD_{\fraks}$.
(In \cite{Gross14}, $\calC^{\fraks}_{\fraks_{i}}$ is called a \emph{cluster chamber}.)
\end{defn}

Let $\bbT_r$ be the \emph{$r$-regular tree\/} wth edges labelled by $1$, \dots, $r$ 
such that all edges emanating  from each vertex $t$ of $\bbT_r$ have different labels.
We abuse the notation $\bbT_r$ for the set of all vertices of $\bbT_r$.
We fix the initial vertex $t_0\in\bbT_r$ arbitrarily.
Then, we attach $\calC^+_{\fraks}$ to $t_0$.
Also, we attach $\calC^{\fraks}_{\fraks_{i}}$ in \eqref{3eq:clusterc1}
to $t_i$, where $t_0$ and  $t_i$ are sequentially connected with edges
labelled by $k_0$, \dots, $k_{i-1}$.
Thus,  the $G$-cones of $\frakD_{\fraks}$ are parametrized by  $\bbT_r$,
though this parametrization is not necessarily one-to-one.
Then,
under the identification $M_{\bbR}\simeq \bbR^r$, $f_i \mapsto \bfe_i$,
 any $G$-cone  $\calC^{\fraks}_{\fraks_{i}}$ is identified with the cone $\sigma(G^{t_0}_{t_i})$ spanned  by the column vectors (the \emph{$g$-vectors}) of the \emph{$G$-matrix\/} $G^{t_0}_{t_i}$ 
in \cite{Fomin07}
for the corresponding cluster pattern.
The collection of $G$-cones $\sigma(G^{t_0}_t)$ ($t\in \bbT_r$) and their faces constitutes a fan
$\Delta(\bfG^{t_0})$ (the \emph{$G$-fan\/}).
Moreover, as a simplicial complex, the $G$-fan $\Delta(\bfG^{t_0})$ is isomorphic to the \emph{cluster complex\/} \index{fan!$G$-}\index{$G$-fan}
in \cite{Fomin03a}.
See Part II for details.

By repeatedly applying the result \eqref{3eq:ggd1}, we have the following property.
\begin{prop}
\label{3prop:cchamber1}

Let $\sigma$ be a face of a $G$-cone $\calC^{\fraks}_{\fraks_i}$
of codimension 1
with the normal vector $n\in N_{\rmpr}^+$.
Then, for any general $z\in \sigma$, we have
\begin{align}
\label{3eq:gxC1}
g_z(\frakD_{\fraks})=\Psi[n]^{\d(n)}.
\end{align}
Moreover, $\d(n)=\d_i$ for some $i\in \{1,\dots, r\}$.
\end{prop}
The above result implies that,
in Theorem \ref{3thm:pos1}, 
for the walls that are  faces of $G$-cones,
only $s=t=1$ occurs.

\begin{ex}
\label{3ex:chamber1}
(a).
Let $\frakD_{\fraks}$ be any CSD of  finite type
with rank 2 or 3
in  Sections \ref{3subsec:rank2} and \ref{3subsec:rank2inf}.
Then, the $G$-cones fill up the entire ambient space $\bbR^r$.
In other words,
$\mathrm{Supp}(\frakD_{\fraks})$
is naturally identified with the $G$-fan $\Delta(\bfG^{t_0})$
of the corresponding cluster pattern.
In particular,
 the property \eqref{3eq:gxC1} holds
 for any wall of $\frakD_{\fraks}$.
 \par
 (b). 
  Let $\frakD_{\fraks}$ be the CSD of   type $A_1^{(1)}$
in  Section \ref{3subsec:rank2inf}.
Then, the $G$-cones fill up  the set $\bbR^2
\setminus \bbR_{> 0}(1,-1)$.
In other words,
$\mathrm{Supp}(\frakD_{\fraks})\setminus \bbR_{> 0}(1,-1)$
in $\bbR^2
\setminus \bbR_{\geq 0}(1,-1)$
is naturally identified with the $G$-fan $\Delta(\bfG^{t_0})$
of the corresponding cluster pattern.
In particular,
 the property \eqref{3eq:gxC1} holds
 for any wall in $\mathrm{Supp}(\frakD_{\fraks})\setminus \bbR_{> 0}(1,-1)$.
 On the other hand,
the wall in \eqref{3eq:a114} is \emph{unreachable\/} \index{wall!unreachable}
by \emph{cluster mutations\/} (i.e., mutations in the sense of cluster algebras).
The  $A_2^{(2)}$ case is similar.
For the non-affine type of rank 2 in  Section \ref{3subsec:rank2inf},  replace  $\bbR_{> 0}(1,-1)$
with the Badlands $ \sigma(\bfv_+, \bfv_-)\setminus \{0\}$ therein,
which is also unreachable by cluster
  mutations.
\end{ex}

\subsection{Rank 3 examples: infinite type}
\label{3subsec:rank32}

Let us continue to give examples of rank 3 CSDs 
based on Construction \ref{3const:CSD1}.
Here we consider the CSDs of \emph{infinite type},
which are further separated into \emph{affine\/}  and   \emph{non-affine\/} types. 
We use the same convention in Section \ref{3subsec:rank31}
 to present rank 3 CSDs.

\begin{figure}
\centering
\includegraphics{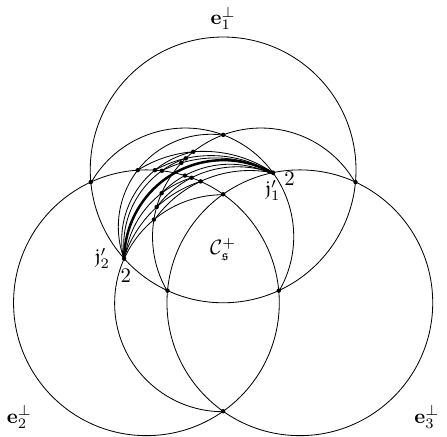}
\caption{CSD of type $A_2^{(1)}$.
 The thick arc is  unreachable by cluster
  mutations.
}
\label{3fig:rank33}
\end{figure}

(a). Type $A_2^{(1)}$. As a prototypical example of affine type, we consider the following case
\begin{align}
B=
\begin{pmatrix}
0 & -1 & -1\\
1 & 0 & -1\\
1 & 1 &0
\end{pmatrix},
\quad
\d_1=\d_2=\d_3=1.
\end{align}
As a new feature compared with the finite type,
we have  perpendicular joints  of type $A_1^{(1)}$ in Section
\ref{3subsec:rank2inf},
which are given by
$\frakj'_1\subset (1,1,0)^{\perp}\cap \bfe_3^{\perp}$ and 
$\frakj'_2 \subset (0,1,1)^{\perp}\cap \bfe_1^{\perp}$ in Figure \ref{3fig:rank33}.
Indeed, we have
\begin{align}
\{ \bfe_3,(1,1,0)\}=
\{(0,1,1),\bfe_1\}=2.
\end{align}
We put this number 2
 in the figure in the position of incoming direction,
so that outgoing walls appear in the opposite direction.
These joints $\frakj'_1$ and $\frakj'_2$ are connected by walls with the common normal vector $(1,1,1)$,
and their wall elements are given by
$\Psi[(1,1,1)]^2$, $\Psi[(2,2,2)]$, $\Psi[(4,4,4)]^{1/2}$, \dots,
based on the result 
the CSD of type  $A_1^{(1)}$
 in Section
\ref{3subsec:rank2inf}.
The union of their supports is the cone $\sigma(\frakj'_1, \frakj'_2)$,
and it
 is  represented by a thick arc in Figure \ref{3fig:rank33}.
It corresponds to the unreachable wall of type $A_1^{(1)}$ with the normal vector $(1,1)$
in Example \ref{3const:CSD1} (b).
Note that $\frakj'_1$ and $\frakj'_2$
are  directly connected only by $\sigma(\frakj'_1, \frakj'_2)$, because 
a hyperplane
containing $\frakj'_1$ and $\frakj'_2$ in $\bbR^3$
is unique.
Consequently, $\frakj'_1$ and $\frakj'_2$ ``alternately'' create new walls,
and they pile up as shown in the figure.
We see that those walls with normal vectors
\begin{align}
(2,1,1),\ (2,2,1),\ (3,2,2),\ (3,3,2),\ \dots
\end{align}
converge to the thick arc \emph{from  above\/} in the figure.
Similarly, 
those walls with normal vectors
\begin{align}
(1,1,2),\ (1,2,2),\ (2,2,3),\ (3,2,2),\ \dots
\end{align}
converges to the thick arc \emph{from below}.
This agrees with the well-known results on 
 the cluster patterns of type $A_2^{(1)}$ by \cite{Fomin02}
 and \cite{Cerulli09}.
 We conclude that the $G$-cones fill up  the region $\bbR^3
\setminus (\sigma(\frakj'_1, \frakj'_2)\setminus \{0\})$.
In other words,
$\mathrm{Supp}(\frakD_{\fraks})\setminus  (\sigma(\frakj'_1, \frakj'_2)\setminus \{0\})$
in $\bbR^3
\setminus  (\sigma(\frakj'_1, \frakj'_2)\setminus \{0\})$
is naturally identified with the $G$-fan $\Delta(\bfG^{t_0})$
of the corresponding cluster pattern.
The walls in $\sigma(\frakj'_1, \frakj'_2)$ is \emph{unreachable\/} by cluster
   mutations.

\begin{figure}
\centering
\includegraphics{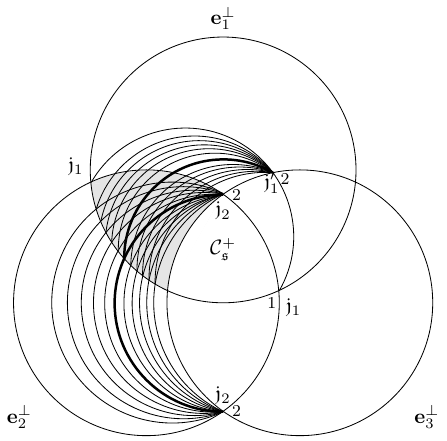}
\caption{Example of rank 3 CSD of non-affine type.
The close-up of the shaded region is given in  Figure \ref{3fig:na2}.}
\label{3fig:na1}
\end{figure}

\begin{figure}
\centering
\includegraphics[scale=0.70]{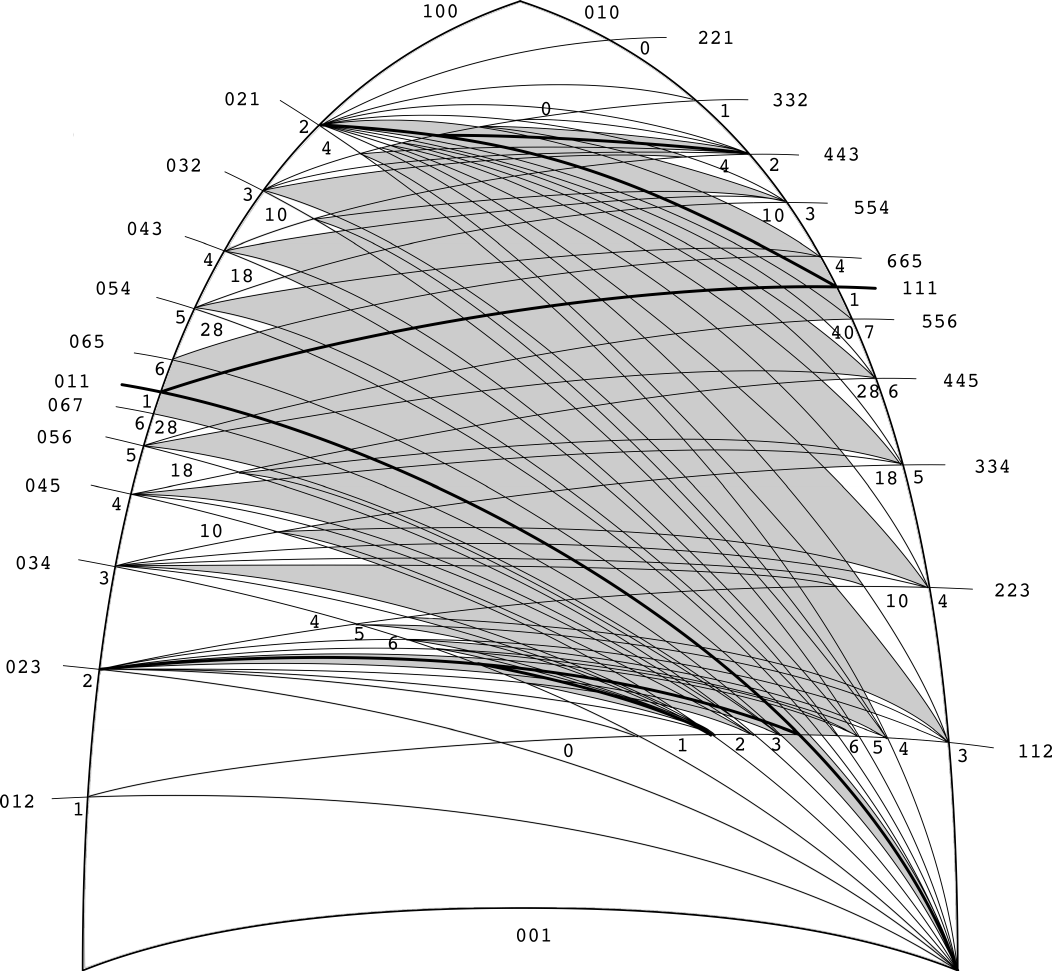}
\caption{Uncompleted picture of the Badlands of rank 3.
This is the close-up of the shaded region in Figure \ref{3fig:na1} }
\label{3fig:na2}
\end{figure}

(b). Non-affine type. As a prototypical example of non-affine type, we consider the following case
\begin{align}
B=
\begin{pmatrix}
0 & -1 & 0\\
1 & 0 & -2\\
0 & 2 &0
\end{pmatrix},
\quad
\d_1=\d_2=\d_3=1.
\end{align}
The CSD has the remarkable \emph{fractal\/} nature.
This
nature is simply originated in the fractal structure of the underlying
$r$-regular tree $\bbT_r$ ($r\geq 3$) for the corresponding cluster complex,
and it has been  observed in some figures in the literature (e.g., \cite{Fock11,Najera11})
for the so-called Markov quiver.
This implies that it is impossible to draw a precise picture in principle.
Thus, we change the strategy  as follows:
\begin{itemize}
\item
We set the goal to obtain a reasonably clear  picture of the \emph{Badlands of rank 3}.
\item
We abandon to faithfully follow Construction \ref{3const:CSD1}
by uniformly increasing the degree $\d(n)$ of a newly attached wall element $\Psi[n]$.
Rather, we do the mixture of Construction \ref{3const:CSD1} and  ``cluster  mutation approach''.
Namely, 
having Construction \ref{3const:CSD1} in mind,
we increase the region of $G$-cones in the scattering diagram.
Since 
the Badlands happens to be the complement of the $G$-fan in this example, 
the shape of the Badlands will gradually emerge by this procedure.
\item
While  proceeding, we carefully choose added walls
so that  both local and global structures of the Badlands are better visible.
(For example, we ignored some wall with  degree 12,
while  we added a new wall of degree 245.)
\item
Due to the complexity we double-checked the result
with a computer program for computing \emph{$c$-vectors\/}
of the  cluster pattern.
More specifically,  we  used SageMath 9.4
 with \emph{Cluster Algebras package\/} by D. Rupel and S. Stella.
\end{itemize}
The resulting CSD is depicted in Figures \ref{3fig:na1} and \ref{3fig:na2}.

Figures \ref{3fig:na1} gives the big picture of the CSD.
Let $\calR=\calH_1^{>0}\cap \calH_2^{>0}\cap\calH_3^{<0}$,
which corresponds to the shaded region in Figure \ref{3fig:na1}.
Outside the region $\calR$, the structure of the CSD is clearly  presented in the figure,
where there are actually infinitely many walls converging to the
unreachable walls attached to joints of type $A_1^{(1)}$
represented by two thick arcs.
There is no more additional wall outside $\calR$.

The close-up of the region $\calR$ is given in Figure \ref{3fig:na2}.
This is where the Badlands and its fractal nature appear.
Here are remarks and observations:
\begin{itemize}
\item
A three-digit number such as 032 represents the normal vector (0,3,2)
of a wall crossing the boundary of the region $\calR$.
\item 
A number such as  4 at a joint in $n_1^{\perp}\cap n_2^{\perp}$
  represents $|\{n_1,n_2\}|=4$.
It is put in the position of incoming direction,
so that outgoing walls appear in the opposite direction
as in the previous case.

\item
A thick arc represents  the supports of walls attached to a joint of type $A_1^{(1)}$
that is unreachable by cluster mutations.
They consist of the ``backbone'' of the total shaded region (explained below).
Six thick arcs are presented here, but in the limit $\ell\rightarrow \infty$ there will be
infinitely many ones added in the top and bottom parts of the back-born.
However, they will be so small and tightly aligned that they are not visible in this page size anyway.
\item
An unshaded triangle is already completed under the algorithm,
and it is a $G$-cone.
Meanwhile, a shaded triangle or quadrilateral is not completed yet.
Thus, when the degree $\ell$ increases, 
the corners are sharpened and new corners appear.
In other words,
the total  shaded region   reduces,
and  the complexity of the boundary increases.

\item
By the above comments and the picture presented in  Figure \ref{3fig:na2},
it is already clear that in the limit $\ell \rightarrow \infty$
the total shaded region (in $\bbR^3$) converges to some connected closed domain, say $\calB$.
It is also characterized as the unreachable region by cluster mutations.
The boundary of $\calB$ is a fractal.
Moreover,
if the density conjecture on the Badlands of rank 2 in Section
\ref{3subsec:rank2inf} holds,
the set $\calB$ is also densely filled with nontrivial walls.
Thus, we call  $\calB$ the \emph{Badlands\/} of this CSD. \index{Badlands}

\end{itemize}

Some other examples of rank 3 CSD of infinite type are found in \cite{Muller15, Reading19}.
It seems that there are several patterns for non-affine infinite type CSDs even for the rank 3 case,
and to classify them is an interesting problem.

\notes
Most results are taken from \cite[\S1.3, Appendix C.2]{Gross14}
with added/modified proofs.
In particular,
the proof of Lemma \ref{3lem:es1}
 is significantly simplified from
the one of \cite[Theorem~1.28]{Gross14}
by 
based on Construction \ref{3const:CSD1}
and  Proposition \ref{3prop:bound1} (b).

\newpage

\section{Broken lines and theta functions}
\label{3sec:broken1}

In this section we introduce broken lines and theta functions for scattering diagrams.
We prove the positivity of coefficients of theta functions for cluster scattering diagrams.
\subsection{Broken lines and theta functions}

In this section we again use the principal $x$-representation for
wall elements.
For notational simplicity, for any $\tilde m_i\in \tilde M^{\circ}=M^{\circ} \oplus N$,
let $m_i$ denote the $M^{\circ}$-component of $\tilde m_i$.
Namely, $\tilde m_i = m_i + n_i$ ($m_i\in M^{\circ}$, $n_i\in N$) for some $n_i$.
We recall that $\langle n, \tilde m_i \rangle = \langle n, m_i \rangle$ holds 
for any $n\in N$ by \eqref{3eq:bibi1}.

\begin{defn}[Broken line] \index{broken line}
\label{3defn:broken1}
Let $\frakD$ be any scattering diagram for a given initial seed $\fraks$.
Let  $\tilde m_0 \in \tilde M^{\circ}\setminus \{0\}$ and  $Q\in M_{\bbR}\setminus
\mathrm{Supp}(\frakD)$.
A \emph{broken line for $\tilde m_0$ with endpoint $Q$} is a piecewise-linear
curve $\gamma:(-\infty, 0]\rightarrow  M_{\bbR}\setminus
\mathrm{Sing}(\frakD)$ satisfying the following properties:
\begin{itemize}
\item[(1).]
The endpoint $\gamma(0)$ is $Q$.
\item[(2).]
The initial velocity  $\gamma'(t)$ in $t\rightarrow -\infty$
is $-m_0$.
\item[(3).]
The curve $\gamma$ bends
and change velocity only when it intersects walls of $\frakD$.
(It can cross walls also \emph{without bending and changing velocity}.)
Moreover,
when $\gamma$ bends, it always \emph{crosses\/}  walls,
not \emph{reflects\/} at them.
\end{itemize}
More conditions follow below.
Suppose that $\gamma$ bends $J\geq 0 $ times.
Let $t_1 < t_2 < \cdots < t_J<0$ be the points (``time'') of the parameter $t$
for $\gamma(t)$ when $\gamma$ bends.
(It may cross multiple walls with a common normal vector at each time $t_j$.)
We set $t_0=-\infty$ and $t_{J+1}=0$.
Thus, $\gamma$
 is linear and with constant velocity in each interval $I_j=(t_{j},t_{j+1})$ ($j=0,\dots,J$).
 As an extreme case, for $J=0$, 
 the curve $\gamma$ is a half line ending at $Q$
 with constant velocity $-m_0$.
 To each  interval $I_j$,
 we attach a monomial $c_j x^{\tilde m_j}\in \bbk[\tilde  M^{\circ}]$ recursively as follows:
\begin{itemize}
 \item[(4).]
To  the initial interval $I_0$, the monomial $x^{\tilde m_0}$ is attached,
where $\tilde m_0$ is the given data. 
\item[(5).] For $j\geq 1$, let $\gamma_j$ be a segment of $\gamma$ for
the interval
$(t_j-\delta,t_j+\delta)$ with sufficiently small $\delta>0$ such that $\gamma_j$ crosses walls
only at $t_j$.
Then, $c_j x^{\tilde m_j}$ is a monomial in
$
\frakp_{\gamma_j, \frakD}(c_{j-1}x^{\tilde m_{j-1}})$,
where $\frakp_{\gamma_j, \frakD}$ acts as the principal $x$-representation.
See \eqref{3eq:pgD1} for more explicit description.
\end{itemize}
Finally, the velocity of $\gamma$  in the interval $I_j$ is given  as follows:
\begin{itemize}
\item[(6).] The velocity $\gamma'(t)$ in the interval $I_j$ is given by $-m_j$,
where $c_jx^{\tilde m_j}$ is the attached monomial.
This agrees with the  conditions (2) and (4).
\end{itemize}
This completes the definition of a broken line.
We then attach a monomial in $\bbk[\tilde M^{\circ}]$ to each broken line $\gamma$ by
\begin{align}
\label{3eq;mono1}
c_{\gamma} x^{\tilde m_\gamma}:=c_J x^{\tilde m_J}.
\end{align}
\end{defn}

Several remarks are in order.

(a).
At  $t_j$, the velocity changes
from $-m_j$ to  $-m_{j+1}=-m_j - k_j p^*(n_j)$, where $n_j$ is the common normal vector
of the walls at $t_j$, and $k_j$ is some positive integer.
Since $p^*(n_j)\in n_j^{\perp}$,
the broken line always \emph{crosses\/}  walls,
\emph{not reflects nor stops} at them as stated in the condition (3).

(b). In particular,  we have $\tilde m_j= \tilde m_0+\sum_{i=1}^r a_i \tilde p^*(e_i)$
$(a_i \in \bbZ_{\geq 0})$.
\par
(c).
For $\tilde m' = \tilde m_0 + \sum_{i=1}^r a_i \tilde p^*(e_i)$,
we define the degree $\deg_{\tilde m_0}(\tilde m')$ by 
\begin{align}
\deg_{\tilde m_0}(\tilde m'):=\sum_{i=1}^r a_i.
\end{align}
Then, by (a), we have
\begin{align}
\label{3eq:dec1}
\deg_{\tilde m_0}(\tilde m_0)=0 < \deg_{\tilde m_0}(\tilde m_1) < \cdots < \deg_{\tilde m_0}(\tilde m_J)=\deg_{\tilde m_0}(\tilde m_{\gamma}).
\end{align}

\par
(d). 
For the segment $\gamma_j$ in the condition (5),
suppose that $\gamma_j$ crosses walls
$[\frakd_{\lambda}, f_{{\lambda}}]_{n_j}$ $(\lambda\in \Lambda_j$)
at $t_j$.
Note that $\gamma'=-m_{j-1}$ just before $\gamma$ crosses these walls.
Then,
the intersection signs 
in \eqref{3eq:factor1}
are given by
\begin{align}
\label{3eq:factor2}
\epsilon_j
=
\begin{cases}
1 & \langle n_{j}, m_{j-1}\rangle>0,\\
-1 & \langle n_{j}, m_{j-1}\rangle<0.\\
\end{cases}
\end{align}
Thus, $\epsilon_j \langle n_{j}, m_{j-1}\rangle =|\langle n_{j}, m_{j-1}\rangle|$.
Therefore,  by Definition \ref{3defn:pop1},
we have
\begin{align}
\label{3eq:pgD1}
\frakp_{\gamma_j, \frakD}(c_{j-1}x^{\tilde m_{j-1}})
=c_{j-1}x^{\tilde m_{j-1}
}
\prod_{\lambda\in \Lambda_j}
f_{{\lambda}}^{| \langle \d(n_j)n_j ,  m_{j-1}\rangle|}.
\end{align}
It is crucial that there is \emph{no division\/}  by $f_{{\lambda}}$ in this expression
for the forthcoming positivity of theta functions.

\begin{defn}[General position] \index{general position}
We say that  $Q\in M_{\bbR}\setminus
\mathrm{Supp}(\frakD)$ is \emph{in general position\/} if
for each $\l>0$ there is some neighborhood $U_{\l}$ of $Q$
such that, for the reduction $\frakD_{\l}$, any broken line for $Q' \in U_{\l}$  
converges to a broken line for $Q$
in the limit $Q'\rightarrow Q$.
(Namely, no broken line crosses $\mathrm{Sing}(\frakD_{\l})$ in the limit.)
\end{defn}

Below we always and implicitly assume that $Q$ is in general position.

\begin{defn}[Theta function] \index{theta function}
Under the same assumption and notations in Definition
\ref{3defn:broken1},
the \emph{theta function $\vartheta_{Q,\tilde m_0}$} for $\tilde  m_0$ with endpoint $Q$  is defined by
\begin{align}
\label{3eq:thetadef1}
\vartheta_{Q,\tilde m_0}
:=\sum_{\gamma\in B(Q,\tilde m_0)}c_{\gamma} x^{\tilde m_\gamma},
\end{align}
where $B(Q,\tilde m_0)$ is the set of  all broken lines for $\tilde m_0$ with endpoint $Q$.
We also set
\begin{align}
\vartheta_{Q,0}=1.
\end{align}
\end{defn}

The function $\vartheta_{Q,\tilde m_0}$ only depends on the equivalence class
of $\frakD$, though the condition for $Q$ to be in general position  
depends on the choice of a representative of $\frakD$.

\begin{lem}[{\cite[Proposition~3.4]{Gross14}}]
For any scattering diagram $\frakD$, we have
\begin{align}
\vartheta_{Q,\tilde m_0}\in x^{\tilde m_0} \bbk[[\tilde P]].
\end{align}
\end{lem}

\begin{proof}
By Remark (b), 
each term $c_{\gamma} x^{\tilde m_\gamma}$ belongs to $x^{\tilde m_0} \bbk[[\tilde P]]$.
Thus, it is enough to show that,
for a given $\tilde m' \in \tilde m_0+\tilde P$,
there are only finitely many $\gamma$ such that
$\tilde m_{\gamma}=\tilde m'$.
We note the following facts:
\begin{itemize}

\item
By the constraint \eqref{3eq:dec1},
such $\gamma$ bends at most $\deg_{\tilde m_0}(\tilde m')$ times.
\item
The walls of $\frakD$ that contribute to $\tilde m'$
belongs to the reduction $\frakD_{\l}$ at $\l=\deg_{\tilde m_0}(\tilde m')$,
which has only finitely many walls.
\item
At each bending there are only finitely many 
possibilities of bending.
\end{itemize}
Therefore, there are only finitely many possibilities of such $\gamma$.
\end{proof}

\begin{prop}[{\cite[Proposition~3.8]{Gross14}}]
\label{3prop:theta2}
Let $\frakD$ be any scattering diagram for a given seed $\fraks$,
and let
$\calC^+_{\fraks}$ be the positive orthant
in \eqref{3eq:chamber1}.
For $\tilde m_0\in \tilde M^{\circ}$,
suppose that $m_0\in \calC^+_{\fraks}\cap M^{\circ}$.
Also, suppose that $Q\in \mathrm{Int}(\calC^+_{\fraks})$.
Then, we have
\begin{align}
\label{3eq:thetain1}
\vartheta_{Q,\tilde m_0}=x^{\tilde m_0}.
\end{align}
\end{prop}

\begin{proof}
We may assume that  $\tilde m_0 \neq 0$.
It is clear that
$Q+\bbR_{\geq 0} m_0$
is a broken line for $\tilde m_0$ with end point $Q$,
and it does not intersect any walls of $\frakD_{\fraks}$;
therefore, the associated monomial is $x^{\tilde m_0}$.
We claim that under the condition there is no  broken line $\gamma$
 with  bending.
Suppose that such $\gamma$  exists.
We use the notation in Definition \ref{3defn:broken1}.
For each break point $t_j$ therein,
let $n_j\in N^+$ be a unique (not necessarily primitive) 
element such that
$\gamma$ intersects  walls in $n_j^{\perp}$ at $t_j$
and that the velocity of $\gamma$ is shifted by $-p^*(n_j)$.
Let $L_j$ be the segment of $\gamma$ for the interval $I_j$.
Then, we claim the following:
\begin{align}
\label{3eq:ine1}
L_j \subset H_j^-:=\biggl\{ z\in M_{\bbR} \biggm | \biggl\langle
\sum_{i=1}^j n_i, z \biggr\rangle \leq 0
\biggr\}
\quad
(j=1,\dots, J).
\end{align}
We show it by the induction on $j$.\
This is true for $j=1$, because the support of the wall crossed by $\gamma$ at $t_1$
is in $n_1^{\perp}$, and $\gamma$ is in $\calC^+_{\fraks}$ for $t \rightarrow -\infty$.
Now suppose that \eqref{3eq:ine1} hold up to $j$.
Then, we have $  \langle
\sum_{i=1}^j n_i, \gamma(t_{j+1}) \rangle \leq0$.
Also, $  \langle n_{j+1}, \gamma(t_{j+1})  \rangle=0$.
Thus, $ \gamma(t_{j+1}) \in H_{j+1}^-$.
Moreover, the new velocity vector $\gamma'_{j+1}=-m_0 -\sum_{i=1}^{j+1} p^*(n_i)$
satisfies the inequality
\begin{align}
\begin{split}
\biggl\langle
\sum_{i=1}^{j+1} n_i, \gamma'_{j+1}  \biggr\rangle
&=
\biggl\langle
\sum_{i=1}^{j+1} n_i, -m_0  \biggr\rangle
+
\biggl\{
\sum_{i=1}^{j+1} n_i, \sum_{i=1}^{j+1} n_i \biggr\}\\
&=
\biggl\langle
\sum_{i=1}^{j+1} n_i, -m_0  \biggr\rangle
\leq 0.
\end{split}
\end{align}
Thus, the fact  \eqref{3eq:ine1} holds for $j+1$.
Since $\mathrm{Int}(\calC^+_{\fraks})\cap H_j^-=\emptyset$, $\gamma$ cannot reach $Q\in \mathrm{Int}(\calC^+_{\fraks})$.
\end{proof}

\subsection{Transitivity on theta functions}
Now we assume that $\frakD$ is consistent.

\begin{prop}[{\cite[Lemma~4.9]{Carl10}}]
\label{prop:theta1}
Let $\frakD$ be a consistent scattering diagram.
Let $\tilde m_0\in\tilde  M^{\circ}$ and $Q,\, Q'\in M_{\bbR}\setminus
\mathrm{Supp}(\frakD)$.
Then, for any admissible curve $\gamma$ from $Q$ to $Q'$,
we have
\begin{align}
\label{3eq:theta3}
\vartheta_{Q',\tilde m_0}=\frakp_{\gamma,\frakD}(\vartheta_{Q,\tilde m_0}).
\end{align}
\end{prop}
The rest of this subsection is devoted to
prove this proposition.
We may assume that $\tilde m_0 \neq 0$.

Let  $J$ be the maximal ideal of 
$\bbk[[\tilde P]]$ generated by $\tilde P\setminus\{0\}$.

\begin{lem}[{\cite[Lemma~4.7]{Carl10}}]
\label{3lem:theta1}
Let  $\frakD_{\l}$ be the reduction of $\frakD$ at $\l$.
Suppose that $Q$ and $Q'$ are in the same chamber of
$\mathrm{Supp}(\frakD_{\l})$.
Then, we have the equality
\begin{align}
\label{3eq:theta2}
\vartheta_{Q,\tilde m_0}
\equiv
\vartheta_{Q',\tilde m_0}
\mod
x^{\tilde m_0} J^{\l+1}.
\end{align}
\end{lem}
\begin{proof}
Since we  consider $\vartheta_{Q,\tilde m_0}$
only modulo
$x^{\tilde m_0} J^{\l+1}$,
one can replace $\frakD$ with 
$\frakD_{\l}$ to calculate $\vartheta_{Q,\tilde m_0}$.
First,
suppose that $Q$ and $Q'$ are connected by a curve $\beta$
in the same chamber
such that any point on $\beta$ is in general position for $\frakD_{\l}$.
Then, each broken line with endpoint $Q$
continuously deformed to a broken line with endpoint $Q'$.
Therefore, the theta function is unchanged.

Suppose  that there is no such curve  between $Q$ and $Q'$.
Let us consider the special case where $Q$ and $Q'$ are connected by a curve $\beta$
in the same chamber
such that there is a unique point $Q_0$ on $\beta$
that is not in general position.
   (The general case  reduces to this case.)
 Then, there is at least one broken line $\gamma$ with endpoint $Q$
 such that, in the limit $Q\rightarrow Q_0$,
 the broken line  intersects with some boundary
 $\partial\frakd$ or joint $\frakj=\frakd_1\cap \frakd_2$ of $\frakD_{\l}$.

\begin{figure}
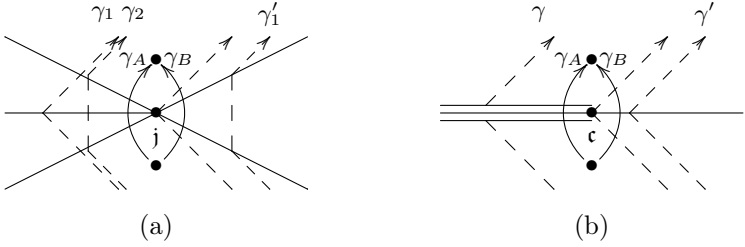

\centering
\leavevmode
\xy
(0,-15)*{\text{(a)}};
(-7,13)*{\text{\small $\gamma_1$}};
(-3,13)*{\text{\small $\gamma_2$}};
(15,13)*{\text{\small $\gamma'_1$}};
(-3,7)*{\text{\small $\gamma_A$}};
(3,7)*{\text{\small $\gamma_B$}};
(0,-3)*{\text{\small $\frakj$}};
(0,0)*+{\bullet};
(0,-7)*+{\bullet};
(0,7)*+{\bullet};
(0,-7)*+{};(0,7)*+{}
   **\crv{(-5,-3)&(-5,3)}
          ?>*\dir{};
       (0,-7)*+{};(0,7)*+{}
       **\crv{(5,-3)&(5,3)}
       ?>*\dir{};
\ar@{>}(-1,6);(-0.7,6.3)
\ar@{>}(1,6);(0.7,6.3)
\ar@{-} (-20,-10); (20,10)
\ar@{-} (20,-10); (-20,10)
\ar@{-} (0,0); (-20,0)
\ar@{--} (-5,-10); (-15,0)
\ar@{-->} (-15,0); (-5,10)
\ar@{--} (-4,-10); (-9,-5)
\ar@{--} (-9,-5); (-9,5)
\ar@{-->} (-9,5); (-4,10)
\ar@{--} (15,-10); (10,-5)
\ar@{--} (10,-5); (10,5)
\ar@{-->} (10,5); (15,10)
\ar@{--} (10,-10); (0,0)
\ar@{-->} (0,0); (10,10)
\endxy
\hskip50pt
\xy
(0,-15)*{\text{(b)}};
(-7,13)*{\text{\small $\gamma$}};
(15,13)*{\text{\small $\gamma'$}};
(-3,7)*{\text{\small $\gamma_A$}};
(3,7)*{\text{\small $\gamma_B$}};
(0,-3)*{\text{\small $\frakc$}};
(0,0)*+{\bullet};
(0,-7)*+{\bullet};
(0,7)*+{\bullet};
(0,-7)*+{};(0,7)*+{}
   **\crv{(-5,-3)&(-5,3)}
          ?>*\dir{};
       (0,-7)*+{};(0,7)*+{}
       **\crv{(5,-3)&(5,3)}
       ?>*\dir{};
\ar@{>}(-1,6);(-0.7,6.3)
\ar@{>}(1,6);(0.7,6.3)
\ar@{-} (-20,-1); (0,-1)
\ar@{-} (-20, 0); (20,0)
\ar@{-} (0,1); (-20,1)
\ar@{--} (-5,-10); (-14,-1)
\ar@{-->} (-14,1); (-5,10)
\ar@{--} (15,-10); (5,0)
\ar@{-->} (5,0); (15,10)
\ar@{--} (10,-10); (0,0)
\ar@{-->} (0,0); (10,10)
\endxy
\caption{Broken lines  in the proof of Lemma 
\ref{3lem:theta1}.}
\label{3fig:broken1}
\end{figure}

First, let us consider the case $\frakj$. 
    See Figure \ref{3fig:broken1} (a).
 Suppose  that $\gamma$ crosses  the walls
 that contain  the  joint $\frakj$
 with the incoming attached monomial $c_{\mathrm{in}} x^{\tilde m_{\mathrm{in}}}$ and the outgoing attached monomial 
 $c_{\mathrm{out}}x^{\tilde m_{\mathrm{out}}}$.
There may be several broken lines  $\gamma_1=\gamma$, \dots, $\gamma_p$
with endpoint $Q$
 such that
 their bending patterns differ from $\gamma$ only 
 inside these walls.
 They have the common
  incoming attached monomial  $c_{\mathrm{in}} x^{\tilde m_{\mathrm{in}}}$
  and the common  exponent $\tilde m_{\mathrm{out}}$ for the outgoing attached monomials.
 Then, all these broken lines have a common limit 
 for $Q\rightarrow Q_0$.
 Let $\gamma_A$ be the half of a small loop around $\frakj$,
 which is on the same side of $\frakj$ with $\gamma_i$'s
 as in Figure \ref{3fig:broken1} (a).
  Then, the sum
  of the outgoing attached monomials 
  for  $\gamma_1$, \dots, $\gamma_p$ 
  is exactly the  terms of the monomials in $\frakp_{\gamma_A, \frakD_{\l}}(c_{\mathrm{in}} x^{\tilde m_{\mathrm{in}}})$
  with exponents $x^{\tilde m_{\mathrm{out}}}$.
    We repeat the same argument to the opposite side of $\frakj$
  with  the opposite half $\gamma_B$ of 
   the small loop around $\frakj$.
Then,  by the consistency of $\frakD$, we have
  \begin{align}
  \frakp_{\gamma_A, \frakD_{\l}}(c_{\mathrm{in}} x^{\tilde m_{\mathrm{in}}})
  \equiv
  \frakp_{\gamma_B, \frakD_{\l}}(c_{\mathrm{in}} x^{\tilde m_{\mathrm{in}}})
  \mod
x^{\tilde m_0} J^{\l+1}.
  \end{align}
Thus,  passing through $Q_0$ does not change 
$\vartheta_{Q,\tilde m_0}$
modulo
$x^{\tilde m_0} J^{\l+1}$.
Therefore, we obtain the equality \eqref{3eq:theta2}.

Next,  we consider the remaining case $\partial\frakd$.
    See Figure \ref{3fig:broken1} (b).
We may concentrate on a face $\frakc$ of $\frakd$ of codimension 2
such that $\frakc$ is not a joint.
Then, $\frakc$ is contained  in some walls  of $\frakD_{\l}$ with the common normal vector $n$.
In this case, for  each broken line for $Q$, there is a
unique broken line for $Q'$  with   the same attached incoming and outgoing monomial
modulo $x^{\tilde m_0} J^{\l+1}$, thanks to the consistency of $\frakD$.
Thus,  passing through $Q_0$ does not change 
$\vartheta_{Q,\tilde m_0}$
modulo
$x^{\tilde m_0} J^{\l+1}$.
\end{proof}

\begin{lem}[{\cite[Lemma~4.8]{Carl10}}]
\label{3lem:theta2}
Suppose that $Q$ and $Q'$ are in different chambers of
$\mathrm{Supp}(\frakD_{\l})$.
Then, for any admissible curve $\gamma$ from $Q$ to $Q'$,
we have
\begin{align}
\label{3eq:theta4}
\vartheta_{Q',\tilde m_0}
\equiv 
\frakp_{\gamma,\frakD_{\l}}(\vartheta_{Q,\tilde m_0})
\mod
x^{\tilde m_0} J^{\l+1}.
\end{align}
\end{lem}

\begin{proof}
We may assume that $Q$ and $Q'$ are in adjacent chambers of
$\mathrm{Supp}(\frakD_{\l})$,
since  the general case follows from this case
by composing the automorphisms $\frakp_{\gamma,\frakD_{\l}}$.
Let $n\in N_{\rmpr}^+$ be the normal vector of the walls
separating $Q$ and $Q'$.
We decompose
\begin{gather}
\vartheta_{Q,\tilde m_0}
=
\sum_{\tilde m\in \tilde M^{\circ}}
c_{\tilde m} x^{\tilde m}
=
\vartheta_+
+
\vartheta_0
+
\vartheta_-,
\\
\vartheta_+
=
\sum_{\scriptstyle \tilde m\in \tilde M^{\circ} \atop  \scriptstyle \langle n, m \rangle >0}
c_{\tilde m} x^{\tilde m},
\quad
\vartheta_0
=
\sum_{\scriptstyle \tilde m\in \tilde M^{\circ} \atop  \scriptstyle \langle n, m \rangle =0}
c_{\tilde m} x^{\tilde m},
\quad
\vartheta_-
=
\sum_{\scriptstyle \tilde m\in \tilde M^{\circ} \atop  \scriptstyle \langle n, m \rangle <0}
c_{\tilde m} x^{\tilde m}.
\end{gather}
Also, we decompose $\vartheta_{Q',\tilde m_0}=
\vartheta'_+
+
\vartheta'_0
+
\vartheta'_-
$
in the same way.
Thanks to Lemma \ref{3lem:theta1},
we may take $Q$ and $Q'$ to be very close to the wall and also move parallel to the wall
inside each chamber as we like.
Without loosing generality, one can assume that $\langle n, Q\rangle >0$ and $\langle n, Q'\rangle <0$.

\begin{figure}
\centering
\leavevmode
\xy
(-20,-3)*{\text{\small $Q'$}};
(-20,3)*{\text{\small $Q$}};
(-10,3)*+{\bullet};
(0,3)*+{\bullet};
(10,3)*+{\bullet};
(-16,-3)*+{\bullet};
\ar@{-->} (-4,9); (-16,-3)
\ar@{--} (6,9); (-3,0)
\ar@{-->} (-3,0); (-16,-3)
\ar@{--} (16,9); (7,0)
\ar@{-->} (7,0); (-16,-3)
\ar@{-} (-25,0); (25,0)
\endxy
\caption{Broken lines in the proof of Lemma 
\ref{3lem:theta2}.}
\label{3fig:broken2}
\end{figure}

\par
(i). $\vartheta_+$ and $\vartheta'_+$.
For any exponent $\tilde m$ of $ \vartheta'_+$,
 $\langle n,m\rangle >0$ holds.
Then,
any broken line for $\vartheta'_+$ 
is obtained by extending a broken line  for $\vartheta_+$.
See Figure \ref{3fig:broken2}.
Thus, we have
\begin{align}
\vartheta_+'
\equiv
\frakp_{\gamma,\frakD_{\l}}(\vartheta_+)
\mod
x^{\tilde m_0} J^{\l+1}.
\end{align}
\par
(ii). $\vartheta_-$ and $\vartheta'_-$.
The situation is opposite, and we have
\begin{align}
\vartheta_-
\equiv
\frakp_{\gamma,\frakD_{\l}}^{-1}(\vartheta'_-)
\mod
x^{\tilde m_0} J^{\l+1}.
\end{align}
\par
(iii). $\vartheta_0$ and $\vartheta'_0$.
The broken lines for both do not cross the walls separating $Q$ and $Q'$
in the final step.
Then,
the proof of Lemma \ref{3lem:theta1} applies 
to $\vartheta_0$ and $\vartheta'_0$,
and
 we have
\begin{align}
\vartheta_0'
\equiv\vartheta_0 \equiv\frakp_{\gamma,\frakD_{\l}}(\vartheta_0)
\mod
x^{\tilde m_0} J^{\l+1}.
\end{align}

We can unify the above three equalities into a single one
\begin{align}
\vartheta_{Q',\tilde m_0}
\equiv 
\frakp_{\gamma,\frakD_{\l}}(\vartheta_{Q,\tilde m_0})
\mod
x^{\tilde m_0} J^{\l+1}
\end{align}
as desired.
\end{proof}

Since $\l$ is arbitrary in \eqref{3eq:theta4},  we have the desired equality \eqref{3eq:theta3}.
This completes the proof of Proposition \ref{prop:theta1}.

\subsection{Positivity  and mutation invariance of theta functions}

Now let us specialize to a CSD $\frakD_{\fraks}$
for a given seed $\fraks$.

We have the following fundamental result on the  positivity of theta functions for a CSD.
\begin{thm}[{Positivity of theta functions \cite[Theorem~1.13 \& Remark 3.2]{Gross14}}]
\label{3thm:positivetheta1} \index{positivity (of theta functions)}
For  a CSD $\frakD_{\fraks}$, 
every theta function $\vartheta_{Q,\tilde m_0}
\in x^{\tilde m_0} \bbk[[\tilde P]]$
has only positive integer coefficients.
\end{thm}

\begin{proof}
We choose a positive realization of $\frakD_{\fraks}$
in Theorem \ref{3thm:pos1},
so that every wall function $f[tn]^s$ of $\frakD_{\fraks}$ has
only positive integer coefficients.
Then, 
the product in the right hand side
of \eqref{3eq:pgD1} has also 
only positive integer coefficients.
Thus, $c_{\gamma}$ in \eqref{3eq;mono1} is a positive integer.
 Therefore, the theta function defined by \eqref{3eq:thetadef1}
 has
only positive integer coefficients.
\end{proof}

In the rest, let us present the mutation invariance
of theta functions  of a CSD
and its consequences.
Let us define a piecewise-linear transformation in parallel to $T_{k}$ in
\eqref{3eq:Tk1},
\begin{align}
\label{3eq:Tk2}
\begin{matrix}
\tilde T_{k}=\tilde T_{k,\fraks}:& \tilde M^{\circ}& \rightarrow &\tilde M^{\circ}\\
& \tilde m & \mapsto& 
\tilde m+ [\langle \d_k e_k,  m\rangle]_+ \tilde p^*(e_k)
=
\begin{cases}
\tilde S_{k}(\tilde m) & m \in  \calH_{k}^+,\\
\tilde m & m \in  \calH_{k}^-,
\end{cases}
\end{matrix}
\end{align}
where $\tilde S_{k}(\tilde m)$ is the one in \eqref{3eq:Sk11}.

\begin{prop}[{\cite[Proposition~3.6]{Gross14}}]
\label{prop:thetamut1}
For a given seed $\fraks$, let $\fraks'=\mu_k(\fraks)$
be the mutation of $\fraks$ in direction $k$.
\par
(a). There is a one-to-one correspondence between
the broken lines for $\tilde m_0$ with endpoint $Q$
with respect to a CSD $\frakD_{\fraks}$
and 
the broken lines for $\tilde T_k(\tilde m_0)$ with endpoint $T_k(Q)$
with respect to a CSD $\frakD_{\fraks'}=T_k(\frakD_{\fraks})$.
The locus of the line is transformed by 
the piecewise-linear map $T_k$ in \eqref{3eq:Tk1}.
\par
(b). (Mutation invariance of theta functions.) \index{mutation invariance!of theta functions}
The theta function $\vartheta^{\fraks}_{Q,\tilde m_0}$ with respect to $\frakD_{\fraks}$
and the one
$\vartheta^{\fraks'}_{Q, \tilde m_0}$ with respect to $\frakD_{\fraks'}$
are related by
\begin{align}
\label{3eq:thetamut1}
\vartheta^{\fraks'}_{T_k(Q),\tilde T_k(\tilde m_0)}
=
\begin{cases}
\tilde S_k(\vartheta_{Q, \tilde m_0}^{\fraks})
& Q\in \calH_{k}^+,
\\
\vartheta_{Q, \tilde m_0}^{\fraks}
&
Q\in \calH_{k}^-.
\end{cases}
\end{align}
where  $\calH_k^{\pm}$ is the one
in \eqref{3eq:Hk1},
and $\tilde S_k$ acts on the exponents of $x$ in  $\vartheta_{Q, \tilde m_0}^{\fraks}$.
\end{prop}
\begin{proof}
(a).
In $\calH_k^{-}$ the walls and broken lines are unchanged,
while
in $\calH_k^{+}$ the walls and broken lines are transformed by the linear map $S_k$.
Thus, the rule for broken lines are preserved in each $\calH_k^{\pm}$.
Therefore, 
we only need to check the rule when a broken line crosses the wall $\bfw_{e_k}=[e_k^{\perp},
1+x^{\tilde p^*(e_k)}]_{e_k}$.
\par
Case (1). Suppose that a broken line $\gamma$   with respect to $\frakD_{\fraks}$ crosses $e_k^{\perp}$ 
from $\calH_k^{-}$ with attached monomial $cx^{\tilde m}$.
Then, $\langle e_k, -  m\rangle >0$.
After crossing $e_k^{\perp}$, the new attached monomial is a term
in $cx^{\tilde m}(1+x^{\tilde p^*(e_k)})^{|\langle \d_k e_k,  m\rangle |}$.
Applying $\tilde  S_k$, it transforms to a term in
\begin{align}
\begin{split}
&\quad\ cx^{\tilde m+\langle \d_ke_k,  m\rangle  \tilde p^*(e_k)}(1+x^{\tilde p^*(e_k)})^{|\langle \d_k e_k,  m\rangle |}
\\
&=
cx^{\tilde m-|\langle \d_ke_k,  m\rangle | \tilde p^*(e_k)}(1+x^{\tilde p^*(e_k)})^{|\langle \d_k e_k, m\rangle |}
\\
&=
cx^{\tilde m}(1+x^{\tilde p^*(e'_k)})^{|\langle \d_k e'_k, m\rangle |}.
\end{split}
\end{align}
This is the rule for a broken line  with respect to $T_k(\frakD_{\fraks})$.
\par
Case (2). Suppose that a broken line  $\gamma$
 with respect to $\frakD_{\fraks}$
 crosses $e_k^{\perp}$ 
from $\calH_k^{+}$ with attached monomial $cx^{\tilde m}$.
Then, $\langle e_k, -  m\rangle  <0$.
After crossing $e_k^{\perp}$, the new attached monomial is a term
in $cx^{\tilde m}(1+x^{\tilde p^*(e_k)})^{|\langle \d_k e_k, m\rangle |}$,
which  is rewritten as
\begin{align}
cx^{\tilde m+\langle \d_ke_k,  m\rangle  \tilde p^*(e_k)}(1+x^{\tilde p^*(e'_k)})^{|\langle \d_k e'_k, m\rangle |}.
\end{align}
Since the attached monomial for the incoming segment is transformed
from $cx^{\tilde m}$ 
to $cx^{\tilde m+\langle \d_ke_k,  m\rangle  \tilde p^*(e_k)}$,
this is  the correct rule for a broken line  with respect to $T_k(\frakD_{\fraks})$.
\par
(b).
Let $\gamma$ be a broken line with respect to $\frakD_{\fraks}$,
and let $\gamma'=T_k(\gamma)$ be
the corresponding  broken line with respect to $\frakD_{\fraks'}$.
Let $c_{\gamma} x^{\tilde m_{\gamma}}$
and $c'_{\gamma} x^{\tilde m'_{\gamma}}$
be the monomials attached with 
$\gamma$ and $\gamma'$, respectively.
Then, by (a), we have
\begin{align}
c'_{\gamma} x^{\tilde m'_{\gamma}}
=
\begin{cases}
c_{\gamma} x^{\tilde S_k(\tilde m_{\gamma})}
& Q\in \calH_{k}^+,
\\
c_{\gamma} x^{\tilde m_{\gamma}}
&
Q\in \calH_{k}^-.
\end{cases}
\end{align}
Thus, we have \eqref{3eq:thetamut1}.
\end{proof}

We have the following corollary of
Propositions \ref{3prop:theta2}
and \ref{prop:thetamut1}.

\begin{cor}[{\cite[Corollary~3.9]{Gross14}}]
\label{3cor:theta2}
Let $\frakD_{\fraks}$ be a CSD,
and let
$\calC^{\fraks}_{\fraks_{i}}$ be
the $G$-cone 
in \eqref{3eq:clusterc1}.
For $\tilde m_0\in \tilde M^{\circ}$,
suppose that $m_0\in \calC^{\fraks}_{\fraks_{i}} \cap M^{\circ}$.
Also, suppose that $Q\in \mathrm{Int}(\calC^{\fraks}_{\fraks_{i}})$.
Then, we have
\begin{align}
\label{3eq:thetain2}
\vartheta_{Q,\tilde m_0}=x^{\tilde m_0}.
\end{align}
\end{cor}
\begin{proof}
By \eqref{3eq:clusterc1},
\begin{align}
\calC^+_{\fraks_{i}}
=
(T_{k_{i-1},\fraks_{i-1}}\circ \cdots \circ T_{k_0,\fraks_{0}})(\calC^{\fraks}_{\fraks_{i}} ).
\end{align}
Let $T:=T_{k_{i-1},\fraks_{i-1}}\circ \cdots \circ T_{k_0,\fraks_{0}}$
and $\tilde T:=\tilde T_{k_{i-1},\fraks_{i-1}}\circ \cdots \circ \tilde T_{k_0,\fraks_{0}}$.
Let $S: M_{\bbR}\rightarrow M_{\bbR}$ be the linear map
that locally coincides with $T$ on the cone $\calC^{\fraks}_{\fraks_{i}}$.
Let $\tilde T$ and $\tilde S$ be the corresponding maps  for $\tilde M^{\circ} \rightarrow \tilde M^{\circ}$.
Then, applying 
\eqref{3eq:thetamut1}
repeatedly, we obtain
\begin{align}
\vartheta^{\fraks_i}_{T(Q),\tilde T(\tilde m_0)}
=
\tilde S(\vartheta_{Q, \tilde m_0}),
\end{align}
where $\tilde S$ acts on the exponents of $x$.
Meanwhile,
$\vartheta^{\fraks_i}_{T(Q),\tilde T(\tilde m_0)}=x^{\tilde T(\tilde m_0)}$
by 
Proposition \ref{3prop:theta2}.
Thus, we have
$\vartheta_{Q,\tilde m_0}=x^{\tilde m_0}$.
\end{proof}

Combining Proposition \ref{prop:theta1} and Corollary \ref{3cor:theta2},
we have another fundamental theorem on theta functions for a CSD.

\begin{thm}[{\cite[Theorem 4.9]{Gross14}}]
\label{3thm:thetaz1}
Let $\frakD_{\fraks}$ be a CSD,
and let
$\calC^{\fraks}_{\fraks_{i}}$ be
the $G$-cone 
in \eqref{3eq:clusterc1}.
For $\tilde m_0\in \tilde M^{\circ}$,
suppose that $m_0\in \calC^{\fraks}_{\fraks_{i}} \cap M^{\circ}$.
Also, suppose that $Q\in \mathrm{Int}(\calC_{\fraks}^+)$.
Then, for any admissible curve $\gamma$ from  any point in $\mathrm{Int}( \calC^{\fraks}_{\fraks_{i}})$
to $Q$, we have
\begin{align}
\label{3eq:thetain3}
\vartheta_{Q,\tilde m_0}=\frakp_{\gamma, \frakD_{\fraks}}(x^{\tilde m_0}).
\end{align}
\end{thm}

To conclude, let us describe
the relation between theta functions
and cluster variables  ($x$-variables)  for the corresponding cluster pattern with principal coefficients
in
\cite{Fomin07}.
See Section 6 of Part II for a detailed account.

As explained briefly in Section \ref{3subsec:chamber1},
each
$G$-cone $\calC^{\fraks}_{\fraks_{i}}$ in  \eqref{3eq:clusterc1}
is identified with a cone spanned by the $g$-vectors $\bfg_{1;t}$, \dots, $\bfg_{r;t}$
of the $G$-matrix $G^{t_0}_t$ for some $t\in \bbT_r$.
Let $\bfx_t=(x_{1;t},\dots,x_{r;t})$ be the  {\em cluster} at $t$
for the corresponding cluster pattern with principal coefficients.
By the celebrated \emph{Laurent phenomenon} by \cite{Fomin02}, any 
cluster variables $x_{i;t}$ ($i=1$, \dots, $r$)
are expressed as a Laurent polynomial in the initial cluster variables $\bfx_{t_0}=\bfx=(x_1,\dots,x_r)$
and the initial principal coefficients $\bfy_{t_0}=\bfy=(y_1, \dots, y_r)$, where negative exponents appear only for variables $x_i$.
We identify $x^{\tilde m}$ ($\tilde m \in \tilde M^{\circ}$) with the above variables 
as $x^{f_i}=x_{i}$ and $x^{e_i}=y_i$.
See \eqref{3eq:cseed1}.
We note that 
\begin{align}
x^{\tilde p^*(e_i)}=y_i\prod_{j=1}^r x_{j}^{b_{ji}}
\end{align}
are the \emph{$\hat{y}$-variables} $\hat{y}_i$ in \cite{Fomin07}.
Under this identification,
it follows from Theorem \ref{3thm:thetaz1}
that the cluster variables $x_{i;t}$  are identified with (a part of)
theta functions as
\begin{align}
x_{i;t}=\vartheta_{Q, \bfg_{i;t}},
\end{align}
where  the RHS is the one given in \eqref{3eq:thetain3}.

\notes
All results are taken from \cite[\S3, \S4]{Gross14}
and \cite[\S4]{Carl10}
adapted to the present formulation of scattering diagrams.

\newpage
\fancyhead[CE]{III. Appendix}
\fancyhead[CO]{Appendix \thesection. \thesectiontitle}

\appendix
\section{Source code of Ordering Algorithm}
\label{3sec:source1}
\index{Ordering Algorithm!source code}

For the reader's convenience we present a source code of Ordering Algorithm (Algorithm \ref{3algo:ord1})
 written for SageMath 9.4 \footnote{Sage~Mathematics~Software, The Sage Development Team, 
  https://www.sagemath.org.}.
It is an almost direct translation of Algorithm \ref{3algo:ord1},
and no effort is made for the efficiency and sophistication.
(If any minor or major bug is found, please let us know.)

The following examples explain 
 how to use it on Sage Notebook.

\begin{ex}
(a).
Let us order the product $\Psi[e_2]^5\Psi[e_1]^3$
modulo $G^{>5}$.
{\small
\begin{verbatim}
In:
order(5,[[0,1,5],[1,0,3]])
Out:
[[1, 0, 3],
 [3, 1, 5],
 [2, 1, 15],
 [3, 2, 125],
 [1, 1, 15],
 [2, 2, 60],
 [2, 3, 270],
 [1, 2, 30],
 [1, 3, 30],
 [1, 4, 15],
 [0, 1, 5]]
\end{verbatim}
}

\par
(b). One can order any (not necessarily anti-ordered) product
with factors of the form $\Psi[tn]^{s/t}$ $(n\in N_{\rmpr}^+, \ s,t \in \bbZ_{>0})$
as specified in Algorithm \ref{3algo:ord1}.
{
\small
\begin{verbatim}
In:
order(6,[[0,1,1],[1,2,1],[2,1,1],[2,2,1/2],[1,0,2]])
Out:
[[1, 0, 2],
 [4, 1, 1],
 [3, 1, 2],
 [2, 1, 2],
 [4, 2, 19],
 [3, 2, 16],
 [1, 1, 2],
 [2, 2, 13/2],
 [3, 3, 33],
 [2, 3, 10],
 [1, 2, 1],
 [2, 4, 9/2],
 [1, 3, 1],
 [0, 1, 1]]
\end{verbatim}
}
\end{ex}

Enjoy exploring the Badlands!

\newpage
\centerline{Source code of Ordering Algorithm for SageMath 9.4}
\begin{lstlisting}
# Ordering Algorithm for dilogarithm products
# written by Tomoki Nakanishi for SageMath 9.4
# ver.2021.10.30.b
# usage: order(L,C)

def bf(a,b):
    # bilinear form {a,b} such that {e_2,e_1}=1
    # a, b: rank 2 integer vectors
    res = a[1]*b[0]-a[0]*b[1]
    return res

def check_ordered(C):
    # check dilog product C is ordered completely
    # C has factors of [n]^{c}, c>0 any, len(C)>1
    l = len(C) # l>1
    for i in range(l-1):
        if bf(C[i][:2],C[i+1][:2]) > 0: # if {n',n}>0
            return False # not ordered
        elif bf(C[i][:2],C[i+1][:2]) == 0: # {n',n}=0
            if C[i][0]+C[i][1] > C[i+1][0]+C[i+1][1]: # deg(n')>deg(n)    
                return False # not ordered
    return True # ordered

def decompose_initial(C):
    # decompose dilog product C for the main routine (1)
    # C has factors [tn]^{s/t}
    res = []
    l = len(C)
    for i in range(l):
        if C[i][2].is_integer() and C[i][2] > 1: # if the power is integer >1
            for j in range(C[i][2]):
                res = res+[[C[i][0],C[i][1],1]] # decompose to [tn]'s
        else:
            res = res+[C[i]]
    return res
            
def decompose(p,C):
    # decompose dilog product C for subroutine (p), p>1
    # p>0: integer
    # C=[n']^{c'}[n]^{c}) is an p-exchangeable pair
    res = []
    c0 = C[0][2] # multiple of 1/p
    c1 = C[1][2] # multiple of 1/p
    for i in range(p*c0):
        res = res+[C[0][:2]+[1/p]] # decompose to [n]^{1/p}
    for i in range(p*c1):
        res = res+[C[1][:2]+[1/p]] # decompose to [n]^{1/p}
    return res

def join(C):
    # join factors with common [n] in dilog product C
    # C has factors of [n]^{c}, len(C)>1
    res = C
    l = len(res) # l>1
    flag_joined = False
    while not flag_joined:
        l = len(res)
        for i in range(l-1):
            if  res[i][:2] == res[i+1][:2]:
                res = res[:i]+[res[i][:2]+[res[i][2]+res[i+1][2]]]+res[i+2:]
                break                
            if i == l-2: # if this is the last pair, it is joined
                    flag_joined = True
    return res

def order_p_partial(L,p,C):
    # order dilog product C up to G^{>L} by p-pentagon (i.e.,{n',n}=p)
    # L,p>0: integer
    # C has factors of [n]^{1/p}, len(C)>1
    res = C
    flag_ordered_p = False
    while not flag_ordered_p: # finish if p-ordered 
        l = len(res) # l>1
        for i in range(l-1):    
            if bf(res[i][:2],res[i+1][:2]) < 0: # if ordered, do nothing
                pass
            elif bf(res[i][:2],res[i+1][:2]) == 0: # if parallel
                if res[i][0]+res[i][1] > res[i+1][0]+res[i+1][1]:
                    # if deg(left)>deg(right), commute it
                    res = res[:i]+[res[i+1]]+[res[i]]+res[i+2:]
                    break
            elif bf(res[i][:2],res[i+1][:2]) > 0: # if anti-ordered
                if res[i][0]+res[i][1]+res[i+1][0]+res[i+1][1] > L:
                    # if deg(n+n')>L, commute it
                    res = res[:i]+[res[i+1]]+[res[i]]+res[i+2:]
                    break
                elif bf(res[i][:2],res[i+1][:2]) == p:
                    # else if {n',n}=p, apply p-pentagon relation
                    if res[i][2] == 1/p and res[i+1][2]== 1/p:
                        res = res[:i]+[res[i+1]]+[[res[i][0]+res[i+1][0],\
                            res[i][1]+res[i+1][1],1/p],res[i]]+res[i+2:]
                        break
            if i == l-2: # if we reach the last pair, it is p-ordered
                flag_ordered_p = True            
    return res

def order_p(L,p,C):
    # subroutine (p); recursively defined
    # L,p>0: integer
    # for p=1, C has factors [tn]^{s/t} (n:primitive vector, s,t>0 integers)
    # for p>1, C has factors of [n]^{1/p}
    # len(C)>1
    res = C
    while not check_ordered(res): # do until completely ordered
        if p == 1: # only in the main routine (1)
            res=decompose_initial(res) # decompose to [n]
        res = order_p_partial(L,p,res) # p-order the product
        l = len(res) # l>1
        for i in range(l-1):
            if bf(res[i][:2],res[i+1][:2]) > 0: # find q-admissible pair
                q = bf(res[i][:2],res[i+1][:2])
                Cin = decompose(q,[res[i],res[i+1]]) # decompose to [n]^{1/q}
                Cout = order_p(L,q,Cin) # go to subroutine (q) (recursively)
                res = res[:i]+Cout+res[i+2:]
                break
    res = join(res)
    return res
  
def order(L,C):
    """    L: integer > 0, C: list of [n1,n2,c]
    n1,n2: integers>=0, c: rational>0 such that c*gcd(n1,n2)=integer>0 """
    # main routine (1)
    # L>0: integer
    # C has factors [tn]^{s/t} (n:primitive positive vector, s,t>0 integers)
    if not (L > 0 and L.is_integer()): # check L
        return print("the first argument is illegal")
    l = len(C) # l>=1
    for i in range(l):
        n1 = C[i][0]
        n2 = C[i][1]
        if not (n1.is_integer() and n2.is_integer()\
                and n1 >= 0 and n2 >= 0 and n1+n2 > 0): # check (n1,n2)
            return print("the second argument is illegal")            
        s = C[i][2]*math.gcd(n1,n2)
        if not (s > 0 and s.is_integer()): # check c
            return print("the second argument is illegal")
    if l == 1: # for a single factor, just return
        return C
    res = order_p(L,1,C) # enter the main routine (1)
    return res     \end{lstlisting}


\newpage
\pagestyle{myheadings}

\bibliography{../../biblist/biblist.bib}

\providecommand{\bysame}{\leavevmode\hbox to3em{\hrulefill}\thinspace}
\providecommand{\MR}{\relax\ifhmode\unskip\space\fi MR }
\providecommand{\MRhref}[2]{%
  \href{http://www.ams.org/mathscinet-getitem?mr=#1}{#2}
}
\providecommand{\href}[2]{#2}
\begin{thebibliography}{GHKK18}

\bibitem[Bri17]{Bridgeland16}
T.~Bridgeland, \emph{Scattering diagrams, {H}all algebras and stability
  conditions}, Algebraic Geometry \textbf{4} (2017), 523--561; arXiv:1603.00416
  [math.AG].

\bibitem[CI12]{Cerulli09}
G~Cerulli~Irelli, \emph{Cluster algebras of type {$A_2^{(1)}$}}, Algebras and
  Representation Theory \textbf{15} (2012), 977--1021; arXiv:0904.2543
  [math.RA].

\bibitem[CPS]{Carl10}
M.~Carl, M.~Pumperla, and B.~Siebert, \emph{A tropical view of
  {L}andau-{G}inzburg models}, preprint, 2010, available at
  https://www.math.uni-hamburg.de/home/siebert/preprints/LGtrop.pdf.

\bibitem[DM21]{Davison19}
B.~Davison and T.~Mandel, \emph{Strong positivity for quantum theta bases of
  quantum cluster algebras}, Invent. Math. (2021), published online;
  arXiv:1910.12915 [math.RT].

\bibitem[FG09]{Fock03}
V.~V. Fock and A.~B. Goncharov, \emph{Cluster ensembles, quantization and the
  dilogarithm}, Annales Sci. de l'\'Ecole Norm. Sup. \textbf{42} (2009),
  865--930; arXiv:math/0311245 [math.AG].

\bibitem[FG16]{Fock11}
V.~Fock and A.~Goncharov, \emph{Cluster poisson varieties at infinity}, Selecta
  Math. (N.S.) \textbf{22} (2016), 2569--2589; arXiv:1104.0407 [math.AG].

\bibitem[FZ02]{Fomin02}
S.~Fomin and A.~Zelevinsky, \emph{Cluster algebras {I}. {F}oundations}, J.
  Amer. Math. Soc. \textbf{15} (2002), 497--529 (electronic);
  arXiv:math/0104151 [math.RT].

\bibitem[FZ03a]{Fomin03a}
\bysame, \emph{Cluster algebras {II}. {F}inite type classification}, Invent.
  Math. \textbf{154} (2003), 63--121; arXiv:math/0208229 [math.RA].

\bibitem[FZ03b]{Fomin03b}
\bysame, \emph{Y-systems and generalized associahedra}, Ann. of Math.
  \textbf{158} (2003), 977--1018; arXiv:hep--th/0111053.

\bibitem[FZ07]{Fomin07}
\bysame, \emph{Cluster algebras {IV}. {C}oefficients}, Compositio Mathematica
  \textbf{143} (2007), 112--164; arXiv:math/0602259 [math.RT].

\bibitem[GHK15]{Gross13}
M.~Gross, P.~Hacking, and S.~Keel, \emph{Birational geometry of cluster
  algebras}, Algebr. Geom. \textbf{2} (2015), 137--175; arXiv:math/1309.2573.

\bibitem[GHKK18]{Gross14}
M.~Gross, P.~Hacking, S.~Keel, and M.~Kontsevich, \emph{Canonical bases for
  cluster algebras}, J. Amer. Math. Soc. \textbf{31} (2018), 497--608;
  arXiv:1411.1394 [math.AG].

\bibitem[GPS10]{Gross09}
M.~Gross, R.~Pandharipande, and B.~Siebert, \emph{The tropical vertex}, Duke
  Math. J. (2010), no.~153, 297 -- 362; arXiv:0902.0779 [math.AG].

\bibitem[Gro11]{Gross11}
M.~Gross, \emph{Tropical geometry and mirror symmetry}, CBMS Regional Conf.
  Ser. in Math., no. 114, Amer. Math. Soc., 2011.

\bibitem[GS11]{Gross07}
M.~Gross and B.~Siebert, \emph{From affine geometry to complex geometry},
  Annals of Math. \textbf{174} (2011), 95--138; arXiv:math/0703822.

\bibitem[IOTW]{Igusa15}
K.~Igusa, K.~Orr, G.~Todorov, and J.~Weyman, \emph{Modulated semi-invariants},
  arXiv:1507.03051 [math.RT].

\bibitem[IT21]{Igusa20}
K~Igusa and G.~Todorov, \emph{Picture groups and maximal green sequences},
  Electronic Research Archive \textbf{29} (2021), 3031--3068; arXiv:2007.14584
  [math.RT].

\bibitem[ITW]{Igusa16}
K.~Igusa, G.~Todorov, and J.~Weyman, \emph{Picture groups of finite type and
  cohomology in type {$A_n$}}, arXiv:1609.02636 [math.RT].

\bibitem[Jac79]{Jacobson79}
N.~Jacobson, \emph{Lie algebras}, Dover Publications, New York, 1979.

\bibitem[KS06]{Kontsevich06}
M.~Kontsevich and Y.~Soibelman, \emph{Affine structures and non-{A}rchimedean
  analytic spaces}, Prog. Math. \textbf{244} (2006), 321--385;
  arXiv:math/0406564.

\bibitem[KS08]{Kontsevich08}
\bysame, \emph{Stability structures, motivic {D}onaldson-{T}homas invariants
  and cluster transformations}, 2008, arXiv:0811.2435 [math.AG].

\bibitem[KS14]{Kontsevich13}
\bysame, \emph{Wall-crossing structures in {D}onaldson-{T}homas invariants,
  integrable systems and mirror symmetry}, Homological mirror symmetry and
  tropical geometry, Lect. Notes Unione Ital., vol.~15, Springer, 2014,
  pp.~197--308; arXiv:1303.3253 [math.AG].

\bibitem[Lew81]{Lewin81}
L.~Lewin, \emph{Polylogarithms and associated functions}, North-Holland,
  Amsterdam, 1981.

\bibitem[Mag54]{Magnus54}
W.~Magnus, \emph{On the exponential solution of differential equations for a
  linear operator}, Commun. Pure Appl. Math. \textbf{7} (1954), 649--673.

\bibitem[Mat21]{Matsushita21}
K.~Matsushita, \emph{Consistency relations of rank 2 cluster scattering
  diagrams of affine type and the pentagon relation}, 2021, arXiv:2112.04743
  [math.QA].

\bibitem[Mul16]{Muller15}
G.~Muller, \emph{The existence of a maximal green sequence is not invariant
  under quiver mutation}, Electron. J. Combin. \textbf{23} (2016), {\#}P2.47,
  23 pages; arXiv:1503.04675.

\bibitem[Nak21]{Nakanishi19}
T.~Nakanishi, \emph{Synchronicity phenomenon in cluster patterns}, J. London
  Math. Soc. \textbf{103} (2021), 1120--1152; arXiv:1906.12036 [math.RA].

\bibitem[NC13]{Najera11}
A.~N\'ajera~Ch\'avez, \emph{On the c-vectors and g-vectors of the {M}arkov
  cluster algebra}, S\'eminaire Lontharingien de Combinatoire \textbf{69}
  (2013), B29d;arXiv:1112.2357.

\bibitem[Rea]{Reading19}
N.~Reading, \emph{Dominance phenomena: mutation, scattering, and cluster
  algebras}, arXiv:1802.10107 [math.CO].

\bibitem[Rea14]{Reading12}
\bysame, \emph{Universal geometric cluster algebras}, Mathematische Zeitschrift
  \textbf{277} (2014), 499--547; arXiv:1209.3987.

\bibitem[Rea20a]{Reading18}
\bysame, \emph{A combinatorial approach to scattering diagrams}, Algebraic
  Combinatorics \textbf{3} (2020), 603--636; arXiv:1806.05094 [math.CO].

\bibitem[Rea20b]{Reading17}
\bysame, \emph{Scattering fans}, Int. Math. Res. Notices \textbf{23} (2020),
  9640--9673; arXiv:1712.06968.

\bibitem[Rei09]{Reineke08}
M.~Reineke, \emph{Poisson automorphisms and quiver moduli}, J. Inst. Math.
  Jussieu \textbf{9} (2009), 653--667, arXiv:0804.3214 [math.RT].

\end{thebibliography}
\fontsize{10pt}{10.5pt}\selectfont

\newpage
\printindex

\end{document}